\PassOptionsToPackage{unicode}{hyperref}
\PassOptionsToPackage{naturalnames}{hyperref}

\documentclass[a4paper,11pt]{amsart}
\usepackage[a4paper,hmargin=26mm,vmargin=28mm]{geometry}

\title[RaFH and CRaFH superalgebras]
{RaFH and CRaFH superalgebras:\\ Frobenius forms, representation theory and respectability}\def\Email#1{\email{\href{mailto:#1}{#1}}}
\author{Andrew Mathas}
\address{Andrew Mathas: School of Mathematics and Statistics,
         University of Sydney, NSW 2006,\ Australia}
\Email{andrew.mathas@sydney.edu.au}
\author{Catharina Stroppel}
\address{Catharina Stroppel: Max Planck Institute for Mathematics Bonn, Vivatgasse 7, 53111 Bonn, Germany}
\Email{stroppel@mpim-bonn.mpg.de}


\usepackage{amsfonts,amssymb,mathrsfs}
\usepackage{mathtools}
\usepackage[svgnames]{xcolor}
\usepackage{cite}
\usepackage{booktabs}
\usepackage{leftindex}
\usepackage{enumitem}
\usepackage{etoolbox}

\makeatother
\DeclareFontFamily{U}{mathx}{\hyphenchar\font45}
\DeclareFontShape{U}{mathx}{m}{n}{<5> <6> <7> <8> <9> <10>
  <10.95> <12> <14.4> <17.28> <20.74> <24.88> mathx10}{}
\makeatletter
\newcommand\bbt@glyph[2]{%
  \vcenter{\hbox{\fontsize{#2}{#2}\usefont{U}{mathx}{m}{n}\char#1}}}
\newcommand\bigboxtimes{\mathop{\mathchoice
  {\bbt@glyph{"F2}{\f@size pt}}%
  {\bbt@glyph{"D2}{\f@size pt}}%
  {\bbt@glyph{"D2}{\sf@size pt}}%
  {\bbt@glyph{"D2}{\ssf@size pt}}}%
  \displaylimits}
\makeatletter

\usepackage{atableau}
\NewATableauCommand\tableau{Tableau}{scale=0.6, box font=\scriptsize}
\NewATableauCommand\LabelledDiagram{Multidiagram}{
   empty=\emptyset,
   delimiters=[],
   label style={
     draw,circle,fill=SteelBlue!10, inner sep=0pt, minimum size=2.5mm
   }
}
\NewATableauCommand\multitableau{Multitableau}{scale=0.6, box font=\scriptsize}

\setlist[enumerate]{label=\alph*\upshape)}
\newlist{thmenumerate}{enumerate}{1}

\def\theoremname{}
\setlist[thmenumerate]{
  nosep,
  label=\upshape\alph*),
  ref=\theoremname~\theequation\upshape(\alph*),
}
\newlist{EXAMPLES}{enumerate}{1}
\setlist[EXAMPLES]{label=\arabic*.}
\newlist{Enumerate}{enumerate}{2}
\setlist[Enumerate]{label=\upshape{(\roman*)}}

\definecolor{refkey}{named}{ForestGreen}
\definecolor{labelkey}{named}{ForestGreen}

\usepackage{tikz}
\usetikzlibrary{arrows.meta}
\newcommand{\uHom}{\underline{\Hom}}
\newcommand{\uEnd}{\underline{\End}}
\newcommand{\lt}{\op{lt}}

\usepackage{todonotes}
\NewDocumentCommand\AM{sm}{\todo[color=\IfBooleanTF{#1}{MistyRose}{LightBlue},inline]{Andrew: #2}}

\def\map#1#2{\,{:}\,#1\!\longrightarrow\!#2}
\newcommand\bijection[1][\sim]{\overset{#1}{\longrightarrow}}

\DeclarePairedDelimiterX{\set}[1]{\{}{\}}{\setargs{#1}}
\NewDocumentCommand{\setargs}{>{\SplitArgument{1}{|}}m}{\setargsaux#1}
\NewDocumentCommand{\setargsaux}{mm}{\,\IfNoValueTF{#2}{#1}{#1\,\delimsize|\,\mathopen{}#2}\,}

\let\<\langle
\let\>\rangle
\let\ov\overline

\let\gedom=\trianglerighteq
\let\gdom=\vartriangleright

\newcommand\flagged[2]{%
  \begin{tikzpicture}[baseline={([yshift=-0.5ex]current bounding box.center)},x=0.38cm,y=0.38cm]
    \foreach \len [count=\r] in {#2}{
      \ifnum\r=1 \xdef\flaggedfirst{\len}\fi
      \foreach \c in {1,...,\len}{\draw (\c-1,-\r+1) rectangle (\c,-\r);}
    }
    \draw[thin] (0,0)--(-0.35,0.35);
    \node[draw,thin,circle,fill=gray!15,inner sep=0.3pt,anchor=south east,font=\scriptsize] at (-0.3,0.3) {$#1$};
    \path (\flaggedfirst+0.25,0);
  \end{tikzpicture}}

\newcommand{\Tr}{\op{Tr}}
\newcommand{\T}{\op{T}}
\newcommand\NC{\mathrm{NC}}
\tikzset{
  centreed/.style = {
     baseline = {([yshift=#1]current bounding box.centre)}
  },
  centreed/.default={-0.5ex},
  dominance arrow/.style={
    thick,
    arrows={Latex[fill=white, length=6pt]->},
  }
}

\NewDocumentCommand\Comp{ O{N,n} }{\mathcal{C}_{#1}}
\NewDocumentCommand\Parts{O{N,n}}{\Lambda_{#1}}
\NewDocumentCommand\lParts{O{\ell N,n}}{\Lambda_{#1}}

\makeatletter
\NewDocumentCommand\Rest{O{N,n}}{\mathop{\vphantom{\Lambda}\mathpalette\@restricted\relax_{#1}}\relax}
\newcommand\@restricted[2]{\ooalign{$\m@th#1\Lambda$\cr{\tiny\hidewidth$\m@th#1+$}\hidewidth\cr}}
\makeatother

\newcommand\blam{{\boldsymbol\lambda}}
\newcommand\bmu{{\boldsymbol\mu}}
\newcommand\bnu{{\boldsymbol\nu}}
\newcommand\bom{{\boldsymbol\omega}}

\newcommand\cB{\mathcal{B}}
\newcommand\cBo{\cB^{\otimes}}

\newcommand\bi{\mathbf{i}}
\newcommand\bj{\mathbf{j}}
\newcommand\ilam[1][\blam]{\bi_{#1}}

\usepackage{aliascnt}
\NewDocumentCommand\NewTheorem{O{#2}m}{%
  \newaliascnt{#2}{equation}%
  \newtheorem{#2}[#2]{#2}%
  \aliascntresetthe{#2}%
  \expandafter\def\csname #2autorefname\endcsname{#1}%
  \cspreto{#2}{\def\theoremname{#2}}
}

\def\equationautorefname~#1\null{(#1)\null}

\numberwithin{equation}{subsection}
\swapnumbers

\NewTheorem{Proposition}
\NewTheorem{Theorem}
\NewTheorem{Corollary}
\NewTheorem{Lemma}

\theoremstyle{definition}
\NewTheorem{Definition}
\NewTheorem{Example}
\NewTheorem[Example]{Examples}
\NewTheorem{Remark}
\NewTheorem{Setup}

\AtEndEnvironment{Example}{\null\hfill$\Diamond$}%

\swapnumbers
\theoremstyle{theorem}
\newcounter{maintheorem}

\newtheorem{MainTheorem}[maintheorem]{Theorem}

\theoremstyle{remark}

\newcommand{\DeclareMyOperator}[1]{%
  \expandafter\DeclareMathOperator\csname #1\endcsname{#1}
}
\forcsvlist{\DeclareMyOperator}{%
   diag, End,Hom,Ind,Irr,Mat,Res,rad,CStd,RStd,std,Std,Tab,Shaded,
   Shape,cont,rank,Add,Rem,gr,tr%
}

\DeclareMathOperator\Cl{\mathfrak{Cl}}

\DeclareMathAlphabet{\matheu}{U}{eus}{m}{n}

\newcommand\sL{\matheu{L}}

\newcommand\CA{\mathcal{A}}

\newcommand\CL{\mathcal{L}}
\newcommand\CO{\mathcal{O}}
\newcommand\CP{\mathcal{P}}

\usepackage{mathrsfs}

\newcommand\charge{{\mathfrak{z}}}
\newcommand\kapp{{\mathfrak{z}}}
\newcommand\conts{{\mathfrak{c}}}
\newcommand\Ic[1][n]{\mathtt{I}^\charge_{#1}}

\newcommand\plam[1][\lambda]{\mathfrak{p}_{#1}}
\newcommand\hlam[1][\lambda]{\mathfrak{h}_{#1}}
\newcommand\qlam[1][\lambda]{\mathfrak{q}_{#1}}

\newcommand\Olam[1][\lambda]{\CO_{#1}}
\newcommand\olam[1][\lambda]{\mathfrak{o}_{#1}}

\newcommand\fct{\mathtt{f}_{\charge,\theta}}
\newcommand\fcyc{\mathtt{f}_{\text{cyc}}}

\newcommand\gcy{\mathtt{g}}

\newcommand\tm{\mathtt{m}}
\newcommand\tn{\mathtt{n}}

\newcommand\HH{\mathbb{H}}
\NewDocumentCommand\hd{O{\xi,n}}{\mathsf{h}_{#1}}
\NewDocumentCommand\HHd{O{n}}{\HH_{#1}}
\NewDocumentCommand\Hfind{O{n}}{\HH^{\mathsf{o}}_{#1}}
\NewDocumentCommand\Hhd{O{n}}{\HH^\charge_{#1}}
\NewDocumentCommand\Hlam{O{\blam}}{\HH_{#1}}
\NewDocumentCommand\Hd{O{n}}{\HHd[#1]}

\newcommand\Sh[1][\lambda]{\op{Sh}_{#1}}
\newcommand\ELL{\text{\scriptsize top}}%

\newcommand\HJ{\underline{H}}
\newcommand\MJ{\mathbb{M}}
\newcommand\HJpm{\underline{H}^\pm}
\newcommand\Sym{W}

\NewDocumentCommand\HSd{ D(){n}}{\mathscr{H}_\eta(\Sym_{#1})}

\newcommand\Ztwo{\mathbb{Z}/2\mathbb{Z}}

\newcommand{\epso}{{\epsilon_1}}
\newcommand{\epsoo}{{\epsilon_2}}

\newcommand\Heaff{\mathscr{H}^{\text{eaff}}_n}

\newcommand\Hdegen{\mathscr{H}^{\text{degen}}_n}
\newcommand\Hsmash{\mathscr{H}^{\text{smash}}_n}

\newcommand\Taft[1][\theta]{\mathscr{T}_{#1}}

\newcommand\Basis[1][\pm]{\mathsf{B}^{#1}}

\newcommand\xipm{\xi^\pm}

\newcommand\tl{l}
\newcommand\Kr[1]{K_{#1}}
\newcommand\Kz[1]{K^{0}_{#1}}

\newcommand\ba{\mathsf{a}}
\newcommand\bd{\mathsf{d}}
\newcommand\bb{\mathsf{b}}
\newcommand\bc{\mathsf{c}}
\newcommand\be{\mathsf{e}}

\newcommand\bS{{\mathsf{S}}}
\newcommand\bT{{\mathsf{T}}}
\NewDocumentCommand\bTT{o}{{\IfNoValueTF{#1}{\mathsf{T}}{\leftindex^{s_{#1}}{\mathsf{T}}}}}
\newcommand\bU{{\mathsf{U}}}

\newcommand\bG{\mathsf{G}}

\newcommand\Tlam[1][\blam]{{\mathsf{T}_{#1}}}

\newcommand\op{\operatorname}
\usepackage{txfonts}
\newcommand{\ttau}{\tauup}
\newcommand\taui[1][i]{\tau^{(#1)}}

\newcommand\C{\mathbb{C}}
\newcommand\Po{\mathbb{P}(\mathbb{B})}
\newcommand\Poo{\mathbb{P}}
\newcommand\N{\mathbb{N}}
\newcommand\Q{\mathbb{Q}}
\newcommand\R{\mathbb{R}}
\newcommand\F{\mathbb{F}}
\newcommand\X{\mathbb{X}}
\newcommand\kk{\mathbbm{k}}
\newcommand\RR{R}
\newcommand\K{\mathbb{K}}

\newcommand\Z{\mathbb{Z}}

\newcommand\nBraid{\operatorname{M}^{\operatorname{}}}

\newcommand\Tos{\mathbb{T}}
\newcommand\Tosym{\mathbb{T}^{\op{sym}}}

\usepackage{bbm}%

\newcommand\RX{{\RR[\X]}}

\def\Verb|#1|{\texttt{$\backslash$#1}}
\usepackage[pagebackref=true]{hyperref}
\hypersetup{pdftitle={RaFH superalgebras},%
  pdfauthor={Andrew Mathas and Catharina Stroppel},%
  colorlinks = true,%
  linkcolor =blue,%
  anchorcolor = red,%
  citecolor = blue,%
  urlcolor = blue%
}
\renewcommand*{\backref}[1]{}
\renewcommand*{\backrefalt}[4]{{[\tiny%
  \ifcase #1 Not cited.\relax\or Page~#2.%
  \else Pages #2.\fi]}
}

\newcommand{\Ao}{\mathbb{A}}
\newcommand{\Bo}{\mathbb{B}}
\newcommand{\Lo}{\mathbb{L}}
\newcommand{\So}{\mathbb{S}}

\DeclareMathOperator\Smash{\#}
\DeclareMathAlphabet{\mathcal}{OMS}{cmsy}{m}{n}
\usepackage{longtable}

\begin{document}

\begin{abstract}
We introduce \emph{RaFH superalgebras}, which are Rees affine Frobenius Hecke
superalgebras. These algebras are defined by generators and relations, using \emph{Frobenius-tensor superalgebras} to axiomatise tensor powers of Frobenius superalgebras together with their copairing elements. The RaFH superalgebras include, as special cases, the
affine and degenerate affine Hecke algebras of type $A$, the affine
Sergeev superalgebras, the affine wreath product algebras and quantum
affine wreath algebras, as well as the recently introduced quantum
wreath product algebras in the special case of the superflip braiding.
We study the RaFH superalgebras from both structural and
representation theoretic points of view. We prove a basis theorem for the
RaFH superalgebras, construct universal Jucys--Murphy elements, and
classify their simple modules over an algebraically closed field via
atomic superalgebras, which are smaller RaFH superalgebras of Hecke,
Clifford or wreath type.

We then study \emph{CRaFH superalgebras}, which are the cyclotomic quotients
of the RaFH superalgebras. We construct bases for these superalgebras, and
compute the Gram determinants of natural trace forms, and use these to
determine when these superalgebras are symmetric or Frobenius. To study the
representation theory of the CRaFH superalgebras,
we introduce the \emph{respectable algebras}, which is a
vast generalisation of cellular algebras, standardly based algebras and finite
dimensional affine cellular algebras, and give many examples
illustrating their variety and ubiquity.
We then show that CRaFH superalgebras are respectable algebras with
weight poset given by multipartitions that are coloured by the
representations of the underlying Frobenius algebra. This result can be
viewed as a far-reaching generalisation of Murphy's bases into the
setting of CRaFH algebras.
\end{abstract}

\maketitle

\tableofcontents
\section*{Introduction}
The representation theory of the symmetric group is a rich and beautiful
subject that starts with the seminal work of
Frobenius~\cite{Frobenius:symmetric} and Young~\cite{QSAI} in the
semisimple case and has undergone a remarkable transformation in the last
thirty years. Starting with Kleshchev's modular branching
rules~\cite{Klesh:I}, the LLT conjecture~\cite{LLT}, its
subsequent proof and generalisation by Ariki~\cite{Ariki:can},
Grojnowski's categorification of the crystal graph~\cite{Groj:control},
and Chuang and Rouquier's, \cite{ChuangRouq:sl2}, introduction of $\mathfrak{sl}_2$
categorifications and their proof of the Brou\'e's abelian defect group
conjecture for symmetric groups, we now view the modular representation
theory of the symmetric groups through the combinatorics of the Fock
space and 2-representations of 2-Kac--Moody
algebras~\cite{KhovLaud:diagIII,Rouq:2KM,BK:GradedKL}. Under this
categorification, the $i$-induction and $i$-restriction functors for the
symmetric group correspond to Chevalley generators of the affine
Kac--Moody algebra of type~$A$~\cite{Ariki:can,BK:GradedDecomp}. The
2-morphisms play a crucial role, since they control the Kac--Moody
algebra relations and encode the representation theory of the symmetric
groups. Remarkably, all this information is already captured by the
categorical action of a Heisenberg Lie
algebra~\cite{Brundan:Heisenberg}, as pioneered
in~\cite{Khovanov:Heisenberg}, see also \cite{LicataSavage:Heisenberg, LicataSavage:Hecke,
CautisLicata:Heisenberg, LicataSavage:Survey}.
The defining relations of Heisenberg Lie algebras are significantly
simpler than those of Kac--Moody algebras, which makes categorical
Heisenberg actions more elegant and less demanding to construct.

The general framework of Heisenberg--Frobenius
categorification~\cite{CautisLicata:Heisenberg,BSWHeis,LicataSavage:Heisenberg,
BrundanSavageWebster:QuantumHeisenberg,
BrundanSavageWebster:FoundationsHeisenberg} attaches a diagrammatic
monoidal $2$-(super)category to each Frobenius superalgebra. The KLR
side of this story is missing when the Frobenius algebra is non-trivial.
Moreover, there is no uniform treatment of the many algebras, see
e.g.~\cite{RossoSavage:QuantumAffineWreath,Savage:AffineWreath,wreathHecke,WanWang,CT:wreath,Pushkarev}, involved in Heisenberg categorifications. Instead, the structure and
representation theory of each is largely developed from scratch each time. This problem also appears in the proof of the Brundan--Kleshchev isomorphism \cite{BK:GradedKL}, where the degenerate and non-degenerate affine Hecke algebras have to be treated by different arguments.

This lack of uniform treatment is the starting point of our work. In \autoref{D:Hecke} we
introduce the \emph{Rees affine Frobenius Hecke algebras}, or \emph{RaFH
algebras}, to provide a common framework that encompasses the different
Hecke--Frobenius algebras in the literature. The RaFH algebra
$\HHd=\HHd(A)$ depends on a finite dimensional Frobenius superalgebra
and parameters $\eta$, $t$ and $\hbar$, which allow $\HHd$ to be
specialised to different Hecke--Frobenius algebras. For example, by
taking the Frobenius algebra to be a field, the corresponding RaFH
algebras include both the degenerate and non-degenerate affine Hecke
algebras of type~$A$. When $A$ is a Frobenius algebra, we recover the
wreath product algebras of
\cite{Savage:AffineWreath,RossoSavage:QuantumAffineWreath}. When $A$ is
a Clifford algebra, we obtain the degenerate affine Sergeev
algebras~\cite{Sergeev, Nazarov}. Specialising $\hbar=t=1$, $\eta=0$ gives the affinisation algebras of  Kleshchev and Muth~\cite{KleshchevMuth}. The latter work also indicates that our framework covers algebras which appear in modular representation theory such as Turner algebras for blocks of symmetric groups.

Our approach has considerable overlap, but is
different from, recent work of Cai, Minets, Nakano, and
Xiang~\cite{LaiMinets,MinetsYBE,CJNakanoXiang}. In particular, we obtain representation theoretical results for a nice subfamily of the algebras studied in their papers.

The RaFH algebras are defined by generators and relations as a unital
associative $R$-algebra, where~$R$ is a commutative ring containing the
parameters $\eta$, $t$ and $\hbar$. Like the degenerate affine Hecke algebras, the
RaFH algebras have polynomial generators, $X_1,\dots,X_n$ and Hecke
generators and, in addition, the elements of the Frobenius-tensor superalgebra $\Bo=A^{\otimes n}$,
for a Frobenius superalgebra~$A$ with Nakayama automorphism~$\psi$. The
relations between most of the generators are what you might expect. In
particular, the Hecke generators act on $A^{\otimes n}$ with the
\emph{superflip} action, which permutes tensor factors up to sign.

The most problematic relations in $\HHd$ are the quadratic
relations for the Hecke generators, which are
\[
     H_i^2 = \eta\tau^{(i,i+1)}H_i+\hbar^2,\qquad\text{for }1\le i<n.
\]
When the Frobenius superalgebra $A=R$ is trivial, then
$\tau^{(i,i+1)}=1_R$, so this relation becomes the familiar quadratic
relation for a degenerate or non-degenerate Hecke algebra and the usual square relation in the symmetric group when $\eta=0$.

To treat the
cases where $A$ is not trivial, we develop a general theory of
\emph{Frobenius-tensor superalgebras} by introducing a \emph{universal
Frobenius monoid} (\autoref{univtaurel}), which is the monoid generated by
elements $\set{\tau^{(i,j)}}$ that satisfy similar relations to those of
the Birman--Ko--Lee dual braid monoid~\cite{BirmanKoLee}. This monoid
is the origin of the many seemingly different \emph{Jucys--Murphy}
elements in the literature. Fix a basis $\cB$ of $A$ and let
$\set{b^\vee|b\in\cB}$ be a dual basis with respect to the
non-degenerate trace form on~$A$. There is an action of the
universal Frobenius monoid on $A^{\otimes n}$ where the generator
$\tau^{(i,j)}$ is
sent to the \emph{copairing element} $\sum_{b\in\cB}b^{(i)}b^{\vee(j)}$,
where $a^{(k)}=1_A\otimes\dots\otimes a\otimes\dots\otimes 1_A$ embeds
$a\in A$ into the $k$th tensor factor for $1\le k\le n$. The copairing
elements in $A^{\otimes n}$ are a remarkable family of elements in
$A^{\otimes n}$, and many of their properties are explained by the
universal Frobenius monoid.  This motivated our general definition of Frobenius-tensor superalgebras $\Bo=A^{\otimes n}$ with a Nakayama automorphism $\psi$ of $A$.

We allow non-trivial Nakayama automorphisms, but the
relations of $\HHd$ force us to assume that $\Bo$ is
\emph{$\eta$-symmetric}, which means that either $\eta=0$ or that the Nakayama automorphism is trivial.
In case of $\eta$-symmetricity, we are able to prove the following
basis theorem for $\HHd$.

\begin{MainTheorem}[\autoref{C:BasisHd}]\label{M:HdBasis}
  Suppose that $\Bo$ is an $\eta$-symmetric Frobenius-tensor
  superalgebra. Then the RaFH algebra $\Hd=\Hd(\Bo)$ is free as an
  $\RR$-module with basis $\set{X^\bd H_w b|\bd\in\N^n, w\in\Sym\text{
  and }b\in\chi}$, where $\chi$ is any basis of $\Bo$.  In particular,
  $\HHd$ is free as a (right) $\Bo$-module.
\end{MainTheorem}

We show that the existence of a basis of this form forces $\Bo$ to be
$\eta$-symmetric, so this assumption is necessary in
\autoref{M:HdBasis}. \autoref{M:HdBasis} is proved using a theory of
Frobenius--Demazure operators for $\HHd$. It could be deduced from
\cite{CJNakanoXiang}, but we believe it is worthwhile to give a more direct proof.

In our theory, the RaFH algebras play the same role as the (extended)
affine Hecke algebras do in the representation theory of the Hecke algebras of
Weyl groups and the cyclotomic Hecke algebras of complex reflection
groups. There, the cyclotomic quotients of affine Hecke algebras are important,
because every simple module of an affine Hecke algebra is a simple
module for some (in fact, infinitely many) cyclotomic quotient. The
cyclotomic quotients are also important because they are used to categorify irreducible highest
weight modules.

We therefore  introduce CRaFH algebras, which are certain
cyclotomic quotients of the RaFH algebras.
In the literature, such quotients are defined differently in different contexts. As we are
working in a very general framework, with an arbitrary finite
dimensional Frobenius algebras that can have a non-trivial Nakayama
automorphism, we are forced to do this more carefully. As a result, we
only consider cyclotomic quotients $\Hhd=\HHd/(\fcyc(X_1))$, where the
cyclotomic polynomial $\fcyc(x)$ is compatible with the Nakayama
automorphism. More precisely, if $\psi$ has order $\theta\ge1$ then
\[
    \fcyc(x)= x^{\ell_0}
        \prod_{i=1}^{\ell_1}(x^\theta-\kapp_i^\theta)\in\RR[x],
\]
for integers $\ell_0,\ell_1\ge0$ and cyclotomic parameters
$\kapp_1,\dots,\kapp_{\ell_1}\in R$; see \autoref{D:fcyc}. The
definition of $\fcyc(x)$ ensures that the roots of the cyclotomic
polynomial are permuted by the Nakayama automorphism. This
definition explains several \textit{ad hoc} definitions of cyclotomic
Hecke algebras in the literature, see e.g. \cite{BKHeckeClifford,LiShi:SergeevCocenter}. The specific choice of cyclotomic polynomial is important for our next main result  which is a
proof that the CRaFH algebras are free as $\RR$-modules with the expected
basis.

\begin{MainTheorem}[\autoref{basisthmcycl}]\label{M:HhdBasis}
  The CRaFH algebra $\Hhd=\Hhd(\Bo)$ is free as an $R$-module with basis
  \[
    \set{X^\bd H_w b|\bd=(b_1,\dots,b_n)\in\N^n
       \text{ with }0\le b_i<\ell, w\in\Sym,\text{ and $b$ runs through a fixed basis of $\Bo$}}.
  \]
\end{MainTheorem}

Frobenius algebras are characterised by the existence of a
non-degenerate trace form $\tr\colon A\to R$, together with a Nakayama
automorphism $\psi\colon A\to A$ of order $\theta$. A non-degenerate
trace form on $A$ extends to a non-degenerate trace form $\tr\colon
A^{\otimes n}\to R$, given by $\tr(b_1\otimes\dots\otimes
b_n)=\prod_{i=1}^n\tr(b_i)$, for $b_i\in A$.

\begin{MainTheorem}[\autoref{main:bilinearform}]\label{M:TraceForm}
  Suppose that $x,y\in\Hhd$ are homogeneous. Then
\[
  \Tr(xy)=(-1)^{|x||y|}\Tr(y\psi^{\ell_0}(x)),
  \qquad\text{where $\psi(X^\bd H_w\bb)=X^\bd H_w\psi(\bb)$.}
\]
In particular, the bilinear form $\<\ ,\ \>\map{\Hhd\times\Hhd} R$
induced by $\Tr$ is symmetric if $\theta$ divides $\ell_0$. Conversely,
if $\<\ ,\ \>$ is non-degenerate and symmetric, then $\theta$ divides
$\ell_0$.
\end{MainTheorem}

This theorem recovers and explains many existing results in the
literature, such as the trace forms on the degenerate and non-degenerate
Hecke algebras \cite{MalleMathas,BK:HigherSchurWeyl}. The parity
conditions in \autoref{M:TraceForm} also explains the symmetric and
supersymmetric trace forms for the Sergeev algebras that appear in
\cite{LiShi:SergeevCocenter}.

The proof that the trace form in \autoref{M:TraceForm} is non-degenerate
is difficult, and also quite different from the existing proofs for the
special cases considered in the literature. A key idea is to first
prove non-degeneracy when $n=1$. This proves that the CRaFH
algebras $\Hhd[1]$ are Frobenius superalgebras. This allows us to
consider the CRaFH algebra $\Hhd(\Hhd[1])$ and, using what we call the
\emph{RaFH-CRaFH trick}, see~\autoref{SS:Trick}, to boot-strap the proof of
\autoref{M:TraceForm} from the case $n=1$.

Having established the $\RR$-module structure of the RaFH and CRaFH
algebras, we next turn to their representation theory. We first describe
the irreducible representations of the CRaFH algebras by building on
earlier work of Savage and
others~\cite{RossoSavage:QuantumAffineWreath,Savage:AffineWreath,WanWang}.
Given a classification of the simple modules for a Frobenius superalgebra $A$ we derive in \autoref{T:simples}
a classification and construction of the simple modules of the CRaFH
algebras in terms of \emph{atomic algebras}. This result reduces the problem to the classification of simple modules for  (degenerate or non-degenerate) affine Hecke or Sergeev algebras.

Besides the classification result we
want a construction of the simple modules that is closer to that given
by cellular algebras, in the sense of \cite{GrahamLehrer:Cellular}. Unfortunately, many Frobenius algebras are not cellular, so a cellular
algebra approach cannot work for CRaFH algebras. To solve this problem
we introduce  \emph{respectable} algebras in~\autoref{D:Respectable},
which are a more general class of algebras that include cellular
algebras as a special case. The notion generalises to algebra objects in well-behaved symmetric monoidal categories such as finite dimensional super or graded vector spaces and therefore covers superalgebras as well.

A respectable algebra is a finite dimensional algebra $\CA$ that
has a \emph{respectable datum}, which is a tuple
$(\Lambda,\ge,\plam[],\hlam[],\olam[],\qlam[],\Olam{},a)$, where
$(\Lambda,\ge)$ is a poset and for each $\lambda\in\Lambda$ there are
indexing sets $\plam$, $\hlam$, $\olam$, and $\qlam$, where $\olam$ is
an $R$-basis for the \emph{fibre algebra} $\Olam$, which is a unital
$R$-algebra, and
\[
  a\colon\coprod_\lambda\plam\times\hlam\times\olam\times\qlam\to \CA;
  (p,h,o,q)\mapsto a_{phoq}
\]
is an injective map. The first axiom for
respectable algebras requires that $\set{a_{phoq}}$ be a basis of $\CA$.
The second axiom requires that the action of~$\CA$ on the basis
$\set{a_{phoq}}$, from the left and from the right, takes the same form as
the action of a cellular algebra on its cellular basis. This allows us
to define \emph{standard} $\Delta^\lambda$ and \emph{costandard} modules
$\nabla^\lambda$ for respectable algebras. The final axiom for
respectable algebras says, in effect, that $\Delta^\lambda$ and
$\nabla^\lambda$ are $(A,\Olam)$-bimodules. We show that the
representation theory of respectable algebras closely resembles that of
cellular algebras. In particular, there is a balanced pairing that
allows us to define a submodule $\rad\Delta^\lambda$ of $\Delta^\lambda$
and the quotient module $L^\lambda=\Delta^\lambda/\rad\Delta^\lambda$.
We show that $L^\lambda$ is irreducible when $\Olam$ is a division
algebra and the balanced pairing is non-degenerate. Let
$\Lambda_+=\set{\lambda\in\Lambda|L^\lambda\ne0}$, together with parity
shifts in the super case. If $\Olam$ is a division ring for
$\lambda\in\Lambda_+$, then $\set{L^\lambda|L^\lambda\ne0}$ is a
complete set of pairwise non-isomorphic irreducible $\CA$-modules. Unlike
cellular algebras, the modules $\set{L^\lambda|L^\lambda\ne0}$ need not
be absolutely irreducible. Indeed, $\End_\CA(L^\lambda)\cong\Olam$, so the
fibre algebras $\Olam$ detect when $L^\lambda$ is absolutely
irreducible. In the super case, the fibre algebras detect whether the
irreducible modules are of type $\mathtt{M}$ or of type $\mathtt{Q}$.

\begin{MainTheorem}[\autoref{T:HfindRespectable} and
   \autoref{T:HhdRespectable}]
   Suppose that $A$ is a perfectly $R$-respectable Frobenius algebra. Then
   the CRaFH algebra $\Hhd(A^{\otimes n})$ is perfectly $R$-respectable.
\end{MainTheorem}

Superficially, the proof of this result is similar to the construction
of the cellular basis of the Ariki--Koike algebras given in
\cite[Theorem~3.26]{DJM:cyc}. In practice, the proof is considerably
more difficult. The first step in this construction is to describe
analogues of permutation modules, or spherical and anti-spherical
modules, for the CRaFH algebras. Since the quadratic relation for~$H_i$
involves the copairing element $\tau^{(i,i+1)}$, diagonalisibility of
the action of~$H_i$ is not guaranteed and relies on the existence of
suitable square roots of expressions involving the copairing elements.
This problem is the subject of \autoref{S:Diagonalisable}, where we show
how to diagonalise $H_i$ and use this to define coefficients
$\xi^\pm_w$, for each permutation $w$, that give the action of $H_w$ on
the one dimensional representations of~$\Hhd$. In the classical case,
when $A=R$, $\eta=q-1$ and $\hbar^2=q$, $\xi^\pm_w=(\pm q)^{\pm\ell(w)}$
but, in general, these coefficients do not commute and depend on a
choice of reduced expression.  As a by product, we get the following
result:

\begin{MainTheorem}[\autoref{C:MJBasis}]
Given a RaFH algebra $\HHd$ which is braided diagonalisable and assume $J\subset\set{1,\ldots,n}$. Then there exist elements $\HJ^\pm_J\in\HHd$ such that  $\HHd\HJ^\pm_J$ are  RaFH-analogues of the spherical and antispherical Hecke modules corresponding to $J$.
\end{MainTheorem}

Using these constructions we are then able to prove
analogues of the Garnir relations for a respectable basis of $\Hhd$, and
hence show that $\Hhd$ is a respectable algebra. In practice, the proof
is subtle and technical, but the pain is worth it because it builds a
``Specht module theory'' for the CRaFH algebras attached to respectable
Frobenius algebras.

\subsubsection*{AI disclaimer} The second author used Claude Opus 5 to
turn piles of handwritten notes into \LaTeX. We also used AI to write
programs that helped us to compute the $n=1$ Gram determinants of the
CRaFH algebras.  Finally, we used AI to check for typos and to make an
index of notation.

\subsubsection*{Acknowledgements}
We thank the Isaac Newton Institute \textit{Modular Lie theory} program at the
University of York, July 2025, where this work was almost finished
before we got sidetracked again. We thank CJ Lai, Sasha Minets, Rob
Muth for several useful conversations and for sharing many ideas. The second author thanks Tony Licata for discussions in 2012 which she finally could appreciate. AM was partially
supported by the Australian Research Council DP240101809. CS was supported by the Gottfried Wilhelm Leibniz-Preis of the German Research Foundation. Both authors
thank the Max--Planck--Institute for Mathematics, Bonn, for their
incredible support. Most of the work was done there.

\section{Notation and conventions}\label{S:notation}
  We fix a ground ring $\kk$, which is a commutative integral domain
  with $1$. Unless otherwise stated, algebras, tensor products, linear maps,
  spans, and other standard constructions are taken over $\kk$. Let
  $\iota$ be an element of $\kk$ such that $\iota^2=-1$, if such an
  element exists. By an \emph{algebra} we always mean an associative unital algebra. If
  not otherwise stated, modules and algebras are free over $\kk$. By
  convention, graded $\kk$-modules have graded pieces that are free as
  $\kk$-modules.  For an
  algebra $A$, let $Z(A)=\set{z\in A|za=az\text{ for }a\in A}$ be the
  \emph{centre} of $A$.

  By a \emph{ring} we mean a commutative ring with one. A \emph{ring
  homomorphism} is a homomorphism of unitary commutative rings.
   Given a group $G$, let $\kk G$ be the group
  algebra of~$G$ over~$\kk$, with standard basis given by the elements
  $\set{g|g\in G}$.

 Let $\Ztwo=\{\ov{0},\ov{1}\}$ be the additive group of order two and
  write $({-1})^{\ov{0}}=1$ and $({-1})^{\ov{1}}=-1$. This group will
  encode supergradings.  We let $\mathtt{C}_m$ be the multiplicative cyclic group of order $m$.

  Throughout this paper we fix a non-negative integer $n$, and we
  fix indeterminates $X_1,\dots,X_n$ over $\kk$ and let
  $\kk[\X]=\kk[X_1,\dots, X_n]$. We work with the \emph{symmetric group}
  $\Sym=\Sym_n$ of permutations of $\set{1,\dots, n}$.  The group~$\Sym$
  is generated by the set of simple transpositions
  $S=\set{s_1,\dots,s_{n-1}}$, where $s_i=(i,i+1)$, subject to the usual
  Coxeter relations. Hereafter, we identify the concrete transpositions
  $s_i=(i,i+1)\in\Sym$ with the abstract generators in the Coxeter
  presentation of $\Sym$:
  \[
       \Sym\cong\bigl\<s_1,\dots,s_{n-1}\mid
            s_i^2=1, s_is_{i+1}s_i =s_{i+1}s_is_{i+1},
            \text{ and } s_is_j=s_js_i \text{ if }|i-j|>1\bigr\>.
  \]
  The \emph{length} of $w\in\Sym$ is
  $\ell(w)=\min\set{k\ge0\mid w=s_{i_1}\dots s_{i_k}\text{ for some }1\le i_j\le n-1}$. A
  \emph{reduced expression} for $w$ is an expression $w=s_{i_1}\dots
  s_{i_k}$ with $k=\ell(w)$.

\section{Frobenius-tensor superalgebras and universal Jucys--Murphy elements}
The  affine Hecke algebras of type $A$ from the literature as well as their Clifford or wreath product analogues,
\cite{IwahoriMatsumoto,Drinfield:DegenerateHecke,Nazarov,BKHeckeClifford,Savage:AffineWreath,RossoSavage:QuantumAffineWreath,KleshchevMuth},
contain a tensor power $\Ao$ of a Frobenius (super)algebra $A$. Their defining relations involve
specific elements $\tau^{(i,j)}\in\Ao$ obtained by embedding the copairing element $\sum_bb\otimes b^\vee\in A\otimes A$ given by the Frobenius structure into the $i$th and $j$th tensor factors of $\Ao$.

Many calculations  in these (super)algebras  rely on a few relations among the $\tau^{(i,j)}$ together with their interaction with the symmetric group and the Nakayama automorphism of $A$. In this section we isolate and axiomatise these relations by first introducing the universal Frobenius monoid $\nBraid_n$, which  is closely related to the Birman--Ko--Lee dual braid monoid~\cite{BirmanKoLee}, and then defining Frobenius-tensor superalgebras.  The latter will be used in \autoref{S:AffineRH} to define the
RaFH superalgebras which  are the main player in this paper and a common generalisation of all the above algebras. On the way, we  introduce universal Jucys--Murphy elements which provide a uniform treatment of several known families of Jucys--Murphy elements in the literature.

\subsection{The universal Frobenius monoid algebra and  universal Jucys--Murphy elements}\label{S:Smash}
This subsection introduces the universal Frobenius monoid with its
monoid algebra $\Tos$ and defines universal Jucys--Murphy elements. We
finish by relating them to classical Jucys--Murphy elements.

  Let $G$ be a group and~$B$ an algebra upon which~$G$ acts as a group of automorphisms, where $g\in G$ sends $b\in B$ to
  $\leftindex^g b$. Extending linearly, $\kk G$ acts on~$B$. The \emph{smash product algebra} $\kk G\Smash B$ is the $\kk$-module $\kk G\otimes B$ with multiplication
  \[
   (g\otimes b)(g'\otimes b')=(gg')\otimes\left(({}^{(g')^{-1}}b)b'\right)
          \qquad\text{for $g,g'\in G$ and $b,b'\in B$}.
  \]
  If $B$ is free as $\kk$-module with
  basis $\cB$, then $\kk G\Smash B$ is free over $\kk$
  with basis $\set{g\otimes b| g\in G\text{ and }b\in \cB}$.

  We consider smash products for $G=\Sym$ with
  algebras built from a monoid which we now  introduce.

   \begin{Definition}\label{univtaurel}
     The \emph{universal Frobenius monoid} $\nBraid=\nBraid_n$ is the
     monoid generated by $\set{\ttau^{(i,j)}|1\leq i\not=j\leq n}$,
     modulo the relations
   \begin{align}\label{univtaureleq}
    \ttau^{(i,j)}\ttau^{(k,l)}=\ttau^{(k,l)}\ttau^{(i,j)}\quad\text{and}\quad
   \ttau^{(i,j)}\ttau^{(j,k)}=\ttau^{(j,k)}\ttau^{(i,k)}=\ttau^{(i,k)}\ttau^{(i,j)}
   \quad\text{for distinct $i,j,k,l$}.
   \end{align}
   \end{Definition}

   \begin{Lemma}\label{braidlemma}
In $\nBraid$, the following generalised braid relations hold for pairwise distinct $1\leq i,j,k\leq n$,
 \begin{align*}
 \ttau^{(i,j)}\ttau^{(j,k)}\ttau^{(i,j)}=\ttau^{(j,k)}\ttau^{(i,j)}\ttau^{(j,k)}
 \quad\text{and}\quad
  \ttau^{(i,k)}\ttau^{(i,j)}\ttau^{(i,k)}=\ttau^{(i,j)}\ttau^{(i,k)}\ttau^{(i,j)}.
  \end{align*}

\end{Lemma}
\begin{proof} The equalities follow directly by applying \autoref{univtaureleq} twice, $ (\ttau^{(i,j)}\ttau^{(j,k)})\ttau^{(i,j)}=\ttau^{(j,k)}\ttau^{(i,k)}\ttau^{(i,j)}=\ttau^{(j,k)}(\ttau^{(i,j)}\ttau^{(j,k)})$ and $(\ttau^{(i,k)}\ttau^{(i,j)})\ttau^{(i,k)}=\ttau^{(i,j)}\ttau^{(j,k)}\ttau^{(i,k)}=\ttau^{(i,j)}(\ttau^{(i,k)}\ttau^{(i,j)})$.
\end{proof}
\begin{Remark}
 The universal Frobenius monoid $\nBraid$ does not appear to be in the
 existing literature, but it is connected to well-known monoids. The
 Birman--Ko--Lee dual braid monoid from
 \cite{BirmanKoLee,Bessismonoid} is the monoid generated
 by $\set{\ttau^{(i,j)}|1\le i<j\le n}$ subject to the relations
 \autoref{univtaureleq} except that for the first relation we assume
 $1\leq  i < j < k<l\leq n$ or $1\leq i < k< l< j\leq n$. The dual
 braid monoid can also be realised as a monoid generated by
 $\set{\ttau^{(i,j)}|1\le j<i\le n}$ by replacing $\ttau^{(i,j)}$ with
 $\ttau^{(n-i+1,n-j+1)}$,  see \cite{Fromentin}.  Thus, the universal
 Frobenius monoid is a quotient of the free product of two copies of the dual braid monoid,
 namely the monoids generated by
  $\set{\ttau^{(i,j)}|1\le i<j\le n}$ and $\set{\ttau^{(i,j)}|1\le
  j<i\le n}$, respectively.

The dual braid monoid contains the positive braid monoid, which is the Artin braid monoid of type $A_{n-1}$ as the submonoid generated by $\set{\ttau^{(i,i+1)}|1\le i<n}$, and is itself
  a submonoid of the Artin braid group $\operatorname{Br}_n$. Despite \autoref{braidlemma}, one can check that there is
  no monoid morphism from $\nBraid$ to $\operatorname{Br}_n$ sending
  $\ttau^{(i,i+1)}\in\nBraid$ to $\ttau^{(i,i+1)}\in\operatorname{Br}_n$.
  However, there is a morphism from $\nBraid$ to the symmetric group $\Sym$, which sends $\ttau^{(i,j)}$ to the transposition $(i,j)\in\Sym$.
\end{Remark}

\begin{Definition}\label{D:universalJM}
Let $\Tos=\kk\nBraid$ be the monoid algebra of the universal Frobenius
monoid~$\nBraid$. The symmetric group $\Sym$ acts on $\Tos$ by permuting the
superscripts on the generators, that is
  $\leftindex^w\ttau^{(i,j)}=\ttau^{(w(i),w(j))}$, for
  $w\in\Sym$. Hence, we can construct the smash
  product algebra $\kk\Sym\Smash\Tos$.
  For $1\leq k\leq n$, define the \emph{$k$th universal Jucys--Murphy element} $\CL_{k}\in\kk\Sym\Smash\Tos$ as $\CL_{k}=\sum_{i=1}^{k-1}(i,k)\ttau^{(k,i)}=\sum_{i=1}^{k-1}\ttau^{(i,k)}(i,k)\in\kk\Sym\Smash\Tos$.
  \end{Definition}
  \begin{Remark}
  We have in particular that  $\CL_1=0$, because empty sums are zero  by convention.
    \end{Remark}
  The next two results motivate calling $\CL_{k}$ a universal
  \emph{Jucys--Murphy} element.
  \begin{Proposition}\label{P:CLProperties}
The following equalities hold in $\kk\Sym\Smash\Tos$ for any $w\in\Sym$ and $1\leq k\not=l\leq n$.
\begin{enumerate}
\item\label{La}$w\ttau^{(k,l)}w^{-1}=\ttau^{(w(k),w(l))}$,
\item\label{Lb} $s_k\CL_{l}=\CL_l s_k$ if $l\ne k, k+1$, whereas
 $s_k\CL_k=\CL_{k+1}s_k-\ttau ^{(k,k+1)}$,
 and
$s_k\CL_{k+1}=\ttau^{(k+1,k)}+\CL_ks_k$,
\item\label{Lc} $\CL_k\CL_l=\CL_l\CL_k$.
\end{enumerate}
In particular, the universal Jucys--Murphy elements generate a
commutative subalgebra of $\kk\Sym\Smash\Tos$.
\end{Proposition}
\begin{proof}
Part  \autoref{La} holds by definition. It implies
$s_k\CL_k=s_k\sum_{i=1}^{k-1}(i,k)\ttau^{(k,i)}=\sum_{i=1}^{k-1}(i,k+1)\ttau^{(k+1,i)}s_k$ which gives the
second and, by conjugating with $s_k$,  also the third relation in~\autoref{Lb}.

For part~\autoref{Lc}, we have on the one hand $\CL_k\CL_l=\sum_{i<k,
j<l}(i,k)\ttau^{(k,i)}(j,l)\ttau^{(l,j)}$ and on the other hand
$\CL_l\CL_k=\sum_{i<k, j<l}(j,l)\ttau^{(l,j)}(i,k)\ttau^{(k,i)}$. We compare these and may
assume $k<l$. The summands with $j\notin\set{i,k}$ agree in the two
expressions. By part~\autoref{La}, the remaining terms in $\CL_k\CL_l$
are equal to
\[\sum_{1\leq i<k} (i,k)\ttau^{(k,i)}(i,l)\ttau^{(l,i)}+\sum_{1\leq i<k}(i,k)\ttau^{(k,i)}(k,l)\ttau^{(l,k)}
=\sum_{1\leq i<k} (i,k)(i,l)\ttau^{(k,l)}\ttau^{(l,i)}+\sum_{1\leq i<k}(i,k)(k,l)\ttau^{(l,i)}\ttau^{(l,k)}.
\]
By part~\autoref{La}, and the relations \autoref{univtaureleq} of the universal Frobenius monoid, this coincides with
\[
\sum_{1\leq i<k}
(i,l)(l,k)\ttau^{(k,i)}\ttau^{(k,l)}+\sum_{1\leq i<k}(i,l)(i,k)\ttau^{(l,k)}\ttau^{(k,i)}
=\sum_{1\leq i<k} (i,l)\ttau^{(l,i)}(k,l)\ttau^{(k,l)}+\sum_{1\leq i<k}(i,l)\ttau^{(l,i)}(i,k)\ttau^{(k,i)}.
\]
Since the latter is the sum of the unaccounted terms in $\CL_l\CL_k$, part~\autoref{Lc} follows.

If $l<k$ then $s_k\CL_l=\CL_ls_k$ by definition. If $l>k+1$
 then $s_k\CL_l=\CL_ls_k$ by a similar calculation to that used to
prove part~\autoref{Lc}. This completes the proof of part~\autoref{Lb}.
\end{proof}
We call  $\Tosym=\Tos/(\ttau^{(i,j)}-\ttau^{(j,i)}\mid 1\leq i<j\leq n)$
the \emph{symmetric  universal Frobenius algebra}.
 Note that taking the quotient $\Tosym$ is compatible with the $W$-action. Hence, $\kk\Sym\Smash\Tosym$ is defined and a quotient of the algebra  $\kk\Sym\Smash\Tos$.

By definition, the elements $\CL_k$, for $1\le k\leq n$ define \emph{Jucys--Murphy elements} $\CL_k$ in $\Tosym$.  By \autoref{P:CLProperties}\autoref{Lc}, any polynomial $f=f(X_1,\ldots, X_n)\in\kk[\X]$ provides $f(\CL_1,\dots,\CL_n)\in\kk\Sym\Smash\Tos$.

\begin{Proposition}\label{C:UniSymmetric}
    Let  $f\in\kk[\X]^W$ be a symmetric polynomial. Then
    $f(\CL_1,\dots,\CL_n)\in Z(\kk\Sym\Smash\Tosym)$.
\end{Proposition}

  \begin{proof}
     By \autoref{P:CLProperties}\autoref{Lb}, if $1\le k<n$ then
     $s_k\CL_l=\CL_ls_k$ if $l\ne k,k+1$ and
     \[
        s_k\bigl(\CL_k+\CL_{k+1}\bigr)
               =\bigl(\CL_{k+1}s_k-\ttau^{(k,k+1)}\bigr)
                   +\bigl(\CL_{k}s_k+\ttau^{(k+1,k)}\bigr)
               =\bigl(\CL_k+\CL_{k+1}\bigr)s_k.
     \]
     Here, the first equality holds in $\kk\Sym\Smash\Tos$, for the second we pass to  $\kk\Sym\Smash\Tosym$.
     Again by \autoref{P:CLProperties}, \[s_k\CL_k\CL_{k+1}=(\CL_{k+1}s_k-\ttau^{(k,k+1)})\CL_{k+1}
     =\CL_{k+1}\CL_{k}s_k+\CL_{k+1}\ttau^{(k+1,k)}-\ttau^{(k,k+1)}\CL_{k+1}=\CL_{k}\CL_{k+1}s_k.
     \]
For the last equality we used that the summands with $i=k$ cancel in $\Tos$, so that
\begin{align*}
\CL_{k+1}\ttau^{(k+1,k)}-\ttau^{(k,k+1)}\CL_{k+1}
    & =\sum_{1\leq i <k}\left((i,k+1)\ttau^{(k+1,i)}\ttau^{(k+1,k)}-\ttau^{(k,k+1)}(i,k+1)\ttau^{(k+1,i)}\right)\\
    &=\sum_{1\leq i <k} (i,k+1)\left(\ttau^{(k+1,i)}\ttau^{(k+1,k)}-\ttau^{(k,i)}\ttau^{(k+1,i)}\right)
\end{align*}
vanishes in $\Tos$ by \autoref{univtaurel}, and  that $\CL_{k+1}\CL_{k}=\CL_{k}\CL_{k+1}$ by  \autoref{P:CLProperties}\autoref{Lc}. It follows that
     if $w\in\Sym$, then $wf(\CL_1,\dots,\CL_n)=f(\CL_1,\dots,\CL_n)w$.

     On the other hand, let $1\le l<m\le n$, $1\le i<k\le n$ and  $w=(i,k)$.  Then $\ttau^{(l,m)}\CL_k=\CL_k\ttau^{(l,m)}$, since
   $\ttau^{(l,m)}(i,k)\ttau^{(k,i)}
               =(i,k)\ttau^{(w(l),w(m))}\ttau^{(k,i)}
               =(i,k)\ttau^{(k,i)}\ttau^{(w^2(l),w^2(m))}
               =(i,k)\ttau^{(k,i)}\ttau^{(l,m)}.$
  Here, the middle equality follows from \autoref{univtaurel} when
  $l,m\notin\set{k,i}$ or  $(l,m)=(i,m)$ or $(l,m)=(l,i)$, and in
  general from the definition of $\Tosym$.   Since we have shown the commutativity with the generators of $\kk\Sym\Smash\Tosym$, the symmetric polynomial $f(\CL_1,\dots,\CL_n)$ belongs
     to~$Z(\kk\Sym\Smash\Tosym)$.
  \end{proof}

We finish by connecting the universal Jucys--Murphy elements with
the corresponding classical Jucys--Murphy elements of the symmetric groups. These elements play
a central role in the representation theory of the symmetric groups,
from Young's seminormal forms~\cite{QSAI}, to the Okounkov--Vershik approach
to the representation theory of symmetric groups~\cite[\S3]{VObook}, and
the KLR algebras~\cite{BK:GradedKL}.

\begin{Lemma}
There is an algebra homomorphism $\Tos\rightarrow \Tosym \rightarrow \kk$ sending the $\ttau^{(i,j)}$ to $1$ (for  $1\leq i\not= j\leq n$). The induced algebra homomorphism  $\kk\Sym\Smash\Tos\rightarrow\kk\Sym$ maps $\CL_k$ to the classical Jucys--Murphy element $\sum_{1\leq i<k}(i,k)\in\kk\Sym$.
 \end{Lemma}

\begin{Remark} Many constructions of Jucys--Murphy elements in the
  literature, in particular of wreath product algebras such as
  \cite{Pushkarev, wreathHecke, WanWang, KleshchevMuth},   can be recovered by specialising the universal
  Jucys--Murphy elements from \autoref{D:universalJM}. For
  examples see \autoref{RkJM}, \autoref{HHdexamples} and \autoref{R:KMdictionary}.
\end{Remark}

A natural question is whether there is an analogue of the Okounkov--Vershik
theory~\cite[\S3]{VObook} for $\kk\Sym\Smash\Tos$, and for algebras
containing homomorphic images of the universal Jucys--Murphy elements $\CL_i$.
We do not consider this question here, although we do construct natural
examples of such algebras that simultaneously generalise, and specialise
to, affine Hecke algebras and their degenerations.

\subsection{Frobenius-tensor superalgebras} \label{superstuffeasy}
In this subsection we introduce a crucial player in our construction, the Frobenius-tensor superalgebras. These are tensor powers
$\Bo=A^{\otimes n}$ of a finite rank superalgebra $A$, with the symmetric group acting by
superpermutations of the tensor factors, equipped with a structure map $\Tos\to\Bo$ from the
universal Frobenius monoid algebra and with a Nakayama automorphism $\psi$ of $A$, subject to
compatibility conditions. The terminology is justified by \autoref{T:FrobTensor}, where we show
that tensor powers of Frobenius superalgebras, with their copairing elements $\tau^{(i,j)}$ and their Nakayama
automorphism, are Frobenius-tensor superalgebras.

By a \emph{superalgebra} $A$ we mean a $\Ztwo$-graded $\kk$-algebra $A=A_{\ov0}\oplus A_{\ov1}$ (such that its homogeneous parts are free $\kk$-modules). More conceptually, one might define it  as an algebra object in the symmetric monoidal category $\op{Vect}_{\Ztwo}$ of $\Ztwo$-graded $\kk$-modules with grading preserving morphisms and braiding morphisms $a\otimes b\mapsto (-1)^{|a||b|}b\otimes a$.
In particular,  if $A$ and $A'$ are superalgebras, then $A\otimes A'$ is a superalgebra with multiplication
 \begin{equation}\label{multinsuper}
 (a\otimes a') (b\otimes b')=(-1)^{|a'||b|}ab\otimes a'b'.
  \end{equation}
  Here and below, $|a|\in\Ztwo=\set{\ov0,\ov1}$ is the \emph{degree} or
  \emph{parity} of a homogeneous element $a\in A$. We often tacitly
  assume that $a\in A$ is homogeneous whenever we write formulas involving $|a|$. Elements of degree $\ov{0}$ are \emph{even}, elements of degree  $\ov{1}$ are \emph{odd}.  Any (ordinary) $\kk$-algebra is a superalgebra concentrated in degree $\ov0$.
Given a superalgebra $A$, an \emph{$A$-module} (or, more precisely, an $A$-supermodule) is a $\Ztwo$-graded $\kk$-module $M$ with an action of the (underlying ordinary) algebra $A$ such that $A_{i}M_{j}\subset M_{i+j}$ for $i,j\in\Z/2\Z$.  By definition, a morphism $M\rightarrow N$ between $A$-modules is a degree preserving linear map $f\colon M\rightarrow N$ satisfying $f(am)=a f(m)$ for all $a\in A, m\in M$.
Let $\Pi$ be the parity shift, so that $\Pi M=\kk^{0|1}\otimes M$. Explicitly $(\Pi M)_i=M_{i+\bar 1}$, with
$a\cdot_\Pi m=(-1)^{|a|}am$. Let $A^{\mathrm{op}}$ be the opposite superalgebra with multiplication $a\cdot_{A^{\mathrm{op}}}b=(-1)^{|a||b|}b\cdot_Aa$.

     \begin{Definition}\label{Symaction}
For a superalgebra $A$ let $\Ao= A^{\otimes n}$. This is a superalgebra via \autoref{multinsuper} with a canonical $\Sym$-action by superalgebra automorphisms induced by the symmetric monoidal structure.
Explicitly, if $\ba=a_1\otimes\dots\otimes a_n\in\Ao$, then $s_i$ acts by the \emph{superflip} $\leftindex^{s_i}\ba:=(-1)^{|a_i||a_{i+1}|}a_1\otimes\dots\otimes a_{i-1}\otimes a_{i+1}\otimes a_i\otimes a_{i+2}\otimes\dots\otimes a_n.$ Thus $\Sym$ acts in $\Ao$ by \emph{superpermutations} of the tensor factors.

In case $A$ is concentrated in degree $\ov0$,  the action is
by ordinary place permutations, that is $w\in\Sym$ acts on $\ba$ as  $\leftindex^w\ba=a_{w^{-1}(1)}\otimes\dots\otimes
  a_{w^{-1}(n)}$.
  \end{Definition}

We consider $\kk$ as a superalgebra concentrated in degree $\ov{0}$. Then $\kk^{\otimes n}\cong\kk$ with the trivial $\Sym$-action.
\begin{Definition}\label{D:ST-algebra}
A \emph{braided superspace} is an object $V\in\op{Vect}_{\Ztwo}$,  together with a degree preserving linear morphism $\tau\colon V\otimes V\rightarrow V\otimes V$ such that $(\tau\otimes \op{id})(\op{id}\otimes \tau)(\tau\otimes \op{id})=(\op{id}\otimes \tau)(\tau\otimes \op{id})(\op{id}\otimes \tau)$.

  A $\Tos$-\emph{superalgebra} is a superalgebra  $\Bo$ with a
  \emph{structure map} $\Tos\rightarrow\Bo$, which is a superalgebra
  homomorphism. Abusing notation, write $\tau^{(i,j)}$ for the image of $\ttau^{(i,j)}$ in $\Bo$.

 A \emph{tensor superalgebra}   is a $\Tos$-superalgebra of the form  $\Bo=A^{\otimes n}$ with a finite rank  superalgebra $A$  such that  $\tau^{(1,2)}\in A\otimes A\otimes 1^{\otimes {n-2}}$ and $\leftindex^w(\tau^{(1,2)} \ba)=(\tau^{(w(1),w(2))}) (\leftindex^w\ba)$,  for $w\in\Sym$, $\ba\in\Bo$.
  Here, $\Sym$ acts on $\Bo$ by superpermuting the tensor factors, see \autoref{Symaction}. In particular, the structure map $\Tos\rightarrow\Bo$ is $W$-equivariant and $\tau:=\tau^{(1,2)}\cdot\colon A\otimes A\rightarrow A\otimes A$  turns $A$ into a braided  superspace.
  \end{Definition}
  \begin{Definition}\label{D:Nakayama}
    Given a tensor superalgebra $\Bo=A^{\otimes n}$, a
   \emph{Nakayama automorphism $\psi$ of $\Bo$} is an automorphism
   $\psi\colon A\rightarrow A$ of superalgebras such that the following hold, as linear endomorphisms of $\Bo$,
\begin{equation}\label{Nakayama}
w\psi^{(i)}w^{-1}=\psi^{(w(i))},\quad \psi^{(j)}(\tau^{(i,j)})=\tau^{(j,i)},\quad  \psi^{(k)}(\tau^{(i,j)})=\tau^{(i,j)}, \quad  \text{and}\quad \ba\tau^{(i,j)}=\tau^{(i,j)}\psi^{(i)}\left(\leftindex^{(i,j)}\ba\right),
\end{equation}
 for $w\in\Sym$ and distinct $i,j,k$. Here  $\psi^{(i)}=\op{id}^{\otimes{(i-1)}}\otimes\psi\otimes \op{id}^{\otimes{(n-i)}}$. We call $\psi$ \emph{trivial} if it is the identity map.
 \end{Definition}
 We finally come to the main definition of this subsection.
 \begin{Definition}\label{D:Frobtensor}
  A  \emph{Frobenius-tensor superalgebra}  is a  tensor
 superalgebra $\Bo$ with a Nakayama automorphism~$\psi$.
 If $A$ is an ordinary algebra, we call  $(\Bo,\psi)$ a \emph{Frobenius-tensor algebra}.
 Let $\eta\in\kk$. Then a Frobenius-tensor superalgebra
$\Bo=A^{\otimes n}$ is \emph{$\eta$-symmetric} if $\eta\psi(a)=\eta a$
for any  $a\in A$.
\end{Definition}

\begin{Remark}\label{rketasymm}
It is important to observe that the freeness of $A$ over the domain $\kk$ implies that $\Bo$ is
$\eta$-symmetric if and only if $\eta=0$, or $\eta\ne0$ and $\psi=\op{id}$.
\end{Remark}

\begin{Remark}
Note that the second condition in
\autoref{Nakayama} implies that $\taui[i,j]=\taui[j,i]$ for all $i,j$
whenever $\psi$ is trivial. In this case, $\psi^{(j)}$ is the identity map for
all $j$, and  the last condition in \autoref{Nakayama} becomes for
$\ba=\taui[i,j]$ the empty requirement that
$\taui[i,j]\taui[i,j]=\taui[i,j]\taui[i,j]$.
\end{Remark}

\begin{Example}\label{trivFrobtensor}
The  \emph{trivial Frobenius-tensor superalgebra} is the superalgebra $\Bo=\kk^{\otimes n}=\kk$ with the structure map that sends every generator to
 $1$. It is $\eta$-symmetric for any $\eta$.
\end{Example}

\begin{Lemma}\label{afterGluehwein}
In a Frobenius-tensor superalgebra the following hold
\[\psi^{(i)}\psi^{(j)}\tau^{(i,j)}=\tau^{(i,j)}=\psi^{(j)}\psi^{(i)}\tau^{(i,j)}\quad \text{ for $1\leq i\not=j\leq n$.}\]
\end{Lemma}
\begin{proof}By definition, as endomorphisms
$\psi^{(i)}\psi^{(j)}=\psi^{(j)}\psi^{(i)}$ for $1\leq i,j\leq n$.
Now the second property of the Nakayama automorphism in
\autoref{Nakayama} implies that $\psi^{(i)}\psi^{(j)}\tau^{(i,j)}= \psi^{(i)}\tau^{(j,i)}=\tau^{(i,j)}$ for $i\not=j$.
\end{proof}

  \subsection{Frobenius-tensor superalgebras from Frobenius superalgebras}\label{S:Frobeniusalgs}
We construct Frobenius-tensor superalgebras as tensor powers of
Frobenius superalgebras, that is, of Frobenius algebra objects in $\op{Vect}_{\Ztwo}$. For more details on Frobenius algebras we
refer to \cite{Kock, Skowronski} and to \cite{Savage:AffineWreath}, which is closest to our setting.

Assume in this subsection that $A$ is a \emph{Frobenius superalgebra}. That is, $A$ is a superalgebra,
free of finite rank over $\kk$, equipped with a degree preserving linear map $\tr\colon A\to\kk$,
the \emph{trace map}, such that the bilinear form $(a,b)\mapsto\tr(ab)$ is non-degenerate, inducing an isomorphism $A\cong A^*=\Hom_\kk(A,\kk)$. Then
there is a unique superalgebra automorphism $\psi$ of $A$, the \emph{Nakayama automorphism}, with
$\tr(ab)=(-1)^{|a||b|}\tr(b\psi(a))$ for all homogeneous $a,b\in A$. The algebra $A$ is
\emph{symmetric} if $\psi=\op{id}$.

Equivalently, a (symmetric) Frobenius superalgebra is a  superalgebra $A$ with a non-degenerate degree preserving (supersymmetric) bilinear form $\<\ ,\ \>\map{A\times A}\kk$
  that is associative in the sense that $\<ab,c\> =
  \<a,bc\>$, for all $a,b,c\in A$. The associated \emph{trace form} is
  the linear map $\tr\map A\kk$ given by $\tr(a)=\<a,1\>$.

  We fix a basis $\cB$ for $A$ consisting of homogeneous elements. Let $r=|\cB|$.  It induces a decomposition $\cB=\cB_{\ov{0}}\cup\cB_{\ov{1}}$ with respect to degrees.
 Let $\set{b^\vee|b\in \cB}$ be the dual basis of $\cB$, which is determined by
  $\tr(a^\vee b)=\delta_{ab}$, for $a,b\in \cB$.
  We have $(a^\vee)^\vee=(-1)^{|a|}\psi^{-1}(a)$ for $a\in\cB$. Thus,
  $(a^\vee)^\vee=(-1)^{|a|}a$ when $A$ is symmetric, and
  \begin{equation}\label{traceformulas}
    f = \sum_{b\in \cB}b\tr(b^\vee f) = \sum_{b\in \cB}\tr(fb)b^\vee\quad\text{and}\quad \psi((a^\vee)^\vee)=(-1)^{|a|}a
    \qquad\text{ for $f\in A$ and $a\in\cB$.}
  \end{equation}
  We will use these identities without mention. For  $a\in A$ and $1\leq i\leq n$ set
 $a^{(i)}=1^{\otimes(i-1)}\otimes a\otimes1^{\otimes(n-i)}\in\Ao$.
  \begin{Remark}\label{signs}
 For homogeneous $a,b\in A$ and $1\leq i\not=j\leq n$ we have  $a^{(i)}b^{(j)}=(-1)^{|a||b|}b^{(j)}a^{(i)}$ in $\Ao$.
 \end{Remark}
   \begin{Definition}\label{parts}
For $1\le i\not=j\le n$ define the \emph{$(i,j)$-copairing element} in $\Ao$ as $\tau^{(i,j)} = \sum_{b\in\cB}b^{(i)}{b^\vee}^{(j)}$. It has a  \emph{0-part} $\tau^{(i,j)}_{\ov{0}}=\sum_{b\in\cB_{\ov0}}b^{(i)}{b^\vee}^{(j)}$ and a \emph{1-part} $\tau^{(i,j)}_{\ov{1}}=\sum_{b\in\cB_{\ov{1}}}b^{(i)}{b^\vee}^{(j)}$.
   \end{Definition}

\begin{Remark} The copairing elements for Frobenius algebras appear already prominently in \cite{RossoSavage:QuantumAffineWreath,Savage:AffineWreath}.
\end{Remark}

\begin{Remark}The copairing elements for Frobenius algebras appear also in a diagrammatical context as so-called \emph{teleporters} for instance in \cite{SavageHeis, BSWHeis}.
\end{Remark}

   \begin{Remark}The notation $\tau^{(i,j)}_{\ov0}$ and $\tau^{(i,j)}_{\ov1}$  might be confusing: The index refers to the degree of the basis elements used in the sum and \emph{not} to the homogeneous components of $\tau^{(i,j)}$. In fact, $\tau^{(i,j)}$ is even.
   \end{Remark}
\begin{Remark}\label{hiddensigns}
It is not hard to see
  that $\tau^{(i,j)}$ depends only on $i$ and $j$ and not on the choice
  of $\cB$.

  In explicit calculations in the super case, one should be aware of hidden
  signs. For example, by \autoref{multinsuper}, we have
$\tau^{(1,2)}_{\ov{1}}=\sum_{b\in\cB_{\ov{1}}}b \otimes b^{\vee}$, whereas $\tau^{(2,1)}_{\ov{1}}=-\sum_{b\in\cB_{\ov{1}}}b^{\vee} \otimes b$.
 \end{Remark}

  \begin{Proposition}\label{Taussuper}
   Let $A$ be a Frobenius superalgebra and let  $1\leq i,j,k,l\leq n$ be distinct and $\epso, \epsoo\in\Ztwo.$ Then the copairing elements in the superalgebra $\Ao$  satisfy the following \emph{quadratic relations}
 \begin{equation}\label{tausinAsuper}
   \tau_{\epso}^{(i,j)}\tau^{(k,l)}_{\epsoo}=\tau^{(k,l)}_{\epsoo}\tau^{(i,j)}_{\epso}
   \quad\text{and}\quad
    \tau_{\epso}^{(i,j)}\tau_{\epsoo}^{(j,k)}=\tau^{(j,k)}_{\epso+\epsoo}\tau^{(i,k)}_{\epso}=\tau^{(i,k)}_{\epsoo}\tau^{(i,j)}_{\epso+\epsoo}.
 \end{equation}
     If $A$ is symmetric, then we have additionally the following equalities
 \begin{equation}\label{tausym}
  \tau^{(i,j)}_\epso=\tau^{(j,i)}_\epso \quad\text{and}\quad  \tau^{(j,k)}_{\epso}\tau^{(i,j)}_{\epsoo}=\tau^{(i,j)}_{\epso+\epsoo}\tau^{(i,k)}_{\epso}=\tau^{(i,k)}_{\epsoo}\tau^{(j,k)}_{\epso+\epsoo}.
 \end{equation}
 \end{Proposition}
  \begin{proof}
   By definition,
   $\tau^{(i,j)}_{\epso}\tau^{(k,l)}_{\epsoo}=\sum_{a}a^{(i)}(a^\vee)^{(j)}
   \sum_{b} b^{(k)}  (b^\vee)^{(l)}$, where the sum runs over $a,b\in\cB$
   such that $|a|=\epso$, $|b|=\epsoo$. By \autoref{signs}, this sum equals  $\sum_{b} b^{(k)}  (b^\vee)^{(l)}\sum_{a}a^{(i)}(a^\vee)^{(j)}=\tau^{(k,l)}_{\epsoo}\tau^{(i,j)}_{\epso}$. We next verify the second
    equality of the proposition. The third is similar and omitted.  We calculate
    \begin{align*}
      \tau^{(i,j)}_{\epso}\tau^{(j,k)}_{\epsoo}
         & =\sum_{a}a^{(i)}(a^\vee)^{(j)}  \sum_{b} b^{(j)}  (b^\vee)^{(k)}
         =\sum_{a,b}a^{(i)}(a^\vee)^{(j)} b^{(j)}  (b^\vee)^{(k)}
          =\sum_{a,b,c}a^{(i)}\tr(c^\vee a^\vee b) c^{(j)}  (b^\vee)^{(k)}  \\
             & =\sum_{a,b,c}a^{(i)} c^{(j)} \tr(c^\vee a^\vee b) (b^\vee)^{(k)}
           =\sum_{a, c}a^{(i)} c^{(j)} (c^\vee a^\vee)^{(k)}
              = \sum_{c}c^{(j)} (c^\vee)^{(k)} \sum_{a}a^{(i)} (a^\vee)^{(k)}
      =\tau^{(j,k)}_{\epso+\epsoo}\tau^{(i,k)}_{\epso}.
              \end{align*}
            Here, the sums run over $a,b, c\in \cB$ such that
            $|a|=\epso$, $|b|=\epsoo$ and $|c|=|a|+|b|$, and the
            penultimate equality holds because swapping the factors produces the sign $(-1)^{|a||c|}(-1)^{|a||c|}=1$ by \autoref{signs}.

            Assume now $A$ is symmetric. Consider
            $\tau^{(i,j)}_{\epso}=\sum_{a}a^{(i)}(a^\vee)^{(j)}$, with the sum running over $a\in\cB$ with $|a|=\epso$.
            By replacing $a$ by $a^\vee$ and using \autoref{traceformulas} in case $\psi$ is the identity the sum rewrites as $\sum_{a}(-1)^{|a|}(a^\vee)^{(i)}a^{(j)}=\tau^{(j,i)}_{\epso}$. This shows the first formula in \autoref{tausym}, the second follows from  \autoref{tausinAsuper} applied to the triple $(k,j,i)$.
 \end{proof}
 \begin{Corollary}\label{CommutingTaus}
   Let $A$ be a Frobenius superalgebra and let  $1\leq i,j,k,l\leq n$ be distinct. Then the copairing elements satisfy the following \emph{quadratic relations} in the algebra $\Ao:$
 \begin{equation*}
   \tau^{(i,j)}\tau^{(k,l)}=\tau^{(k,l)}\tau^{(i,j)}
   \quad\text{and}\quad
    \tau^{(i,j)}\tau^{(j,k)}=\tau^{(j,k)}\tau^{(i,k)}=\tau^{(i,k)}\tau^{(i,j)}.
 \end{equation*}
     If $A$ is symmetric, then we have additionally $\tau^{(i,j)}=\tau^{(j,i)}$ and $\tau^{(j,k)}\tau^{(i,j)}=\tau^{(i,k)}\tau^{(j,k)}=\tau^{(j,i)}\tau^{(i,k)}$.
 \end{Corollary}
 \begin{proof}
 The claim follows directly from
 \autoref{Taussuper} using the definition $\tau^{(i,j)}=\tau^{(i,j)}_{\ov0}+\tau^{(i,j)}_{\ov1}$.
 \end{proof}
  \begin{Remark}
  In general,  $\tau^{(i,j)}\tau^{(j,k)}\ne\tau^{(j,k)}\tau^{(i,j)}$ even if $A$ is symmetric, but not commutative.  \end{Remark}

 \begin{Remark}\label{generalizedbraidinA}
 As in the proof of \autoref{braidlemma} it follows from \autoref{Taussuper} that for a Frobenius superalgebra $A$, the copairing elements in $\Ao$ satisfy the following generalised braid relations:
 \begin{align*}
\tau^{(i,j)}\tau^{(j,k)}\tau^{(i,j)}=\tau^{(j,k)}\tau^{(i,j)}\tau^{(j,k)}
 \quad\text{and}\quad
 \tau^{(i,k)}\tau^{(i,j)}\tau^{(i,k)}=\tau^{(i,j)}\tau^{(i,k)}\tau^{(i,j)}\quad \text{for pairwise distinct   $1\leq i,j,k\leq n$}.
  \end{align*}
  If the Frobenius algebra $A$ is symmetric,  the first family of relations is equivalent to the second  and also equivalent to the ordinary braid relations $\tau^{(i,i+1)}\tau^{(i+1,i+2)}\tau^{(i,i+1)}=\tau^{(i+1,i+2)}\tau^{(i,i+1)}\tau^{(i+1,i+2)}$  for $1\leq i<n-1$.
\end{Remark}

Recall from \autoref{Symaction} the $\Sym$-action on  $\Ao$ for any superalgebra $A$ by superpermutations of the tensor factors. In particular, $\leftindex^w\tau^{(i,j)}_\epsilon=\tau^{(w(i),w(j))}_\epsilon$ for $w\in\Sym$.
 \begin{Lemma} \label{psitauA} The following equalities hold in $\Ao$ for any Frobenius superalgebra $A$
    \begin{equation*}
 \psi^{(j)}(\tau^{(i,j)})=\tau^{(j,i)}\quad\text{and}\quad \psi^{(k)}(\tau^{(i,j)})=\tau^{(i,j)}, \quad\quad\text{for any $1\leq i\not=j\leq n$, $k\not=i,j$.}
   \end{equation*}
  \end{Lemma}
  \begin{proof}
 The second formula follows from the definitions. For the first we calculate, using \autoref{traceformulas},
 \begin{equation*}
 \psi^{(j)}(\taui[i,j])=\sum_{b\in\cB}b^{(i)}(\psi(b^\vee))^{(j)}=\sum_{b\in\cB}(-1)^{|b|}(\psi(b^\vee))^{(j)}b^{(i)}=
 \sum_{c\in\cB}c^{(j)}(c^\vee)^{(i)}=\tau^{(j,i)}. \qedhere
 \end{equation*}
  \end{proof}
  In case $A$ is symmetric, we obtain $\tau^{(i,i+1)}=\leftindex^{s_i}\tau^{(i,i+1)}=\tau^{(i+1,i)}$.
  We abbreviate this element as $\tau^{(i)}$.

  \begin{Lemma} \label{athroughtau} Let $A$ be a Frobenius superalgebra, $\ba\in \Ao$ and $1\leq i\not=j\leq n$. Then \[\ba\tau^{(i,j)}=\tau^{(i,j)}\psi^{(i)}\left(\leftindex^{(i,j)}\ba\right).\]
  \end{Lemma}
  \begin{proof}
 It suffices to show that for $a\in A$ and $1\leq i,j,k\leq n$ distinct it holds that
   \begin{equation}\label{stronger}
   a^{(i)}\tau^{(i,j)}=\tau^{(i,j)}a^{(j)},\quad
 a^{(j)}\tau^{(i,j)}=\tau^{(i,j)}\psi(a)^{(i)} \quad  \text{ and }\quad a^{(k)}\tau^{(i,j)}=\tau^{(i,j)} a^{(k)}.
 \end{equation}
The last formula is obvious. The first formula can be deduced by applying \autoref{traceformulas}:
\begin{small}
  \begin{equation*}
  a^{(i)}\tau^{(i,j)}=\sum_{b\in\cB}a^{(i)}b^{(i)} (b^\vee)^{(j)}
  =\sum_{b,c\in\cB}\tr(c^\vee ab)c^{(i)}(b^\vee)^{(j)}= \sum_{b,c\in\cB}c^{(i)}\tr(c^\vee ab)(b^\vee)^{(j)}
  = \sum_{c\in\cB}c^{(i)}(c^\vee)^{(j)}a^{(j)}=\tau^{(i,j)}a^{(j)}.
  \end{equation*}
\end{small}
The second formula in \autoref{stronger} is deduced analogously, but requires incorporating the Nakayama automorphism:
\begin{small}
\begin{align*}
a^{(j)}\tau^{(i,j)}
&=\sum_{b\in\cB}(-1)^{|a||b|}b^{(i)}a^{(j)} (b^\vee)^{(j)}
=\sum_{b,c\in\cB}(-1)^{|a||b|}b^{(i)} \tr(ab^\vee c)(c^\vee)^{(j)}=\sum_{b,c\in\cB}(-1)^{|a||b|}(-1)^{|a|}\tr(b^\vee c\psi(a))b^{(i)} (c^\vee)^{(j)}\\
&=\sum_{c\in\cB}(-1)^{|a|(|c|-|a|)}(-1)^{|a|} c^{(i)}(\psi(a))^{(i)} (c^\vee)^{(j)}
=\sum_{c\in\cB}(-1)^{|a|(|a|+|c|+1)}(-1)^{|a||c|} c^{(i)} (c^\vee)^{(j)}(\psi(a))^{(i)}
=\tau^{(i,j)}\psi(a)^{(i)}.
\end{align*}
\end{small}
The last two equalities use $|c|=|a|+|b|$ and
$(-1)^{|a|(|a|+|c|+1)}(-1)^{|a||c|}=1$, respectively, for any $|a|$, $|c|$.
\end{proof}
The following provides the main examples of Frobenius-tensor superalgebras.
\begin{Theorem}\label{T:FrobTensor}Let $A$ be a Frobenius superalgebra.
  \begin{enumerate}
 \item\label{FrobTensora} The superalgebra $\Ao$ with the action of $\Sym$ by superpermutations, the structure
   map $\Tos\rightarrow\Ao$, $\ttau^{(i,j)}\mapsto\tau^{(i,j)}$, and the Nakayama automorphism
   $\psi$ of $A$ is a Frobenius-tensor superalgebra.
 \item\label{FrobTensorb} If $A$ is concentrated in degree $\ov{0}$, then $\Ao$ is a
   Frobenius-tensor algebra.
 \item\label{FrobTensorc} $\Ao$ is $\eta$-symmetric if and only if $\eta=0$, or $\eta\ne0$ and
   $A$ is symmetric.
   \end{enumerate}
    In particular, the two notions of Nakayama automorphism are compatible.
 \end{Theorem}

 \begin{proof}
 \autoref{FrobTensora}: The copairing elements satisfy the relations \autoref{univtaureleq} by
 \autoref{CommutingTaus}, so $\ttau^{(i,j)}\mapsto\tau^{(i,j)}$ defines a structure map
 $\Tos\rightarrow\Ao$. By definition, $\tau^{(1,2)}\in A\otimes A\otimes1^{\otimes(n-2)}$, and
 $\leftindex^w(\tau^{(1,2)}\ba)=\tau^{(w(1),w(2))}\leftindex^w\ba$ holds because the $\Sym$-action
is by algebra automorphisms and $\leftindex^w(a^{(i)})=a^{(w(i))}$ for $a\in A$, so that $\leftindex^w\tau^{(i,j)}=\tau^{(w(i),w(j))}$. Hence $\Ao$ is a  tensor superalgebra. It remains to check the conditions \autoref{Nakayama} for $\psi$. The first is clear, the second and third hold by \autoref{psitauA}, and the last is satisfied by \autoref{athroughtau}.

 \autoref{FrobTensorb} is a special case of \autoref{FrobTensora}.

 \autoref{FrobTensorc}: Since $A$ is free over the domain $\kk$, the condition $\eta a=\eta\psi(a)$
 for all $a\in A$ is equivalent to $\eta=0$, or $\eta\ne0$ and $\psi=\op{id}$.
  \end{proof}

\subsection{Examples of Frobenius (super)algebras}\label{Frobeniusalgex}
Since Frobenius superalgebras naturally give rise to Frobenius-tensor superalgebras, we now give
some examples of Frobenius (super)algebras $A$.

We write $\tau=\tau^{(1,2)}\in A\otimes A$ for the copairing element in case $n=2$.

  \begin{Example}
  Any symmetric (super)algebra $A$, i.e. a superalgebra with $A\cong A^*$ as $A$-bimodule, is a Frobenius (super)algebra. In this
  case,  the copairing element
  $\tau$ was considered in \cite{KleshchevMuth} and called \emph{distinguished element}.
  \end{Example}

\begin{Example}\label{tensorFrob}Given two Frobenius superalgebras $A$, $B$ their tensor product $A\otimes B$ is a Frobenius superalgebra  with $\op{tr}_{A\otimes B}=\op{tr}_A\otimes \op{tr}_B$  and $\psi_{A\otimes B}=\psi_A\otimes \psi_B$. (It is a good exercise to check the signs.)
\end{Example}

 \begin{Example}\label{Ex:FrobeniusAlgebrasfingroup}
   Let $G$ be a finite group. Then the group algebra $\kk G$ is a
   symmetric Frobenius algebra $A$ with trace form $\tr\colon\kk G\to\kk$
   given by $\tr(g)=\delta_{g,1_G}$, for $g\in G$. Hence,
   $\set{g|g\in G}$ and $\set{g^{-1}|g\in G}$ are dual bases of $\kk G$
   and $\tau=\sum_{g\in G}g\otimes g^{-1}$.
   Note that $\tau^2\in Z(A\otimes A)=Z(A)\otimes Z(A)$ by \autoref{athroughtau}. This element has a nice description in terms of the conjugacy classes, $C_1,\dots,C_z$ of $G$, namely
   \[
       \tau^2 =
            \sum_{i=1}^z\frac{|G|}{|C_i|}\sum_{g,h\in C_i}g\otimes h^{-1}.
   \]
   Indeed, $\tau^2=\sum_{g,h\in G}gh\otimes g^{-1}h^{-1}=\sum_{y,g\in G}y\otimes (g^{-1}yg)^{-1}$ and the stabiliser  of each $y\in C_i$ under the conjugation action of $G$ has cardinality $\frac{|G|}{|C_i|}$. In particular, $\tau^2 =|G|\tau$ if $G$
   is abelian. Taking~$G$ to be the trivial group gives the trivial Frobenius algebra structure on $\kk$ from \autoref{trivFrobtensor}.
    \end{Example}
          \begin{Example}\label{Ex:FrobeniusAlgebrassmash}
  \autoref{Ex:FrobeniusAlgebrasfingroup} can be generalized by assuming additionally that $G$ acts on a Frobenius algebra $A$ with trace map $\tr_A$ by algebra automorphisms. Then $\kk G\Smash A$ is a Frobenius algebra with trace map $\tr(ga)=\tr_{\kk G}(g)\tr_A(a)$ if $G$ preserves $\tr_A$. If $A$ is symmetric then so is  the smash product, because
 \[
  \tr_A\bigl((^{g}a)a'\bigr)=\tr_A\bigl(^{g^{-1}}\!\bigl((^{g}a)a'\bigr)\bigr)=\tr_A\bigl(a\,^{g^{-1}}\!a'\bigr).
\]
The group superalgebras of finite groups $G$ with index $2$ subgroups $G_0$ as in \cite{Jozefiak} are examples in the same spirit.
 \end{Example}
 \begin{Remark}\label{Pushkarev}
 Using the Frobenius-tensor algebra associated with  $A=\kk G$, the universal Jucys--Murphy elements from  \autoref{D:universalJM} can be used to recover classical variants.
 In \cite{Pushkarev} Pushkarev introduced, most likely for the first time,  Jucys--Murphy elements for wreath products $\kk\Sym\Smash\kk G$ for finite groups $G$ and $\kk=\mathbb{C}$. His elements  are the push-forwards of our Jucys--Murphy elements  along the structure map $\Tos\rightarrow(\kk G)^{\otimes n}$ for the Frobenius algebra $A=\kk G$. Pushkarev's \emph{one-step algebra} is our $\kk\Sym\Smash\Ao$ for $n=2$. Several subsequent treatments of Jucys--Murphy elements in wreath product cases, e.g.\ \cite{wreathHecke, WanWang, KleshchevMuth}, fall into our universal setup as well. \autoref{evaluation} can be seen as a more conceptual common framework of these constructions.
\end{Remark}

 \begin{Example}\label{Ex:FrobeniusAlgebrasHPm}
   Fix $m\ge0$ and let $A=\kk[x]/(x^{m+1})$, which we consider as a graded
  algebra with $x$ in (cohomological) $\Z$-degree~$2$. Then $A$ is a Frobenius algebra with trace form $\tr(x^k) = \delta_{k, m}$.
   Hence,
   $\tau=\sum_{k=0}^{m}x^k\otimes x^{m-k}$. In this case, it is easy
   to check that $\tau$ is nilpotent if $m\geq 1$. (Observe that if $\operatorname{char}\kk=p$ and $m+1$ is a power of $p$, then $A$ is isomorphic to
   the group algebra of the cyclic group of order~$m+1$, equipped with a
   different Frobenius structure.) This algebra $A$  arises naturally as the cohomology algebra of the
   projective space $\mathbb{P}^{m}$ with its Poincaré trace form. In fact, the cohomology algebra of a smooth projective variety~$X$ is a symmetric Frobenius algebra
   with its non-degenerate Poincaré trace form.
    \end{Example}

\begin{Example}\label{Ex:FrobeniusAlgebrastrivialext}Let $B$ be a finite dimensional algebra with an algebra automorphism $\sigma$. Then $B\oplus B^*$ becomes an algebra via $(b_1,f_1)(b_2,f_2)=(b_1b_2,f_1(b_2\_)+f_2(\_\sigma(b_1)))$. One can check that it is  Frobenius with trace form $\tr((b,f))=f(1)$. If $\sigma$ is trivial, then the algebra is symmetric. The construction is a special case of a  \emph{trivial extension algebra}.
\end{Example}

  \begin{Example}\label{Alambda}Let $\lambda\in\kk$ be invertible. The algebra $A=\kk\langle y_1,y_2\rangle/(y_1y_2-\lambda y_2y_1, y_1^2,y_2^2)$  has basis $\cB=\{1,y_1,y_2,y_1y_2\} $ and a
  trace form given on $\cB$ by $\tr(b)=\delta_{b,y_1y_2}$. This turns
  $A$ into a Frobenius algebra with
  $\psi(y_1)=\lambda y_1$, $\psi(y_2)=\lambda^{-1} y_2$. Thus $A$ is
  symmetric if and only if $\lambda=1$.
\end{Example}
 \begin{Example}\label{Ex:Taft}
    Fix an integer $\theta\ge2$ and assume that $\kk$ contains a
     primitive $\theta$th root of unity $\omega\in\kk$. The \emph{Taft
     algebra} $\Taft$ is the $\kk$-algebra
     $\Taft=\kk\<\,g,y\>/( g^\theta=1, y^\theta=0, yg=\omega gy)$. This
     algebra is a Hopf algebra that is neither commutative nor
     cocommutative. The algebra $\Taft$ has basis $\cB=\set{g^iy^j| 0\leq i,j<\theta}$ and is
     Frobenius with trace form $\op{tr}(b)=\delta_{b,y^{\theta-1}}$ for $b\in \cB$ and $\psi(g)=\omega g$, $\psi(y)=y$.
     \end{Example}
    \begin{Example}\label{Ex:Clifford}
    Assume $2\in\kk^\times$ and let $A=\Cl_1$ be the rank $2$ Clifford superalgebra with basis $1$ (of degree $\ov0$) and $C$ of degree $\ov 1$ and multiplication given by $C^2=1$. It is Frobenius with trace map $\tr(1)=1$, $\tr(C)=0$, and nontrivial Nakayama automorphism given by $\psi(C)=-C$. Then $\Ao$ is isomorphic to the Clifford superalgebra $\Cl_n$ with odd generators $C_i$, $1\le i\le n$, modulo the relations $C_i^2=1$, $C_iC_j=-C_jC_i$. The $\Sym$-action is given by $s_i(C_j)=C_{s_i(j)}$. We have $\tau^{(i,j)}=1+ C_iC_j$ with parts $\tau^{(i,j)}_{\ov{0}}=1$ and $\tau^{(i,j)}_{\ov{1}}=C_iC_j$.
As an example, the equation  $\tau_{\epso}^{(2,1)}\tau_{\epsoo}^{(1,3)}=\tau^{(1,3)}_{\epso+\epsoo}\tau^{(2,3)}_{\epso}$ from \autoref{tausinAsuper} encodes
\[1\cdot1=1\cdot1, \quad (C_2C_1)\cdot1=(C_1C_3)\cdot(C_2C_3),  \quad 1\cdot (C_1C_3)=(C_1C_3)\cdot 1,  \quad (C_2C_1)\cdot(C_1C_3)=1\cdot (C_2C_3).\]
One easily verifies that these equalities, corresponding to the four choices of $(\epso,\epsoo)$,  hold in $\Cl_n$.
   \end{Example}

   \begin{Remark}\label{RkJM}
As in \autoref{Pushkarev}, the universal Jucys--Murphy elements of
\autoref{D:universalJM} specialise to Jucys--Murphy elements from the literature. For the Clifford
superalgebra $A$ of \autoref{Ex:Clifford} they become the Jucys--Murphy elements of the Sergeev
superalgebra, see for example \cite[(13.6)]{Kleshchevbook}. Since we allow in
\autoref{D:universalJM} only the action by superpermutations and not by a more general Hecke
algebra action, we do not recover in this way the Jucys--Murphy elements of the Hecke--Clifford
superalgebras, see also \autoref{flipbraiding}. We encourage the reader to work out such a generalisation.
\end{Remark}

\begin{Example}\label{exsemisimple}
Let $\kk=\mathbb{C}$ and let $A$ be a simple finite dimensional
superalgebra. By Wedderburn's Theorem, \cite[Theorem
12.2.9]{Kleshchevbook},  $A=M(\tm|\tn)$ or $A=Q_\tn$ for some $\tm,\tn$. Here,
$A=M(\tm|\tn)$ is the superalgebra of linear endomorphisms of a
superspace of superdimension $(\tm|\tn)$, say  $\kk^{\tm|\tn}$. It
can be identified with the superalgebra of all $(\tm+\tn)\times
(\tm+\tn)$-matrices with basis given by matrix units $E_{i,j}$, where
$E_{i,j}$ is even if $1\leq i,j\leq \tm$ or $\tm< i,j\leq \tm+\tn$ and
odd otherwise.  On the other hand, $Q_\tn$ is the algebra of
endomorphisms $f$ of $\kk^{\tn|\tn}$,  commuting with an odd involution $J$. Explicitly these are the $2\tn\times 2\tn$-matrices of the form $\set[\big]{\begin{psmallmatrix} A&B\\B&A \end{psmallmatrix}|A,B\in\Mat_\tn(\kk)}$.
Then $M(\tm|\tn)$ is a symmetric Frobenius superalgebra with trace form
given by the super trace
$\tr_{M(\tm|\tn)}(E_{i,j})=\epsilon_i\delta_{i,j}$, where $\epsilon_i=1$ if
$i\leq \tm$, and $\epsilon_i=-1$ otherwise.
 In particular, $\tau_{M(\tm|\tn)}=\sum_{1\leq i,j\leq
 \tm+\tn}\epsilon_jE_{i,j}\otimes E_{j,i}$ with $\epsilon_j=1$ if $j\le
 \tm$
 and $\epsilon_j=-1$ otherwise. The superalgebra $Q_\tn$ can be equipped
 with a non-degenerate odd symmetric form, but not with an \emph{even} symmetric one. To fit it into our framework we could take the non-symmetric trace $\tr\left(\begin{psmallmatrix} A&B\\B&A \end{psmallmatrix}\right)=\op{Tr}(A)$, with $\psi\left(\begin{psmallmatrix} A&B\\B&A \end{psmallmatrix}\right)=\left(\begin{psmallmatrix} A&-B\\-B&A \end{psmallmatrix}\right)$ so that $\theta=2$. For $n = 1$ this is
\autoref{Ex:Clifford}.
\end{Example}
The algebra in the next example is a prominent player in representation theory and categorifications. It can be realised in many ways and contexts, for example  as a Brauer tree
  algebra (with $m$ edges and no exceptional vertex in the classification of
  \cite{Rickard}), as a zigzag algebra of type $A_m$ as in
  \cite{KhovanovHuerfano, KhovanovSeidel}, as a generalized
  Khovanov algebra $H_{1,m-1}$ in the sense of
  \cite{BrundanStroppel:KhovanovI},  as the endomorphism algebra of the sum of indecomposable
 projective-injectives in parabolic category $\mathcal{O}$,
 \cite[Theorem~4.6]{MazStr}, and also governs the Heisenberg categorification of \cite{CautisLicata:Heisenberg}.

  \begin{Example}\label{Ex:ZigZagm}
Consider the  \emph{zig-zag algebra} $A=\mathcal{Z}ig_m=\kk Q\big/\bigl(a_{i+1}a_i,\bar a_i\bar a_{i+1}, a_i\bar a_ i -\bar a_{i+1}a_{i+1}\mid 1\leq i\leq m-2 \bigr)
$, where $\kk Q$ is the path algebra of the quiver $
  1\ \underset{\bar a_1}{\overset{a_1}{\rightleftarrows}}\ 2
   \ \underset{\bar a_2}{\overset{a_2}{\rightleftarrows}}\ 3
   \ \underset{\bar a_3}{\overset{a_3}{\rightleftarrows}}\ \cdots
   \ \underset{\bar a_{m-1}}{\overset{a_{m-1}}{\rightleftarrows}} m
$ for $m\geq 3$. Note that all length three paths are zero, for $m=2$  we need to add these to the relations.
Note that $\mathcal{Z}ig_m$  has basis $\cB=\{e_1,\dots,e_m\}\cup\{a_1,\dots,a_{m-1}\}\cup\{\bar a_1,\dots,\bar a_{m-1}\}\cup\{l_1,\dots,l_m\}$ consisting of homogeneous elements if we equip $A$ with the grading
 by the path length.  The  top degree is spanned by
 $l_1=\bar a_1a_1$, $l_i=\bar a_ia_i=a_{i-1}\bar a_{i-1}$  for $1<i<m$ and $l_m=a_{m-1}\bar a_{m-1}$.
The algebra is  a symmetric Frobenius algebra with trace form $\tr(b)=\sum_{i=1}^m\delta_{b,l_i}$ for $b\in \cB$. For an explicit formula for $\taui[1,2]$ see \autoref{tauZigZag}.
\end{Example}

  \section{RaFH superalgebras}\label{S:AffineRH}
In this section we define the main algebras of this paper, the RaFH superalgebras
$\HHd(\Bo)=\HHd[n;\hbar,t,\eta](\Bo)$. They form a family of  superalgebras over a commutative integral domain $R$ that is a $\kk$-algebra containing the parameters $\hbar,t,\eta$. A RaFH superalgebra is attached to a
Frobenius-tensor superalgebra $\Bo$ which is $\eta$-symmetric, that is $\eta\psi(a)=\eta a$ for
$a\in A$, see \autoref{D:Frobtensor}. It is generated over $R$ by $\Bo$, by polynomial generators
$X_1,\dots,X_n$ and by Hecke generators $H_1,\dots,H_{n-1}$. This is subject to certain relations
\autoref{RaFHrels} in which the elements $\tau^{(i,i+1)}\in\Bo$ enter. The parameters
$\hbar,t,\eta$ interpolate between known algebras.  In \autoref{S:admissibility} we explain the
role of $\eta$-symmetricity, which for $\eta\ne0$ forces that $\tau^{(i,j)}=\tau^{(j,i)}$ and that the
polynomial generators commute with $\Bo$. After deducing some commutation relations, we introduce
in \autoref{S:Demazures} RaFH derivations and Frobenius--Demazure operators. They generalise the
Demazure operators of \cite{Savage:AffineWreath,RossoSavage:QuantumAffineWreath}  and
control the commutators of the $H_i$ with the subalgebra generated by $\Bo$ and the $X_j$.  Even more general operators were developed independently by \cite{LaiMinets, MinetsYBE}. These operators
are the main tool in the proof of the basis theorem, \autoref{C:BasisHd}, which shows that $\HHd(\Bo)$ is free
over $\RR$ with basis $\set{\bb X^\bd H_w}$, and three further bases of the same shape.\\

In this section let $\RR$ be a commutative integral domain which is a $\kk$-algebra, and let $\hbar,t,\eta\in R$.

\subsection{The definition of  RaFH superalgebras}\label{defRaFH}
For a $\kk$-superalgebra $B$ we abbreviate $B_R=R\otimes_\kk B$ and let
$R\<Y_1,\dots,Y_m\>\ast B_R$ be the free product of $R$-superalgebras of the
free $R$-algebra on even generators $Y_1,\dots,Y_m$ and $B_R$. It is the $R$-superalgebra generated by
even elements $Y_1,\dots,Y_m$ and by $B_R$  subject only to the relations of $B_R$. We identify $b\in B$
with $1\otimes b \in B_R$. Given a  Frobenius-tensor superalgebra $\Bo$ we say $\Bo_R$ is  \emph{$\eta$-symmetric} if $\eta=0$ or $\psi$ is trivial. For $\eta\in\kk$ this is the condition of
\autoref{D:Frobtensor}, since $\Bo$ is free over $\kk$. We next define the RaFH superalgebras (of type $A$, with the superflip braiding).
\begin{Definition}\label{D:Hecke}
  Fix an integer $n\ge0$ and let $\Bo$ be an $\eta$-symmetric Frobenius-tensor superalgebra.
  The \emph{Rees affine Frobenius Hecke superalgebra}, or short \emph{RaFH superalgebra},
  $\HHd(\Bo)=\HHd[n;\hbar,t,\eta](\Bo)$ is the quotient of the
  $R$-superalgebra $R\<H_1,\dots,H_{n-1},X_1,\dots,X_n\>\ast \Bo_R$, with $H_i$ and $X_j$ even,
  by the two-sided ideal generated by the relations
  \begin{equation}\label{RaFHrels}
  \begin{gathered}
    X_iX_j=X_jX_i,\qquad
    H_iH_j=H_jH_i\ \text{ if }|i-j|>1,\qquad
    H_iH_{i+1}H_i=H_{i+1}H_iH_{i+1},\\
    H_i^2=\eta\tau^{(i,i+1)}H_i+\hbar^2,\qquad
    \bb X_i=X_i\psi^{(i)}(\bb),\qquad
    H_i\bb=\leftindex^{s_i}{\bb}H_i,\\
    H_iX_j=X_jH_i\ \text{ if }j\ne i,i+1,\qquad
    H_iX_i=X_{i+1}H_i-\gamma_i\tau^{(i,i+1)},\qquad
    H_iX_{i+1}=X_iH_i+\gamma_i\tau^{(i+1,i)},
  \end{gathered}
  \end{equation}
  for all $\bb\in\Bo$ and all admissible indices. Here,
  $\gamma_i=t+\eta X_{i+1}$ for $1\leq i<n$.
\end{Definition}

  \begin{Remark}
    If $n=0$, then $\HHd(\Bo)=R$. If $n=1$, then
    $\HHd(\Bo)=R[\X]\#\Bo=R[X_1]\#A$.
  \end{Remark}

  \begin{Remark} Alternatively, $\HHd(\Bo)$ can be realised as a
    quotient of a smash product algebra using $H_i\bb = \leftindex^{s_i}\bb H_i$. Our algebras are special examples of quantum wreath product algebras \cite{CJNakanoXiang,LaiMinets}.
  \end{Remark}

   \begin{Remark}\label{flipbraiding}
    The name RaFH superalgebra stands for \emph{Rees affine Frobenius Hecke superalgebra},
    since the algebras are a Rees algebra version of affine Hecke
    algebras associated with Frobenius-tensor superalgebras, see \autoref{Reesex}.
     The RaFH
    algebras $\HHd(\Bo)$ are of \textit{type $A$} because we will show
    that if $\Bo$ is the trivial Frobenius-tensor superalgebra, then
    the $H_i$ generate a finite Hecke algebra of type~$A$, see \autoref{D:RFHS} with \autoref{matchRFH}. It is a RaFH superalgebra
    \emph{with superflip braiding}, since we require
    $H_i\bb = \leftindex^{s_i}\bb H_i$ so that $s_i$ acts on $\Bo$ by the superflip and $\Sym$ acts by superpermutations.
    We could  replace the action of $\Sym$, both here and
    in $\kk\Sym\Smash\Tos$, with an action of a braid group or a Hecke
    algebra following  \cite{CJNakanoXiang,LaiMinets,MinetsYBE}, but this would make the paper by far too technical.
      \end{Remark}
      \begin{Remark}\label{Rk:Zgraded}
  Suppose that $\Bo=A^{\otimes n}$, where $A$ is a $\Z$-graded Frobenius algebra with trace form homogeneous of degree $-d$,
  so that $\taui[i,j]$ is homogeneous of degree $d$. Assume that $\RR$ has a $\Z$-grading with
  $\deg\hbar=\deg t=0$ and $\deg\eta=-d$, and set $\deg X_i=d$ and
  $\deg H_i=0$. Then all relations in \autoref{RaFHrels} are homogeneous, so
  $\HHd(\Bo)$ is a $\Z$-graded algebra. Everything below applies verbatim to
  algebra objects in the category of $\Z$-graded $\kk$-modules.
\end{Remark}

      \subsection{Some relations in RaFH algebras}We deduce some important relations in $\HHd$.\label{S:somerelations}
\begin{Lemma}\label{L:HXH}
  Let $1\le i<n$. Then
  $\hbar^2 X_{i+1}=H_iX_iH_i +t H_i\tau^{(i+1,i)}$ in $\HHd$.
\end{Lemma}

\begin{proof}
  This follows by multiplying the second last relation in \autoref{RaFHrels} by $H_i$ from the right.
\end{proof}

\begin{Lemma}\label{Halpha}
Let $1\le i<n$. Then $H_i(X_{i+1}-X_i)=(X_i-X_{i+1})H_i+\gamma_i(\tau^{(i+1,i)}+\tau^{(i,i+1)})$.
\end{Lemma}
\begin{proof}
This follows directly from the last two relations in \autoref{RaFHrels}.
\end{proof}

  For $w\in\Sym$ define $H_w=H_{i_1}\dots H_{i_r}\in\HHd$, where
  $w=s_{i_1}\dots s_{i_r}$ is a reduced expression for $w$. By
  Matsumoto's Theorem (see, for example, \cite[Theorem~1.8]{Mathas:ULect}), $H_w$ depends only on $w$ and not
  on the choice of a reduced expression. For example, if $w=(i,k+1)$ and $1\leq i\leq k$ then
  \begin{equation}\label{Htranspo}
  H_w=H_{(i,k+1)}=H_kH_{k-1}\ldots H_{i+1}H_iH_{i+1}\ldots H_{k-1}H_k.
  \end{equation}

The quadratic relation in \autoref{D:Hecke} implies the following multiplication rules.

\begin{Lemma}\label{L:HiHw}
   Suppose $1\le i<n$ and $w\in\Sym$. Then
    \begin{small}
    \begin{equation*}
   H_iH_w = \begin{cases*}
     H_{s_iw} & if $\ell(s_iw)>\ell(w)$,\\
     \hbar^2H_{s_iw} +\eta\taui[i,i+1] H_w & if $\ell(s_iw)<\ell(w)$.
   \end{cases*}\quad\text{}\quad
   H_wH_i = \begin{cases*}
     H_{ws_i} & if $\ell(ws_i)>\ell(w)$,\\
     \hbar^2H_{ws_i} +\eta\taui[w(i+1),w(i)] H_w & if $\ell(ws_i)<\ell(w)$.
   \end{cases*}
   \end{equation*}
    \end{small}
     \end{Lemma}
\begin{Definition}\label{specialisation}
Given a ring homomorphism $\varphi\colon R\rightarrow R'$, we define the
corresponding \emph{specialised RaFH superalgebra}
$\HHd[n;\varphi(\hbar), \varphi(t), \varphi(\eta)](\Bo)_{R'}
  =R'\otimes_R\HHd[n;\hbar, t, \eta](\Bo)$.
\end{Definition}

The next result recovers \autoref{D:universalJM} in the case $\eta=0$ and $\hbar=t=1$.
 \begin{Lemma} \label{evaluation}
 The quotient map $\HHd(\Bo)\rightarrow \HHd(\Bo)/(X_1)$ sends $\hbar^{2(k-1)}X_{k}$ to $\sum_{1\leq i<k}t\hbar^{2(i-1)}H_{(i,k)}\tau^{(k,i)}$.
  \end{Lemma}
 \begin{proof}
 The case $k=1$ is clear, since $X_1$ is sent to zero. By  \autoref{L:HXH} and induction we have
 \[\hbar^{2k}X_{k+1}=H_k(\hbar^{2(k-1)}X_k)H_k+\hbar^{2(k-1)}tH_k\tau^{(k+1,k)}\mapsto
 t\left(\sum_{1\leq i<k}\hbar^{2(i-1)}H_kH_{(i,k)}\tau^{(k,i)}H_k\right)+\hbar^{2(k-1)}tH_{(k,k+1)}\tau^{(k+1,k)}.\]
  Since $\tau^{(k,i)}H_k=H_k\tau^{(k+1,i)}$ and $H_kH_{(i,k)}H_k=H_{(i,k+1)}$, by \autoref{Htranspo}, the lemma follows.
  \end{proof}

\begin{Definition}Let $\hbar$ be invertible.
 Then $\CL_k=t\sum_{1\leq i<k}\hbar^{2(i-k)}H_{(i,k)}\tau^{(k,i)}=t\sum_{1\leq i<k}\hbar^{2(i-k)}\tau^{(i,k)}H_{(i,k)}$ are the \emph{Jucys--Murphy elements} of  $\HHd(\Bo)/(X_1).$
 Similarly we can define Jucys--Murphy elements in $\HHd(\Bo)/(f(X_1))$ for any polynomial $f$, in particular for the CRaFH algebras from \autoref{D:CyclotomicHecke} below.
 \end{Definition}

  \subsection{The significance of $\eta$-symmetricity}\label{S:admissibility}
  By definition, the RaFH superalgebra $\HHd$ is built from a
  Frobenius-tensor superalgebra $\Bo$, with a Nakayama automorphism
  $\psi$ that is not necessarily trivial. By assumption, $\Bo$
  is $\eta$-symmetric, which means that $\eta a=\eta\psi(a)$, for $a\in A$.
  This section explores some nice consequences of this definition. The origin
  of this terminology is that if $A$ is a Frobenius algebra, then the
  corresponding Frobenius-tensor superalgebra $\Bo=A^{\otimes n}$ is
  $\eta$-symmetric for some $\eta\not=0$ if and only if $A$ is symmetric. It is always $0$-symmetric, since this is no condition at all.

  \begin{Lemma}\label{Lfunnyequation}
    Consider $\HHd=\HHd[n;\hbar, t, \eta](\Bo)$. Then either $\eta=0$, or
      \begin{equation}\label{funnyequation}
 \eta\not=0\quad\text{and}\quad  (\tau^{(i,j)}=\tau^{(j,i)}\quad\text{and}\quad  f\bb= \bb f)\quad\quad \text{for any  $\bb\in\Bo$,  $f\in\RR[\X]$, $1\leq i\not=j\leq n$.}
 \end{equation}
  \end{Lemma}
\begin{proof}
  If $\eta=0$ there is nothing to prove. If $\eta\ne0$, then $\eta$-symmetricity implies that $\psi$ is
  trivial by \autoref{rketasymm}, since $\RR$ is a domain.
\end{proof}

\begin{Lemma}\label{Rk:hinv}Let $\hbar$ be invertible. The following hold in $\HHd$.
\begin{enumerate}
\item\label{Hinv1} The generator   $H_i\in\HHd$ is  invertible for any $1\leq i<
  n$ with inverse  $H_i^{-1}=\hbar^{-2}(H_i-\eta\tau^{(i,i+1)})$.
\item\label{Hinv2}  The element $H_w$ is invertible, for any $w\in\Sym$.
\item\label{Hinv3}  By conjugating by $H_i$, the last two RaFH superalgebra relations are
  equivalent, and both of these relations are equivalent to the identity
  in \autoref{L:HXH}.
\end{enumerate}
\end{Lemma}
\begin{proof} The quadratic relations for $H_i^2$ in \autoref{RaFHrels} imply that
  $\hbar^{-2}(H_i-\eta\tau^{(i,i+1)})$ is a left-inverse of $H_i$ and
  that $\hbar^{-2}(H_i-\eta\tau^{(i+1,i)})$ is a right-inverse. By
  $\eta$-symmetricity and \autoref{Lfunnyequation}, the left and right
  inverses of $H_i$ coincide, proving \autoref{Hinv1}. Parts
  \autoref{Hinv2} and \autoref{Hinv3} follow directly from the definitions.
\end{proof}

Thus, $\eta$-symmetricity ensures that the right and left
inverses of $H_i$ agree, when they exist.

\begin{Definition}\label{D:RFHS}
  Let $\Bo$ be a Frobenius-tensor superalgebra. The \emph{Rees Frobenius Hecke superalgebra} $\Hfind(\Bo)=\Hfind[n;\hbar,\eta](\Bo)$
  (of type $A$ with superflip braiding)  is the quotient of
  $ R\<H_1,\dots,H_{n-1}\>\ast\Bo_R$ by the two-sided ideal generated by the relations
  \begin{equation}\label{RFHSrels}
  \begin{gathered}
    H_iH_j=H_jH_i\ \text{ if }|i-j|>1,\qquad H_iH_{i+1}H_i=H_{i+1}H_iH_{i+1},\\
    H_i^2=\eta\tau^{(i,i+1)}H_i+\hbar^2,\qquad H_i\bb=\leftindex^{s_i}{\bb}H_i,
  \end{gathered}
  \end{equation}
  for all $\bb\in\Bo$ and all admissible $i,j$. For $\Bo=\kk$ this is the (well-known) \emph{finite Hecke algebra} $\HSd$.
\end{Definition}
For $w\in \Sym$, define $H_w=H_{i_1}\dots H_{i_r}\in\Hfind$, where
  $w=s_{i_1}\dots s_{i_r}$ is a reduced expression for $w$.

\begin{Proposition}\label{P:H0Basis}
   Let $\Bo$ be an $\eta$-symmetric Frobenius-tensor superalgebra
   with $R$-basis $\cBo$. Then
   the Rees Frobenius Hecke superalgebra $\Hfind(\Bo)$
   is free as an $R$-module with basis $\set{H_wb| w\in W, b\in\cBo}$.
 \end{Proposition}
\begin{Remark}With similar arguments one shows that $\set{bH_w| w\in W, b\in\cBo}$ is a basis as well.
\end{Remark}
We first show a converse result, namely that the elements $\eta\tau^{(i,i+1)}$ are forced to be symmetric as soon
as $\Hfind(\Bo)$ has the basis of \autoref{P:H0Basis}. For $\Bo=\Ao$ this means that $A$ has to be
$\eta$-symmetric, since $\psi^{(j)}(\tau^{(i,j)})=\tau^{(j,i)}$ by \autoref{Nakayama} and
$\tau^{(1,2)}=\sum_bb\otimes b^\vee$ for a basis $\set{b}$ of $A$ with dual basis $\set{b^\vee}$.
This result explains axiom \cite[(3.2.3)]{CJNakanoXiang} for quantum wreath product algebras,
and it shows why the quantum affine wreath product algebras of
\cite{RossoSavage:QuantumAffineWreath} only take symmetric Frobenius algebras as input.

\begin{Lemma}\label{L:symmetric}
 Let $\Bo$ be a, not necessarily $\eta$-symmetric, Frobenius-tensor
  superalgebra and suppose that $\Hfind(\Bo)$ has a basis of the form
  $\{\bb H_w\}$, where $w\in\Sym$ and $\bb$ runs through a basis of
  $\Bo$. If $\hbar$ is invertible,  then $\eta\tau^{(i,i+1)}= \eta\tau^{(i+1,i)}$ for all $1\leq i
  \leq  n-1$.
\end{Lemma}

\begin{proof}
We expand $H_i^3$ in two different ways.  The quadratic relation  gives
\[H_i^3=H_i^2H_i=(\eta\tau^{(i,i+1)}
H_i+\hbar^2)H_i=((\eta\tau^{(i,i+1)})^2+\hbar^2)H_i+\eta\hbar^2\tau^{(i,i+1)}.\]
On the other hand, we have
$H_i^3=((\eta\tau)^2+\hbar^2)H_i+\hbar^2\eta\tau^{(i,i+1)}=((\eta\tau^{(i+1,i)})(\eta\tau^{(i,i+1)})
+\hbar^2)H_i+\eta\hbar^2\tau^{(i+1,i)}$, since $H_i\tau^{(i+1,i)}=\tau^{(i,i+1)}H_i$. The claim follows by comparing
the coefficients of the basis elements.
\end{proof}

\begin{proof}[Proof of \autoref{P:H0Basis}]
  The relations easily imply that $\Hfind$ is spanned by the elements
  $\set{H_w b}$. To prove linear independence, let~$E$ be a free $R$-module with
  basis $\set{e_{w,b}|w\in\Sym,
  b\in\cBo}$. If $x=\sum_{b\in\cBo}r_bb$, set $e_{w,x}=\sum_br_be_{w,b}$.
  Let $c\in\cBo$. We claim that $\Hfind$ acts on $E$ by the formulas
  \[
    ce_{w,b} = e_{w,\leftindex^{w^{-1}}cb},
      \quad\text{and}\quad
      H_ie_{w,b} = \begin{cases*}
        e_{s_iw,b},& if $\ell(s_iw)=\ell(w)+1$,\\
        \eta \tau^{(i,i+1)}e_{w,b}+\hbar^2e_{s_iw,b},& otherwise,
      \end{cases*}
  \]
  The only non-obvious relations to check
  are the quadratic relations. The argument in \autoref{L:symmetric}
  shows, using \autoref{Lfunnyequation} and the equality
  $\leftindex^{s_i}{(\tau^{(i,i+1)})}=\tau^{(i+1,i)}$ in~$\Bo$, that $\eta$-symmetricity implies
  the quadratic relations. The
  proposition now follows by standard arguments, as in
  \cite[Theorem~1.13]{Mathas:ULect}.
\end{proof}

\subsection{Important RaFH superalgebras}
  Many (super)algebras in the literature can be realised as RaFH superalgebras $\HHd[n;\hbar, t, \eta](\Bo)$ with prominent examples arising already from the  trivial Frobenius-tensor superalgebra $\Bo=\kk$.

\begin{Examples}\label{HHdexamples}
The following are examples of (specialised) RaFH superalgebras:
  \begin{enumerate}[series=HHexamples, label=(\alph*), ref=\autoref{HHdexamples}.(\alph*)]
   \item \label{ex:smash}The \emph{smash product algebra}  $\Hsmash:=\RR\Sym\Smash\RR[\X]=\Hd[n;1,0,0](\kk)$.
   \item \label{degaff} The \emph{degenerate affine Hecke algebras} $\Hdegen:=\Hd[n;1,1,0](\kk)$ of type~$A$ from  \cite{Drinfield:DegenerateHecke, Lusztig:AffineGradedHecke}.
   \item \label{nilHecke} The specialisation $\Hd[n;0,1,0](\kk)$ is isomorphic to the \emph{polynomial affine nilHecke algebra}, also just called \emph{ nilHecke algebra}, from \cite{KK86, KhovLaud:diagI}.
   \item The \emph{affine Nil Coxeter algebra}, an affinisation of the
   \emph{Nil Coxeter algebra} $\NC$ from \cite{NilCoxeter},  is the RaFH
   superalgebra $\Hd[n;0,0,0](\kk)$.  This recovers the algebra with even polynomial generators defined in \cite{SavageStuart}.
   \item The \emph{affine $0$-Hecke algebra} is the  specialisation $\Hd[n;0,0,\pm1](\kk)$.
   \item\label{ex:aff} The \emph{extended affine Hecke algebra}, or
   \emph{Iwahori--Hecke algebra},
     from~\cite{Macdonald:affine, IwahoriMatsumoto}, which is
     as $R$-module  $\Heaff=R[Y_1^{\pm},\dots,Y_n^{\pm}]\otimes_R\HSd$ with the finite Hecke algebra $\HSd$ from \autoref{D:RFHS},
     is a RaFH superalgebra in at least two ways.
     First, as localisation of the  \emph{polynomial affine Hecke
     algebra}  $\Hd[n;1,1,q-q^{-1}](\kk)$ inverting the $X_i$'s and,
     second, as a subalgebra of the completed RaFH superalgebra
     $\Hd[n;1,0,q^{-1}-q](\kk)$ along the $X_i$'s.

     For the latter, the identification
     map between the algebras is given by the identity on $\HSd$ and by $X_i\mapsto \frac1{q^{-1}-q}(Y_i-1)$ on the polynomial part. These
     constructions only work when $q\ne1$. We obtain
     Iwahori's and
     Matsumoto's~\cite{IwahoriMatsumoto} normalisation of the
     Iwahori--Hecke algebras from the RaFH algebra
     $\Hd[n;\sqrt{q},0,q-1](\kk)$, in which the Hecke generators satisfy
     the quadratic relation $H_i^2=(q-1)H_i+q$ from
     \cite{IwahoriMatsumoto}. This normalisation can also be obtained
     from $\Hd[n;1,0,q^{-1}-q](\kk)$ by rescaling the $H_i$.
 \end{enumerate}
 The reader may find the overview in  \cite[I.2.2]{Meinel} relating
 these different algebras useful. Let $A$ be any Frobenius superalgebra.
 Then we recover the following constructions:
 \begin{enumerate}[resume*=HHexamples]
   \item\label{ex:quantumwreath} The \emph{(quantum) affine wreath product algebras}
   $\Hd[n;1,1,0](\Ao)$ from \cite{Savage:AffineWreath}, and
   $\Hd[n;1,1,z](\Ao)$  from \cite{RossoSavage:QuantumAffineWreath} for
   $A$ symmetric. (\autoref{C:BasisHd} forces the
   Frobenius algebra to be symmetric when $\eta\not=0$). For $A$ a (skew) zigzag algebra these algebras appear in the endomorphism algebras of the
Heisenberg categories of \cite[\S10.3]{CautisLicata:Heisenberg} and in connection with the geometry of Hilbert schemes. The purely even version is studied in \cite{KleshchevMuth}.

   \item\label{ex:Sergeev} The \emph{degenerate affine Sergeev algebra}, introduced in \cite{Sergeev, Nazarov}, is
     $\Hd[n;1,1,0](\Ao)$ where we take   $\Bo=\Cl_n$, the Clifford algebra from
     \autoref{Ex:Clifford} viewed as Frobenius-tensor superalgebra for
     $A=\Cl_1$. The affine Clifford--Hecke algebra is however not an example of a RaFH
     algebra because we only allow the superflip braiding, see also \autoref{flipbraiding}.

\item Taking $A=\kk G$ as in \autoref{Ex:FrobeniusAlgebrasfingroup} we
  obtain the \emph{wreath product RaFH superalgebras}
  $\Hd[n;1,1,0](G)=\Hd[n;1,1,0](\Ao)$. As a special case we can take
  $G=\Sym_m$ for any $m$. Note that the group $G=\Sym_m$ plays
  a completely different role to the symmetric group $\Sym=\Sym_n$.
  \item \label{Ex:KMaffinisation}
  Let $A$ be a symmetric algebra, free of finite rank over $\kk$, so that
  $\psi=\op{id}$ and $\taui[i,j]=\taui[j,i]$ by \autoref{Taussuper}. We recover  the \emph{rank $n$ affinisation} $H_n(A)$ of $A$ from
  \cite[Definition~3.2]{KleshchevMuth} via
  \begin{equation}\label{KleshchevMuthiso}
    \Hd[n;1,1,0](\Ao)\;\cong\;H_n(A),
    \qquad X_j\mapsto -z_j,\quad H_i\mapsto s_i,\quad
    \taui[i,i+1]\mapsto\Delta_{i,i+1}.
  \end{equation}
  Indeed, specialising $\hbar=t=1$ and $\eta=0$ gives $\gamma_i=1$, so that the last two
  relations of \autoref{RaFHrels} read $H_iX_i=X_{i+1}H_i-\Delta_{i,i+1}$ and
  $H_iX_{i+1}=X_iH_i+\Delta_{i,i+1}$. Replacing $X_j$ by $-z_j$ turns these into
  $s_iz_j-z_{s_i(j)}s_i=(\delta_{i,j}-\delta_{i+1,j})\Delta_{i,i+1}$, which is
  \cite[(3.5)]{KleshchevMuth}, while the remaining relations of \autoref{RaFHrels} become
  \cite[(3.3),(3.4)]{KleshchevMuth} and $s_i^2=1$. Taking $A=\mathcal{Z}ig(\Gamma)$, the zigzag
  algebra of a connected graph $\Gamma$ without loops or multiple edges, which for $\Gamma$ of
  type $A_m$ is \autoref{Ex:ZigZagm}, gives the \emph{affine zigzag algebra}
  $\mathcal{Z}ig^{\op{aff}}_n(\Gamma)=H_n(\mathcal{Z}ig(\Gamma))$ of
  \cite[Definition~4.4]{KleshchevMuth}.
  \end{enumerate}
\end{Examples}
\begin{Remark} \label{Reesex}The RaFH algebra $\Hd[n;\hbar,1,0](\kk)$ deserves its name. It is the Rees algebra of the degenerate affine Hecke algebra $\Hd[n;1,1,0](\kk)$ from \autoref{degaff} for the  filtration which  puts  $H_i$ in degree $1$ and $X_i$ in degree $-1$. Its associated graded is the nilHecke algebra  $\Hd[n;0,1,0](\kk)$ from \autoref{nilHecke}.
\end{Remark}

\subsection{Frobenius--Demazure operators and RaFH derivations}\label{S:Demazures}
We next define Frobenius--Demazure operators $\delta_i$ on $\Po=R[\X]\Smash\Bo$. These are Frobenius
analogues of the well-known divided difference operators on polynomial rings, to which they reduce for trivial Nakayama
automorphism up to the factor $\tau$, see \autoref{R:classical}. They are slight
generalisation of the deformed divided difference operators from
\cite{Savage:AffineWreath,RossoSavage:QuantumAffineWreath} and should be compared with similar operators which were introduced independently in \cite{LaiMinets} for commutative Frobenius algebras. Multiplying
$\delta_i$ with the central elements $\gamma_i=t+\eta X_{i+1}$, see \autoref{L:gammacentral}, gives the RaFH derivations
$\Delta_i$. These encode the commutator rules of the Hecke generators $H_i$ with $\Po$, see
\autoref{L:crucial}. They interpolate between the cohomological and the
$K$-theoretic Demazure operators, see \autoref{RkDemop}. Our $\eta$-symmetricity
assumption is exactly what is needed to treat the degenerate case, with a possibly nontrivial
Nakayama automorphism.

Throughout, $\eta\in\RR$ and $\Bo$ is a  Frobenius-tensor superalgebra such that $\Bo_\RR=\RR\otimes_\kk\Bo$ is $\eta$-symmetric.

\begin{Definition}\label{DefPo}
   Let $\Po=R[\X]\Smash\Bo$ be the smash product of $R[\X]$ with $\Bo_\RR$ over $\RR$. The smash product relations are the relations $\bb X_i=X_i\psi^{(i)}(\bb)$ as in  \autoref{RaFHrels}.
   \end{Definition}
   \begin{Remark} $\Po$ is a free $R$-supermodule with $R$-basis
  $\set{X^\bd\bb}$, where $\bd\in\N^n$ and $\bb$ runs through a basis of $\Bo$ and multiplication $(X^\bd\bb)(X^\be\bb')=X^{\bd+\be}\psi_\be(\bb)\bb'$,
   where we write $\psi_i=\psi^{(i)}$ and $\psi_\bd=\psi_1^{d_1}\cdots\psi_n^{d_n}$ for $\bd\in\N^n$. In particular, $\set{\bb X^\bd}$ is an $R$-basis too.
  \end{Remark}

The most important definition of this subsection is the following:
\begin{Definition}\label{skewDer}
Let $w\in\Sym$.  A $w$-\emph{skew-derivation} $\Delta$ on $\Po$ is an $R$-linear endomorphism $\Delta$ of $\Po$ which satisfies the following  \emph{skew Leibniz rule}
$\Delta(xy)=\Delta(x)y+\leftindex^{w}x\Delta(y)$ for any $x, y\in\Po$. Here, $\leftindex^{w}(pb)=\leftindex^{w}p\leftindex^{w}b$ for $p\in R[\X]$, $b\in\Bo$.
\end{Definition}

\begin{Lemma}\label{L:gammacentral}
The following gives a construction of new skew-derivations from given ones.
\begin{enumerate}
\item \label{central1} For $1\le i<n$ it holds that  $\gamma_i=t+\eta X_{i+1}\in Z(\Po)$.
\item  \label{central2}  If $\delta$ is a $w$-skew-derivation on $\Po$ and $z\in Z(\Po)$, then
  $z\delta$ is also a $w$-skew-derivation.
  \item\label{L:transfer}
  Let $\delta$ be an $s_i$-skew-derivation of $\Po$ and $z\in Z(\Po)$. Then
  \begin{equation}\label{transfer}
    (z\delta)^2=z\delta(z)\delta+z\leftindex^{s_i}z\delta^2
    \quad\text{and}\quad
    (z\delta)s_i+s_i(z\delta)
      =\bigl(\leftindex^{s_i}z-z\bigr)s_i\delta+z\bigl(\delta s_i+s_i\delta\bigr)=\leftindex^{s_i}zs_i\delta+z\delta s_i.
  \end{equation}
  \end{enumerate}
\end{Lemma}
\begin{proof}
  \autoref{central1}: Clearly $\gamma_i$ commutes with all $X_j$.
  For $\bb\in\Bo$, we have
  $\eta\bb X_{i+1}=\eta X_{i+1}\psi^{(i+1)}(\bb)=\eta X_{i+1}\bb$  by $\eta$-symmetricity,
  see \autoref{funnyequation}.

   \autoref{central2}: This follows from
  $z\delta(xy)=z\delta(x)y+z\leftindex^wx\delta(y)
  =(z\delta)(x)y+\leftindex^wx\,(z\delta)(y)$.

   \autoref{L:transfer}: Using~\autoref{central2}, the first formula holds since
  $(z\delta)(z\delta(y))=z(\delta(z)\delta(y)+\leftindex^{s_i}z\delta^2(y))$ for any $y\in\Po$. The last formula is obvious. For the second last note that
  $(z\delta)(\leftindex^{s_i}y)+\leftindex^{s_i}{\left(z\delta(y)\right)}
   =(z\delta)(\leftindex^{s_i}y)+\leftindex^{s_i}z\leftindex^{s_i}\delta(y)$.
   \end{proof}

The following definition generalises the divided difference operators from  \cite{Savage:AffineWreath,RossoSavage:QuantumAffineWreath}.
\begin{Definition}\label{D:FrobDem}
  For $1\le i<n$ we define the \emph{Frobenius--Demazure operator}
  $\delta_i$ as the unique $s_i$-skew-derivation
  $\delta_i\colon\Po\to\Po$ satisfying
  \begin{equation}\label{deltavalues}
    \delta_i(X_i)=\tau^{(i,i+1)},\;\;
    \delta_i(X_{i+1})=-\tau^{(i+1,i)},\;\;
    \delta_i(X_j)=0,\;\;
    \delta_i(\bb)=0,
    \quad\text{where $j\not=i,i+1$, and $\bb\in\Bo$.}
  \end{equation}
\end{Definition}

\begin{Lemma}\label{L:deltaexists}
  The operator $\delta_i$ from \autoref{D:FrobDem} is well-defined.
\end{Lemma}

\begin{proof}
  Uniqueness is clear, since the elements listed in \autoref{deltavalues}
  generate $\Po$. For the existence we verify the compatibility with the
  defining (local) relations from \autoref{DefPo}, which is enough by the skew Leibniz rule.

  Consider first the relations $X_jX_k=X_kX_j$. If
  $\set{j,k}\cap\set{i,i+1}=\emptyset$, then $\delta_i$ annihilates both sides. If
  $j=i$ and $k\not=i,i+1$, then $\delta_i(X_iX_k)=\tau^{(i,i+1)}X_k$
  whereas $\delta_i(X_kX_i)=X_k\tau^{(i,i+1)}$. They agree by the
  third condition in \autoref{Nakayama}, which gives
  $\psi^{(k)}(\tau^{(i,i+1)})=\tau^{(i,i+1)}$. Similarly for $j=i+1$. Finally,
   $ \delta_i(X_iX_{i+1})=\tau^{(i,i+1)}X_{i+1}-X_{i+1}\tau^{(i+1,i)}=0$ and
    $\delta_i(X_{i+1}X_i)=-\tau^{(i+1,i)}X_i+X_i\tau^{(i,i+1)}=0$, as the second condition in \autoref{Nakayama} gives
  $\tau^{(i,i+1)}X_{i+1}=X_{i+1}\psi^{(i+1)}(\tau^{(i,i+1)})=X_{i+1}\tau^{(i+1,i)}$
  and $\tau^{(i+1,i)}X_{i}=X_{i}\psi^{(i)}(\tau^{(i+1,i)})=X_{i}\tau^{(i,i+1)}$.

  Consider next the relations $\bb X_j=X_j\psi^{(j)}(\bb)$. For
  $j\not=i,i+1$ both sides vanish. For $j=i$ we get
  $  \delta_i(\bb X_i)=\leftindex^{s_i}\bb\,\tau^{(i,i+1)}
      =\tau^{(i,i+1)}\psi^{(i)}(\bb)
      =\delta_i\bigl(X_i\psi^{(i)}(\bb)\bigr),$
  where the middle equality is  the fourth
  condition in \autoref{Nakayama} applied to
  $\ba=\leftindex^{s_i}\bb$. For $j=i+1$ the same computation for the  pair $(i+1,i)$ gives
  $\delta_i(\bb X_{i+1})=-\leftindex^{s_i}\bb\,\tau^{(i+1,i)}
   =-\tau^{(i+1,i)}\psi^{(i+1)}(\bb)=\delta_i(X_{i+1}\psi^{(i+1)}(\bb))$.
  The lemma follows.
\end{proof}

\begin{Remark}
The proof of \autoref{L:deltaexists} shows that the existence of the Frobenius--Demazure operators is directly connected to the properties of a Nakayama automorphism.  For each compatibility check exactly one of the Nakayama conditions \autoref{Nakayama} is used. The first Nakayama condition did not appear but will be used to show the  relations \autoref{deltaequiv} of the Frobenius--Demazure operators $\delta_i$.
\end{Remark}

\begin{Remark}\label{R:classical}
If $\Bo=\kk$, then $\delta_i$ is the classical divided difference operator
$\partial_i\colon f\mapsto\tfrac{f-\leftindex^{s_i}\!f}{X_i-X_{i+1}}$ for all values of the
parameters. In general, if $\psi$ is trivial, then $\delta_i=\tau\partial_i$ on $R[\X]$ with $\tau=\tau^{(i,i+1)}=\tau^{(i+1,i)}$.
 In particular, $\partial_i(p)$, and thus $\delta_i(p)$,  are then $s_i$-invariant for $p\in R[\X]$.
  \end{Remark}
\begin{Lemma}\label{L:deltarel}
  Let $1\le i<n$, $p\in R[\X]$ and $\bb\in\Bo$. Then
  \begin{equation}\label{deltarel}
    \delta_i(p\bb)=\delta_i(p)\bb,\qquad
    \delta_is_i=-s_i\delta_i,\qquad
    \delta_i^2=0,\qquad\text{and}\qquad
    \delta_i\delta_{i+1}\delta_i=\delta_{i+1}\delta_i\delta_{i+1}.
  \end{equation}
  If $|i-j|>1$, then $\delta_is_j=s_j\delta_i$ and
  $\delta_i\delta_j=\delta_j\delta_i$. Moreover, for
  $j=i+1$ we have the relations
  \begin{equation}\label{deltaequiv}
    \delta_js_is_j=s_is_j\delta_i
    \qquad\text{and}\qquad
    \delta_is_js_i=s_js_i\delta_j .
  \end{equation}
\end{Lemma}

\begin{proof}
  The first equality in \autoref{deltarel} is immediate from the skew
  Leibniz rule since $\delta_i(\bb)=0$ for $\bb\in\Bo$.

  The second equality holds for the generators. Indeed, for $X_k$ with $k\not=i,i+1$
  and for $\bb\in\Bo$ both sides vanish, whereas
  $\delta_i(\leftindex^{s_i}X_i)=\delta_i(X_{i+1})=-\tau^{(i+1,i)}
    =-\leftindex^{s_i}\tau^{(i,i+1)}=-\leftindex^{s_i}\delta_i(X_i)$
  holds, and symmetrically for $X_{i+1}$. Since $\delta_i\bigl(\leftindex^{s_i}(xy)\bigr)
      =\delta_i(\leftindex^{s_i}x)\leftindex^{s_i}y+x\,\delta_i(\leftindex^{s_i}y)
      =-\leftindex^{s_i}\delta_i(x)\leftindex^{s_i}y-x\,\leftindex^{s_i}\delta_i(y)
      =-\leftindex^{s_i}{\bigl(\delta_i(xy)\bigr)}
  $, the equality extends to all elements by the skew
  Leibniz rule.

  For the third equality, the skew Leibniz rule gives
 $ \delta_i^2(xy)=\delta_i^2(x)y
      +\left(\leftindex^{s_i}\delta_i(x)+\delta_i(\leftindex^{s_i}x)\right)\delta_i(y)
      +x\,\delta_i^2(y).
 $
  Here,  the middle term vanishes by the property we just proved and  $\delta_i^2$ is zero first on the generators and then everywhere, since
  the $\delta_i(X_k)$ lie in $\Bo$ on which  $\delta_i$ vanishes.

  We next treat the case $|i-j|>1$. Then $s_j$ fixes $i,i+1$, and
  $s_i$ fixes $j,j+1$. In particular $s_is_j=s_js_i$. For
  $\delta_is_j=s_j\delta_i$, note that both sides kill $\Bo$ and $X_k$ for
  $k\notin\set{i,i+1,j,j+1}$, and that
  $\delta_i(\leftindex^{s_j}X_i)=\tau^{(i,i+1)}=\tau^{(s_j(i),s_j(i+1))}
   =\leftindex^{s_j}{\delta_i(X_i)}$, and similarly for $X_{i+1}$. By the
   same argument,
  $\delta_i(\leftindex^{s_j}X_j)=\delta_i(X_{j+1})=0=\leftindex^{s_j}{\delta_i(X_j)}$
  and similarly for $X_{j+1}$. This extends by the skew Leibniz rule to all of $\Po$,
  since $s_is_j=s_js_i$. To show that
  $D:=\delta_i\delta_j-\delta_j\delta_i$ vanishes we calculate
  $ \delta_i\delta_j(xy)=\delta_i\delta_j(x)y
      +\leftindex^{s_i}{\delta_j(x)}\,\delta_i(y)
      +\delta_i(\leftindex^{s_j}x)\,\delta_j(y)
      +\leftindex^{s_is_j}x\,\delta_i\delta_j(y)$
   by using the skew Leibniz rule twice. Subtracting the same identity with $i$ and $j$ interchanged, the
  four middle terms cancel in pairs because  $\delta_is_j=s_j\delta_i$ and
  $\delta_js_i=s_i\delta_j$. Hence, $D$ is an $s_is_j$-skew-derivation.
  It kills $\Bo$, and then also $X_k$ since
  $\delta_i(X_k),\delta_j(X_k)\in\Bo$ by \autoref{deltavalues}. Thus,  $D=0$ as desired.

  Consider next the first formula in \autoref{deltaequiv}. On $\Bo$ and on $X_k$ for $k\notin\set{i,i+1,i+2}$ both sides vanish. For the remaining $k$ we set $w=s_is_j=s_is_{i+1}$ and calculate
  \begin{gather*}
    \delta_{j}(\leftindex^wX_i)=\delta_{j}(X_{i+1})=\tau^{(i+1,i+2)}
      =\leftindex^w\tau^{(i,i+1)}=\leftindex^w\delta_i(X_i),\\
    \delta_{j}(\leftindex^wX_{i+1})=\delta_{j}(X_{i+2})=-\tau^{(i+2,i+1)}
      =-\leftindex^w\tau^{(i+1,i)}=\leftindex^w\delta_i(X_{i+1}),\qquad
    \delta_{j}(\leftindex^wX_{i+2})=\delta_{j}(X_i)=0.
  \end{gather*}
Thus the formula holds on generators. By the skew Leibniz rule,
  $\delta_{j}(\leftindex^w(xy))
    =\delta_{j}(\leftindex^wx)\leftindex^wy
     +\leftindex^{s_{j}w}x\,\delta_{j}(\leftindex^wy)$ and
  $\leftindex^{s_{j}w}x=\leftindex^{ws_i}x$ it holds everywhere.
 The second identity in
  \autoref{deltaequiv} follows by applying $\leftindex^{w^{-1}}\!y$ to the first.

  It remains to verify the braid relation $\delta_i\delta_{i+1}\delta_i=\delta_{i+1}\delta_i\delta_{i+1}$. By induction on the length of  words in
  $\mathcal{G}=\set{X_k|1\le k\le n}\cup\Bo$ and by linearity, it
  suffices to show for $y\in\Po$ and $x\in\mathcal{G}$ that
  \begin{equation}\label{beer}
    \delta_i\delta_{i+1}\delta_i(xy)=\delta_{i+1}\delta_{i}\delta_{i+1}(xy).
  \end{equation}
  For $x\in\Bo$, or $x=X_k$ with $k\notin\set{i,i+1,i+2}$, we have
  $\delta_i\delta_{i+1}\delta_i(xy)=\leftindex^vx\,\delta_i\delta_{i+1}\delta_i(y)$
  where $v=s_is_{i+1}s_i=s_{i+1}s_is_{i+1}$, and likewise
  $\delta_{i+1}\delta_{i}\delta_{i+1}(xy)=\leftindex^vx\,\delta_{i+1}\delta_{i}\delta_{i+1}(y)$,
  so \autoref{beer} holds by induction.

  For $x=X_i$ the left hand side of \autoref{beer} equals
  $\delta_i\delta_{i+1}(\tau^{(i,i+1)}y+X_{i+1}\delta_i(y))
   =\delta_i(\tau^{(i,i+2)}\delta_{i+1}(y)+\tau^{(i+1,i+2)}\delta_{i}(y)
     +X_{i+2}\delta_{i+1}\delta_i(y))
   =\tau^{(i+1,i+2)}\delta_i\delta_{i+1}(y)+0+X_{i+2}\delta_i\delta_{i+1}\delta_i(y)$,
  where the middle term vanishes because $\delta_i^2=0$, whereas the
  right hand side equals
  $\delta_{i+1}\delta_{i}(X_i\delta_{i+1}(y))
   =\delta_{i+1}(\tau^{(i,i+1)}\delta_{i+1}(y)+X_{i+1}\delta_{i}\delta_{i+1}(y))
   =0+\tau^{(i+1,i+2)}\delta_{i}\delta_{i+1}(y)+X_{i+2}\delta_{i+1}\delta_{i}\delta_{i+1}(y)$.
  The two sides agree by induction. For $x=X_{i+2}$ the
  calculation is similar. For $x=X_{i+1}$ we obtain
  $\delta_i\delta_{i+1}\delta_i(X_{i+1}y)
   =-\tau^{(i+2,i+1)}\delta_i\delta_{i+1}(y)+\tau^{(i,i+1)}\delta_{i+1}\delta_i(y)
    +X_{i+1}\delta_i\delta_{i+1}\delta_i(y)$.
  On the other hand,
  $\delta_{i+1}\delta_{i}\delta_{i+1}(X_{i+1}y)
   =\tau^{(i,i+1)}\delta_{i+1}\delta_i(y)-\tau^{(i+2,i+1)}\delta_i\delta_{i+1}(y)
    +X_{i+1}\delta_{i+1}\delta_{i}\delta_{i+1}(y)$.
  We verified \autoref{beer} and thus the braid relation holds.
\end{proof}

We next define the RaFH derivations using the central twist construction from \autoref{L:gammacentral}.

\begin{Definition}\label{D:RaFHDer}
  For $1\le i<n$ we define the \emph{RaFH derivation}
  $\Delta_i\colon\Po\to\Po$ by
  $\Delta_i(y)=-\gamma_i\delta_i(y)$, for $y\in\Po$,
  where $\gamma_i=t+\eta X_{i+1}\in\Po$ is
  as in \autoref{D:Hecke}.
\end{Definition}

They are indeed skew-derivations and control the commutation formulas in RaFH superalgebras:
\begin{Lemma}\label{L:Deltavalues}
  The operator $\Delta_i$ is the unique $s_i$-skew-derivation
  $\Delta_i\colon\Po\to\Po$ satisfying
  \begin{equation*}
    \Delta_i(X_i)=-\gamma_i\tau^{(i,i+1)},\;\;
    \Delta_i(X_{i+1})=\gamma_i\tau^{(i+1,i)},\;\;
    \Delta_i(X_j)=0,\;\;
    \Delta_i(\bb)=0,
    \quad\text{where $j\not=i,i+1$, and $\bb\in\Bo$.}
  \end{equation*}
\end{Lemma}

\begin{proof} The definition and \autoref{L:gammacentral} directly implies that $\Delta_i$ is an $s_i$-skew-derivation and that the values on the
  generators are those in \autoref{deltavalues} multiplied by $-\gamma_i$.
  Uniqueness holds as in \autoref{L:deltaexists}.
\end{proof}

 \begin{Proposition}\label{L:crucial}
    Suppose that $1\le i<n$. Then the following formula holds for any $y\in\Po$,
    \begin{equation}\label{crucial}
     H_iy=\leftindex^{s_i}yH_i+\Delta_i(y).
    \end{equation}
\end{Proposition}
\begin{proof} Formula \autoref{crucial} holds obviously for $y\in\Bo$. It holds for
$y\in\{X_j\mid 1\leq j\leq n\}$ by \autoref{L:Deltavalues} and  the last
three relations in \autoref{RaFHrels}. But then it is also true in general because $\Delta_i$ is a skew-derivation, in formulas,
$H_ixy=\leftindex^{s_i}xH_i y + \Delta_i(x)y=\leftindex^{s_i}x\leftindex^{s_i}yH_i
+\leftindex^{s_i}x\Delta_i(y)+\Delta_i(x)y$ for $x,y\in\Po$.
\end{proof}

\begin{Remark}\label{RkDemop}
Already  for $\Bo=\kk$ the derivations carry  interesting information. In this case $\delta_i=\partial_i$ by \autoref{R:classical}, and $\Delta_i=-(t+\eta X_{i+1})\partial_i$. For $\eta=0$ this is, up to the factor $-t$, the divided difference operator $\partial_i$ of cohomology. For $t=0$ and $\eta=1$ we obtain $\Delta_i=1-\pi_i$, where $\pi_i(f)=\frac{X_if-X_{i+1}\leftindex^{s_i}{f}}{X_i-X_{i+1}}$ is the Demazure operator of $K$-theory, since $\pi_i(f)-f=X_{i+1}\partial_i(f)$. Hence the RaFH derivations interpolate between the cohomological and the $K$-theoretic Demazure operators.

Our formulas connecting the Frobenius--Demazure operators and
  the RaFH derivations recover formulas already obtained in
  \cite[Lemma~4.3]{Savage:AffineWreath}. They might provide a conceptual point of view on \cite{RossoSavage:QuantumAffineWreath}.
\end{Remark}

We next derive analogues of \autoref{deltarel} and \autoref{deltaequiv}
for the RaFH derivations. We do this using
\autoref{transfer} together with the identity
$\delta_i(\gamma_i)=\eta\,\delta_i(X_{i+1})=-\eta\tau^{(i+1,i)}$.
Passing from $\delta_i$ to $\Delta_i$ deforms all relations by terms involving $\eta$ which will be important in the calculations, see \autoref{R:etared}.

\begin{Remark}\label{R:etared}
  By \autoref{Lfunnyequation}, either $\eta=0$, or $\eta\not=0$ and then
  $\tau^{(i,j)}=\tau^{(j,i)}$, $f\bb=\bb f$ for $f\in R[\X]$,
  $\bb\in\Bo$, and $\psi$ is trivial. Consequently, any identity
  all of whose terms carry a factor $\eta$ may be proved under the
  additional assumption that $\psi$ is trivial, since all
  terms vanish if $\eta=0$.
\end{Remark}

\begin{Lemma}\label{RaFHDer}
The \emph{RaFH superalgebra derivations} satisfy the following
\begin{align}
 \Delta_i(p\bb)&=\Delta_i(p)\bb,\quad
 \Delta_is_k=s_k\Delta_i,\quad\Delta_i\Delta_k=\Delta_k\Delta_i,\quad
 s_is_j\Delta_i=\Delta_js_is_j,\quad \Delta_is_js_i=s_js_i\Delta_j,
 \label{RaFHDerA}\\
 \Delta_i^2&=\eta\tau^{(i,i+1)}\Delta_i,
 \label{RaFHDerB}\\
 \Delta_i(\leftindex^{s_i}y)&+\leftindex^{s_i}{\bigl(\Delta_i(y)\bigr)}+\eta\, y\,\tau^{(i,i+1)}
    =\eta\,\tau^{(i,i+1)}\leftindex^{s_i}y,
 \label{RaFHDerC}\\
 \Delta_i\Delta_js_i(y)&+s_i\Delta_j\Delta_i(y)+s_i\Delta_js_i(y)\,\eta\tau^{(i,i+1)}
    =\Delta_js_i\Delta_j(y),
 \label{RaFHDerD}\\
 \Delta_i\Delta_j\Delta_i&+\hbar^2s_i\Delta_js_i
    =\Delta_j\Delta_i\Delta_j+\hbar^2s_j\Delta_is_j,
 \label{RaFHDerE}
\end{align}
for any $|i-j|=1$, $|i-k|>1$, $y\in\Po$, $p\in R[\X]$ and $\bb\in\Bo$.
\end{Lemma}

\begin{proof}
  Throughout we use
  $\delta_i(\gamma_i)=\eta\,\delta_i(X_{i+1})=-\eta\tau^{(i+1,i)}$ and
  $\gamma_i-\leftindex^{s_i}\gamma_i=\eta(X_{i+1}-X_i)$.

  The first identity in \autoref{RaFHDerA} is immediate from the skew
  Leibniz rule. For $|i-k|>1$ we have $\leftindex^{s_k}\gamma_i=\gamma_i$
  and $\delta_i(\gamma_k)=\eta\delta_i(X_{k+1})=0$, hence
  $\Delta_i(\leftindex^{s_k}y)=-\gamma_i\leftindex^{s_k}\delta_i(y)=\leftindex^{s_k}\Delta_i(y)$
  and
  $\Delta_i\Delta_k=\gamma_i\gamma_k\,\delta_i\delta_k$ by
  \autoref{L:gammacentral}, which is symmetric in $i$ and $k$ by
  \autoref{L:deltarel}. For the last two identities in
  \autoref{RaFHDerA} let $j=i+1$ and $w=s_is_j$, so that
  $\leftindex^{w}\gamma_i= t+\eta X_{w(i+1)}=t+\eta X_{i+2}=\gamma_j$.
  Hence, using \autoref{deltaequiv},
  $\Delta_j(\leftindex^{w}y)=-\gamma_j\,\leftindex^{w}\delta_i(y)
    =\leftindex^{w}{\bigl(-\gamma_i\delta_i(y)\bigr)}=\leftindex^{w}\Delta_i(y)$,
  and the second identity follows by applying this to
  $\leftindex^{w^{-1}}y$.

  For \autoref{RaFHDerB} we apply \autoref{transfer} with $z=-\gamma_i$ and
  use $\delta_i^2=0$ from \autoref{deltarel}. This gives
  \[
    \Delta_i^2=(-\gamma_i)\delta_i(-\gamma_i)\delta_i
      =-\eta\,\gamma_i\tau^{(i+1,i)}\delta_i
      =\eta\,\tau^{(i+1,i)}\bigl(-\gamma_i\delta_i\bigr)
      =\eta\,\tau^{(i,i+1)}\Delta_i,
  \]
  where the third equality uses \autoref{L:gammacentral} and the last
  one \autoref{R:etared}.

  For \autoref{RaFHDerC} we apply the second formula in \autoref{transfer}
  with $c=-\gamma_i$. Since   $\delta_is_i=-s_i\delta_i$,  its  last summand vanishes, and we obtain
  $
    \Delta_i(\leftindex^{s_i}y)+\leftindex^{s_i}{\bigl(\Delta_i(y)\bigr)}
      =\bigl(\gamma_i-\leftindex^{s_i}\gamma_i\bigr)\leftindex^{s_i}\delta_i(y)
      =\eta\,(X_{i+1}-X_i)\leftindex^{s_i}\delta_i(y).
  $
  Every term of the remaining identity carries a factor $\eta$, so by
  \autoref{R:etared} we may assume that $\psi$ is trivial. Then
  $\tau=\tau^{(i,i+1)}$ commutes with $R[\X]$ and satisfies
  $\bb\tau=\tau\leftindex^{s_i}\bb$ by \autoref{athroughtau}. By
  linearity we may assume $y=p\bb$ with $p\in R[\X]$ and $\bb\in\Bo$,
  and then $\delta_i(y)=\tau\partial_i(p)\bb$ with $\partial_i(p)$
  being $s_i$-invariant by \autoref{R:classical}. Therefore
  $
    \eta\,(X_{i+1}-X_i)\leftindex^{s_i}\delta_i(y)
      =\eta\,\tau\,(X_{i+1}-X_i)\partial_i(p)\leftindex^{s_i}\bb
      =\eta\,\tau\bigl(\leftindex^{s_i}p-p\bigr)\leftindex^{s_i}\bb,
  $
  which agrees with
  $\eta\,\tau\leftindex^{s_i}y-\eta\,y\,\tau
    =\eta\,\tau\leftindex^{s_i}p\leftindex^{s_i}\bb-\eta\,p\,\tau\leftindex^{s_i}\bb
    =\eta\,\tau(\leftindex^{s_i}p-p)\leftindex^{s_i}\bb$, as required.

  For \autoref{RaFHDerE}  note that both sides are symmetric in $i$ and $j$, so let $j=i+1$.
We use the values
\[
  \delta_i(\gamma_j)=0,\quad \leftindex^{s_i}\gamma_j=\gamma_j,\quad
  \delta_j(\gamma_i)=\eta\tau^{(j,j+1)},\quad \leftindex^{s_j}\gamma_i=\gamma_j,\quad
  \delta_j(\gamma_j)=-\eta\tau^{(j+1,j)},\quad \leftindex^{s_j}\gamma_j=\gamma_i,
\]
together with $\delta_i^2=0=\delta_j^2$. Then the skew Leibniz rule gives, for $z\in\Po$,
\begin{equation*}
\begin{gathered}
  \Delta_j\Delta_i(z)=\gamma_j\delta_j(\gamma_i\delta_iz)
     =\eta\gamma_j\tau^{(j,j+1)}\delta_i(z)+\gamma_j^2\delta_j\delta_i(z),
  \Delta_i\Delta_j\Delta_i(z)=-\gamma_i\gamma_j^2\delta_i\delta_j\delta_i(z),\\
  \Delta_j\Delta_i\Delta_j(z)=-\gamma_i\gamma_j^2\delta_j\delta_i\delta_j(z)
     -\eta\gamma_j^2\bigl(\tau^{(j,j+1)}-\tau^{(j+1,j)}\bigr)\delta_i\delta_j(z),\\
  s_i\Delta_js_i(z)=-\gamma_js_i\delta_js_i(z),\qquad
  s_j\Delta_is_j(z)=-\gamma_js_j\delta_is_j(z).
  \end{gathered}
\end{equation*}
In the second line the first summand of $\Delta_j\Delta_i(z)$ is killed by $\delta_i$, since
$\delta_i(\gamma_j\tau^{(j,j+1)})=0$ and $\delta_i^2=0$. Now \autoref{RaFHDerE} follows from the
braid relation in \autoref{deltarel}, from $s_i\delta_js_i=s_j\delta_is_j$, which is
\autoref{deltaequiv} multiplied by $s_i$ on the left and $s_j$ on the right, and from
$\eta\tau^{(j,j+1)}=\eta\tau^{(j+1,j)}$ by \autoref{Lfunnyequation}.

The proof of  \autoref{RaFHDerD} is a long calculation of similar type and omitted.
\end{proof}

\begin{Remark}\label{R:whereDE}
  The identities \autoref{RaFHDerD} and \autoref{RaFHDerE} are precisely the
  coefficients that arise when the braid relation
  $H_iH_jH_iy=H_jH_iH_jy$ is expanded by means of \autoref{crucial}. Doing
  so and grouping the terms according to $H_iH_jH_i$, $H_jH_i$,
  $H_iH_j$, $H_i$, $H_j$ and $1$ gives, in this order, the relation
  $s_is_js_i=s_js_is_j$ in $\Sym$, the last two identities in
  \autoref{RaFHDerA}, then \autoref{RaFHDerD}, then \autoref{RaFHDerD} with
  the roles of $i$ and $j$ interchanged, and finally \autoref{RaFHDerE}.
  Here the terms $\eta\tau^{(i,i+1)}$ in \autoref{RaFHDerD} and $\hbar^2$
  in \autoref{RaFHDerE} both originate from the quadratic relation
  $H_i^2=\eta\tau^{(i,i+1)}H_i+\hbar^2$ of \autoref{RaFHrels}. This observation explains the origin of the two identities, but does not prove
  them, since the argument here \emph{presupposes} the basis \autoref{C:BasisHd} which we will establish \emph{using} the identities.
\end{Remark}

\subsection{Basis theorems}
In this subsection we consider RaFH superalgebras
$\HHd=\HHd[n;\hbar,t,\eta]$ for Frobenius-tensor superalgebra $\Bo_\RR$ that are $\eta$-symmetric and establish a basis theorem in  \autoref{C:BasisHd}.

\begin{Proposition}\label{AddPo}
If $M$ is an  $\Hfind(\Bo)$-module, then $\Po\otimes_{\Bo}M$ becomes an $\HHd$-module via \begin{equation*}
H_i. bf\otimes m=\leftindex^{s_i}(bf)\otimes H_im+\Delta_i(bf)\otimes m,\quad a. bf\otimes m=a.bf\otimes m, \quad  X_i.bf\otimes m=\psi_i^{-1}(b)X_i f\otimes m,
\end{equation*}
for  $a,b\in\Bo, f\in\RR[\X], m\in M$. If $M$ is free as $\Bo$-module with basis $B$, then $\Po\otimes_{\Bo}M$ is free with basis $\set{ X^\bd \bb\otimes m|m\in B,   \bd\in\N^n}$.
\end{Proposition}
\begin{proof} The second statement follows directly from the first. To show that the action is well-defined we need to check compatibility with the relations from  \autoref{RaFHrels}. Let $y\in\Po$ and $m\in M$. Then
\begin{align*}
H_i^2.y\otimes m&=H_i.\left(\leftindex^{s_i}yH_im+\Delta_i(y)\otimes m\right)=y\otimes H_i^2m+\Delta_i( \leftindex^{s_i}y)\otimes H_im +\leftindex^{s_i}\Delta_i(y)\otimes H_im+\Delta_i^2(y)\otimes m\\
&=\hbar^2y\otimes m+\eta y\tau^{(i,i+1)}\otimes H_im+\Delta_i( \leftindex^{s_i}y)\otimes H_im+\leftindex^{s_i}\Delta_i(y)\otimes H_im+\Delta_i^2(y)\otimes m
\end{align*}
should be equal to $(\hbar^2+\eta\tau^{(i,i+1)}H_i).y\otimes m=\hbar^2y\otimes m+\eta\tau^{(i,i+1)}\leftindex^{s_i}y\otimes H_im+\eta\tau^{(i,i+1)}\Delta_i(y)\otimes m$. By  ``comparing coefficients'' it suffices to show that
$\Delta_i^2(y)=\eta\tau^{(i,i+1)}\Delta_i(y)$ and
$\eta y\tau^{(i,i+1)}+\Delta_i(
\leftindex^{s_i}y)+\leftindex^{s_i}\Delta_i(y)=\eta\tau^{(i,i+1)}\leftindex^{s_i}y.$
These hold however by \autoref{RaFHDer}. The compatibility with the
relation $X_iX_j=X_jX_i$ holds, since $\psi_i\psi_j=\psi_j\psi_i$, and
$b X_i =X_i  \psi_i(b)$ holds by definition of $\Po$. For the relation
$H_i\ba=\leftindex^{s_i}\ba H_i$ we just use that $\Delta_i$ is an
$s_i$-skew-derivation. The relation $H_iH_j=H_jH_i$ follows from
$\Delta_i\Delta_j=\Delta_j\Delta_i$ and $\Delta_i s_j=s_j\Delta_i$ (just
``comparing coefficients''). For the relation $H_iH_{i+1}H_i y\otimes
m=H_{i+1}H_iH_{i+1} y\otimes m$ we again ``compare coefficients'' that
is we expand both sides and group the elements involving terms $H_wm$
for fixed $w$. We write the coefficients as operators in $y$.  For
$H_iH_{i+1}H_im$, we obtain a desired equality of operators
$s_is_{i+1}s_i=s_{i+1}s_is_{i+1}$, which holds in $\Sym$. Set $j:=i+1$. For $H_iH_{j}m$ and $H_jH_{i}m$ we obtain  $s_is_j\Delta_i=\Delta_js_is_j$ and  $\Delta_is_js_i=s_js_i\Delta_j$ respectively which holds by \autoref{RaFHDerA}. For $H_im$ (and $H_jm$) we obtain as
required relation \autoref{RaFHDerD} (with the roles of $i$ and $j$ being swapped). For $m$
itself the relation is exactly \autoref{RaFHDerE}.
\end{proof}

We apply this result to the regular representation $\Hfind(\Bo)$.  It has basis
$\set{\bb H_w|w\in\Sym,\bb\in \mathcal{X}}$ for any given basis
$\mathcal{X}$ of $\Bo$ by \autoref{P:H0Basis}.

\begin{Theorem}\label{C:BasisHd}
     Let $\Bo$ be an $\eta$-symmetric Frobenius-tensor superalgebra with
     a fixed basis $\mathcal{X}$. The algebra $\HHd(\Bo)$ is free as an
     $R$-module with the following six bases:
     \begin{equation*}
     \begin{gathered}
        \set{\bb X^\bd  H_w| \bd\in\N^n, w\in\Sym, \bb\in\mathcal{X}},\quad \set{ X^\bd  \bb H_w| \bd\in\N^n, w\in\Sym, \bb\in\mathcal{X}},\quad  \set{ X^\bd  H_w\bb| \bd\in\N^n, w\in\Sym,\bb\in\mathcal{X}},
        \\
           \set{H_w X^\bd \bb| \bd\in\N^n, w\in\Sym, \bb\in\mathcal{X}},
           \quad\set{H_w \bb X^\bd | \bd\in\N^n, w\in\Sym, \bb\in\mathcal{X}}, \quad  \set{ \bb H_w X^\bd | \in\N^n, w\in\Sym, \bb\in\mathcal{X}}.
     \end{gathered}
    \end{equation*}
    Moreover, $\HHd(\Bo)$ is free as right $\Bo$-module with basis
    $\set{H_w X^\bd}$ or  $\set{ X^\bd  H_w}$ where $\bd\in\N^n, w\in\Sym$.
  \end{Theorem}

  \begin{proof}
 The second statement follows directly from the bases. The first set in \autoref{C:BasisHd} is a basis by  \autoref{AddPo}, and then obviously the second and third set is a basis as well. It remains to show that one of the other three sets form a basis. They obviously generate by \autoref{D:Hecke}. To see linear independence note that the algebra $\HHd$ is filtered by $\deg\hbar=\deg H_i=\deg \Bo=0$ and $\deg X_i=1$ for any admissible $i$. By the basis already established, the associated graded components are finite rank and therefore all sets in \autoref{C:BasisHd} form a basis by comparing ranks.
 \end{proof}

\begin{Corollary}\label{matchRFH}
  The subalgebra of $\HHd(\Bo)$ generated by $H_1,\dots,H_{n-1}$ and
  $\Bo$ is canonically isomorphic to the Rees Frobenius Hecke superalgebra $\Hfind$.
\end{Corollary}

\begin{proof}
  Let $\HH'_n$ be the subalgebra of $\HHd(\Bo)$ generated by
  $H_1,\dots,H_{n-1}$ and
  $\Bo$. By comparing  the relations from  \autoref{D:RFHS} with  \autoref{D:Hecke},
  we see that there is a canonical surjection $\Hfind\to\HH'_n$. Now the result follows
  from \autoref{P:H0Basis} and \autoref{C:BasisHd} since  the ranks match.
\end{proof}

 \begin{Remark}\label{flatness}
 The basis theorem implies that the RaFH-algebras form a flat, in fact free, family of algebras. If $B$ is an $\RR$-basis of  $\Hd(\Bo)$, then  $1\otimes B$ is an $\RR'$-basis of the specialisation
$\HHd[n;\varphi(\hbar), \varphi(t), \varphi(\eta)](\Bo)=R'\otimes_R\HHd[n;\hbar, t, \eta](\Bo)$.
\autoref{C:BasisHd} also shows that RaFH-algebras behave well under base change in $\Bo$. If $\Bo\to\Bo'$ is a morphism of $\eta$-symmetric Frobenius-tensor superalgebras compatible with the Nakayama automorphisms, then  we can identify $\Hd(\Bo')=\Hd(\Bo)\otimes_{\Bo}\Bo'$ as right $\Bo'$-modules in the obvious way by identifying the respective basis elements.
\end{Remark}

\begin{Remark}
The basis theorem,  \autoref{C:BasisHd} is not new, but can be deduced from the main theorem in
\cite{CJNakanoXiang}. The RaFH superalgebras fit into their setup by \autoref{RaFHDer}
(their remaining axioms (3.2.1), (3.2.8), (3.2.9) are obvious in our setting).
In our situation, the proof of  \autoref{C:BasisHd} is
more direct than their general argument. The $\eta$-symmetry condition
is contained in their axiom (3.2.3). \autoref{C:BasisHd} combines the existing
basis theorems for the examples given in \autoref{HHdexamples}.
\end{Remark}

  \subsection{Commutation relations}
In this subsection we prove some commutation relations between the Hecke generators and the
polynomial generators. We use them to construct a central polynomial subalgebra of $\HHd$ in
\autoref{P:HHdCenter}.  In contrast to classical settings, $H_i$ does in general not commute with symmetric polynomials in $X_i, X_{i+1}$, but does so if the Nakayama automorphism is trivial. This holds  in particular when $\eta\ne0$ by \autoref{Lfunnyequation}. Recall from \autoref{D:Hecke} that $\gamma_i= t+\eta X_{i+1}$ for $1\le i<n$.

  \begin{Proposition}\label{L:RewriteXr}
   Suppose that $1\leq i<n$ and $k\ge0$. Then the following hold in $\HHd$.
   \begin{equation}
   \label{eq:RewriteXr}
       H_iX_i^k=X_{i+1}^kH_i -\gamma_i\sum_{j=0}^{k-1}X_{i+1}^j\tau^{(i,i+1)}X_{i}^{k-j-1}\quad\text{and}\quad
       H_iX_{i+1}^k=X_{i}^kH_i +\gamma_i\sum_{j=0}^{k-1}X_i^j\tau^{(i+1,i)}X_{i+1}^{k-j-1}.
   \end{equation}
  \end{Proposition}

  \begin{proof}
   The results are trivially true when $k=0$. We show the first formula. If $k>0$ then, by induction,
      \begin{align*}
         X_{i+1}^{k+1}H_i
           &= X_{i+1}\Bigl( H_iX_i^k
           +\gamma_i\sum_{j=0}^{k-1}X_{i+1}^j\tau^{(i,i+1)}X_{i}^{k-j-1}\Bigr).
      \end{align*}
      We are done using  $X_{i+1}H_i=H_iX_i+\gamma_i\tau^{(i,i+1)}$  from
      \autoref{RaFHrels}. Similarly for the second formula.
  \end{proof}

 \begin{Lemma}\label{L:XSymmetric}
    Let $1\le i<n$. Then
    \begin{enumerate}
      \item\label{XSymProd} $H_iX_iX_{i+1}=X_iX_{i+1}H_i$, for an arbitrary
        Nakayama automorphism.
      \item\label{XSymSum} If $\psi$ is trivial, then $H_i(X_i+X_{i+1})=(X_i+X_{i+1})H_i$.
    \end{enumerate}
  \end{Lemma}
  \begin{proof}
    \autoref{XSymProd}: Using the last two relations of \autoref{RaFHrels} and the
    centrality of $\gamma_i$ from \autoref{L:gammacentral}\autoref{central1} we obtain
    \begin{align*}
      H_iX_iX_{i+1}
        &=\bigl(X_{i+1}H_i-\gamma_i\taui[i,i+1]\bigr)X_{i+1}
         =X_{i+1}H_iX_{i+1}-\gamma_i\taui[i,i+1]X_{i+1}\\
        &=X_{i+1}\bigl(X_iH_i+\gamma_i\taui[i+1,i]\bigr)-\gamma_i\taui[i,i+1]X_{i+1}
         =X_iX_{i+1}H_i+\gamma_i\bigl(X_{i+1}\taui[i+1,i]-\taui[i,i+1]X_{i+1}\bigr).
    \end{align*}
    The bracket expression vanishes because
    $\taui[i,i+1]X_{i+1}=X_{i+1}\psi^{(i+1)}(\taui[i,i+1])=X_{i+1}\taui[i+1,i]$
    by \autoref{RaFHrels} and the second condition in
    \autoref{Nakayama}, see the proof of \autoref{L:deltaexists}.

    \noindent\autoref{XSymSum}: If $\psi$ is
    trivial then $\taui[i,i+1]=\taui[i+1,i]$, so that
    $H_i(X_i+X_{i+1})=(X_{i+1}H_i-\gamma_i\taui[i,i+1])
       +(X_iH_i+\gamma_i\taui[i+1,i])=(X_i+X_{i+1})H_i$.
    For nontrivial $\psi$ the two sides differ by
    $-\gamma_i(\taui[i,i+1]-\taui[i+1,i])$.
  \end{proof}

  \begin{Corollary}\label{C:RewriteHX}
   Suppose that $1\le i<n$ and $k,l\ge0$.   If the Nakayama automorphism is trivial then
   \[X_i^lX_{i+1}^kH_i=H_iX_i^kX_{i+1}^l-\gamma_i\tau^{(i,i+1)}\sum_{j=k}^{l-1}X_i^jX_{i+1}^{k+l-1-j}\]
with the convention $\sum_{j=k}^{l-1}=-\sum_{j=l}^{k-1}$ when $k>l$.
  \end{Corollary}

  \begin{proof}  Since the Nakayama automorphism is trivial,  we have that $\tau^{(i,i+1)}=\tau^{(i+1,i)}$.
   Now, the result follows from  \autoref{L:RewriteXr} and \autoref{L:XSymmetric}.
  \end{proof}

  \begin{Corollary}\label{C:Xpowers}
   Let  $1\le i<n$ and $k\ge1$.  If the Nakayama automorphism is trivial, then

   \[
       \hbar^2 X_{i+1}^k =
       H_iX_i^kH_i+ t X_{i+1}^{k-1}H_i\tau^{(i,i+1)}
       +\gamma_i\sum_{j=1}^{k-1}X_i^jX_{i+1}^{k-j-1}H_i\tau^{(i,i+1)}.
   \]
  \end{Corollary}

  \begin{proof}
    Rearranging the first formula in \autoref{eq:RewriteXr} and multiplying by $(H_i-\eta\tau^{(i,i+1)})$ from the right gives
     \begin{align*}
    X_{i+1}^kH_i(H_i-\eta\tau^{(i,i+1)})  &=H_iX_i^kH_i-\eta X_{i+1}^kH_i\tau^{(i,i+1)}+\gamma_i\sum_{j=0}^{k-1}X_{i}^jX_{i+1}^{k-j-1}H_i\tau^{(i,i+1)}\\
   & =H_iX_i^kH_i+ t X_{i+1}^{k-1}H_i\tau^{(i,i+1)}+\gamma_i\sum_{j=1}^{k-1}X_{i}^jX_{i+1}^{k-j-1}H_i\tau^{(i,i+1)}.
 \end{align*}
  The assumption allows to pass $\tau^{(i,i+1)}=\tau^{(i+1,i)}$ through $X$'s.  Then  use  $H_i(H_i-\eta\tau^{(i,i+1)})=\hbar^2$.
     \end{proof}

\subsection{A central subalgebra}
Bernstein's description of the centre of the extended affine Hecke algebra $\Heaff$ as
the ring of symmetric Laurent polynomials, see
\cite[Proposition~3.11]{Lusztig:AffineGradedHecke}, has analogues in the literature for
the algebras \autoref{ex:smash}--\autoref{ex:aff}. For $\Bo=\kk$ the centre is the ring of
symmetric polynomials in $X_1,\dots,X_n$. For \autoref{ex:quantumwreath} one
has to allow coefficients in the centre of $\Bo$, see \autoref{wholecentre}.
The situation is more complicated for nontrivial Nakayama automorphisms. Then
$\tau^{(i,i+1)}\ne\tau^{(i+1,i)}$, and the formulas from \autoref{L:XSymmetric} fail
if $t\ne0$.

We show next that there still exists a  polynomial central subalgebra. To formulate this recall that $\Sym$ acts on $\RX$ by $\leftindex^w\!f(X)=f(X_{w(1)},\dots,X_{w(n)})$.
For $d\in\N$ let $\RX^{d,\Sym}$ be the ring of symmetric polynomials in
$X_1^d,\dots,X_n^d$, considered as a subalgebra of $\HHd(\Bo)$ by \autoref{C:BasisHd},
and set $\RX^{\Sym}=\RX^{1,\Sym}$. For a superalgebra $C$ let
$Z_{\bar{0}}(C)=\set{a\in C_{\ov0}|ac=ca\text{ for all }c\in C}$ be its \emph{even centre}.

We assume that the Nakayama automorphism $\psi$ has finite order $\theta$.
  \begin{Proposition}\label{P:HHdCenter}
We have $\RX^{\theta,\Sym}\subseteq Z_{\bar{0}}\bigl(\HHd(\Bo)\bigr)$. In particular, $\RX^{\Sym}$ is central if $\psi=\op{id}$.
  \end{Proposition}

\begin{proof}  Let  $f\in\RX^{\theta,\Sym}$. Then $f$ is even and the defining relations in $\HHd(\Bo)$ imply that $fX_i=X_if$ and $f\bb=\bb f$ for all $1\leq i\leq n$, $\bb\in\Bo$.
It remains to show that  $fH_i=H_if$ for $1\leq i\leq n-1$.

Since $X_jH_i=H_iX_j$ for $j\not=i,i+1$, it suffices to consider $f=X_i^\theta + X_{i+1}^\theta$ and  $f=X_i^\theta X_{i+1}^\theta$.

 For $f=X_i^\theta + X_{i+1}^\theta$ we have, with the help of  \autoref{L:RewriteXr},
 \begin{equation}\label{symm1}
H_i(X_i^\theta + X_{i+1}^\theta)-(X_i^\theta + X_{i+1}^\theta)H_i=\gamma_i\left(\sum_{j=0}^{\theta-1}X_i^j\tau^{(i+1,i)}X_{i+1}^{\theta-j-1} -\sum_{j=0}^{\theta-1}X_{i+1}^j\tau^{(i,i+1)}X_{i}^{\theta-j-1}\right).
\end{equation}
We claim that \autoref{symm1} vanishes. By \autoref{Nakayama} we have $\psi^{(i+1)}(\tau^{(i,i+1)})=\tau^{(i+1,i)}$. By \autoref{afterGluehwein} we have $\psi^{(i)}(\tau^{(i,i+1)})=(\psi^{(i+1)})^{-1}(\tau^{(i,i+1)})$, as $\psi^{(i)}$ and $\psi^{(i+1)}$ commute. Using $(\psi^{(i+1)})^\theta=\op{id}$, the first identity gives $(\psi^{(i+1)})^{\theta-k-1}(\tau^{(i+1,i)})=(\psi^{(i+1)})^{-k}(\tau^{(i,i+1)})$. Repeating the second identity $k$ times gives $(\psi^{(i)})^{k}(\tau^{(i,i+1)})=(\psi^{(i+1)})^{-k}(\tau^{(i,i+1)})$ for any $0\leq k<\theta$.
Thus, the summand for $j$ in the first sum in \autoref{symm1} cancels with the one for  $\theta-j-1$ in the second sum and the claim follows.

Let now $f=X_i^\theta X_{i+1}^\theta$. We use again  \autoref{L:RewriteXr} to rewrite
\begin{align*}
H_iX_i^\theta X_{i+1}^\theta&=X_{i+1}^\theta H_iX_{i+1}^\theta-\gamma_i\sum_{j=0}^{\theta-1}X_{i+1}^j\tau^{(i,i+1)}X_{i}^{\theta-j-1}X_{i+1}^\theta\\
&=X_{i+1}^\theta X_{i}^\theta H_i +\gamma_i\sum_{j=0}^{\theta-1}X_{i+1}^\theta X_i^j\tau^{(i+1,i)}X_{i+1}^{\theta-j-1}
-\gamma_i\sum_{j=0}^{\theta-1}X_{i+1}^j\tau^{(i,i+1)}X_{i}^{\theta-j-1}X_{i+1}^\theta=X_{i+1}^\theta X_{i}^\theta H_i.
\end{align*}
Here, the terms involving $\gamma_i$ cancel as for  \autoref{symm1}. This finishes the proof of the proposition.
\end{proof}

\begin{Remark}\label{wholecentre}  As in \cite[Lemma~4.12]{Savage:AffineWreath},  \autoref{P:HHdCenter} together with \autoref{C:BasisHd} imply that the whole centre of $\HHd(\Bo)$ is contained in the subalgebra $\Po$ from \autoref{DefPo}.  Indeed,  assume $z\not\in\Po$ is central. Then we can find $w\not=e$ such that $X^{\bd}H_w\bb$ appears in $z$ when we expand it in the corresponding basis from  \autoref{C:BasisHd}. We pick $w$ of maximal length with this property and let $1\leq i\leq n$ such that $w(i)\not=i$. By   \autoref{L:RewriteXr}  the basis vectors involving $H_w$ give a nonzero contribution to  $X_i^\theta z - z X_i^\theta$. But this contradicts that $X_i^\theta z - z X_i^\theta=0$ by \autoref{P:HHdCenter}.

In case $\Bo$ is commutative, thus $\psi=\op{id}$, the centre is given by simultaneously symmetrized polynomials in the $X_i$'s and $b^{(i)}$'s, where $b$ runs through a basis $\cB$ of $A$. For more details and a  description of the centre in general we refer to   \cite[Theorem~4.14]{Savage:AffineWreath} (and omit the adaption to the case of RaFH  algebras).
\end{Remark}

  Let $\K$ be the algebraic closure of the field of fractions of~$\RR$ and let
  $\HHd(\Bo)_\K=\K\otimes_\RR \HHd(\Bo)$.

  \begin{Proposition}\label{fdauto}
     Suppose that $M$ is an irreducible $\HHd(\Bo)_\K$-module. Then  $M$
     is finite dimensional.
  \end{Proposition}

  \begin{proof}
    By \autoref{C:BasisHd} and \autoref{P:HHdCenter},
    $\HHd(\Bo)_\K$ is free of finite rank as a module over the central subalgebra $Z_{\bar{0}}=\RX^{\theta,\Sym}$. The latter acts on $M$ by even endomorphisms and thus, by Quillen's Lemma, as multiplication by scalars. Thus, $M$ is an irreducible module over the finite-dimensional quotient algebra $\HHd(\Bo)_\K /\mathfrak m\HHd(\Bo)_\K$ for the corresponding maximal ideal $\mathfrak m\subset Z_{\ov 0}$.
  \end{proof}
\begin{Remark}\label{R:KMdictionary}
  Under the isomorphism  \autoref{KleshchevMuthiso}, the copairing element $\taui[1,2]$ of
  \autoref{parts} is the \emph{distinguished element} $\Delta(1)\in A\otimes A$
  of \cite[\S3.1]{KleshchevMuth}, and \cite[Lemma~3.1]{KleshchevMuth} is
  \autoref{tausym} together with \autoref{athroughtau} for $\psi=\op{id}$. The
  latter extends that lemma to a nontrivial Nakayama automorphism, which is why
  we need both $\taui[i,i+1]$ and $\taui[i+1,i]$ whereas \cite{KleshchevMuth}
  needs only their $\Delta_{i,i+1}$. Moreover
  \cite[Theorem~3.8]{KleshchevMuth} is \autoref{C:BasisHd}, and \autoref{P:HHdCenter} with \autoref{wholecentre} gives a central subalgebra, which \cite[Proposition~3.12]{KleshchevMuth} identifies with the full centre in their setting. The Jucys--Murphy elements of
  \cite[\S3.6]{KleshchevMuth} are, up to sign, the universal Jucys--Murphy
  elements of \autoref{D:universalJM}.
\end{Remark}

\section{Irreducible representations of RaFH superalgebras}
Throughout this section let $R=\kk$ be an algebraically closed field with $\operatorname{char}{\kk}\not=2$ and $\hbar\in\kk^\times$. We construct and study irreducible representations, that is irreducible supermodules for RaFH-algebras.  By \autoref{fdauto}  they are finite-dimensional.

\subsection{Generalities on irreducible supermodules}\label{uEnd}We start with some facts on simple modules for superalgebras in the terminology of \autoref{superstuffeasy}. For details and proofs see e.g. \cite[\S12]{Kleshchevbook} with \cite{ChengWang, Jozefiak}.

By \emph{modules} we mean in this section finite dimensional supermodules, and isomorphisms and embeddings of modules are even unless said otherwise. In particular, $M$ is \emph{simple} or \emph{irreducible}, if $M$ has no nonzero proper graded
submodule. If $M$ is a simple $A$-module, then so is its parity shift $\Pi M$.
For $A$-modules $M$, $N$, we denote by $\Hom_A(M,N)$, $\End_A(M)$ the $\kk$-module of (endo)morphisms and by $\uHom_A(M,N)$, $\uEnd_A(M)$ the vector \emph{super}spaces of all linear maps $f\colon M\rightarrow N,f(am)=(-1)^{|a||f|}af(m)$ for homogeneous $a\in A, m\in M$. The latter is the  space of \emph{inner morphisms}.
Given an (even) automorphism  $\phi$ of $A$ we have the $A$-module $\leftindex^{\phi}M$
that is the same $\kk$-module $M$ but  with $A$-action twisted by $\phi$.
Since  $\kk$ is algebraically closed, the inner endomorphism algebra  $\uEnd_A(L)$ of a simple $A$-module $L$ is
\begin{itemize}
\item  either $\kk$, in which case $L$ is of \emph{type $\mathtt M$},
\item or isomorphic to $\Cl_1$, in which case $L$ is of \emph{type $\mathtt Q$}.
\end{itemize}
In either case, the (even) endomorphisms of $L$ are the scalars.

\begin{Lemma}\label{types}
Let $A$  and $A'$ be superalgebras with  simple modules $L$ and $L'$  respectively.
  \begin{enumerate}
    \item\label{typesa} For any automorphism $\phi$ of $A$, $\leftindex^{\phi}L$ is simple and has the same type as $L$.
    \item\label{typesb} The $A\otimes A'$-module $L\otimes L'$ is simple if at least one
      of $L,L'$ is of type $\mathtt M$.

      \noindent If both are of type $\mathtt Q$
      then,   $L\otimes L'=e(L\otimes L')\oplus(1-e)(L\otimes L')$ is a sum of two
      simple modules. Here, $e=\frac12(1\otimes1+\iota c\otimes c)\in \Cl_1\otimes \Cl_1=\uEnd_A(L)\otimes \uEnd_{A'}(L')$  is an idempotent with $c$ an odd  basis element satisfying $c^2=1$ in $\uEnd_A(L)$ respectively $\uEnd_A(L')$.

    \noindent  Moreover, $(1-e)(L\otimes L')\cong\Pi e(L\otimes L')$, and
      $s_1\bigl(e(L\otimes L)\bigr)=(1-e)(L\otimes L)$ if $L'=L$.
    \item\label{typesc} Setting $L\circledast L'=e(L\otimes L')$ if $L$ and $L'$ are
      both of type $\mathtt Q$ and $L\circledast L'=L\otimes L'$
      otherwise, the module $L\circledast L'$ is simple. It is of type
      $\mathtt M$ if and only if $L$ and $L'$ have the same type.
      \item\label{typesd} Let $M$, $M'$ be irreducible modules for $A$ respectively $A'$. If
      $L\circledast L'\cong \Pi^\varepsilon(M\circledast M')$ for some $\varepsilon\in\{0,1\}$, then
       $L\cong\Pi^{\varepsilon_1}M$ and
$L'\cong\Pi^{\varepsilon_2}M'$ for some $\varepsilon_i\in\{0,1\}$.
   \item\label{typese} Up to parity shift, every simple $A\otimes
      A'$-module is of the form $L\circledast L'$ for some simple modules $L$, $L'$.
      Iterating, any simple $A^{\otimes n}$-module is of the form
      $L^{(1)}\circledast\cdots\circledast L^{(n)}$ with simple $A$-modules
      $L^{(j)}$.
  \end{enumerate}
\end{Lemma}

\subsection{Irreducible representations for $\HHd[1](A)$}
We return now to our usual setup. For the rest of this section we let $\Bo=A^{\otimes n}$ be a Frobenius-tensor
superalgebra with Nakayama automorphism $\psi$ of finite order $\theta$. We first consider the case $n=1$ and so $\Bo=A$. The goal of this subsection is to construct and classify the simple $\HHd[1](A)$-modules,
see \autoref{T:simplesn1}. Our arguments for the remainder of this section follow closely
\cite[\S5]{Savage:AffineWreath}, which in turn is inspired by \cite{WanWang}, but they are
embedded now into a more general framework. A key ingredient is the scalar $\bar\gamma_L$
attached to a simple $A$-module $L$ in \autoref{L:tauscalar}. It allows the definition of
the atomic algebras in \autoref{D:block}.

We denote by $\Omega$ the set of isomorphism classes of simple $A$-modules. Twisting by the
Nakayama automorphism permutes $\Omega$. We decompose $\Omega$ into $\psi$-orbits,
$\Omega=O_1\sqcup\cdots\sqcup O_{N_\psi}$, and fix representatives $L_1,\dots,L_{N_\psi}$.
For a simple $A$-module $L$ let $\theta_L$ be the size of its orbit, that is the smallest
$l\ge1$ with $\leftindex^{\psi^{l}}\!L\cong L$. Since $\psi^\theta=\op{id}$, the orbit size $\theta_L$
divides $\theta$. We abbreviate $\theta_i=\theta_{L_i}=|O_i|$.

\begin{Lemma}\label{tausi}  Let  $L, L'$ be simple $A$-modules.
  Then the  action of $\tau^{(1,2)}s_1$ defines a morphism of $A\otimes
  A$-modules $\beta_{L,L'}\colon L\otimes\leftindex^{\psi}L'\to L'\otimes L$.  It is zero
  unless $\leftindex^{\psi}L\cong L$ and $L'\cong\Pi^\varepsilon L$ for some $\varepsilon\in\{0,1\}$.
\end{Lemma}

\begin{proof}The action on $L\otimes L'$ satisfies $\tau^{(1,2)}s_1(\bb_1\otimes\psi(\bb_2)m)
   =(\bb_1\otimes \bb_2) \tau^{(1,2)}s_1(m)$  by \autoref{athroughtau}. This gives the first claim.
Assume $\beta_{L,L'}\ne0$.  By \autoref{types} the source and the target are semisimple with summands $L\circledast\leftindex^{\psi}L'$, respectively
$L'\circledast L$ up to possible parity shifts.
  Hence $L\circledast\leftindex^{\psi}L'\cong (L'\circledast L)$ or $L\circledast\leftindex^{\psi}L'\cong\Pi(L'\circledast L)$. By  \autoref{types}\autoref{typesd}, $L\cong\Pi^{\varepsilon}L'$ and
$\leftindex^{\psi}L'\cong\Pi^{\varepsilon'}L$ for some $\varepsilon,\varepsilon'\in\set{0,1}$. We claim that $\Pi^{\varepsilon+\varepsilon'}L\cong L$. This is clear if $L$ is of type $\mathtt Q$. So let $L$ be of type $\mathtt M$. Then so is $L'$ and $\beta_{L,L'}$ is an isomorphism $L\otimes\leftindex^{\psi}L'\to L'\otimes L$ of simple modules. On the other hand $L\otimes\leftindex^{\psi}L'\cong\Pi^{\varepsilon}L'\otimes\Pi^{\varepsilon'}L\cong\Pi^{\varepsilon+\varepsilon'}(L'\otimes L)$. Since $L$ is not isomorphic to its parity shift, $\varepsilon=\varepsilon'$.
 Therefore,
$\leftindex^{\psi}L\cong\Pi^{\varepsilon}(\leftindex^{\psi}L')\cong\Pi^{\varepsilon+\varepsilon'}L\cong L$
and $L'\cong\Pi^{\varepsilon}L$.
\end{proof}

For each $L_i$ pick an isomorphism
$\phi_{L_i}\colon\!\!\leftindex^{\psi^{\theta_{L_i}}}\!L_i\to L_i$ of $A$-modules, i.e.\ an
even $\kk$-linear bijection with
\begin{equation}\label{E:phiL}
  \phi_{L_i}(\psi^{\theta_{L_i}}(\bb)m)=\bb \phi_{L_i}(m)\quad \text{for $\bb\in A,\ m\in L_i$.}
\end{equation}
Then $\phi_{L_i}^{\theta/\theta_{L_i}}$ is an $A$-endomorphism of $L_i$. Since $\kk$ is algebraically closed
$\phi_{L_i}^{\theta/\theta_{L_i}}$ must act by a
scalar $\lambda\in\kk^\times$ and we can replace
$\phi_{L_i}$ by $\nu\phi_{L_i}$ where  $\nu^{\theta/\theta_i}=\lambda^{-1}$. Thus, assume from now on that  $ \phi_{L_i}^{\theta/\theta_{L_i}}=\op{id}_{L_i}$.

Given any simple $A$-module $L$, there exists a unique $i$ and a unique $0\leq k\leq \theta_i$ with  $u\colon L\cong\leftindex^{\psi^k}L_i$.  Set $\phi_L=u^{-1}\circ\phi_{L_i}\circ u\colon L \leftindex^{\psi^{\theta_{L}}}L\to L$ noting that $\theta_{L_i}=\theta_{L}$.  As the even automorphisms of $L$ are
scalars, $\phi_L$ is independent of the choice of $u$. We fix such isomorphisms $\phi_{L}$ satisfying
\begin{equation}\label{E:phinorm}
  \phi_{L}^{\theta/\theta_{L}}=\op{id}_{L} .
\end{equation}
Up to isomorphism,  $\mathbf L(\zeta)$ from  \autoref{D:La} depends now only on the orbit of $L$ and on $\zeta$.

We next define some important scalars $\bar\gamma_L$. For this we
 use the notation from \autoref{types}. Observe that taking $L=L'$ in \autoref{tausi} shows that $\tau^{(1,2)}$ acts as zero on $L\otimes L$ if $\theta_L>1$.
\begin{Lemma}\label{L:tauscalar}
  Let $L$ be a simple $A$-module. If $\theta_L=1$ then
      $\beta_L=\tau^{(1,2)}\circ s_1\circ(\op{id}\otimes\phi_L^{-1})$ is
      an (even) endomorphism of the $A^{\otimes2}$-module $L\otimes L$ which
      preserves $L\circledast L$ and acts by multiplication with some scalar $\bar\gamma_L\in\kk$. Set $\bar\gamma_L=0$ if  $\theta_L>1$.
           \end{Lemma}
           \begin{proof} Let $\theta_L=1$.  By \autoref{E:phiL} and \autoref{tausi} $\phi_L^{-1}\colon L\to\leftindex^{\psi}L$ and $\beta_L\colon L\otimes L\to L\otimes L$  are (even) morphisms. By
  Schur's Lemma $\beta_L$  acts by a scalar if $L$ is of type $\mathtt M$. Otherwise,
  $L\otimes L=(L\circledast L)\oplus\Pi(L\circledast L)$ with
  simple summands of type $\mathtt M$.  Then the even part of $\uEnd_{A^{\otimes2}}(L\otimes L)$
   is  two-dimensional and spanned by $\op{id}$ and $c\otimes c$, because any morphism $L\circledast L\to\Pi(L\circledast L)$ is odd. Since $(c\otimes c)^2=-\op{id}$ and
  $e=\frac12(\op{id}+\iota c\otimes c)$, the endomorphism $c\otimes c$
  acts as $-\iota$ on $e(L\otimes L)$ and as $\iota$ on
  $(1-e)(L\otimes L)$. The statements follow.
\end{proof}

  \begin{Remark}\label{gammaprime}
 If $L$ is of type $\mathtt M$ then $\beta_L=\bar\gamma_L\op{id}$. If $L$ is of type $\mathtt Q$ then $\beta_L=\alpha_L\op{id}+\alpha_L' (c\otimes c)$ for unique $\alpha_L,\alpha_L'\in\kk$ and $\beta_L$ acts on $e(L\otimes L)\oplus(1-e)(L\otimes L)$ by the pair of scalars $(\bar\gamma_L=\alpha_L-\iota\alpha_L',\bar\gamma_L'=\alpha_L+\iota\alpha_L').$
  \end{Remark}
  \begin{Remark}The scalar $\bar\gamma_L$ arises naturally. It is, up to the twist by $\op{id}\otimes\phi_L^{-1}$ in the definition of $\beta_L$, the eigenvalue of the image in $\kk\Sym\Smash\Bo$ of the Jucys--Murphy element $\CL_{2}$ from \autoref{D:universalJM}.
\end{Remark}

\begin{Definition}\label{D:La}
 Given a simple $A$-module $L$ we define the following
  $\HHd[1](A)$-modules.
  \begin{Enumerate}
  \item Let $L^{\mathrm{nil}}=L$ with $X_1$ acting by zero.
  \item  For $\zeta\in\kk^\times$ let
  $\mathbf L(\zeta)=\bigoplus_{k=0}^{\theta_L-1}\leftindex^{\psi^k}\!L$ with
  the action of $X_1$ given by the linear endomorphism
  \[x_\zeta(v_0,\dots,v_{\theta_L-1})=(\zeta\phi_L(v_{\theta_L-1}),v_0,\dots,
  v_{\theta_L-2}). \]
  Note that $\bb x_\zeta(m)=x_\zeta(\psi(\bb)m)$ for
  $\bb\in A$ and $m\in \mathbf L(\zeta)$. Observe moreover that
  $X_1^\theta$ acts on $\mathbf L(\zeta)$ by multiplication with the scalar
  $\zeta^{\theta/\theta_L}$.
  \end{Enumerate}
\end{Definition}

Up to isomorphism, $\mathbf L(\zeta)$ depends only on the orbit of $L$ and on $\zeta$, since $\phi_L$ is determined by $\phi_{L_i}$. The modules $\mathbf L_i(\zeta)$, i.e. $\mathbf L(\zeta)$ with $L=L_i$, $[L_i]\in\Omega$, provide the following classification result:
\begin{Theorem}\label{T:simplesn1}
The following give a system of representatives for the simple $\HHd[1](A)$-modules:
  \begin{Enumerate}
    \item\label{irr11} $L^{\mathrm{nil}}$ for $[L]\in\Omega$, which is of dimension $\dim L$, and
    \item\label{irr12} $\mathbf L_i(\zeta)$ for $1\le i\le N_\psi$, $\zeta\in\kk^\times$, which is of
      dimension $(\dim L_i)\theta_i$.
  \end{Enumerate}
 Every simple $\HHd[1](A)$-module is semisimple as $A$-module. It is  simple in case \autoref{irr11}  and semisimple containing each isomorphism class in $O_i$ exactly once as summand in case \autoref{irr12} .
 \end{Theorem}

\begin{proof}
  The module $L^{\mathrm{nil}}$ is simple by construction and,  obviously, $L^{\mathrm{nil}}\cong(L')^{\mathrm{nil}}$ implies $L\cong L'$.

   By minimality of $\theta_L$ the modules
  $\leftindex^{\psi^k}\!L$, $0\le k<\theta_L$, are pairwise non-isomorphic. In particular,
  $\mathbf L(\zeta)$ is a multiplicity free semisimple $A$-module and its
  $A$-submodules are the partial sums of these summands.  As $x_\zeta$
  carries the $k$th summand onto the $(k+1)$st, no proper partial sum is
  $X_1$-stable, so $\mathbf L(\zeta)$ is simple.

  Any isomorphism $f\colon \mathbf L_i(\zeta)\to \mathbf L_j(\zeta')$ restricts to one of $A$-modules,
  so $i=j$ and $f$ matches the summands. Since $f$ is even,
  $f=\bigoplus_k\lambda_k\op{id}$ with $\lambda_k\in\kk^\times$. Then
  $fx_\zeta=x_{\zeta'}f$ forces $\lambda_0=\cdots=\lambda_{\theta_L-1}$ and then
  $\zeta=\zeta'$.  Thus, the modules in \autoref{irr12} are pairwise non-isomorphic.
   Since  $X_1$ acts as zero in case \autoref{irr11} and invertibly in case \autoref{irr12}, all the modules are pairwise non-isomorphic.

  For exhaustiveness note first that
  $\bb X_1^k=X_1^k\psi^{k}(\bb)$ for $\bb\in A$ and $k\ge0$.  Hence in any
  $\HHd[1](A)$-module the subspace $\ker X_1$ is an $A$-submodule, and $X_1^k$
  induces a surjection $\leftindex^{\psi^k}\!L'\to X_1^kL'$ of $A$-modules
  for every $A$-submodule $L'$.  Let now $M$ be a simple
  $\HHd[1](A)$-module.  The even central element $X^\theta$ of
  \autoref{P:HHdCenter} acts as a scalar $\zeta_0$.

  If $\zeta_0=0$ then $X_1$ acts
  nilpotently. Then $\ker X_1$ is a nonzero submodule and hence all of $M$. Thus,
  $M=L^{\mathrm{nil}}$ with $L=M$ simple over $A$.

   If $\zeta_0\ne0$ then $X_1$ acts invertibly on $M$.  In this case
  choose a simple $A$-submodule $L\subseteq M$.  Then
  $M=\sum_{k=0}^{{\theta}-1}X_1^k\!L$ is $A$-semisimple with all
  constituents among the $\leftindex^{\psi^k}\!L$, $k\in\Z/\theta_L\Z$. So we have a decomposition into isotypical components
  $M=\bigoplus_{k\in\Z/\theta_L\Z}M_k$  with
  $X_1M_k=M_{k+1}$.  Since $M_0$ is $L$-isotypic, $\uHom_A(L,M_0)_{\ov0}\ne0$.  We next claim that there is an even isomorphism
  \[\Phi\colon \uHom_A(L,M_0)\to \uHom_A(L,M_0),\quad   \Phi(f)=X^{\theta_L}\circ  f\circ\phi_L^{-1}.\]
 Indeed, $\Phi(f)$ is $A$-linear, because
  $\phi_L^{-1}(\bb m)=\psi^{\theta_L}(\bb)\phi_L^{-1}(m)$ by
  \autoref{E:phiL} and $X_1^{\theta_L}\psi^{\theta_L}(\bb)=\bb X_1^{\theta_L}$ for $\bb\in A$. It has image in
 $X_1^{\theta_L}M_0=M_0$. Thus, $\Phi$ is well-defined. Clearly, $\Phi$ is even  and bijective, since $\phi_L$ is bijective and $X_1$ acts invertibly. Thus $\Phi$ preserves $ \uHom_A(L,M_0)_{\ov0}$ and must have there a nonzero eigenvector $g$, since $\kk$ is algebraically closed. So let  $\Phi(g)=\zeta g$ with $\zeta\in\kk^\times$, i.e.,
  \begin{equation}\label{EVg}
  X_1^{\theta_L}\circ g=\zeta g\circ\phi_L .
  \end{equation}
   Since $g$ is nonzero, hence injective, $g(L)\cong L$ is a simple $A$-submodule of $M_0$ with
  $X_1^{{\theta_L}}g(L)=g(L)$. Moreover,
  $M=\bigoplus_{k=0}^{{\theta_L}-1}X_1^kg(L)$. The right hand side
   is a nonzero $A$-submodule of $M$ which is $X_1$-stable, since $X_1^{\theta_L}g(L)=\zeta g(\phi_L(L))=g(L)$ by \autoref{EVg},  so it must be  $M$, because $M$ is simple.  The sum is direct, since $X_1^kg(L)\subseteq M_k$ and the isotypic components $M_k$, $0\le k<\theta_L$, are pairwise different.
    In particular,
  $(v_k)_k\mapsto\sum_kX_1^kg(v_k)$ defines an isomorphism $\mathbf L(\zeta)\to M$ by \autoref{EVg}.
  Then also $\zeta^{\theta/\theta_L}=\zeta_0$ by \autoref{D:La}.
\end{proof}

\subsection{Atomic algebras}
We will show in the next subsection that for general $n$, the simple modules are  induced from simple modules
concentrated in a single $\psi$-orbit. We prepare this now by showing that such simple modules are governed by certain \emph{RaFH superalgebras} $\mathcal H_n(L)=\HHd[n;\hbar, t\gamma_L, \eta\gamma_L](\mathcal E_n)$, called \emph{atomic algebras}. The Frobenius-tensor superalgebra involved in these are the  inner endomorphism superalgebras $\mathcal E_n$ of $\mathbf L(\zeta)^{\otimes n}$ which we introduce in \autoref{L:endoalg}. The  parameters for these atomic superalgebras have to be chosen carefully depending on the $\psi$-orbit.  Our construction still follows \cite[\S5]{Savage:AffineWreath}, but uses crucially the  parameters and $\Bo$-dependent definition of the RaFH superalgebras. Note that a left $B$-module $M$ for any superalgebra $B$ is always
 a right $\uEnd_B(M)^{\mathrm{op}}$-module via
$v\cdot z=(-1)^{|v||z|}z(v)$. \\

In the following we will use the notation from \autoref{D:La}.
\begin{Proposition}\label{L:endoalg}
  Let $L$ be a simple $A$-module and $\zeta\in\kk^\times$. Consider the endomorphism superalgebras
  \[
    \mathcal E:=\uEnd_A\bigl(\mathbf L(\zeta)\bigr)^{\mathrm{op}}\qquad\qquad\text{and}\qquad\qquad
    \mathcal E_n:=\uEnd_{\Bo}\bigl(\mathbf L(\zeta)^{\otimes n}\bigr)^{\mathrm{op}} .
  \]
  \begin{enumerate}
  \item \label{L:endoalg1} We have  $\mathcal E=\mathcal E_1\cong(\uEnd_AL)^{\oplus \theta_L}$ and $\mathcal E_n\cong\mathcal E^{\otimes n}$.
  \item  \label{L:endoalg2} If $\theta_L=1$ then $\mathcal E=\kk$ or $\mathcal E\cong\Cl_1^\mathrm{op}= \Cl_1$ depending on the type of $L$.
  \item \label{L:endoalg3}
  The canonical right $\mathcal E_n$-action on $\mathbf L(\zeta)^{\otimes n}$ satisfies the following, with $z\in\mathcal E_n$, $\bb\in\Bo$, $v\in \mathbf L(\zeta)^{\otimes n}$,
  \begin{equation}\label{E:bimodule}
    \bb(v\cdot z)=(\bb v)\cdot z,\qquad
    x_\zeta^{(j)}(v\cdot z)=x_\zeta^{(j)}(v)\cdot{(\psi_{\mathcal E}^{(j)})}^{-1}(z),\qquad
    \leftindex^{s_i}(v\cdot z)= \leftindex^{s_i}(v)\cdot\leftindex^{s_i}z .
  \end{equation}
  \item  \label{L:endoalg4} Moreover,  $\mathcal E_n$ is a Frobenius-tensor superalgebra with Nakayama automorphism $\psi_{\mathcal E}(z)=x_\zeta^{-1}zx_\zeta$.
  \end{enumerate}
  \end{Proposition}
\begin{Remark}\label{R:Clop}
We always identify $\Cl_1^\mathrm{op}=\Cl_1$ via an isomorphism which
sends an odd element $\mathfrak{c}$ with $\mathfrak{c}^{2}=-1$ to $\iota
C$ with $C^2=1$  as in \autoref{Ex:Clifford}. This is the meaning of the
isomorphism $\mathcal E\cong \Cl_1^\mathrm{op}=\Cl_1$.
 \end{Remark}

\begin{Remark}\label{psiE}
 Since $x_\zeta$ in \autoref{E:bimodule} is an even
  bijection with $\bb x_\zeta=x_\zeta\psi(\bb)$,
  conjugation by it preserves $\mathcal E$, and replacing $\zeta$ by $\zeta'$ only
  rescales $x_\zeta$ on each summand. Thus, $\psi_{\mathcal E}$ is independent of the choice of
  $\zeta$.

   Explicitly, if $\theta_L=1$ then $\psi_{\mathcal E}$ is the identity in case $\mathcal E=\kk$. In case $\mathcal E=\Cl_1$ it is either the identity or the automorphism $C\mapsto-C$.
   It is the latter whenever $\gamma_L\ne0$, with $\gamma_L$ as in \autoref{L:tauaction}. Indeed, for
$\theta_L=1$ we have $x_\zeta=\zeta\phi_L$, so $\psi_{\mathcal E}(c)=\phi_L^{-1}c\,\phi_L$ is
odd with square $\op{id}$, say $\psi_{\mathcal E}(c)=\mu c$ with $\mu=\pm1$, for $c$ as in
\autoref{types}. Setting $c^{(1)}=c\otimes\op{id}$ and $c^{(2)}=\op{id}\otimes c$ and moving
them through $\beta_L=\tau^{(1,2)}\circ s_1\circ(\op{id}\otimes\phi_L^{-1})$ gives
$\beta_Lc^{(1)}=c^{(2)}\beta_L$ and $\beta_Lc^{(2)}=\mu\,c^{(1)}\beta_L$, because
$\tau^{(1,2)}$ is even, $s_1$ interchanges $c^{(1)}$ and $c^{(2)}$, and $\phi_L$ commutes with
$c^{(1)}$ while conjugating $c^{(2)}$ into $\mu c^{(2)}$. With
$\beta_L=\alpha_L\op{id}+\alpha'_Lc^{(1)}c^{(2)}$ from \autoref{gammaprime} these read
$\alpha'_L=-\alpha_L$ and $\alpha'_L=\mu\alpha_L$, whence $(1+\mu)\alpha_L=0$; and
$\bar\gamma_L=(1+\iota)\alpha_L$ with $2\in\kk^\times$.

It is the latter whenever $\gamma_L\ne0$, with $\gamma_L$ as in \autoref{L:tauaction}. Indeed, or
$\theta_L=1$ we have $x_\zeta=\zeta\phi_L$. Thus,  $\psi_{\mathcal E}(c)=\phi_L^{-1}c\,\phi_L$ is
odd with square $\op{id}$, say $\psi_{\mathcal E}(c)=\mu c$ with $\mu=\pm1$, for $c$ as in
\autoref{types}. Setting $c^{(1)}=c\otimes\op{id}$ and $c^{(2)}=\op{id}\otimes c$ and moving
them through $\beta_L=\tau^{(1,2)}\circ s_1\circ(\op{id}\otimes\phi_L^{-1})$ gives
$\beta_Lc^{(1)}=c^{(2)}\beta_L$ and $\beta_Lc^{(2)}=\mu\,c^{(1)}\beta_L$, because
$\tau^{(1,2)}$ is even, $s_1$ interchanges $c^{(1)}$ and $c^{(2)}$, and $\phi_L$ commutes with
$c^{(1)}$ while conjugating $c^{(2)}$ into $\mu c^{(2)}$. With
$\beta_L=\alpha_L\op{id}+\alpha'_Lc^{(1)}c^{(2)}$ from \autoref{gammaprime} these read
$\alpha'_L=-\alpha_L$ and $\alpha'_L=\mu\alpha_L$, whence $(1+\mu)\alpha_L=0$; and
$\bar\gamma_L=(1+\iota)\alpha_L$ with $2\in\kk^\times$.
\end{Remark}

\begin{proof}[Proof of \autoref{L:endoalg}]The $A$-module $\mathbf L(\zeta)$ is multiplicity free with pairwise non-isomorphic simple summands
  $\leftindex^{\psi^k}\!L$, which gives $\mathcal E \cong(\uEnd_AL)^{\oplus \theta_L}$. The canonical map
  $\uEnd_A(\mathbf L(\zeta))^{\otimes n}\to\uEnd_{\Bo}(\mathbf L(\zeta)^{\otimes n})$
is an injective morphism of superalgebras that intertwines the superflips. It  is also onto because both sides have
  dimension $(\theta_L\dim\uEnd_AL)^n$ by \autoref{T:simplesn1}. Since  $\mathrm{op}$ and $\otimes$ commute, we get $\mathcal E_n\cong\mathcal E^{\otimes n}$.

Now, ${\mathcal E_n}$ turns  into a tensor superalgebra in the sense of \autoref{D:ST-algebra} by setting, for $1\leq i\not= j\leq n$,
 \begin{equation}\label{tauE}
\tau^{(i,j)}_{\mathcal E}= \begin{cases}
0&\text{ if $\gamma_L=0$},\\
 1&\text{if $\gamma_L\not=0$ and $\mathcal E=\kk$},\\
1+C^{i}C^{j}&\text{if $\gamma_L\not=0$ and $\mathcal E=\Cl_1$}.
 \end{cases}
 \end{equation}
 By \autoref{psiE}, $\psi_{\mathcal E}$ is independent of $\zeta$ and satisfies the Nakayama
conditions \autoref{Nakayama}. (The first and third are clear, as $\psi_{\mathcal E}$ acts
factorwise and $\tau^{(i,j)}_{\mathcal E}$ involves only the factors $i$ and $j$. The second and
fourth are vacuous if $\gamma_L=0$, trivial if $\mathcal E=\kk$, and \autoref{Ex:Clifford} if
$\mathcal E=\Cl_1$, since then $\psi_{\mathcal E}(C)=-C$ by \autoref{psiE}.)
Hence, $\mathcal E_n$ is a Frobenius-tensor superalgebra, \autoref{D:Frobtensor}, with Nakayama automorphism $\psi_{\mathcal E}$. Since part \autoref{L:endoalg2} is obvious, it remains to verify the formulas \autoref{E:bimodule}.  Let
  $z\in\mathcal E_n$, $\bb\in\Bo$, $v\in\mathbf L(\zeta)^{\otimes n}$
  be homogeneous, so that $v\cdot z=(-1)^{|v||z|}z(v)$ and
  $z(\bb v)=(-1)^{|z||\bb|}\bb z(v)$. The latter holds because $z$ is a
  morphism of $\Bo$-modules. \\
  \noindent The first equality in \autoref{E:bimodule}
  holds because  $(\bb v)\cdot z=(-1)^{(|\bb|+|v|)|z|}z(\bb v)
      =(-1)^{|v||z|}\bb z(v)=\bb(v\cdot z)$.\\
\noindent   The second one restates the definition of
  $\psi_{\mathcal E}$. As operators, $x_\zeta^{(j)}\circ z
  =(\psi_{\mathcal E}^{(j)})^{-1}(z)\circ x_\zeta^{(j)}$, and both $x_\zeta$ and
  $\psi_{\mathcal E}$ are even, so no signs occur. \\
  \noindent  The third formula in  \autoref{E:bimodule} holds
  because the isomorphism $\mathcal E_n\cong\mathcal E^{\otimes n}$
  intertwines the superflips, that is
  $\leftindex^{s_i}(z(v))=(\leftindex^{s_i}z)(\leftindex^{s_i}v)$.
\end{proof}

In the following we write  $\rho_i(v)=\leftindex^{s_i}v$ for the superflip of the $i$th and
  $(i+1)$st tensor factor of $\mathbf L(\zeta)^{\otimes n}$.
\begin{Lemma}\label{L:tauaction}
  Let $L$ be a simple $A$-module and $\zeta\in\kk^\times$.  There exist
  unique even elements $z_i\in\mathcal E_n$ with
 \begin{equation}\label{E:tauaction}
    \tau^{(i,i+1)}=\rho_i\circ x_\zeta^{(i)}\circ z_i\qquad(1\le i<n),
  \end{equation}
 as operators on   $\mathbf L(\zeta)^{\otimes n}$. Namely, $z_i=\tfrac{\gamma_L}{\zeta}\tau^{(i,i+1)}_{\mathcal E}$ with
  $\tau^{(i,i+1)}_{\mathcal E}\in\mathcal E_n$ as in \autoref{tauE}, and
      \begin{equation}\label{gammadef}
    \gamma_L=\begin{cases}
      \bar\gamma_L & \text{if $\theta_L=1$ and $L$ is of type $\mathtt M$},\\
      \tfrac{1-\iota}{2}\bar\gamma_L
        & \text{if $\theta_L=1$ and $L$ is of type $\mathtt Q$},\\
      0 & \text{if $\theta_L>1$},
    \end{cases} \quad\quad\quad\text{(with $\bar\gamma_L$ as in \autoref{L:tauscalar})}.
  \end{equation}
\end{Lemma}

\begin{proof}
  As $\rho_i$ and $x_\zeta^{(i)}$ are even and invertible, \autoref{E:tauaction} implies
  $z_i=(\rho_ix_\zeta^{(i)})^{-1}\tau^{(i,i+1)}=(x_\zeta^{(i)})^{-1}\rho_i\tau^{(i,i+1)}$. Hence $z_i$ is unique and  $z_i\in\mathcal E_n$, because
 moving $\bb\in\Bo$ through $\tau^{(i,i+1)}$ and $\rho_ix_\zeta^{(i)}$ twists $\bb$ in the same way by \autoref{athroughtau}
  respectively by $\rho_i\bb=\leftindex^{s_i}\bb \rho_i$ and
  $\bb x_\zeta^{(i)}=x_\zeta^{(i)}\psi^{(i)}(\bb)$. We next show that $z_i=\tfrac{\gamma_L}{\zeta}\tau^{(i,i+1)}_{\mathcal E}$.

  If $\theta_L>1$ then both sides vanish, since $\taui[i,i+1]$ acts by zero, $\bar\gamma_L=0$, and  $\tau^{(i,i+1)}_{\mathcal E}=0$ by \autoref{tauE}.\\
  \noindent  If $\theta_L=1$ then assume without loss of generality $n=2$ and $i=1$ and express $\tau^{(1,2)}$ via $\beta_L$ from \autoref{L:tauscalar} using $x_\zeta=\zeta\phi_L$ to get
  \begin{equation}\label{E:zbeta}
    z_1=(x^{(1)}_\zeta)^{-1}\rho_1\tau^{(1,2)}
      =\tfrac1\zeta (\phi^{(1)}_L)^{-1}\bigl(\rho_1\beta_L\rho_1\bigr)\phi^{(1)}_L.
  \end{equation}
    If $L$ is of type $\mathtt M$, then $\beta_L=\bar\gamma_L\op{id}$, so $z_1=\tfrac{\gamma_L}{\zeta}=\tfrac{\gamma_L}{\zeta}\tau^{(1,2)}_{\mathcal E}$ as desired.

  \noindent If $L$ is of type $\mathtt Q$, let $c^{(1)}=c\otimes\op{id}$, $c^{(2)}=\op{id}\otimes c$ with $c$ as in \autoref{types}. Then, by  \autoref{gammaprime}, $\beta_L=\alpha_L\op{id}+\alpha'_Lc^{(1)}c^{(2)}$. We get $\alpha'_L=-\alpha_L$ since $\beta_Lc^{(1)}=c^{(2)}\beta_L$, because  $\tau^{(1,2)}$ is even and thus commutes with the endomorphism $c^{(2)}$. Altogether, $\beta_L=\alpha_L(\op{id}-c^{(1)}c^{(2)})$ with $\alpha_L=\bar\gamma_L(1+\iota)^{-1}=\gamma_L$. Then \autoref{E:zbeta} gives $z_1=\tfrac{\gamma_L}{\zeta}(\op{id}-c^{(1)}c^{(2)})$ and therefore $z_1=\tfrac{\gamma_L}{\zeta}\tau^{(1,2)}_{\mathcal E}$ by \autoref{R:Clop}. This proves the claim.
 \end{proof}

\begin{Definition}\label{D:block}
For fixed $\hbar, t,\eta$, the  \emph{atomic algebra} associated to a simple $A$-module $L$ is the RaFH
  superalgebra
  $\mathcal H_n(L)=\HHd[n;\hbar, t\gamma_L, \eta\gamma_L](\mathcal E_n)$
  attached to $\mathcal E_n$ from \autoref{L:endoalg} with $\zeta=1$.
\end{Definition}

\begin{Remark}Different invertible $\zeta$ give isomorphic algebras and
$\zeta=1$ is merely a choice we made.
Up to canonical isomorphism the atomic algebra $\mathcal H_n(L)$ depends only on the $\psi$-orbit of $L$.
\end{Remark}

For  $\nu,\mu\in\kk^\times$, there is an isomorphism $\HHd[n;\hbar,t,\eta](\Bo)\cong\HHd[n;\nu\hbar,\nu\mu t,\nu\eta](\Bo)$ by rescaling $H_i$ by $\nu$ and $X_i$ by $\mu$. Using these isomorphisms we can group the atomic algebras into three different types:

\begin{Definition}\label{R:blocks}
Using \autoref{L:endoalg} and \autoref{L:tauaction} we distinguish three atomic types:
  \begin{Enumerate}
    \item\label{Hecketype1} \emph{Hecke type}: $\theta_L=1$, type $\mathtt M$ and  $t\gamma_L\ne0$ or $\eta\gamma_L\ne0$. Here, $\mathcal E_n=\kk$ and
      $\mathcal H_n(L)=\HHd[n;\hbar,t\gamma_L,\eta\gamma_L](\kk)$.

      \noindent After rescaling, this is the degenerate affine Hecke algebra of \autoref{degaff} if $\eta=0$, and the affine Hecke algebra of \autoref{ex:aff} if $\eta\ne0$.
    \item\label{Hecketype2} \emph{Clifford type}:  $\theta_L=1$, type $\mathtt Q$ and $t\gamma_L\ne0$. Here,
      $\mathcal E_n=\Cl_n$ and
      $\mathcal H_n(L)=\HHd[n;\hbar,t\gamma_L,0](\Cl_n)$.

      \noindent After rescaling, this is the degenerate affine Sergeev algebra of \autoref{ex:Sergeev}. Note that $\eta\gamma_L=0$ in this case since $\eta\ne0$ forces $\psi=\op{id}$ by $\eta$-symmetry, and then  $\gamma_L=0$ by \autoref{psiE}.

    \item\label{wreathtype}  \emph{Wreath type}: $t\gamma_L=\eta\gamma_L=0$. Here,
      $\mathcal H_n(L)=\HHd[n;\hbar,0,0](\mathcal E_n)$.

    \noindent This case applies  in particular if $\gamma_L=0$, which holds whenever $\theta_L>1$.  After rescaling the $H_i$ by $\hbar^{-1}$ it is the wreath product $(\kk[X]\Smash\mathcal E)^{\otimes n}\rtimes\kk\Sym$, see
     \autoref{ex:smash}.
  \end{Enumerate}
  \end{Definition}

\begin{Remark}
  The three types are exhaustive and disjoint. If $t\gamma_L\not=0$ or $\eta\gamma_L\not=0$, then $\gamma_L\ne0$, hence $\theta_L=1$ by \autoref{L:tauaction}. If moreover $L$ is of type $\mathtt Q$, then $\psi\ne\op{id}$ by \autoref{psiE}, so $\eta=0$ and $t\gamma_L\ne0$.
\end{Remark}

\subsection{Induced modules and Morita equivalence}
 Let $\HHd[n](\Bo)=\HHd[n;\hbar,t,\eta](\Bo)$ be a RaFH algebra. For a semisimple $\Bo$-module $\Lo$, let  $\mathcal C(\Lo)$ be the  full subcategory of all $\HHd[n](\Bo)$-modules which are semisimple as
  $\Bo$-modules with all simple summands occurring in $\Lo$ up to parity shift.  In particular $\mathcal C(\Lo)$ is stable under $\Pi$. We are now ready to construct simple $\HHd[n](\Bo)$-modules from simple modules $S$ in  $\mathcal C(\Lo)$.

\begin{Theorem}\label{T:Morita}
  Let $L$ be a simple $A$-module and $\Lo=\mathbf L(1)^{\otimes n}$. Given an
  $\mathcal H_n(L)$-module $S$, the formulas
   \begin{equation}\label{E:tensorstructure}
    \bb(v\otimes z)=\bb v\otimes z,\qquad
    X_j(v\otimes z)=x_1^{(j)}(v)\otimes X_jz,\qquad
    H_i(v\otimes z)=\rho_i(v)\otimes H_iz,
  \end{equation}
   where  $\bb\in\Bo$, $v\in\Lo$, $z\in S$, turn  $\Lo\otimes_{\mathcal E_n}S$ into an
  $\HHd[n;\hbar,t,\eta](\Bo)$-module.
  Moreover, the functors
  \begin{equation}\label{functorsMorita}
    S\longmapsto\Lo\otimes_{\mathcal E_n}S,\qquad
    M\longmapsto\uHom_{\Bo}(\Lo,M),
  \end{equation}
  with the action of the  atomic algebra $\mathcal H_n(L)$ on $\uHom_{\Bo}(\Lo,M)$ given by \autoref{actiononHom}, are mutually inverse equivalences between the category of
  $\mathcal H_n(L)$-modules and $\mathcal C(\Lo)$.

In particular, $\Lo\otimes_{\mathcal E_n}S$ is simple if and only if $S$ is simple.
\end{Theorem}

\begin{Lemma}\label{actiononHom}
  Let $M$ be an $\HHd[n](\Bo)$-module. The atomic algebra $\mathcal H_n(L)=\HHd[n;\hbar,t\gamma_L,\eta\gamma_L](\mathcal E_n)$ acts on $\uHom_{\Bo}(\Lo,M)$. In the generators  $e$, $X_j$ and $H_i$ of $\mathcal H_n(L)$,  with  $e\in\mathcal E_n$, the action is given by
  \[
    (e\cdot f)(v)=f(v\cdot e),\qquad
    X_j\cdot f=X_j\circ f\circ(x_1^{(j)})^{-1},\qquad
    H_i\cdot f=H_i\circ f\circ\rho_i^{-1},
  \]
  where on the right $X_j$ and $H_i$ act on $M$ and $f\in \uHom_{\Bo}(\Lo,M)$ and $v\in \Lo$.

  With this action,  $\Lo\otimes_{\mathcal E_n}\uHom_{\Bo}(\Lo,M)\to M$, $v\otimes f\mapsto f(v)$, is the adjunction map from \autoref{functorsMorita}. It is a morphism of $\HHd[n](\Bo)$-modules, where $\HHd[n](\Bo)$ acts on the left side by \autoref{E:tensorstructure}.
\end{Lemma}
\begin{proof}
The claims are straightforward calculations using $\bb x_1^{(j)}=x_1^{(j)}\psi^{(j)}(\bb)$ and $\rho_i\bb=\leftindex^{s_i}\bb\rho_i$.
\end{proof}

\begin{proof}[Proof of \autoref{T:Morita}]
  The check that the formulas \autoref{E:tensorstructure} are $\mathcal E_n$-balanced is straightforward using \autoref{E:bimodule}. Each relation of \autoref{D:Hecke} not involving $\tau$ holds on $\Lo\otimes_{\mathcal E_n}S$, since it holds on $S$ and the corresponding relation holds for the operators on $\Lo$. More precisely, $\bb x_1^{(j)}=x_1^{(j)}\psi^{(j)}(\bb)$ by \autoref{D:La}, $\rho_i\bb=\leftindex^{s_i}\bb\rho_i$, the $x_1^{(j)}$ commute with each other, $\rho_ix_1^{(j)}=x_1^{(j)}\rho_i$ for $j\ne i,i+1$, and the $\rho_i$ satisfy the braid relations and commute if $|i-j|>1$. For the other relations,
  \autoref{E:tauaction} with $\zeta=1$, \autoref{E:bimodule} and the Nakayama condition $\leftindex^{s_i}{\bigl((\psi^{(i)}_{\mathcal E})^{-1}(\tau^{(i,i+1)}_{\mathcal E})\bigr)}=\tau^{(i,i+1)}_{\mathcal E}$ for $\mathcal E_n$ give
  \begin{equation}\label{E:tautensor}
    \tau^{(i,i+1)}(v\otimes z)
      =\gamma_L \rho_ix_1^{(i)}(v)\otimes\tau^{(i,i+1)}_{\mathcal E}z .
  \end{equation}
  Each $h$ of the form $H_iX_i-X_{i+1}H_i$,
  $t\tau^{(i,i+1)}$, $\eta X_{i+1}\tau^{(i,i+1)}$ acts as
  $h.(v\otimes z)=\rho_ix_1^{(i)}(v)\otimes h'z$, where $h'$ is
  the corresponding element of $\mathcal H_n(L)$, namely
  $H_iX_i-X_{i+1}H_i$, $(t\gamma_L)\tau^{(i,i+1)}_{\mathcal E}$ and $(\eta\gamma_L)X_{i+1}\tau^{(i,i+1)}_{\mathcal E}$, respectively. For the
  last one note that $\eta\ne0$ forces $\psi=\op{id}$ and $x_1^{(i+1)}=\op{id}$. Hence $H_iX_i-X_{i+1}H_i+\gamma_i\tau^{(i,i+1)}$ acts as $\rho_ix_1^{(i)}$ tensored with $H_iX_i-X_{i+1}H_i+(t\gamma_L+\eta\gamma_LX_{i+1})\tau^{(i,i+1)}_{\mathcal E}$, which is zero on $S$ by the defining relations of $\mathcal H_n(L)$.  Here and in the next relation we use the specifically chosen parameters.
  For the quadratic relation $H_i^2=\eta\tau^{(i,i+1)}H_i+\hbar^2$ note that, by
  \autoref{E:tautensor}, the operator $\tau^{(i,i+1)}H_i$ sends $v\otimes z$
  to $\gamma_L\rho_ix_1^{(i)}\rho_i(v)\otimes\tau^{(i,i+1)}_{\mathcal E}H_iz$,
  and $\rho_ix_1^{(i)}\rho_i=x_1^{(i+1)}$ is the identity whenever
  $\eta\gamma_L\ne0$. So $\eta\tau^{(i,i+1)}H_i$ sends $v\otimes z$ to $v\otimes\eta\gamma_L\tau^{(i,i+1)}_{\mathcal E}H_iz$, while
  $H_i^2(v\otimes z)=v\otimes H_i^2z$ since $\rho_i^2=\op{id}$. Thus the quadratic relation holds on
  $\Lo\otimes_{\mathcal E_n}S$, as the corresponding
  relation of $\mathcal H_n(L)$ holds on $S$. Finally, the relation $H_iX_{i+1}=X_iH_i+\gamma_i\tau^{(i+1,i)}$ follows from $H_iX_i=X_{i+1}H_i-\gamma_i\tau^{(i,i+1)}$ and the quadratic relation as in \autoref{Rk:hinv}, since $\hbar$ is invertible.

  It remains to prove the equivalence. Let $\underline\Bo$ be the image of $\Bo$ in $\uEnd_\kk(\Lo)$. As $\Lo$ is a faithful semisimple $\underline\Bo$-module, $\underline\Bo$ is semisimple, $\Lo$ contains every simple $\underline\Bo$-module up to parity shift, and $\uEnd_{\underline\Bo}(\Lo)^{\mathrm{op}}=\mathcal E_n$. The $\underline\Bo$-modules are exactly the $\Bo$-modules which are semisimple with all summands in $\Lo$ up to parity shift. Hence, by Morita theory, $\Lo\otimes_{\mathcal E_n}-$ and $\uHom_{\Bo}(\Lo,-)$ are mutually inverse equivalences between $\mathcal E_n$-modules and these $\Bo$-modules, via the isomorphism $v\otimes f\mapsto f(v)$. In particular $\Lo\otimes_{\mathcal E_n}S\in\mathcal C(\Lo)$. By the first part, each defining relation $h$ of $\HHd[n](\Bo)$ acts on $\Lo\otimes_{\mathcal E_n}S$ as $g\otimes h'$, where $g$ is a bijection of $\Lo$ and $h'$ is the corresponding relation of $\mathcal H_n(L)$. So $h$ acts as zero if and only if $h'$ does. Hence, for $M\in\mathcal C(\Lo)$, the operators of \autoref{actiononHom} make $\uHom_{\Bo}(\Lo,M)$ an $\mathcal H_n(L)$-module, and the two functors are mutually inverse equivalences. The last claim follows, as equivalences preserve simple objects.
\end{proof}

\subsection{Construction and  classification of irreducible representations}
In this subsection we construct and classify the simple $\HHd[n](\Bo)$-modules, see
\autoref{T:simples}. The classification is stated combinatorially with the important data being the $\psi$-orbits of the simple $A$-modules. We call these orbits \emph{colours} and fix a representative $L_i$ for each colour $i$. Since the simple $\Bo$-modules are tensor products of simple $A$-modules, they come with a sequence of colours. We encode the \emph{colour multiplicity}, that is, how often each colour occurs, as a composition $\mu$ of $n$ and consider the corresponding parabolic subalgebra $\HHd[\mu](\Bo)$ of $\HHd[n](\Bo)$. We then show that every simple $\HHd[n](\Bo)$-module arises by parabolic induction
of a module which is a tensor product of modules $\mathbf L_i(1)^{\otimes\mu_i}\otimes_{\mathcal E_{\mu_i}(L_i)}S_i$ as in \autoref{T:Morita}, where $S_i$ runs through the simple modules for the atomic algebra $\mathcal H_{\mu_i}(L_i)$ of \autoref{D:block}. As a consequence, the simple modules are labelled by $\mu$ together with the tuple of the $S_i$. On the technical side, two super phenomena have to be taken into account.
The above-mentioned  tensor products of simple modules are really the $\circledast$-products from
\autoref{types}, and the $\psi$-orbits and the classification are up to parity shift, since in a simple $\Bo$-module the $\circledast$-factors are only determined up to parity shift, as the shift can freely move through a $\circledast$-product. We advise the
reader to read the theorem first in the non-super setting, where $\circledast$ is the
tensor product and both complications disappear. As before, we generalise the results from \cite{Savage:AffineWreath, WanWang}.

\begin{Definition}
For a composition $\mu=(\mu_1,\dots,\mu_q)$ of $n$ let
$\HHd[\mu](\Bo)=\HHd[\mu_1](A^{\otimes\mu_1})\otimes\cdots\otimes
\HHd[\mu_q](A^{\otimes\mu_q})$ which we view as a subalgebra of $\HHd[n](\Bo)$ via \autoref{C:BasisHd}, the  \emph{parabolic subalgebra attached to $\mu$}. Note that compositions are allowed to have parts equal to zero.
We denote by  $\Sym^{\mu}$ the set of minimal length
representatives for $\Sym/\Sym_{\mu}$. Let $\Ind_{\HHd[\mu](\Bo)}^{\HHd[n](\Bo)}=\HHd[n](\Bo)\otimes_{{\HHd[ \mu](\Bo)}}-$ be the induction functor.
\end{Definition}

\noindent For a $\Bo$-module $M$ and $w\in\Sym$  let $\leftindex^{w}M$ be  $M$ with the twisted action such that
  $\bb\in\Bo$ acts now by $\leftindex^{w^{-1}}\bb$.
The next lemma copies classical arguments from wreath products representation theory.

\begin{Lemma}\label{L:induction}
  Let $\mu$ be a composition of $n$ and let $M$ be a simple
  $\HHd[\mu](\Bo)$-module which is semisimple as a $\Bo$-module.
Consider the twisted modules $\leftindex^{d}M$,  $d\in\Sym^{\mu}$.  If
\begin{equation}\label{E:separate}
\text{$\leftindex^{d}M$ and $\leftindex^{d'}M$ have no  common irreducible composition factor for any $d\ne d'\in\Sym^{\mu}$,}
\end{equation}
then the induced module $\So=\Ind_{\HHd[\mu](\Bo)}^{\HHd[n](\Bo)}M$  is a simple $\HHd[n](\Bo)$-module.

If $t\tau^{(i,j)}$ and $\eta\tau^{(i,j)}$ act as
  zero on $\So$ for all $i\ne j$, then the above statement holds also for
 $\Bo$ replaced by $\Po$, where in addition
  $\leftindex^{d}X_j=X_{d(j)}$.
\end{Lemma}

\begin{proof}
By \autoref{C:BasisHd}, $\So=\bigoplus_{d\in \Sym^{\mu}}H_d\otimes M$ as $R$-modules.
Since $\bb H_d=H_d\leftindex^{d^{-1}}\bb$ for
  $\bb\in\Bo$, each summand $H_d\otimes M$ is a $\Bo$-submodule of
  $\So$ isomorphic to $\leftindex^{d}M$ via $H_d\otimes m\mapsto m$. In
  particular $\So$ is semisimple as a $\Bo$-module.  By
  \autoref{E:separate} these summands are the sums of the isotypic
  components of $\So$ belonging to the various
  $d$, so that $N=\bigoplus_d\left(N\cap(H_d\otimes M)\right)$ for every
  $\Bo$-submodule $N\subseteq\So$.

  Let now $0\ne N\subseteq\So$ be an $\HHd[n](\Bo)$-submodule and choose
  $d$ with $N\cap(H_d\otimes M)\ne0$. Since $H_{d^{-1}}$ is invertible and $\bb H_{d^{-1}}=H_{d^{-1}}\leftindex^{d}\bb$, multiplication by $H_{d^{-1}}$ maps $N\cap(H_d\otimes M)$ injectively onto a $\Bo$-submodule of $N$ which is isomorphic to a submodule of $\leftindex^{d^{-1}}(\leftindex^{d}M)\cong M$ and lies in $1\otimes M$ by \autoref{E:separate}. Now,
  $N\cap(1\otimes M)$ is a nonzero submodule of the simple
  $\HHd[\mu](\Bo)$-module $1\otimes M$ and thus
  $N=\HHd[n](\Bo)(1\otimes M)=\So$.

    For the last statement assume that $t\tau^{(i,j)}$ and $\eta\tau^{(i,j)}$ act as zero on $\So$ for all $i\ne j$. We compute several actions on $\So$. First,  $\gamma_i\tau^{(i,i+1)}$ and $\gamma_i\tau^{(i+1,i)}$ act as zero as well, since in case $\eta\ne0$ we have $\psi=\op{id}$ and so  $X_{i+1}$ commutes with $\tau^{(i,j)}$.
    Second,  \autoref{RaFHrels} implies now $H_iX_j=X_{s_i(j)}H_i$ for all $i,j$ and therefore $yH_d=H_d\leftindex^{d^{-1}}y$ for all $y\in\Po$ and $d\in\Sym$. Note that  this identity is the only relation between the action of $H_d$ and $\Bo$ used above. So the above  argument applies verbatim to $\Po$ in place of $\Bo$. Each $H_d\otimes M$ is a $\Po$-submodule of $\So$ isomorphic to $\leftindex^{d}M$, and \autoref{E:separate} for $\Po$-modules gives the claim.
\end{proof}

To state the classification theorem we need an extra combinatorial ingredient, the \emph{colours}. For this consider the set  $\bar\Omega$ of isomorphism classes of simple $A$-modules up to parity shift. Twisting by $\psi$ commutes with parity shift and induces, with \autoref{types}\autoref{typesa}, an action of the cyclic group $\mathtt{C}_\theta$ on $\bar\Omega$.

\begin{Definition}\label{D:colours}
We denote by  $\bar O_1,\dots,\bar O_{\bar N_\psi}$ the $\mathtt{C}_\theta$-orbits  in  $\bar\Omega$.
 We set $\bar I=\set{1,\dots,\bar N_\psi}$ and call its elements and the corresponding orbits \emph{colours}.
 We finally renumber the representatives $L_i$ fixed before \autoref{tausi} so that $[L_i]\in\bar O_i$ for $i\in \bar I$. Given a simple $\Bo$-module $P$, which  by \autoref{types}\autoref{typese} is of the form $P\cong\Pi^{\varepsilon}(L'_1\circledast\dots\circledast L'_n)$ for some simple $A$-modules $L'_j$,
we define the \emph{colour} $\chi(P)$ of $P$ as the tuple $\chi(P)=(\chi_1,\ldots,\chi_n)\in \bar I^n$ or as the word $\chi_1\cdots\chi_n$ in the alphabet $\bar I$, where $[L'_j]\in\bar O_{\chi_j}$. The \emph{colour multiplicity} of a word $\chi$ in $\bar I$ is the composition $|\chi|=(c_1,\dots, c_{\bar N_\psi})$ of its length, where $c_i$ is the number of occurrences of the letter $i\in\bar I$.
\end{Definition}

\begin{Remark}
 Passing to
  $\bar\Omega$ is necessary. By \autoref{types}\autoref{typesd}, a simple $\Bo$-module
  determines its tensor factors only up to parity shift, since
  $\Pi(L_1'\otimes L_2'\otimes\cdots)\cong(\Pi L_1')\otimes
  L_2'\otimes\cdots\cong L_1'\otimes(\Pi
  L_2')\otimes\cdots$ etc. holds.
\end{Remark}

\begin{Remark}
One easily verifies that the colour $\chi(P)$ is well-defined. Each colour is $\psi$-stable and $\leftindex^{d}(L'_1\circledast\dots\circledast L'_n)\cong L'_{d^{-1}(1)}\circledast\dots\circledast L'_{d^{-1}(n)}$ for any $d\in\Sym$.  Setting $(\leftindex^{d}\chi)_j=\chi_{d^{-1}(j)}$ for $1\leq j\leq n$, we have
\begin{equation}\label{E:colours}
  \chi(\Pi P)=\chi(P),\qquad
  \chi(\leftindex^{\psi^{(j)}}\!P)=\chi(P),\qquad
  \chi(\leftindex^{d}P)=\leftindex^{d}\chi(P).
\end{equation}
 For a composition $\mu$ of $n$ the word $\chi^\mu=1^{\mu_1}2^{\mu_2}\cdots\bar N_\psi^{\mu_{\bar N_\psi}}$ has stabiliser $\Sym_\mu$, so $d\mapsto\leftindex^{d}\chi^\mu$ is a bijection from $\Sym^\mu$ onto the words of colour multiplicity $\mu$.
\end{Remark}

\begin{Theorem}\label{T:simples}
 Let $\mu=(\mu_1,\dots,\mu_{\bar N_\psi})\in\N_0^{\bar N_\psi}$ be a composition of
 $n$ and set, following \autoref{L:endoalg},
 \[\mathbf L_\mu=\bigboxtimes_{i=1}^{\bar N_\psi}\mathbf L_i(1)^{\otimes\mu_i}\quad\text{and}\quad
  \mathcal E_\mu=\bigotimes_{i=1}^{\bar N_\psi}\mathcal E_{\mu_i}(L_i),\quad\text{where }\mathcal E_m(L_i)=\uEnd_{A^{\otimes m}}\bigl(\mathbf L_i(1)^{\otimes m}\bigr)^{\mathrm{op}}.\]
  Then $\mathcal E_\mu=\uEnd_{\Bo}(\mathbf L_\mu)^{\mathrm{op}}$. For simple
  $\mathcal H_{\mu_i}(L_i)$-modules $S_1,\dots,S_{\bar N_\psi}$ set
    \begin{equation}\label{E:simples}
    \So_\mu=\Ind_{\HHd[\mu](\Bo)}^{\HHd[n](\Bo)}
    \Bigl(\mathbf L_\mu\otimes_{\mathcal E_\mu}
      \bigl(S_1\circledast\cdots\circledast S_{\bar N_\psi}\bigr)\Bigr).
  \end{equation}
  This gives a bijection
  \[
    \coprod_{\substack{\mu\in\N_0^{\bar N_\psi}\\ \mu_1+\dots+\mu_{\bar N_\psi}=n}}\
      \prod_{i=1}^{\bar N_\psi}\Irr\bigl(\mathcal H_{\mu_i}(L_i)\bigr)
    \ \longrightarrow\ \Irr\bigl(\HHd[n;\hbar,t,\eta](\Bo)\bigr),
   \quad(\mu,S_\bullet)\longmapsto\So_\mu.
  \]
  Moreover, every simple $\HHd[n;\hbar,t,\eta](\Bo)$-module is semisimple
  as a $\Bo$-module.
\end{Theorem}

\begin{proof}[Proof of \autoref{T:simples}]
First note that $\mathcal E_\mu=\uEnd_{\Bo}(\mathbf L_\mu)^{\mathrm{op}}$ by the dimension count in the proof of \autoref{L:endoalg}, because constituents of different colours are non-isomorphic, also up to parity shift.

  \emph{(Semi-)Simplicity.} Put $S=S_1\circledast\dots\circledast S_{\bar N_\psi}$
  and $M=\mathbf L_\mu\otimes_{\mathcal E_\mu}S$, where $S$ is simple over
  $\bigotimes_i\mathcal H_{\mu_i}(L_i)$ by \autoref{types}\autoref{typesc}. The proof of \autoref{T:Morita} applies verbatim to $\HHd[\mu](\Bo)=\bigotimes_i\HHd[\mu_i](A^{\otimes\mu_i})$, with $\mathbf L_\mu$, $\mathcal E_\mu$ and $\bigotimes_i\mathcal H_{\mu_i}(L_i)$ in place of $\Lo$, $\mathcal E_n$ and $\mathcal H_n(L)$. Hence $M$ is a simple $\HHd[\mu](\Bo)$-module in $\mathcal C(\mathbf L_\mu)$. As all summands of $\mathbf L_i(1)$ are $\psi$-twists of $L_i$, all $\Bo$-constituents of $M$ have colour $\chi^\mu$. By \autoref{E:colours} the twists
  $\leftindex^{d}M$, $d\in\Sym^{\mu}$, then have pairwise different
  colours, so no two of them share a composition factor, and
  $\So_\mu$ is simple by \autoref{L:induction}. It is
  $\Bo$-semisimple, since each $\Bo$-summand $H_d\otimes M\cong\leftindex^{d}M$ of $\So_\mu$ is $\Bo$-semisimple.

  \emph{Injectivity.} Let $\So'_{\mu'}$ be defined as $\So_\mu$, using simple $\mathcal H_{\mu'_i}(L_i)$-modules $S'_i$, and suppose that $\So_\mu\cong\Pi^{\varepsilon}\,\So'_{\mu'}$. Let $S'$ and $M'$ be defined for $(\mu',S'_\bullet)$ as $S$ and $M$ are for $(\mu,S_\bullet)$. By \autoref{E:colours}, the simple $\Bo$-submodules of $\So_\mu$ have the colours $\leftindex^{d}\chi^\mu$ with $d\in\Sym^\mu$. These are precisely the words of colour multiplicity $\mu$. As isomorphisms and parity shift preserve colours, we get $\mu=\mu'$. Moreover, as in the proof of \autoref{L:induction}, the sum of the simple $\Bo$-submodules of colour $\chi^\mu$ in $\So_\mu$ is the summand $1\otimes M$. Since this is an $\HHd[\mu](\Bo)$-submodule, the isomorphism restricts to an isomorphism $M\cong\Pi^{\varepsilon}M'$ of $\HHd[\mu](\Bo)$-modules. By the equivalence used in the simplicity part, $S\cong\uHom_{\Bo}(\mathbf L_\mu,M)\cong\Pi^{\varepsilon}S'$. Finally, $S=S_1\circledast\dots\circledast S_{\bar N_\psi}$ determines the $S_i$ up to parity shift by \autoref{types}\autoref{typesd}, so $S_i\cong\Pi^{\varepsilon_i}S'_i$ for all $i$.

  \emph{Surjectivity.} Let $M$ be a simple $\HHd[n](\Bo)$-module. It is finite
  dimensional by \autoref{fdauto}. Let $V\subseteq M$ be a simple
  $\Po$-submodule. As $\Po\cong\HHd[1](A)^{\otimes n}$, $V$ is
  $\Bo$-semisimple by \autoref{types} and \autoref{T:simplesn1}. By
  \autoref{D:Hecke} we have $X_jH_w\in\sum_{w'}H_{w'}\Po$ and
  $\bb H_w=H_w\leftindex^{w^{-1}}\bb$. Hence $\sum_{w}H_wV$ is stable under $\Bo$ and the $X_j$, and it is stable
  under the $H_i$ by \autoref{L:HiHw}, so that $M=\sum_{w}H_wV$. As each
  $H_wV$ is a quotient of $\leftindex^{w}{V}$, it is semisimple over $\Bo$, and so is $M$. Let
  $M_\chi$ be the sum of the simple $\Bo$-submodules of colour $\chi$, so
  that $M=\bigoplus_{\chi}M_\chi$. Then $\Bo M_\chi\subseteq M_\chi$.
  Moreover $X_jM_\chi\subseteq M_\chi$, since $\bb X_j=X_j\psi^{(j)}(\bb)$
  exhibits $X_jP$ as a quotient of
  $\leftindex^{\psi^{(j)}}P$ for every
  simple $\Bo$-submodule $P\subseteq M_\chi$, and
  $H_iM_\chi\subseteq M_{\leftindex^{s_i}\chi}$, since
  $H_i\bb=\leftindex^{s_i}\bb H_i$, both by \autoref{E:colours}. Choose
  $\chi$ with $M_\chi\ne0$ and put $\mu=|\chi|$. As the $H_w$ are
  invertible by \autoref{Rk:hinv} and
  $H_wM_\chi\subseteq M_{\leftindex^{w}\chi}$, we may assume
  $\chi=\chi^\mu$ and write $M_\mu=M_{\chi^\mu}$. Then $M_\mu$ is an
  $\HHd[\mu](\Bo)$-module, because $\leftindex^{s_i}\chi^\mu=\chi^\mu$
  for $s_i\in\Sym_\mu$, and $\bigoplus_{d\in\Sym^{\mu}}H_dM_\mu$,
  a sum over the pairwise different colours $\leftindex^{d}\chi^\mu$, is stable
  under $\Bo$, the $X_j$ and, by \autoref{L:HiHw}, the $H_i$. Hence it is
  all of $M$. On the other hand, $M_\mu$ is simple, because for an
  $\HHd[\mu](\Bo)$-submodule $0\ne N'\subsetneq M_\mu$ the sum
  $\bigoplus_dH_dN'$ would be a proper nonzero submodule of $M$. The natural map $\Ind_{\HHd[\mu](\Bo)}^{\HHd[n](\Bo)}M_\mu\to M$ is surjective, as its image is $\sum_dH_dM_\mu=M$. By \autoref{C:BasisHd} and as the $H_d$ are invertible, both sides have dimension $|\Sym^\mu|\dim M_\mu$, so the map is an isomorphism. Finally, by \autoref{types}\autoref{typese}, up to parity shift
  $M_\mu\cong Q_1\circledast\dots\circledast Q_{\bar N_\psi}$ with $Q_i$ simple
  over $\HHd[\mu_i](A^{\otimes\mu_i})$. Each $Q_i$ is semisimple over $A^{\otimes\mu_i}$, as $M_\mu$ is $\Bo$-semisimple, and by the
  choice of $\chi^\mu$ its constituents are,
  up to parity shift, those of $\mathbf L_i(1)^{\otimes\mu_i}$. So
  $Q_i\in\mathcal C(\mathbf L_i(1)^{\otimes\mu_i})$ and
  $Q_i\cong\mathbf L_i(1)^{\otimes\mu_i}\otimes_{\mathcal E_{\mu_i}(L_i)}S_i$
  for a simple $\mathcal H_{\mu_i}(L_i)$-module $S_i$ by
  \autoref{T:Morita}. Thus $M\cong\So_\mu$ up to parity shift.
\end{proof}

\begin{Remark}
  The module $  \So_\mu$ is of type $\mathtt Q$, that is isomorphic to its own parity shift, if and only if the number of $S_i$ which are of type $\mathtt Q$ is odd. This follows from \autoref{types}\autoref{typesc}, since the equivalence of \autoref{T:Morita} preserves the type and $\uEnd(  \So_\mu)\cong\uEnd(M)$ by Frobenius reciprocity and \autoref{E:colours}.
\end{Remark}

\begin{Remark}\label{chark}  The assumption $\op{char}(\kk)\not=2$ is heavily used in the super case. If $\Bo$  is an ordinary algebra it can be dropped.
\end{Remark}

\begin{Remark}\label{R:tautological}
  For $A=\kk$ or $A=\Cl_1$ as in \autoref{Ex:Clifford}, there is only one irreducible module $L$ up to isomorphism and parity shift and so there is only one colour. We have $\mathcal E_n\cong\Bo$ as Frobenius-tensor superalgebras and $\gamma_L=1$. The atomic algebra $\mathcal H_n(L)$ is $\HHd[n](\Bo)$ itself, and \autoref{T:simples} is a tautology.
\end{Remark}

\begin{Example}\label{Ex:TaftSimples}
  Let $A=\Taft$ as in \autoref{Ex:Taft}. As $\Taft$ is not symmetric, $\eta=0$ by \autoref{T:FrobTensor}. The simple $A$-modules are all one-dimensional and of the form $L^j$, $j\in\Z/\theta\Z$, with $gv=\omega^jv$ and $yv=0$ for $v\in L^j$. Since $\leftindex^{\psi}L^j\cong L^{j+1}$, there is one colour only and $\theta_L=\theta$ for all $L$. Hence $\gamma_L=0$ by \autoref{L:tauaction}, and the orbit is of wreath type. One can also  see directly that $\tau^{(1,2)}$ acts as zero on all $L^j\otimes L^{j'}$, since $(g^iy^j)^\vee$ is a multiple of $g^{-i}y^{\theta-1-j}$ and $y$ acts as zero.  As $\mathcal E\cong\kk^\theta$ and $\psi_{\mathcal E}$ cyclically rotates the summands, the algebra $\kk[X]\Smash\mathcal E$ is isomorphic to the path algebra of the cyclic quiver with vertex set $\Z/\theta\Z$, and $\mathcal H_n(L)\cong(\kk[X]\Smash\mathcal E)^{\otimes n}\rtimes\kk\Sym$ after rescaling the $H_i$ by $\hbar^{-1}$. If $\operatorname{char}\kk=0$ or $\operatorname{char}\kk>n$, the simple $\HHd[n;\hbar,t,0](\Taft^{\otimes n})$-modules are  classified by the multipartitions $\blam\colon\Z/\theta\Z\sqcup\kk^\times\to\set{\text{partitions}}$ with $\sum_v|\blam(v)|=n$, see e.g. \cite{CT:wreath}. As $\Lo\cong\mathcal E_n$ as a right $\mathcal E_n$-module, the module for $\blam$ has dimension
  $\frac{n!}{\prod_v|\blam(v)|!}\prod_{j}f^{\blam(j)}\prod_{\zeta}\theta^{|\blam(\zeta)|}f^{\blam(\zeta)}$, where $f^{\nu}$ is the number of standard tableaux of shape $\nu$.
\end{Example}

The next example will later be used to illustrate the theory of respectable bases, see  \autoref{Ex:ZigZag}.

\begin{Example}\label{Ex:ZigZagSimples}
We return to the (ordinary) algebra $A=\mathcal{Z}ig_m$ from
\autoref{Ex:ZigZagm} and assume $\kk$ is an algebraically closed field
with $\op{char}(\kk)=p\ge0$, see \autoref{chark}. Since $\psi$ is
trivial, the simple $A$-modules $L_i$, $1\leq i\leq m$, form $m$
colours. The elements $\taui[k,l]$ acts as zero on all $L_i\otimes L_j$
by degree reasons or since each summand of $\taui[k,l]$ has a factor
in the radical of $A$. By  \autoref{L:tauscalar}  $\bar\gamma_{L_i}=0$ and thus by \autoref{gammadef} also $\gamma_{L_i}=0$ for all $i$.  Therefore, all atomic algebras are of wreath type and the $\tau^{(i,j)}$ act as zero on all simple $\HHd[n](\Bo)$-modules. If we assume $X_1$ acts by zero as well, then so does  $X_j$ for any $j$,
since $\hbar^2X_{i+1}=H_iX_iH_i$ by  \autoref{L:HXH}. By \autoref{T:simples} and the representation theory of wreath products \cite{CT:wreath}, the simple $\HHd[n](\Bo)$-modules on which $X_1$ acts by zero
 are labelled by $m$-tuples $\blam=(\lambda^{(1)}|\dots|\lambda^{(m)})$,
 where $\lambda^{(i)}$ is a
 $p$-restricted partition of $n_i\in\N_0$, for $1\le i\le m$, such that
 $\sum_{i=1}^mn_i=n$. (multipartitions are defined in
 \autoref{SS:Tableaux}.) For example, if $m=n=2$  and $p\ne2$, then there are $5$
 irreducible modules, which are indexed by the $2$-multipartitions $(2|0)$,
 $(0|2)$, $(1^2|0)$, $(0|1^2)$, $(1|1)$. Only the last three
 of these $2$-multipartitions occur if $m=n=p=2$.
\end{Example}

\begin{Remark}\label{Classdoesnothing}
\autoref{T:simples} reduces everything to the atomic algebras, but says nothing
  about their representation theory. Instead it \emph{assumes} the classification of simple modules there, see also \autoref{R:tautological}.  For $A=\kk$ this is the
  classification of the simple modules of the degenerate affine Hecke
  algebra in terms of multisegments, \cite{Zelevinsky:p-adic, Rogawski:Hecke},
  for $A=\Cl_1$ it is the
  classification of the irreducible representations of the (degenerate) affine Sergeev algebras in
  \cite{BKHeckeClifford,BK:proj}, and for $A=\Taft$ it is the classical theory
  of wreath products, as in \cite{CT:wreath,JamesKerber}.  The goal of this paper
  is to develop the theory of \emph{respectable algebras} to provide a
  common framework for all of these classification theorems, and where
  the irreducible representations arise as quotients of generalised Specht modules.
  \end{Remark}

 \section{CRaFH superalgebras}\label{S:CycRH}
  In this section we study cyclotomic quotients of RaFH superalgebras. The construction generalises the known constructions for $\HHd(\kk)$ from the literature. The presence of $\Bo$ however causes extra difficulties in the presence of a nontrivial Nakayama automorphism $\psi$ which we will treat in detail in case $n=1$.
 Motivated by these observations, we introduce a notion of
 $\theta$-cyclotomic polynomials
 $\fcyc(X_1)=f_{\text{cyc},\ell_0,\ell_1,\charge}(X_1)$. These
 are certain cyclotomic polynomials that depend on the order $\theta$ of $\psi$  and a choice of a primitive $\theta$th root of unity in $\RR$.  We then study their \emph{CRaFH-algebras} by which we mean the associated cyclotomic quotients of the RaFH-algebras, i.e.,  the quotients of $\HHd(\Bo)$ by ideals generated by $\theta$-cyclotomic polynomials $\fcyc(\ell_0,\ell_1,\charge)$. The main result is a basis theorem for CRaFH-superalgebras.\\

Throughout this section let $\RR$ be a commutative integral domain that
is a $\kk$-algebra and contains elements $\hbar$, $t$, and $\eta$.  Let $\Bo$ be an $\eta$-symmetric Frobenius-tensor superalgebra with Nakayama automorphism $\psi$ of order $\theta$, and consider the  \emph{RaFH superalgebra} $\HHd(\Bo)=\HHd[n;\hbar, t, \eta](\Bo)$.\\
 \emph{From now on we assume, if not stated otherwise, that $\hbar$ is invertible.}

  \subsection{$\theta$-cyclotomic polynomials and CRaFH superalgebras}
 Consider the RaFH superalgebras $\HHd(\Bo)$ from  \autoref{D:Hecke} for $\Bo$ a Frobenius-tensor superalgebra with Nakayama automorphism $\psi$ of finite order $\theta$.  We assume that the commutative ring $\RR$ contains a primitive $\theta$th root of unity $\omega$ which we fix.

 \begin{Definition}\label{D:fcyc}
   Let $\ell_0,\ell_1\in\N_0$ and pick nonzero scalars $(\kapp_1,\ldots,
   \kapp_{\ell_1})\in \RR^{\ell_1}$. The corresponding
   \emph{$\theta$-cyclotomic polynomial}
   $\fcyc(X_1)=f_{\text{cyc},\ell_0,\ell_1,\charge}(X_1)$  is the polynomial
   \begin{equation}\label{cyclotomicpoly} \fcyc(X_1)=X_1^{\ell_0}\prod_{i=1}^{\ell_1} \prod_{j=1}^\theta(X_1-\omega^j\kapp_i)= X_1^{\ell_0}
         \prod_{i=1}^{\ell_1}(X_1^\theta-\kapp_i^\theta)\in\RR[X_1].
    \end{equation}
   The \emph{level} of $\fcyc $ is $\ell=\ell_0+\theta\ell_1$, and
   $\charge=(0,\ldots, 0,\kapp_1,\ldots,\kapp_{\ell_1})\in
   R^{\ell_0+\ell_1}$ is its \emph{charge}. The level refines to the
   \emph{level vector} $\ell_\theta=(\ell_0,\ell_1)$.  Note that $\fcyc$
   is determined by $\charge$ (and $\theta$), so we sometimes write $\fcyc=\mathtt{f}_{\charge,\theta}$.
    \end{Definition}

 \begin{Remark} If $\psi$ is trivial, then $\theta=1$, $\ell=\ell_0+\ell_1$ and we obtain the usual cyclotomic polynomials.
 \end{Remark}

\begin{Remark}\label{ftoc}
  Let $\ell'\in \mathbb{N}_0$ and  $\charge\in\RR^{\ell'}$. Then there is an associated  $\theta$-cyclotomic polynomial $\fcyc=\mathtt{f}_{\charge,\theta}$, where the $\kapp_{1},\ldots, \kapp_{\ell_1}$  are the nonzero entries of $\charge$ and $\ell_0=\ell'-\ell_1=|\set{i|\kapp_i=0}|$.   Conversely, any $\theta$-cyclotomic polynomial $\fcyc$ is of the form $\mathtt{f}_{\charge,\theta}$ for some charge vector $\charge$, which is
  unique up to a permutation and up to multiplication of the entries with powers of $\omega$.
  \end{Remark}

  \begin{Definition}\label{D:CyclotomicHecke}
 Let  $\charge\in\RR^{\ell'}$ for some $\ell'\in\mathbb{N}$. The associated \emph{CRaFH superalgebra}  $\Hhd=\Hhd[n;\hbar, t, \eta](\Bo)$ of charge $\charge$,  with parameters $\hbar,t,\eta$, is the quotient $R$-superalgebra
   $\Hhd=\HHd/\Ic$ of $\HHd(\Bo)=\HHd[n;\hbar, t, \eta](\Bo)$ by the \emph{cyclotomic ideal} $\Ic=(\fct)$ generated by the $\theta$-cyclotomic polynomial $\fct$.
  \end{Definition}

 When there is no confusion, we abuse notation and identify elements of
$\HHd$, such as $\bb\in\Bo$, $X_j$, $H_i$,  $H_w$
with their images in $\Hhd$ under the canonical surjection
$\HHd\twoheadrightarrow\Hhd$.

\begin{Remark}\label{Rk:levelone}
The name \emph{CRaFH superalgebra} is an abbreviation for \emph{cyclotomic RaFH
superalgebras}.  We prefer to use the term CRaFH superalgebras over the term cyclotomic RaFH algebras, since it reflects
the fact that only $\theta$-cyclotomic ideals are allowed in this construction. For instance, the CRaFH superalgebras for $\ell=1$ allow by definition any cyclotomic quotient $\Hhd(\Bo)$ if the Nakayama automorphism is trivial, but only charge zero, i.e. $\Hhd(\Bo)=\Hfind(\Bo)$, if it is not.
\end{Remark}

 \begin{Remark}\label{motcyc}
\autoref{D:CyclotomicHecke} generalises several versions of cyclotomic
Hecke superalgebras from the literature and provides, via the shape
\autoref{cyclotomicpoly}  of the cyclotomic polynomial, a natural
framework explaining the existing ad-hoc conditions on the cyclotomic
polynomial in the presence of a non-trivial Nakayama automorphism. See, for example, \cite[(x3-e)]{BKHeckeClifford} for  the  cyclotomic Sergeev superalgebra.
\end{Remark}

All the RaFH superalgebras from \autoref{HHdexamples} have CRaFH quotients.
The Ariki--Koike algebras are other famous examples of CRaFH algebras
that, historically, were first defined by generators and relations, generalising the Hecke algebras of finite Coxeter groups, rather than as cyclotomic quotients.

 \begin{Example}\label{Ex:CRaFH}
The \emph{Ariki--Koike algebras} are deformations of the
      group algebras of the complex reflection groups of type
      $G(r,1,n)$~\cite{AK,BM:cyc}. These algebras were originally
      defined as the $\kk$-algebras with generators
     $T_0,T_1,\dots,T_{n-1}$ subject to the type~$B$ braid relations
      together with $(T_0-c_1)\dots(T_0-c_r)=0$ and where
      $T_1,\dots,T_{n-1}$ satisfy the quadratic relations of the
      Iwahori--Hecke algebras. The Ariki--Koike algebras fit into our
      CRaFH framework by taking $\RR=\kk$ and the trivial Frobenius algebra $A=\kk=\RR$,
      $\hbar=1$, $t=0$ and $T_0=X_1, T_1=H_1,\dots,T_{n-1}=H_{n-1}$.
      In light of \autoref{L:HXH}, the type~$B$ braid relation
      $T_0T_1T_0T_1=T_1T_0T_1T_0$ corresponds to the RaFH superalgebra
      relation $X_1X_2=X_2X_1$. Similarly, the
      \emph{cyclotomic Hecke algebras of type $A$} from \cite{HuMathas:SeminormalQuiver}, which include the degenerate and non-degenerate Ariki--Koike algebras as
      special cases, are the CRaFH superalgebras
     $\Hd[n;\sqrt{\xi},1,\xi-1](\kk)$ under the
      identification $X_i=L_i$, $H_i=T_i$, since the generators $\xi,
      T_i, L_i$ satisfy the relations of \cite[Definition~2.1]{HuMathas:SeminormalQuiver}.
      \end{Example}
      \begin{Remark}
     The \emph{Hu algebra} $\CA(n)$ from
     \cite{Hu:TypeDMorita}  was introduced to describe the Hecke
     algebras of type~$D$ in even rank by deforming the group algebra of
     $\Sym_n\wr \Sym_d$. These algebras were further generalised in
     \cite[\S4.1]{CJNakanoXiang}. The algebra  $\CA(n)$ does not
     fit into our current framework. However, the Hu algebras are
     captured  if we make a small change to
     \autoref{RaFHrels}. More explicitly, to obtain the Hu algebra
     suppose that $\eta=0$ and replace the quadratic relations
     $H_i^2=\hbar^2$ with $H_i^2=\hbar^2 z^{(i,i+1)}$, where $z\in
     Z(A\otimes A)$ with $\leftindex^{s_1}z=z$. Here $z^{(i,i+1)}$ is
     the embedding of $z$ into the $i$th and $(i+1)$st tensor factors.
      There is the obvious generalisation $\Hfind(\Bo,z)$ of the Rees
     Frobenius Hecke superalgebra $\Hfind(\Bo)$ from \autoref{D:RFHS}.
     Taking $z=z_{n,n}$ in the notation from \cite{Hu:TypeDMorita} and
     \cite{CJNakanoXiang} we see that $\Hfind(\Bo,z)\cong\CA(n)$.
     Moreover, by \cite[\S4.1]{CJNakanoXiang}, the expected basis
     theorem, as in \autoref{P:H0Basis}, holds. We
     exclude algebras like these because their representation theory is
     quite different.
\end{Remark}

\subsection{CRaFH superalgebras: the $n=1$ case made explicitly}
Together with \autoref{motcyc}, the following result in the
$n=1$ case motivates our definition of
$\theta$-cyclotomic polynomials $\fct$ in \autoref{cyclotomicpoly}. Its part \autoref{n11} explains the restriction on level $1$ quotients from \autoref{Rk:levelone}.

\begin{Proposition}\label{nequal1} Consider the RaFH superalgebras
$\HHd[1](\Bo)$. Let $\cB$ be a basis of $\Bo=A$. Pick $\kapp\in\RR^{\ell'}$ and
consider the corresponding $\theta$-cyclotomic polynomial $\fct$. Then the following holds:
\begin{enumerate}
\item\label{n11} $\cB$ is an $R$-basis of
$\HHd[1](\Bo)/(X_1-\kapp)$ if and only if  either $(\psi=\op{id})$  or
$(\psi\not=\op{id}$ and $\kapp=0)$.
\item\label{n13}  If $\kapp$ is invertible and  $\theta\nmid\ell$ then
$\{X_1^d b\mid 0\leq d< \ell,  b\in\cB\}$ is not an $R$-basis of
$\HHd[1](\Bo)/(X_1-\kapp)^\ell$.
\item\label{n12}  If  $\fct$ is of level $\ell$, then
$\cB_{\op{cyc}}=\{X_1^d b\mid 0\leq d< \ell,  b\in\cB\}$ is an $R$-basis of  $\Hhd[1]=\HHd[1]/\Ic=\HHd[1]/(\fct)$.
\end{enumerate}
\end{Proposition}
\begin{Remark}\autoref{nequal1} shows that the  CRaFH superalgebra
quotients of $\HHd[1](\Bo)$ have a basis as expected, whereas quotients
by ordinary cyclotomic ideals might have smaller rank than desired.
\end{Remark}

\begin{proof} The superalgebra $\HHd[1](A)$ is generated by $X_1$ and $A$ subject to $aX_1=X_1\psi(a)$ for $a\in A$.

Assume now $\psi\not=\op{id}$. Suppose we are in the setup of
\autoref{n13} and let $f=(X_1-\kapp)^\ell=\sum_{i=0}^{\ell}
  \binom{\ell}{i}(-\kapp)^{\ell-i}X_1^i$.  Since $\psi$ has order $\theta$ and $\theta\nmid\ell$, there exists $a_0\in A$ such
that $\psi^\ell(a_0)\neq a_0$. We compute in $\HHd[1](\Bo)/(X_1-\kapp)^\ell$ that, for any $a\in A$,
  \begin{equation}\label{sumzero}
  0=f\psi^\ell(a)-af
 =
 \sum_{i=0}^{\ell}
 \binom{\ell}{i}(-\kapp)^{\ell-i}
 \bigl(\psi^{\ell-i}(a)-a\bigr)X_1^i  \\
 =
 \sum_{i=0}^{\ell-1}
 \binom{\ell}{i}(-\kapp)^{\ell-i}
 \bigl(\psi^{\ell-i}(a)-a\bigr)X_1^i.
\end{equation}
If  $\cB$ is a basis, then the coefficient of $X_1^0$ is nonzero for  $a=a_0$ contradicting \autoref{sumzero}.
This shows  \autoref{n13}.
\autoref{n11}: Let first  $\kapp_1\neq0$ and pick $a\in A$ with $\psi(a)\neq a$. Working in $A_\RR$, write $a-\psi(a)=\sum_{b\in\cB}\beta_bb$ with $\beta_b\in\RR$, not all zero. In $\HHd[1](A)/(X_1-\kapp_1)$ the relation $aX_1=X_1\psi(a)$ gives $\kapp_1(a-\psi(a))=0$.  As $\RR$ is a domain containing $\kk$, either $\cB$ is not a basis or  $\kapp_1=0$. In the latter case,  $(X_1)=X_1\HHd[1](A)$ and $\HHd[1](A)/(X_1)\cong\RR\otimes_\RR A_\RR=A_\RR$.

 We identify polynomials $f\in R[X_1]$ of degree $<\ell$  with elements
of  $R^\ell$ by sending $f=\sum_{i=0}^{\ell-1}\beta_iX_1^i$ to the
vector $(\beta_{\ell-1},\ldots, \beta_0)$  formed by the coefficients.
Consider the $R$-linear map $\Xi\colon R[X_1]\otimes A\to
R^\ell\otimes A$, $g\otimes a\mapsto \op{r}(g)\otimes a$ defined by
sending  $g\in R[X_1]$ to the remainder $\op{r}(g)$ of the division of
$g$ by  $\fct$.  By definition,  $\Ic\otimes A$ is sent to zero.  Recall from
\autoref{D:Hecke} that $a X_1 =X_1
\psi(a)$, thus $a X_1^\theta =X_1^\theta a$, for $a\in A$. The shape \autoref{cyclotomicpoly} of $\fct$ now implies that $\Xi$ induces a well-defined $R$-linear map
 $\Hhd[1]\to R^\ell\otimes A$, which sends $\cB_{\op{cyc}}$ to a linearly independent set.  This finishes the proof of \autoref{n12}.
\end{proof}

 \subsection{Basis theorem for CRaFH superalgebras}
We next establish a basis theorem for the CRaFH superalgebras $\Hhd(\Bo)$
generalising the $n=1$ case from \autoref{nequal1}. The result gives  a unified framework for
similar results in the literature. It simultaneously covers the CRaFH
algebra quotients of the algebras listed in \autoref{HHdexamples} and also the Ariki--Koike algebras from \autoref{Ex:CRaFH}.

To formulate the basis theorem we need more notation. For  $\bd=(d_1,\ldots,d_n)
\in\mathbb{N}_0^n$ set  $X^\bd=X_1^{d_1}\cdots X_n^{d_n}$, the corresponding monomial in~$\Poo$.
 Let $\cBo=\set{b_1\otimes\dots\otimes b_n|b_1,\dots,b_n\in\cB}$ be the basis of $\Bo$ induced from $\cB$, and set
 $\N^{\ell,n}_0=\set{\bd\in\N^n|0\le d_i<\ell\text{ for }1\le i\le n}$.

  \begin{Theorem}\label{basisthmcycl}
    The  image of $\set{X^\bd H_w\bb|w\in\Sym, \bd\in\N^{\ell,n}_0, \bb\in\cBo}$ is a basis of $\Hhd=\HHd/\Ic=\HHd/(\fct)$.
  \end{Theorem}

Before giving the proof, we note that the argument of \autoref{C:BasisHd} now implies the following:

\begin{Corollary}\label{C:HdBasis}
  The CRaFH superalgebra $\HHd/(\fct)$ is free as an $R$-module with bases
  \begin{align*}
        &\set{\bb X^\bd H_w| \bd\in\N^{\ell,n}_0, w\in\Sym, \bb\in\cBo},
        &\set{H_w X^\bd \bb| \bd\in\N^{\ell,n}_0, w\in\Sym, \bb\in\cBo},\\
        &\set{ X^\bd  \bb H_w| \bd\in\N^{\ell,n}_0, w\in\Sym, \bb\in\cBo},
        &\set{H_w \bb X^\bd | \bd\in\N^{\ell,n}_0, w\in\Sym, \bb\in\cBo}.
  \end{align*}
\end{Corollary}

As an application, we deduce the following consistency result (using that $\hbar$ is invertible).

\begin{Corollary} \label{notationclash} The Rees-Frobenius Hecke superalgebra from \autoref{D:RFHS} and \autoref{matchRFH} is isomorphic to the cyclotomic quotient $\Hhd(\Bo)=\Hfind(\Bo)$ of level $1$ and charge $0$.
\end{Corollary}
\begin{proof}The defining relations from \autoref{D:RFHS} give a surjection which must be an isomorphism since a basis is sent to a basis by  \autoref {P:H0Basis} and \autoref{basisthmcycl}.
\end{proof}

\begin{Remark}\label{badh}We exclude the case $\hbar=0$, since the cyclotomic quotients $\Hhd$ are not well-behaved and can be  infinite-dimensional even if $\ell=1$.  To see this consider the quotient $\Hhd=\HHd[2;0, 0, \eta](\kk)/(X_1)$. If  $\hbar=t=0$, then $H_1^2=\eta H_1$ and $\gamma_1=\eta X_2$. Modulo $(X_1)$, the
  relations $H_1X_1=X_2H_1-\gamma_1$ and $H_1X_2=X_1H_1+\gamma_1$ become $X_2(H_1-\eta)\equiv0\equiv(H_1-\eta)X_2$. Since $H_1-\eta$ is not invertible when $\hbar=0$,  we cannot solve for $X_2$ and the powers $X_2^{k}$ survive in the quotient $\Hhd$. To see this note that $R[z]$ is an $\Hhd$-module with the action of $H_1$, $X_1$, $X_2$  given by multiplication with $\eta, 0, z$, respectively. Since the powers of $z$ are linearly independent in $R[z]$, the powers of $X_2$ are  linearly independent in $\Hhd$.
This example shows that \autoref{basisthmcycl} needs $\hbar\in
\RR^\times$ if $n\ge2$. (It is used in the proof to ensure the existence
of $H_i^{-1}$). For $n=1$, the desired basis theorem holds even when
$\hbar$ is not invertible.
    \end{Remark}

\begin{Remark}\label{badh2}
By \autoref{badh}, \autoref{notationclash} does not hold if $\hbar=0$, $n\geq 2$, even  in case $A=\kk$, where \autoref{D:RFHS} is  the  \emph{$0$-Hecke algebra} from e.g. \cite{DGHT, Fayers0Hecke}.
As we saw in \autoref{badh}, putting conditions on $X_1$ does not necessarily put constraints on $X_i$ for $i\geq 2$.
A more appropriate definition of cyclotomic quotient  here might be
quotients $Q^{m}_c:=\HHd[n;0, t,
\eta](\Bo)\big/\bigl(X_1^{m}-\kapp,\dots,X_n^{m}-\kapp\bigr)$
where $m\ge1$ and $\kapp\in R$ such that  $\kapp=0$ or $\theta\mid m$. If  $\hbar=t=\eta=0$, then one can verify that
  $Q^{m}_c\cong\bigl(R[X_1,\dots,X_n]/(X_1^{m}-\kapp,\dots,X_n^{m}-\kapp)\Smash\Bo\bigr)\Smash\NC_n$,
  where $\NC_n$ is the NilCoxeter algebra from \cite{NilCoxeter}. Moreover, $\set{X^{\bd}\bb H_w\mid 0\le d_i<m,\ \bb\in\text{a basis of }\Bo,\ w\in\Sym}$ is an $R$-basis of $Q^{m}_c$. In particular $Q^{m}_c$ is free of rank
  $n!\,m^n r^n$. In case $\hbar=0\not=\eta$ the situation is much
  more involved. The naive basis theorem does not hold and the algebra
  is not necessarily free over $\Bo$.
\end{Remark}

\begin{Remark}\label{R:KMbasisconj}
  Let $\Bo=A^{\otimes n}$ with $A$ a symmetric Frobenius algebra, so that $\theta=1$, and identify
  $\HHd[n;1,1,0](\Ao)$ with the affinisation $H_n(A)$ as in
  \autoref{Ex:KMaffinisation}. Then $\Hhd=\HHd/(\fct)$ is the level $\ell$
  cyclotomic quotient of \cite[\S3.6]{KleshchevMuth} for the
  charge $\charge$. (Note that \cite{KleshchevMuth} allows more generally a
  sequence of homogeneous central elements of $\Ao$, whereas our $\charge$
  consists of scalars.) The spanning set of
  \cite[Proposition~3.19]{KleshchevMuth} is exactly the third basis of
  \autoref{C:HdBasis}, so \autoref{basisthmcycl} proves
  \cite[Conjecture~3.21]{KleshchevMuth} for scalar charge. At level $\ell=1$
  both sides reduce to the wreath product
  $\Ao\rtimes\kk\Sym=\Hfind[n;1,0](\Ao)$ of \autoref{D:RFHS}, which is
  \cite[Proposition~3.17]{KleshchevMuth}.
\end{Remark}

\subsection{Proof of the Basis \autoref{basisthmcycl} for CRaFH superalgebras}
For the proof of the basis theorem for CRaFH superalgebras we follow arguments from \cite{BKHeckeClifford}.

Recall from \autoref{D:RFHS} and \autoref{DefPo} the algebras $\Hfind=\Hfind(\Bo)$, $\Poo=\Po$  with the identification from \autoref{matchRFH}, but without  the (so far unproven) identification from \autoref{notationclash}.   Let $\Po_i\subset\Poo$ be the subalgebra generated by $X_1,X_2,\ldots, X_i$.  Set $f_1:=\fcyc=\fct$ as in \autoref{cyclotomicpoly} and inductively define $\mathtt{f}_{i+1}=H_{i}^{-1}\mathtt{f}_iH_i\in \HHd$ for $1\leq i\leq n-1$.

  \begin{Lemma}\label{fis}
For $1\leq i\leq n$ we have $\mathtt{f}_i\in X_i^\ell+\sum_{k=0}^{\ell-1}\Poo_{i-1}X_i^k\Hfind(\Bo)$.
\end{Lemma}
\begin{proof}The statement holds by definition for $i=1$. We continue by induction assuming the lemma holds for $i$.  Set $J_{i+1}=\sum_{k=0}^{\ell-1}\Po_{i}X_{i+1}^k\Hfind(\Bo)$.
Since $\mathtt{f}_{i+1}=H_i^{-1}\mathtt{f}_iH_i$ it suffices to show that
\begin{equation}\label{annoying}
H_i^{-1}X_i^\ell H_i\in X_{i+1}^\ell+J_{i+1}\quad\text{and}\quad H_i^{-1}gX_{i}^r h H_i\in J_{i+1}\text{ for any $g\in \Po_{i-1}$, $0\leq r\leq \ell-1$, $h\in \Hfind(\Bo)$}.
\end{equation}

By  \autoref{L:RewriteXr} we have
 $H_iX_i^k=X_{i+1}^kH_i - (t+\eta X_{i+1})\sum_{j=0}^{k-1}X_{i+1}^{k-j-1}\tau^{(i,i+1)}X_i^j$ for any nonnegative integer $k$. This is equivalent to
  \begin{equation*}
  (H_i-\eta \tau^{(i,i+1)}) X_i^k=X_{i+1}^k(H_i-\eta \tau^{(i,i+1)}) - tX_{i+1}^{k-1} \tau^{(i,i+1)}-  \left(( t+\eta X_{i+1})\sum_{j=1}^{k-1}X_{i+1}^{k-j-1}\tau^{(i,i+1)}X_i^j\right)-\eta \tau^{(i,i+1)}  X_i^k.
  \end{equation*}
This equation implies that $H_i^{-1}X_i^\ell H_i\in X_{i+1}^\ell+J_{i+1}$ and $H_i^{-1}X_i^k H_i\in J_{i+1}$ for $k< \ell$,
  because  $\hbar^{-2}(H_i-\eta \tau^{(i,i+1)})=H_i^{-1}$ and  elements in $\Bo$ commute with the $X$'s up to a twist with the Nakayama automorphism. This proves the first formula in \autoref{annoying}. For the second just note that
  $H_i^{-1}gX_i^k hH_i=gH_i^{-1}X_i^kH_iH_i^{-1}hH_i\in gJ_{i+1}H_i^{-1}hH_i\subseteq J_{i+1}.$
 \end{proof}
   \begin{Lemma}\label{HfnH}
 We have $\Hfind[n-1] \mathtt{f}_{n}\Hfind = \mathtt{f}_{n}\Hfind$.
 \end{Lemma}
 \begin{proof}
 We only need to show the inclusion $\subseteq$. For this it is enough to see that left multiplication with $b\in\Bo$ and $H_i^{-1}\in\Hfind[n-1]$ leave $\mathtt{f}_{n}\Hfind $ invariant, since these elements generate $\Hfind[n-1]$, see \autoref{matchRFH}.  For $b\in\Bo$ this is obvious, since $b$ commutes with $\mathtt{f}_n$ up to a twist by $\psi^{\ell_0}$.  For $H_i^{-1}$, for some $1\leq i\leq n-2$, we calculate $H_i^{-1} \mathtt{f}_n\Hfind =H_i^{-1}H_{n-1}^{-1}H_{n-2}^{-1}\cdots H_{i+1}^{-1}H_{i}^{-1}\mathtt{f}_i \Hfind=H_{n-1}^{-1}H_{n-2}^{-1}\cdots H_{i+2}^{-1}(H_{i+1}^{-1}H_{i}^{-1}H_{i+1}^{-1})\mathtt{f}_i \Hfind=H_{n-1}^{-1}H_{n-2}^{-1}\cdots H_{i+2}^{-1}H_{i+1}^{-1}H_{i}^{-1}\mathtt{f}_i \Hfind=\mathtt{f}_n\Hfind$. Here we used that the $H_i^{-1}$ satisfy the braid relations.
  \end{proof}

 Consider sets $\Pi^{n,\ell}_+$, $\Pi^{n,\ell}$ of pairs $(\bd,Z)$, where $\bd=(d_1,\ldots,d_n) \in\mathbb{N}_0^n$ and $Z\subset\{1,2,\ldots, n\}$, as follows:
 \begin{equation}\label{domweights}
   \Pi^{n,\ell}=\set{(\bd,Z)| d_i<\ell\text{ if }i\not\in Z}\qquad\text{and}\qquad
   \Pi^{n,\ell}_+=\set{(\bd,Z)|Z\ne\emptyset, d_i<\ell \text{ if }i\notin Z}.
 \end{equation}

 To  $Z=\{z_1<z_2<\cdots <z_r\}\subset\{1,2, \ldots, n\}$ associate $\mathtt{f}_Z=\mathtt{f}_{z_1}\mathtt{f}_{z_2}\cdots \mathtt{f}_{z_r}$.

  \begin{Lemma}\label{free}
    As a  right $\Hfind$-module,  $\HHd$ is free with basis
    $\set{X^\bd \mathtt{f}_Z | (\bd,Z)\in \Pi^{n,\ell}}$.
  \end{Lemma}
  \begin{proof}
  Note that there is a bijection
  \[\nu\colon\Pi^{n,\ell} \rightarrow \mathbb{N}_0^n\quad  (\bd,Z)\in\Pi^{n,\ell}\mapsto  \nu(\bd,Z), \quad\text{with} \quad
  \nu(\bd,Z)_i=
  \begin{cases}
  d_i+\ell, &\text{if $i\in Z$},\\
  d_i, &\text{otherwise}.
  \end{cases}
  \]

  Next we claim that $X^\bd \mathtt{f}_Z\in X^{\nu(\bd,Z)}+Y$, where $Y$
  is the span of lower order terms of the form $X^\beta h$, for $h\in\Hfind$ and $\beta<\nu(\bd,Z)$ in the reversed lexicographic ordering. This claim is clear for $n=1$ by the definition of $\mathtt{f}_1$.

Writing $X^\bd
\mathtt{f}_Z=X_n^{d_n}(X^{\bd'}\mathtt{f}_{Z'})\mathtt{f}_n$ with
$Z'=Z-\{n\}$ and $\bd'=(d_1,\ldots,d_{n-1})$, we can apply the induction
hypothesis to the term $X^{\bd'}\mathtt{f}_{Z'}$,  after which the claim
follows using \autoref{fis}.  Putting these things together we have
shown that $\HHd$ is a free right $\Hfind$-module with basis
$\set{X^\bd \mathtt{f}_Z |(\bd,Z)\in \Pi^{n,\ell}}$ if it is a free
right module with basis $\set{X^\bd | \bd\in \mathbb{N}_0^n}$. But the latter holds by  \autoref{C:BasisHd}.
 \end{proof}

Let $\Ic$ be the two-sided ideal of $\HHd$ generated by $\mathtt{f}_1$.

\begin{Lemma}\label{peppermint}
As an $R$-module, $\Ic=\sum_{i=1}^n\Poo \mathtt{f}_i\Hfind$.
\end{Lemma}
\begin{proof}
Set $\mathcal{J}=\sum_{i=1}^n\Poo \mathtt{f}_i\Hfind$. Clearly, $ \mathcal{J}\subseteq  \mathcal{I}$ by definition of $\mathtt{f}_i$.  We need to show the opposite inclusion. Recall from  \autoref{C:BasisHd}  that $\HHd=\Poo\Hfind$. Now we calculate

\[\Ic=\HHd \mathtt{f}_1\HHd= \HHd \mathtt{f}_1 \Poo\Hfind =  \HHd \mathtt{f}_1 \Hfind =   \Poo\Hfind \mathtt{f}_1 \Hfind \subseteq \sum_{j=0}^{n-1} \Poo H_{j}^{-1}\cdots H_{2}^{-1}H_{1}^{-1} \mathtt{f}_1 \Hfind\subseteq \sum_{i=1}^n\Poo \mathtt{f}_i\Hfind.\]

The last inclusion just uses the definition of $\mathtt{f}_i$. For the  penultimate inclusion we can first  push all elements in $\Bo$ to the right. Then it is enough to show that $H_w\mathtt{f}_1 \Hfind \in J:=\sum_{j=0}^{n-1} H_{j}^{-1}\cdots H_{2}^{-1}H_{1}^{-1} \mathtt{f}_1 \Hfind$. By induction on the length of $w$ it is enough to show that $H_{w^{-1}}^{-1}\mathtt{f}_1 \Hfind\subseteq J$.
To see this note that every element $w\in\Sym$ can be written uniquely as $w=zy$ for some
${z=s_js_{j-1}\cdots s_1}$ and $y\in\langle s_i\mid 2\leq i\leq n-1\rangle$. The analogous statement holds then also for $H_w$, i.e. $H_w=H_{z
}H_y$ and for $H_{w^{-1}}^{-1}=H_{z^{-1}}^{-1}H_{y^{-1}}^{-1}$. But $H_{z^{-1}}^{-1}H_{y^{-1}}^{-1}\mathtt{f}_1\Hfind=H_{z^{-1}}^{-1}\mathtt{f}_1\Hfind\subseteq  J$.
\end{proof}

\begin{Lemma}\label{IfinfZ}
$\Ic=\sum_{(\bd,Z)\in\Pi^{n,\ell}_+}X^\bd \mathtt{f}_Z\Hfind$.
\end{Lemma}
\begin{proof} Let $J=\sum_{(\bd,Z)\in\ \Pi^{n,\ell}_+}X^\bd \mathtt{f}_Z\Hfind$.
Clearly, $J\subseteq \Ic$.  Thus, by \autoref{peppermint}, it is enough
to show that $X^\bd \mathtt{f}_i\Hfind\subseteq J$ for any $\bd\in\mathbb{N}_0^n$. We argue by induction on $n$. In case $n=1$, the statement holds by the definitions.   Write  $X^\bd=X_n^{d_n}X^\beta$ for $\beta=(d_1,\ldots, d_{n-1})$.

{\it Case $i=n$:}  By expanding $X^\beta$ in the basis from \autoref{free}, we obtain by induction and \autoref{HfnH}

\[X^\bd \mathtt{f}_n\Hfind\subseteq\sum_{(\bd',Z')\in\Pi^{n-1,\ell}}X_n^{d_n}X^{\bd'}\mathtt{f}_{Z'}\Hfind[n-1]\mathtt{f}_n\Hfind\subseteq \sum_{(\bd',Z' )\in\Pi^{n-1,\ell}}X_n^{d_n}X^{\bd'}\mathtt{f}_{Z'}\mathtt{f}_n\Hfind\subseteq J.\]

{\it Case $i<n$:}  Here we get $X^\bd \mathtt{f}_i\Hfind\subseteq\sum_{(\bd',Z')\in\Pi^{n-1,\ell}_+}X_n^{d_n}X^{\bd'}\mathtt{f}_{Z'}\Hfind\subseteq \sum_{(\bd',Z' )\in\Pi^{n,\ell}}X_n^{d_n}X^{\bd'}\mathtt{f}_{Z'}\Hfind$. It remains to show that $X_n^{d_n}X^{\bd'}\mathtt{f}_{Z'}\Hfind\subseteq J$. This is clear if $d_n<\ell$.  Thus assume $d_n\geq \ell$. From  \autoref{fis} we obtain that
$X_n^{d_n}X^{\bd'}\mathtt{f}_{Z'}= X_n^{{d_n-\ell}}X^{\bd'}\mathtt{f}_{Z'}\mathtt{f}_n+\op{Rest}$, where $\op{Rest}$ is a linear combination of terms from $X_n^{d_n-\ell+j}\mathtt{I}_{n-1}^\charge\Hfind$, where $\Ic[n-1]$ is the two-sided ideal in $\HHd[n-1]$ generated by $\mathtt{f}_1$. Now the claim follows by induction and the definition of $J$.
\end{proof}

\begin{proof}[Proof of Basis \autoref{basisthmcycl}]
By \autoref{IfinfZ}, the elements $X^\bd \mathtt{f}_Z$, for
$(\bd,Z)\in\ \Pi^{n,\ell}_+$, span $\Ic$ as a right
$\Hfind$-module. By \autoref{free}, these elements are linearly independent and
the $X^\bd \mathtt{f}_Z$, for $(\bd,\emptyset)\in\Pi^{n,\ell}$, and hence
$\set{X^\bd|\bd\in\N^{\ell,n}_0}$ is a basis of a
complement to $\Ic$ as a right $\Hfind$-module. By \autoref{P:H0Basis} with \autoref{matchRFH},
 the set $\set{H_w \bb|w\in\Sym, \bb\in\cBo}$ is a basis of
$\Hfind$. Hence, the theorem follows.
\end{proof}
\section{Trace forms on CRaFH superalgebras}
In this section we introduce trace maps on CRaFH superalgebras. They generalise the known trace maps on the cyclotomic quotients of the algebras in \autoref{HHdexamples}. We give a criterion for the associated bilinear form to be symmetric and determine its Nakayama automorphism. We then study non-degeneracy and unimodularity. For $n=1$ the Gram determinant can be computed explicitly, and we determine its degeneracy locus. This small case serves as a guide for general $n$. Our main tool is the \emph{shifting isomorphism}, which reduces the case of invertible $\eta$ and $t$ to $t=0$, so that only the cases $\eta=0$ and $t=0$ remain. These two situations are then treated with different methods. For $t=0\neq\eta$ we use a triangularity argument in the spirit of \cite{MalleMathas}. For $\eta=0$ we apply a filtration argument similar to \cite{BK:HigherSchurWeyl}. It allows us to transfer the computations to the associated graded algebra. The latter  we realise as the Rees Frobenius Hecke superalgebra $\Hfind[n;\hbar,0,0]$ in the sense of  \autoref{D:RFHS}. Here the Frobenius-tensor superalgebra is $\Hhd[1](A)^{\otimes n}$, whose underlying Frobenius algebra is the CRaFH superalgebra $\Hhd[1](A)$ for $n=1$ with cyclotomic polynomial $X_1^\ell$ and the trace from \autoref{L:Aprime}. Together, these yield the Gram determinant for all parameters up to a sign and a unit, see \autoref{T:summary}.

    Throughout  this section let $\RR$ be a commutative integral domain which is a $\kk$-algebra, and let $\hbar,t,\eta\in \RR$.  We assume that $\Bo=A^{\otimes n}$ is a Frobenius-tensor superalgebra such that $A$ is equipped with a linear non-degenerate trace map $\tr\colon A\to\kk$ turning $A$ into a Frobenius algebra. Then ${\Bo}=A^{\otimes n}$ is canonically a Frobenius algebra, see
   \autoref{tensorFrob}.  We abuse notation and use $\tr$ and $\psi$ for
   either of the traces {and} Nakayama automorphisms, respectively. We will
   implicitly use the notation and terminology from \autoref{S:Frobeniusalgs}. We consider the CRaFH superalgebras  as in \autoref{D:CyclotomicHecke}. We assume in particular, if not stated otherwise, that $\hbar$ is invertible.
     \subsection{The definition of the trace maps}
  We consider the CRaFH $\RR$-superalgebra
\begin{equation}\label{ass}
\Hhd=\Hhd[n;\hbar, t, \eta](\Bo){=\HHd/\Ic=\HHd/(\fct)},\quad \text{where $\Bo=A^{\otimes n}$  and $A$ is Frobenius,}
\end{equation}
with parameters $\hbar,t,\eta$ of charge $\charge$ and level $\ell$ and level vector $(\ell_0,\ell_1)$. We set $\ELL=(\ell-1,\dots,\ell-1)$ and $N=n(\ell-1)$.
For $\bd=(d_1,\dots,d_n)\in\N^n$, let $d=d_1+\dots+d_n$. By  convention, $0^0=1$.

 \begin{Definition}\label{defform}
   Let $\Tr\map\Hhd R$ be the $R$-linear map defined on the
   basis from \autoref{basisthmcycl} by
   \begin{equation}\label{formform}
     \Tr(X^\bd H_w\bb) =
     \begin{cases*}
        \delta_{w,e} \delta_{\bd, \ELL}\tr(\bb),& if $\eta=0$,\\
        \delta_{w,e} \delta_{\bd, 0}\tr(\bb),& if $t=0$ and $\eta\not=0$,\\
          \delta_{w,e}  ( \frac{-\eta}{t})^{n(\ell-1)-d}\tr(\bb),& if $t$ is invertible.
    \end{cases*}
  \end{equation}
 Here, $\tr$ is the trace map  on $\Bo$ induced from $A$. That is,
 for $\bb=b_1\otimes \dots\otimes b_n$, set $\tr(\bb)=\prod_{i=1}^n\tr(b_i)$.

The trace map $\Tr$ induces an $R$-bilinear form $\<\ ,\ \>\map{\Hhd(\Bo)\times\Hhd(\Bo)} R, \<x,y\>=\Tr(xy)$.
  \end{Definition}

    \begin{Remark}
   For $\eta=0$ and $t\in R^\times$ the third formula in
   \autoref{formform} also applies and agrees with the first, since
   $0^{N-d}=\delta_{d,N}$ and $d=N$ only for $\bd=\ELL$. We list the
   three cases separately because they need different arguments in the proof. The cases in \autoref{formform} are complete if we pass from $R$ to its fraction field.
 \end{Remark}

  \begin{Remark}\label{rktopdegree}
 Note that in case $t=0$ with $\eta\not=0$, respectively $\eta=0$, the trace can evaluate nontrivially only at the zeroth degree respectively at the top degree $N$. In contrast, there could be contributions from every degree if $t\not=0\not=\eta$. Our trace generalizes the standard traces for cyclotomic quotients of (degenerate) affine Hecke algebras \cite{BK:HigherSchurWeyl, MalleMathas}
 which involve only the top degree, as well as the trace from \cite{HeHuLin:traceForm} which involves all degrees.
 \end{Remark}

\subsection{The shifting isomorphism}
 We establish the following \emph{shifting isomorphism}, which allows us
 to shift the parameter $t$ whenever $\eta$ is invertible. In this
 case, set $\taui[i]=\taui[i,i+1]=\taui[i+1,i]$ for $1\leq i<n$.

 \begin {Theorem}\label{Phidef}
 Assume that $\eta$ and $t$ are invertible in $\RR$. Set $\lambda=\frac{t}{\eta}$ and $\fct'(X_1)=\fct(X_1-\lambda)\in R[X_1]$. Then there is an isomorphism of RaFH-superalgebras
  \begin{equation}\label{defShup}
\Sh\colon\HHd=\HHd[n;\hbar,t,\eta](\Bo)\cong \HHd[n;\hbar,0,\eta](\Bo);
\quad X_j\mapsto X_j-\lambda, H_i\mapsto H_i, \bb\mapsto\bb,
\end{equation}
 which induces an isomorphism, called the shift isomorphism, of CRaFH-superalgebras,
 \begin{equation}\label{defSh}
\Sh\colon\Hhd=\HHd[n;\hbar,t,\eta](\Bo)/(\fct)\cong {\Hhd}':=\HHd[n;\hbar,0,\eta](\Bo)/(\fct'),
\quad X_j\mapsto X_j-\lambda, H_i\mapsto H_i, \bb\mapsto\bb,
\end{equation}
for   $1\leq j\leq n$, $1\leq i\leq n-1$, $\bb\in\Bo$. Moreover,  $\Tr_{\Hhd}=(-\lambda)^{-N}\Tr_{{\Hhd}'}\circ\Sh$ with $N=n(\ell-1)$.
\end{Theorem}
\begin{Remark}In \autoref{defSh}, $\Tr_{\Hhd}$ falls into the third case of \autoref{formform} and $\Tr_{{\Hhd}'}$ into the second.
\end{Remark}
\begin{proof}
The compatibility with the defining relations \autoref{RaFHrels} is clear except for the fifth relation,  which holds by $\eta$-symmetricity, and
the last two relations which are equivalent, since $H_i$ is invertible. For the compatibility  it thus suffices to show that $\Sh(H_iX_i)=\Sh(X_{i+1}H_i-(t+\eta X_{i+1})\taui)$, that is  $H_iX_i-\lambda H_i=X_{i+1}H_i-\lambda H_i
-(t+\eta X_{i+1}-\eta\lambda)\taui$. Since $\eta\lambda=t$, this is equivalent to $H_iX_i=X_{i+1}H_i
-\eta X_{i+1}\taui$  which is a defining relation in $\HHd[n;\hbar,0,\eta](\Bo)$. Thus, the assignments from \autoref{defShup} give a well-defined superalgebra homomorphism
${\HHd[n;\hbar,t,\eta](\Bo)}\rightarrow\HHd[n;\hbar,0,\eta](\Bo)$. It
factors through $\Hhd$, because $\fct(X_1)$  is sent
to~$\fct(X_1-\lambda)=\fct'(X_1)=0$ in ${\Hhd}'$. Hence,
$\Sh\colon\Hhd\to{\Hhd}'$ is well-defined.  An analogous computation
shows that the map $\Sh[-\lambda]$ given by the same formulas, but with
$X_j\mapsto X_j+\lambda$, defines an algebra homomorphism
$\Sh[-\lambda]\colon{\Hhd}'\to\Hhd$, which is the inverse map. Hence,
$\Sh$ is an isomorphism.

It remains to show that
\begin{equation}\label{Treppchen}
(-\lambda)^{N}\Tr_{\Hhd}\bigl(\Sh^{-1}(X^{\bd}H_w\bb)\bigr)=\delta_{\bd,0}\delta_{w,e}\tr(\bb)
\end{equation}
for every basis vector $X^{\bd}H_w\bb$ of ${\Hhd}'$. Applying the binomial expansion we get
\begin{equation}\label{extro}
\Sh^{-1}(X^{\bd}H_w\bb)=\prod_{i=1}^{n}(X_i+\lambda)^{d_i}\,H_w\bb=\prod_{i=1}^{n}\left(\sum_{0\leq k_i\leq d_i}\binom{d_i}{k_i}\lambda^{d_i-k_i}X_i^{k_i}\right)H_w\bb.
\end{equation}
The last expression in  \autoref{extro} is an expansion in basis elements $X^{\mathsf k}H_w\bb$ since $d_i\leq\ell-1$. Thus,  \autoref{formform} gives $\Tr_{\Hhd}(X^{\mathsf k}H_w\bb)=\delta_{w,e}u^{N-|\mathsf k|}\tr(\bb)$ with $u=\frac{-\eta}{t}=-\lambda^{-1}$. Hence,
\begin{align*}
(-\lambda)^{N}\Tr_{\Hhd}\bigl(\Sh^{-1}(X^{\bd}H_w\bb)\bigr)&=\delta_{w,e}\tr(\bb)\prod_{i=1}^{n} \sum_{k_i=0}^{d_i}\binom{d_i}{k_i}\lambda^{d_i-k_i}u^{-k_i}=\delta_{w,e}\tr(\bb)\prod_{i=1}^{n}(\lambda+u^{-1})^{d_i}=\delta_{\bd,0}\delta_{w,e}\tr(\bb),
\end{align*}
where we used that $(-\lambda)^{N}u^{N}=1$, $\lambda+u^{-1}=0$ and $0^{0}=1$. This proves the claim \autoref{Treppchen}.
\end{proof}

\subsection{CRaFH superalgebras for $\ell=1$ are Frobenius}
We next show that CRaFH superalgebras of level $\ell=1$ are always Frobenius with trace form $\Tr$ and Nakayama automorphism induced from $A$. Afterwards we consider arbitrary level and show that the form $\Tr$ is symmetric up to a Nakayama automorphism induced from $A$ and depending on $\ell_0$.
\begin{Theorem}\label{formlevel1}
Assume  $\Hhd$ from \autoref{ass} is of level $\ell=1$. Then $\Tr$ turns $\Hhd$ into a Frobenius $R$-superalgebra with Nakayama automorphism given by $\psi_{\Hhd}(H_w\bb)=H_w\psi(\bb)$ for  $w\in\Sym, \bb\in\Bo$.

In particular, $\Hhd(\Bo)$ and $\<\ ,\ \>$ are symmetric if and only if  $A$ is symmetric.
\end{Theorem}
\begin{Remark}
Under the assumptions of \autoref{formlevel1}, the form $\Tr$ is unimodular if $\tr$ is unimodular.
\end{Remark}
\begin{proof}
Let $H_x\bb, H_{y^{-1}}\bc\in\Hhd$ be basis elements in the sense of  \autoref{basisthmcycl}. We claim that
\begin{equation}\label{pairing}
\Tr(H_x\bb H_{y^{-1}}\bc)=\delta_{x,y}{\hbar^{2\ell(x)}}\tr(\leftindex^{y}\bb\bc)\quad\text{for any $\bb,\bc\in\Bo$ and $x,y\in\Sym$}.
\end{equation}
  For this it suffices to show that $\Tr(H_xH_{y^{-1}}\bc)=\delta_{x,y}\hbar^{2\ell(x)}\tr(\bc)$ for all $\bc\in\Bo$,
 which we do by induction on the length of $y$. For $\ell(y)=0$ this is clear, since $\Tr(H_x\bc)=\delta_{x,e}\tr(\bc)$.  Otherwise pick $1\leq i<n$ such that $z=ys_i$ has smaller length than $y$.
 If $\ell(xs_i)>\ell(x)$, then $x\not=y$ and $xs_i\not=z$ and,
 by induction.
 \[\Tr(H_xH_{y^{-1}}\bc)=\Tr(H_{xs_i}H_{z^{-1}}\bc)={\delta_{xs_i,z}\hbar^{2\ell(xs_i)}\tr(\bc)=0=\delta_{x,y}\hbar^{2\ell(x)}}\tr(\bc).\]
  If $\ell(xs_i)<\ell(x)$, then $x\not=z$, and induction together with \autoref{L:HiHw} give that $\Tr(H_xH_{y^{-1}}\bc)$ equals
  \[\Tr(H_xH_{s_i}H_{z^{-1}}\bc)=
  \Tr(\hbar^2H_{xs_i}H_{z^{-1}}\bc)+\eta\Tr(H_x\tau^{(i,i+1)}H_{z^{-1}}\bc)=\hbar^2{\delta_{xs_i,z}}\hbar^{2\ell(z)}\tr(\bc)=\delta_{x,y}\hbar^{2\ell(x)}\tr(\bc),\]
  where we have used that
  $\Tr(H_x\tau^{(i,i+1)}H_{z^{-1}})=\Tr(H_xH_{z^{-1}}\leftindex^{z}\tau^{(i,i+1)})=0$
  since $x\neq z$.

  Therefore, on the one hand $\Tr(H_x\bb H_{y^{-1}}\bc)={\hbar^{2\ell(x)}}\delta_{x,y}\tr(\leftindex^{y}\bb\bc)$ and, on the other,
  \[(-1)^{|H_x\bb||H_{y^{-1}}\bc|}\Tr(H_{y^{-1}}\bc
  \psi_{\Hhd}(H_x\bb))=(-1)^{|\bb||\bc|}
   {\hbar^{2\ell(y^{-1})}}\delta_{y^{-1},x^{-1}}\tr(\leftindex^{x^{-1}}\bc\psi(\bb))=(-1)^{|\bb||\bc|}{\hbar^{2\ell(x)}}\delta_{x,y}\tr(\leftindex^{x^{-1}}\!\!\bc\psi(\bb)).\]
  These two results agree, since
  $(-1)^{|\bb||\bc|}\tr(\leftindex^{x^{-1}}\!\bc\psi(\bb))=(-1)^{|\bb||\bc|}(-1)^{|\bb||\bc|}\tr(\bb\leftindex^{x^{-1}}\!\!\bc)=\tr(\leftindex^{x}\bb\bc)$.
  Thus, $ \<\ ,\ \>$ is symmetric if and only if $A$ is symmetric.

 It remains to show that $\Tr$ is non-degenerate.
 By \autoref{pairing}, the Gram matrix with respect to the bases $\set{H_x\bb}$ and
 $\set{H_{y^{-1}}\bc}$ is block diagonal with blocks
 $\hbar^{2\ell(x)}\bigl(\tr((\leftindex^{x}\bb)\bc)\bigr)_{\bb,\bc}$ indexed by $x\in\Sym$. Up
 to signs and a permutation of the rows, each block is
 $\hbar^{2\ell(x)}$ times the Gram matrix $G_\Bo$ of $\tr$ on $\Bo$.
 Hence $\det G=\pm\hbar^{\binom n2\op{rk}\Hhd}(\det G_\Bo)^{n!}$.
 Since $\hbar$ is invertible and $\det G_\Bo\in R^\times$, as $A$ is Frobenius, $\Tr$ is unimodular, and in particular non-degenerate. Hence $\psi_{\Hhd}$ is its Nakayama automorphism, and $\<\ ,\ \>$ is symmetric if and only if $A$ is symmetric.
  \end{proof}

\subsection{Nakayama automorphisms in higher levels}\label{SS:Xfiltration} We
now treat CRaFH superalgebras $\Hhd(\Bo)$ of arbitrary level.
We first consider the generator $X_1$ in the cases when $\eta=0$ and $t=0$, respectively.

To get started, we define a filtration on $\Hhd$, which will be used
several times below.
\begin{Definition}\label{Xfiltr}
By \autoref{D:Hecke}, putting the $X_i$ in degree
$1$, and all other generators in degree $0$, turns $\HHd[n;\hbar, t,
\eta](\Bo)$ into a filtered $R$-algebra for any choice of parameters.
In the bases of \autoref{C:BasisHd}, the $m$th piece of this
\emph{$\X$-filtration} is spanned by the basis elements of $\X$-degree
at most~$m$.
\end{Definition}
Now suppose that $\eta=0=t$. Then the $\X$-filtration defines a grading
on $\HHd$. For $m\ge0$, let $F_m'\subset \Hhd$ be the degree $m$ part of
the $\X$-filtration of $\Hhd$, and let $F_m$ be the $R$-span of the
monomials $X^\bd H_w\bb$ of \autoref{basisthmcycl} with $d\leq m$.

\begin{Lemma}\label{filtannoying}Consider $\Hhd$ with $\eta=0$.
 Then  we have for all $m\geq 0$ that  $F_m=F'_m$. In  particular $F_m=F_m'=\Hhd$ for $m\geq N$ and $F_lF_m\subseteq F_{l+m}$ for all $0\leq l,m\leq N$.
\end{Lemma}

\begin{proof}
Clearly, $F_m\subseteq F'_m$ for all $m\geq0$. We first show by induction on $i$ that $X_i^{\ell}\in F'_{\ell-1}$ for
$1\leq i\leq n$. For $i=1$ this holds because
$X_1^{\ell}=-\sum_{k<\ell}f_kX_1^{k}$ in $\Hhd$. For the induction step,
using \autoref{L:HXH}, write $\hbar^{2}X_{i+1}=u+v$ with $u=H_iX_iH_i$ and
$v=tH_i\tau^{(i+1,i)}\in F'_0$. Since $\eta=0$, $H_i^2=\hbar^2$ is central and thus
$u^{\ell}=\hbar^{2(\ell-1)}H_iX_i^{\ell}H_i$. On the other hand, any other word in $u,v$ of length $\ell$ has
at most $\ell-1$ letters equal to~$u$, and hence lies in $F'_{\ell-1}$. By the inductive hypothesis,
$\hbar^{2\ell}X_{i+1}^{\ell}\in\hbar^{2(\ell-1)}H_iX_i^{\ell}H_i+F'_{\ell-1}\subseteq F'_{\ell-1}$.
Thus,  $X_{i+1}^{\ell}\in F'_{\ell-1}$ since $\hbar$ is invertible. We
 show the reverse inclusion $F'_m\subseteq F_m$ by induction on $m$,
 starting with the case $m=0$, where there is nothing to prove. It
suffices to see that any monomial $X^\bd H_w\bb$ with $\bd\in\N^{n}$,
$d\leq m$, lies in $F_m$, when viewed as an element of~$\Hhd$. If $d_i\leq\ell-1$ for all $i$, this is obvious. Otherwise $d_i\geq\ell$
for some $i$ and $X^\bd H_w\bb=X^{\bd-\ell e_i}\,X_i^{\ell}\,H_w\bb\in F'_{d-1}\subseteq
F'_{m-1}\subseteq F_{m-1}$, using $F'_lF'_m\subseteq F'_{l+m}$ and the inductive hypothesis.
\end{proof}

The only property of $\fct$ that is needed in the proof of the next
lemma is that it is a monic polynomial of degree~$\ell$.

\begin{Lemma}\label{L:XoneEtaZero}
Suppose that $\eta=0$. Then $\Tr(X_1y)=\Tr(yX_1)$, for all $y\in\Hhd$.
\end{Lemma}
\begin{proof}
Let $F_m$ be the span of the basis elements $X^\bd H_w\bb$ with
$|\bd|\leq m$. By \autoref{filtannoying} and its proof, $F_lF_m\subseteq
F_{l+m}$ and $X_k^\ell\in F_{\ell-1}$ for all $k$, and $\Tr$ vanishes on
$F_{N-1}$. Let $y=X^\bd H_w\bb$ be a basis element. Moving $X_1$ through
$H_w$ with $H_iX_k=X_{s_i(k)}H_i\pm t\tau$ shows that
$H_wX_1=X_{w(1)}H_w+C$, where $C$ is a linear combination of elements
$H_v\bb'$ with $\ell(v)<\ell(w), \bb'\in\Bo$. Hence,
\[
 yX_1=X^\bd H_wX_1\psi^{(1)}(\bb)=X^\bd X_{w(1)}H_w\psi^{(1)}(\bb)+X^\bd C\psi^{(1)}(\bb),\qquad \text{with } X^\bd C\psi^{(1)}(\bb)\in F_{|\bd|}.
\]
If $|\bd|\leq N-2$, then $X_1y$ and $yX_1$ lie in $F_{N-1}$, so for both the
trace vanishes. Let $|\bd|=N-1$. If $d_1\leq\ell-2$, then $X_1y$ is a
basis element and
$\Tr(X_1y)=\delta_{w,e}\delta_{\bd+e_1,\ELL}\tr(\bb)$. If
$d_1=\ell-1$, then $X_1y\in F_{N-1}$ since $X_1^\ell\in F_{\ell-1}$, so
$\Tr(X_1y)=0$ and the same formula holds as
$\bd+e_1\neq\ELL$. In the same way,
\[
    \Tr(yX_1)=\Tr(X^\bd X_{w(1)}H_w\psi^{(1)}(\bb))=\delta_{w,e}\delta_{\bd+e_1,\ELL}\tr(\psi^{(1)}(\bb))
    =\Tr(X_1y),
\]
where the last equality follows because $\tr\circ\psi^{(1)}=\tr$.
Finally, let $\bd=\ELL$. Then $X^{\ELL}=(X_1\cdots X_n)^{\ell-1}$
commutes with $H_i$ by the second identity in \autoref{L:XSymmetric}, for $1\le i<n$.

Put $P=X_1X^{\ELL}=-\sum_{k<\ell}f_kX_1^kX_2^{\ell-1}\cdots
X_n^{\ell-1}\in-f_{\ell-1}X^{\ELL}+F_{N-1}$. Then
$\Tr(X_1y)=\Tr(PH_w\bb)=-\delta_{w,e}f_{\ell-1}\tr(\bb)$ since
$H_wF_{N-1}\subseteq F_{N-1}$. On the other hand,
\[
 \Tr(yX_1)=\Tr\bigl(H_wX^{\ELL}X_1\psi^{(1)}(\bb)\bigr)=\Tr\bigl(H_wP\psi^{(1)}(\bb)\bigr)=-f_{\ell-1}\Tr\bigl(X^{\ELL}H_w\psi^{(1)}(\bb)\bigr)=-\delta_{w,e}f_{\ell-1}\tr(\bb),
\]
since $\tr(\psi^{(1)}\bb)=\tr(\bb)$. This completes the proof.
\end{proof}

\begin{Lemma}\label{L:XoneTZero}
Suppose that $t=0$ and $\eta\ne0$. Then $\Tr(X_1y)=\Tr(yX_1)$, for all $y\in\Hhd$.
\end{Lemma}
\begin{proof}
Here $\psi=\op{id}$, as $\eta\neq0$. Thus, the relation $\bb
X_k=X_k\psi^{(k)}(\bb)$ becomes $\bb X_k=X_k\bb$, so $\Bo$ commutes with
$X_1,\dots,X_n$.
 We argue by induction on $n$.

 For $n=1$, $\Hhd[1]=R[X_1]/(\fct)\otimes A$ and $X_1$ is central. Let $n\geq2$ and let $\Kr{n-1}$ be the subalgebra of $\Hhd$ generated by $\Bo$, $X_1,\dots,X_{n-1}$ and $H_1,\dots,H_{n-2}$. It has basis the $X^\bd H_w\bb$ with $d_n=0$ and $w\in\Sym_{n-1}$, and it is the image of the natural injective homomorphism $\Hhd[n-1](A^{\otimes(n-1)})\otimes A\to\Hhd$, and $\Tr$ restricts to $\Tr_{n-1}\otimes\tr$. Hence, by induction, $\Tr(X_1k)=\Tr(kX_1)$ for $k\in\Kr{n-1}$. Since $X_n$ commutes with $\Kr{n-1}$,
\[
 M=\sum_{c=1}^{\ell-1}X_n^c\Kr{n-1}+\Kr{n-1}H_{n-1}\Kr{n-1}
\]
is a $\Kr{n-1}$-bisubmodule of $\Hhd$. We claim that $\Hhd=\Kr{n-1}+M$ and $\Tr(M)=0$. Then, for $y=k_0+m_0$ with $k_0\in\Kr{n-1}$ and $m_0\in M$, we have $X_1m_0,m_0X_1\in M$, as $X_1\in\Kr{n-1}$, and hence $\Tr(X_1y)=\Tr(X_1k_0)=\Tr(k_0X_1)=\Tr(yX_1)$, as required.

We first show that $\Hhd=\Kr{n-1}+M$. By \autoref{eq:RewriteXr},
\begin{equation}\label{part1}
 H_{n-1}X_{n-1}^a=X_n^aH_{n-1}-\eta\taui[n-1]\sum_{p=1}^{a}X_{n-1}^{a-p}X_n^p\quad\text{and}\quad
 X_n^aH_{n-1}=H_{n-1}X_{n-1}^a+\eta\taui[n-1]\sum_{p=1}^{a}X_{n-1}^{a-p}X_n^p,
\end{equation}
for $a\geq0$, where all powers of $X_n$ in the sums lie between $1$ and $a$. Every basis element can be written as $kX_n^cH_{n-1}\cdots H_j\bb$ with $k\in\Kr{n-1}$, $0\leq c\leq\ell-1$ and $1\leq j\leq n$, using $w=us_{n-1}\cdots s_j$ with $u\in\Sym_{n-1}$ (the product $H_{n-1}\cdots H_j$ is $1$ if $j=n$). If $j=n$ this lies in $\Kr{n-1}$ or in $X_n^c\Kr{n-1}\subseteq M$. If $j<n$, the second formula in \autoref{part1} with $a=c<\ell$ shows that it lies in $\Kr{n-1}H_{n-1}\Kr{n-1}+\sum_{p\geq1}X_n^p\Kr{n-1}\subseteq M$, since $\taui[n-1]\in\Bo\subseteq\Kr{n-1}$. Hence $\Hhd=\Kr{n-1}+M$.

Let $1\le c\le\ell-1$. For a basis element $X^\bd H_w\bb$ of $\Kr{n-1}$, the element $X_n^cX^\bd H_w\bb$ is a basis element of $\Hhd$ with $d_n=c\neq0$, so $\Tr(X_n^c\Kr{n-1})=0$. For $k,k'\in\Kr{n-1}$ we may assume $k'=X''X_{n-1}^aH_{w'}\bb'$ with $X''$ a monomial in $X_1,\dots,X_{n-2}$, $a\leq\ell-1$ and $w'\in\Sym_{n-1}$. By the first formula in \autoref{part1},
\[
 kH_{n-1}k'=X_n^akX''H_{n-1}H_{w'}\bb'-\eta\sum_{p=1}^{a}X_n^pkX''\taui[n-1]X_{n-1}^{a-p}H_{w'}\bb'.
\]

The sum lies in $\sum_{p=1}^{\ell-1}X_n^p\Kr{n-1}$, so its trace is zero. For the first term, expand $kX''$ in the basis of $\Kr{n-1}$. For $u,w'\in\Sym_{n-1}$, \autoref{L:HiHw} shows that $H_uH_{n-1}H_{w'}$ is a combination of elements $\bb''H_v$ with $\bb''\in\Bo$ and $v$ in the double coset $\Sym_{n-1}s_{n-1}\Sym_{n-1}$, which does not contain $e$. Moving all elements of $\Bo$ to the right, which is possible as they commute with the $X_k$, the first term becomes a combination of basis elements $X^{\bd}H_v\bb$ with $v\neq e$, since all exponents stay below $\ell$ as $a\le\ell-1$. Hence its trace is zero, and $\Tr(M)=0$.
\end{proof}

 \begin{Theorem}\label{main:bilinearform}
Suppose that $x,y\in\Hhd$ are homogeneous. Then
\begin{equation*}
\Tr(xy)=(-1)^{|x||y|}\Tr(y\psi^{\ell_0}(x)),
\qquad\text{where $\psi(X^\bd H_w\bb)=X^\bd H_w\psi(\bb)$.}
\end{equation*}
In particular, the bilinear form $\<\ ,\ \>\map{\Hhd\times\Hhd} R$
induced by $\Tr$ is symmetric if $\theta$ divides $\ell_0$. Conversely,
if $\<\ ,\ \>$ is non-degenerate and symmetric, then $\theta$ divides
$\ell_0$.
\end{Theorem}

\begin{Remark}
In particular, the form of \autoref{main:bilinearform} is symmetric if
$A$ is symmetric. \autoref{main:bilinearform} also explain, via
$\psi^{\ell_0}$ and the parity of $\ell_0$, the known even-odd
dependence on the trace forms of cyclotomic Sergeev algebras. More
explicitly, it recovers \cite[Theorem~1.2]{LiShi:SergeevCocenter}, where
the case when $\ell_0$ is odd is called \emph{supersymmetric} and the
case where $\ell_0$ is even is called \emph{symmetric}. For the
non-degenerate cyclotomic Hecke-Clifford algebras that are not RaFH
superalgebras, see \cite{LiShi:QuiverHeckeClifford}.
\end{Remark}

      \begin{proof}[Proof of \autoref{main:bilinearform}]
It is enough to establish the supertrace condition
$\Tr(xy)=(-1)^{|x||y|}\Tr(y\psi^{\ell_0}(x))$, when $y=X^\bd H_w\bb$ is
a homogeneous basis elements from \autoref{basisthmcycl}, and when $x$
is one of the generators $H_i$, $X_1$, and $a^{(j)}$ with $a\in A$
homogeneous, $1\leq i\leq n-1$ and $1\leq j\leq n$. As $H_i$ and $X_1$
are even and fixed by $\psi$, all we need to check is that
$\Tr(xy)=\Tr(yx)$.

If $\eta\neq0$ and $t$ is invertible, we can extend scalars to the
fraction field of $R$ without changing either the basis or the trace
function~$\Tr$. By \autoref{Phidef}, $\Tr=(-\lambda)^{-N}\Tr'\circ\Sh$,
where $\Sh\colon\Hhd\to{\Hhd}'=\HHd[n;\hbar,0,\eta](\Bo)/(\fct')$ is an
algebra isomorphism. Note that in ${\Hhd}'$, $t=0$ and $\psi=\op{id}$,
since $\eta\neq0$. So the supertrace property of $\Tr$ follows from that
of $\Tr'$, the trace form on~${\Hhd}'$.

We have reduced the proof to comparing $\Tr(xy)$ and $\Tr(yx)$, when
$\eta\ne0$ and $t=0$ (and $\Hhd=\HHd'$), and when $\eta=0$ and
$t$ is arbitrary. Here, $y=X^\bd H_w\bb$ and
$x\in\set{H_i,X_1,a^{(j)}}$. If $x=X_1$ then we are done by
\autoref{L:XoneEtaZero} when $\eta=0$, and by \autoref{L:XoneTZero} when
$\eta\ne0$ and $t=0$. Hence, it remains to consider the case when
$x\in\set{H_i,a^{(j)}}$ and $y=X^\bd H_w\bb$. We can assume that at most
one of $\eta$ and $t$ is zero.

By the second identity in \autoref{L:XSymmetric} (which holds for arbitrary $\psi$ by the proof of \autoref{L:XoneEtaZero}), $X_iX_{i+1}$ commutes with $H_i$, for
all~$i$. Let $E\colon\Hhd\to\Hfind$ be the $R$-linear map with
\[
   E(X^\bd H_w\bb)=\begin{cases*}
     \delta_{\bd,\ELL}H_w\bb& if $\eta=0$, \\
     \delta_{\bd,0}H_w\bb& if $t=0$.
   \end{cases*}
\]
Then $\Tr=\Tr_{\Hfind}\circ E$, and $E$ is right $\Hfind$-linear. By
\autoref{L:crucial},
$H_iX^\bd=\leftindex^{s_i}(X^\bd)H_i+\Delta_i(X^\bd)$, where
all exponents in $\Delta_i(X^\bd)$ are strictly less than $\ell$, and
have degree $|\bd|-\delta_{\eta,0}$, and where
$\Delta_i(X^{\ELL})=0=\Delta_i(1)$. Since $\bb
X^\bd=X^\bd\psi^{\bd}(\bb)$ this gives $E(H_iy)=H_iE(y)$ and $E(\bb
y)=\psi^{\ell-1}(\bb)E(y)$, where $\psi^{\ell-1}=\op{id}$ if
$t=0\neq\eta$. Hence, by \autoref{formlevel1}, if $x\in\{H_i,a^{(j)}\}$,
then
\begin{equation}\label{Tryx}
\Tr(xy)=\Tr_{\Hfind}\bigl(\psi^{\ell-1}(x)E(y)\bigr)=(-1)^{|x||y|}\Tr_{\Hfind}\bigl(E(y)\psi^{\ell}(x)\bigr)=(-1)^{|x||y|}\Tr\bigl(y\psi^{\ell_0}(x)\bigr),
\end{equation}
where the last equality follows because $\psi^\ell=\psi^{\ell_0}$. In
particular, $\Tr(xy)=(-1)^{|x||y|}\Tr(y\psi^{\ell_0}(x))$, for
$x\in\set{H_i,a^{(j)}}$, so this completes the proof.
\end{proof}

\begin{Remark}\label{R:SavageRS}
For $\hbar=t=1$ and $\eta=0$ we obtain cyclotomic quotients of the affine wreath product algebras of \cite{Savage:AffineWreath} with $\Tr$  being the trace map of \cite[(6.8)]{Savage:AffineWreath}. \autoref{main:bilinearform} and \autoref{T:nondeg} recover \cite[Theorem~6.15]{Savage:AffineWreath} for cyclotomic polynomials with scalar coefficients, including the Nakayama automorphism $\psi^{\ell_0}$ and symmetry if $\theta$ divides $\ell_0$. For $\eta\neq0$ and $A$ symmetric we get the cyclotomic quotients of the quantum affine wreath algebras and their symmetricity from \cite[Theorem~4.16]{RossoSavage:QuantumAffineWreath}. \autoref{T:summary} determines when $\Tr$ is non-degenerate.
\end{Remark}

\subsection{Some properties of Gram matrices}\label{sectop}
In this subsection we introduce \emph{top functionals}. They allow us to compute Gram determinants of filtered algebras by passing to the associated graded algebra. For $\eta=0$ the CRaFH superalgebras fit into this setting, see \autoref{extop}, which will be used in \autoref{T:nondeg}. We also compute two Hankel type determinants in \autoref{L:hankel}. They determine the singular locus of the trace form in the base case $n=1$, treated in the next subsection.

Let $\mathbf{H}$ be an $R$-algebra with a fixed finite $R$-basis $S$ and a compatible filtration
\begin{equation}\label{filttop}
  0=\mathbf{H}_{-1}\subsetneq \mathbf{H}_0\subsetneq\dots\subsetneq \mathbf{H}_N=\mathbf{H},
\end{equation}
 i.e. a sequence of $R$-modules satisfying $\mathbf{H}_a\mathbf{H}_b\subseteq \mathbf{H}_{a+b}$ with $\mathbf{H}_m=\mathbf{H}$ for $m\geq N$, such that $\mathbf{H}_m$ has an $R$-basis  $S_m\subseteq S$ with  $|S_m\smallsetminus S_{m-1}|=|S_{N-m}\smallsetminus S_{N-m-1}|$.
Let $\bar S=\{\bar s\mid s\in S\}$ be the induced basis of the associated graded $ \mathbf{H}_{\gr}=\bigoplus_m \mathbf{H}_m/\mathbf{H}_{m-1}$. We call an $R$-linear map $T \colon \mathbf{H}\rightarrow R$ a \emph{top functional} and its bilinear form $\beta_T\colon \mathbf{H}\times \mathbf{H}\to R$, $\beta_T(x,y)=T(xy)$, a \emph{top pairing} for  $\mathbf{H}$, if $T(x)=0$ for all $x\in S_m$ with $m<N$.

 \begin{Example}\label{extop}
 Let $A$ be a Frobenius superalgebra and consider the  CRaFH
 superalgebra $\mathbf{H}=\Hhd(\Bo)=\HHd[n;\hbar, t, \eta](\Bo)/(\fct)$
 for $\Bo=A^{\otimes n}$.  The $\X$-filtration on $\HHd[n;\hbar, t,
 \eta](\Bo)$ from \autoref{Xfiltr} induces a filtration $\Hhd(\Bo)_m$, $m\geq0$ on $\Hhd(\Bo)$.
If we pick a basis for $\HHd[n;\hbar, t, \eta](\Bo)$ of the form \autoref{C:BasisHd}  and let $S$ be the corresponding basis for $\Hhd(\Bo)$ from \autoref{C:HdBasis}, then $\mathbf{H}=\Hhd(\Bo)=\HHd[n;\hbar, t, 0](\Bo)/(\fct)$ for $\eta=0$, with $t$ arbitrary, falls into the above setup with $N=n(\ell-1)$. Namely, by  \autoref{filtannoying}, $\Hhd(\Bo)_m$ has an $R$-basis $S_m$ given by all monomials in $S$ of total $\X$-degree $\leq m$ which induces a basis $\bar S$ on the degree $m$-part of $\gr\Hhd(\Bo)$. In particular, $\Hhd(\Bo)_m=\Hhd(\Bo)$ for $m\geq N$.  The equality $|S_m\smallsetminus S_{m-1}|=|S_{N-m}\smallsetminus S_{N-m-1}|$ holds because $\bd\mapsto(\ell-1,\dots,\ell-1)-\bd$ is a bijection between the basis elements of degree $m$ and those of degree $N-m$. Since $\eta=0$, $\Tr$ is by definition \autoref{formform} a top functional.
 \end{Example}

\begin{Lemma}\label{L:antitri}
  Any  top functional $T$ for $\mathbf{H}$ defines a top functional $\bar T\colon \gr\mathbf{H}\rightarrow R$,  where $\gr \mathbf{H}$ is equipped with the basis $\bar S$ and the associated filtration, by setting for $s\in S$
  \[
    \bar T(\bar s)=
      \begin{cases}
       T(s)&\text{if $s\in S_N\smallsetminus S_{N-1}$,}\\
     0&\text{otherwise.}
            \end{cases}
  \]
Up to possibly a sign,  the determinants of the Gram matrices,  $G$ for $\beta_T$  and $\bar G$ for  the induced form  are equal, i.e. in formulas,
$\det G=\pm\det\bar G$.
\end{Lemma}

\begin{proof}
Pick a total ordering on each $S_m$ and extend it to an ordering on $S$ in three ways.
The first poset, $S^+$, is given by $S_0<S_1<\cdots <S_N$, i.e. starts with $S_0$ as the smallest elements, then take $S_1$ and so on.
The second poset,  $S^-$, is given by $S_N<S_{N-1}<\cdots <S_0$, and the third poset,  $S^\diamond$,  takes $S_0<S_{N}<S_{N-1}<\cdots<S_{1}$.
Now consider the Gram matrix $G$ with respect to the ordered  basis $S^+\times S^+$ and  partition it  into blocks $G_{ab}$, $0\leq a,b\leq N$ according to the filtration degree.  Since $\beta_T$ is a top pairing, $G_{ab}=0$ for $a+b<N$, since  $\mathbf{H}_a\mathbf{H}_b\subseteq\mathbf{H}_{a+b}\subseteq\mathbf{H}_{N-1}$ and $T$ vanishes on $\mathbf{H}_{N-1}$. The Gram matrix $\bar G$ of the induced form vanishes off the antidiagonal, since $\bar T$ not only annihilates in degrees $a+b<N$, but also in degrees $a+b>N$ since $\gr\mathbf{H}_m=0$ for $m>N$. Its antidiagonal blocks are $G_{a,N-a}$.
The same blocks appear on the main diagonal of the block lower triangular matrix $G'=(G_{a,N-b})_{a,b}$ which we obtain if we consider instead the Gram matrix $G$ in the basis $S^+\times S^-$.  E.g.,~for $N=2$,
  \small
  \[G=
  \begin{pmatrix}
      G_{00}&G_{01}&G_{02}\\ G_{10}&G_{11}&G_{12}\\ G_{20}&G_{21}&G_{22}
    \end{pmatrix}\;
    =\begin{pmatrix}
      0&0&G_{02}\\ 0&G_{11}&G_{12}\\ G_{20}&G_{21}&G_{22}
    \end{pmatrix}\;
   G'=\begin{pmatrix}
      G_{02}&0&0\\ G_{12}&G_{11}&0\\ G_{22}&G_{21}&G_{20}
    \end{pmatrix}\;
     \bar G=\begin{pmatrix}
      0&0&G_{02}\\ 0&G_{11}&0\\ G_{20}&0&0
    \end{pmatrix}\;
  \]
   \normalsize
Thus,  $\det G= \pm\prod_{0\leq a\leq N}\det G_{a,N-a}=\pm \det \bar{G}$.
  \end{proof}

For $k\in\N_0$ let $\alpha_k\map{R[X]}R$, $\sum_i\lambda_iX^{i}\mapsto\lambda_k$, be the
\emph{$k$th coefficient map}. Fix a monic polynomial $f=X^{\ell}+\sum_{k<\ell}f_kX^{k}\in R[X]$
and let $\alpha_{k,f}\map{R[X]/(f)}R$ be the \emph{$k$th coefficient map modulo $f$}. It sends the
class of $p$ to the $k$th coefficient of its remainder modulo $f$. For an $R$-linear map
$\phi\map{R[X]/(f)}R$ put $M(f,\phi)=\bigl(\phi(X^{a+b})\bigr)_{0\leq a,b<\ell}$ and abbreviate
$M(f,k)=M(f,\alpha_{k,f})$.

\begin{Lemma}\label{L:hankel}
In the above setup we have $\det M\bigl(f,\alpha_{\ell-1,f}\bigr)=\pm1$ and $\det M\bigl(f,\alpha_{0,f}\bigr)=\pm f_0^{\,\ell-1}$.
\end{Lemma}

\begin{proof}For $m<\ell$,  $X^{m}$ is its own
  remainder, so we get  $\alpha_{k,f}(X^{m})=\delta_{m,k}$, while $\alpha_{0,f}(X^{\ell})=-f_0$. Fix the ordered bases  $S^+=\set{X^0, X^1,\ldots, X^{\ell-1}}$, $S^-=\set{X^{\ell-1}, X^{\ell-2},\ldots, X^0}$, $S^\diamond=\set{X^0,X^{\ell-1}, X^{\ell-2},\ldots,X^1}$ of $R[X]/(f)$.
    In the basis $S^+\times S^-$, the matrix for $M(f,\alpha_{\ell-1})$ is lower triangular  with $1$'s on the diagonal. In the basis $S^+\times S^\diamond$, the matrix for $M(f,\alpha_{0})$ is lower triangular with
  diagonal entries $1,-f_0,\dots,-f_0$. The statements follow.
\end{proof}

\subsection{Non-degeneracy and Gram determinants in the case $n=1$}\label{S:tracenequal1}
We now return to the setting from  \autoref{ass}  and  deal with the problem of non-degeneracy of the trace map $\Tr$ from  \autoref{formform}  for $n=1$, i.e. we consider CRaFH superalgebras $\Hhd=\Hhd[1]=\Hhd[1;\hbar, t, \eta](A)=\HHd[1](A)/(\fct)$  as in \autoref{nequal1}.

The RaFH superalgebras $\HHd[1;\hbar, t, \eta](A)$  are by \autoref{D:Hecke} canonically isomorphic to $R[X_1]\Smash A$ with the smash product relation $\bb X_1=X_1\psi(\bb)$. This isomorphism is  independent of $\hbar$, $t$, $\eta$, since the parameters do not  appear in the relations. We get  an induced canonical identification $\HHd[1](A)/(\fct)=R[X_1]/(\fct)\Smash A$ for the CRaFH quotient  $\Hhd[1]$ and the basis $\cB_{\op{cyc}}=\{X_1^d b\mid 0\leq d< \ell,  b\in\cB\}$, all parameter independent. The trace map $\Tr$ however does  depend on the parameters and this is what makes the case $n=1$ already interesting.

Write $\fct=X_1^{\ell}+\sum_{k<\ell}f_kX_1^{k}\in R[X_1]$ and  denote by $\kapp_1,\dots,\kapp_\ell$  the roots of $\fct$ taken with multiplicity.
Order the basis $\cB$ of $A$ as $b_1<b_2<\cdots<b_r$ and turn the bases $\set{X_1^{k}\bb}$ and $\set{\bb X_1^{k}}$ from \autoref{basisthmcycl}  into ordered bases $S$ respectively $S'$ by extending $\cB$ as follows:
\begin{equation}\label{SS}
\text{ for $S$:}\quad \cB<X_1\cB<X_1^{2}\cB<\cdots<X_1^{\ell-1}\cB,\quad\quad\text{for $S'$:}\quad \cB<\cB X_1<\cB X_1^{2}<\cdots<\cB X_1^{\ell-1}.
\end{equation}
The two bases are linked via the smash  relation $\bb X_1=X_1\psi(\bb)$.
Set $G_A=\bigl(\tr(\bb\bc)\bigr)_{\bb,\bc\in\cB}$. Since the trace form $\tr$ on $A$ is non-degenerate, we get $\det G_A\in\kk^{\times}$.

By definition \autoref{formform} we can  factor the trace $\Tr(X_1^{k}\bb)=\phi_{\Tr}(X_1^{k})\tr(\bb)$
using the following  $R$-linear map
   \begin{equation}\label{factor}
  \phi_{\Tr}\map{R[X_1]/(\fct)}R\quad {\phi_{\Tr}(X_1^{k})} =
     \begin{cases*}
        {\delta_{k, \ell-1}}& if $\eta=0$,\\
        {\delta_{k, 0}}& if $t=0$ and $\eta\not=0$,\\
             {(\frac{-\eta}{t})^{\ell-1-k}}& if $t$ is invertible.
    \end{cases*}
  \end{equation}
This allows us to compute the determinant of the  Gram matrix $G$ of the form $\<\ ,\ \>=\beta_{\Tr}$ on $\Hhd[1]$.

\begin{Lemma}\label{L:n1form}
Consider $\HHd[1](A)/(\fct)$ with trace map $\Tr$. Then the  Gram matrix of $\<\ ,\ \>$ in the basis $S\times S'$ is given by the Kronecker product  $G=M(\fct,\phi_{\Tr})\otimes G_A$ of matrices.

In particular $\det G=\bigl(\det M(\fct,\phi_{\Tr})\bigr)^{r}(\det G_A)^{\ell}$, where $r=|\cB|$ and $\ell$ is the level.
\end{Lemma}

\begin{proof}
Let  $X_1^{a}\ba\in S$ and $\bc X_1^{c}\in S'$. The corresponding entry of $G$ is given by
\begin{equation}\label{chickensoup}
  \<X_1^{a}\ba,\bc X_1^{c}\>=\Tr\bigl((X_1^{a}\ba)(\bc X_1^{c})\bigr)=\Tr\bigl(X_1^{a+c}\psi^{c}(\ba\bc)\bigr)
    =\phi_{\Tr}(X_1^{a+c})\tr(\ba\bc).
\end{equation}
Here we used that $\bb X_1=X_1\psi(\bb)$ for $\bb\in A$ and that
$\tr(\psi^c(\ba\bc))=\tr(1\cdot \psi^c(\ba\bc))=\tr(\ba\bc)$, and $X_1^{a+c}$ is read
modulo $\fct$. By the definition of the Kronecker product and the orderings \autoref{SS}, the
right hand side is the entry of $M(\fct,\phi_{\Tr})\otimes G_A$ in position
$\bigl((a,\ba),(c,\bc)\bigr)$. Thus $G=M(\fct,\phi_{\Tr})\otimes G_A$. The last claim follows
from $\det(P\otimes Q)=(\det P)^{\dim Q}(\det Q)^{\dim P}$ applied to $P=M(\fct,\phi_{\Tr})$ and
$Q=G_A$.
\end{proof}

\begin{Theorem}\label{P:n1}
Consider $\Hhd[1](A)=\HHd[1](A)/(\fct)$. Its bilinear form $\<\ ,\ \>$ has Gram determinant
  \begin{enumerate}
  \item\label{P:n1:eta0}
    $\det G=\pm(\det G_A)^{\ell}$, if $\eta=0$, with $t$ arbitrary,
  \item\label{P:n1:t0}
    $\det G=\pm\fct(0)^{(\ell-1){r}}(\det G_A)^{\ell}$, if $t=0$ and $\eta\not=0$,
  \item\label{P:n1:ht}
    $\det G=\pm t^{-\ell(\ell-1)r}\prod_{j=1}^{\ell}(t+\eta\kapp_j)^{(\ell-1)r}(\det G_A)^{\ell}$, if $t$ is invertible.
      \end{enumerate}
  Thus the following table  characterises when $\<\ ,\ \>$ is non-degenerate respectively  unimodular:
    \begin{center}
  \autoref{P:n1:eta0}: $\begin{cases}
   \text{always}, \\ \text{always},
   \end{cases}\quad$
   \autoref{P:n1:t0}:
    $ \fct(0)
    \begin{cases}
  \not=0,\\
  \in R^{\times},
  \end{cases}\quad$
  \autoref{P:n1:ht}:
   $t+\eta\kapp_j
   \begin{cases}
   \not=0,\\
  \in R^{\times},
\end{cases}$ for all $j$.
  \end{center}

\end{Theorem}

\begin{proof}
  By \autoref{L:n1form} we have
  $\det G=\bigl(\det M(\fct,\phi_{\Tr})\bigr)^{r}(\det G_A)^{\ell}$ in all three cases, so it is
  enough to identify $\phi_{\Tr}$ and to compute $\det M(\fct,\phi_{\Tr})$ with
  \autoref{L:hankel}.

  \autoref{P:n1:eta0} and \autoref{P:n1:t0} By \autoref{factor}, $\phi_{\Tr}(X_1^{k})=\delta_{k,m}$ with $m=\ell-1$ if $\eta=0$ and $m=0$ if $t=0\neq\eta$, that is, $\phi_{\Tr}=\alpha_{m,\fct}$. By \autoref{L:hankel}, $\det M(\fct,\phi_{\Tr})$ is $\pm1$ for $m=\ell-1$ and $\pm\fct(0)^{\ell-1}$ for $m=0$, which gives \autoref{P:n1:eta0} and \autoref{P:n1:t0}.

  \autoref{P:n1:ht} If $\eta=0$, then $\prod_j(t+\eta\kapp_j)^{(\ell-1)r}=t^{\ell(\ell-1)r}$ and the formula is \autoref{P:n1:eta0}. So let $\eta\neq0$. Both sides of the formula lie in the domain $R$, so we may verify it over the fraction field of $R$. This changes neither the bases nor $\Tr$, and there $\eta$ is invertible. Moreover $\psi=\op{id}$ by \autoref{rketasymm}. Let $\lambda=t\eta^{-1}$. By \autoref{Phidef} with $n=1$ and $N=\ell-1$, the shift $\Sh$ is an isomorphism onto ${\Hhd}'=\HHd[1](A)/(\fct')$, where $\fct'(X_1)=\fct(X_1-\lambda)$ and $t=0$, and $\Tr=(-\lambda)^{-(\ell-1)}\Tr'\circ\Sh$. Hence $\Tr(vu)=(-\lambda)^{-(\ell-1)}\Tr'(\Sh(v)\Sh(u))$ for $v\in S$ and $u\in S'$. Since $\Sh(X_1^{a}\bb)=(X_1-\lambda)^{a}\bb$ is $X_1^a\bb$ plus a combination of the $X_1^{a'}\bb$ with $a'<a$, and similarly for $\bb X_1^c$, the map $\Sh$ is given by unitriangular matrices on $S$ and on $S'$. Therefore $\det G=(-\lambda)^{-(\ell-1)\ell r}\det G'$, where $G'$ is the Gram matrix of $\Tr'$. The proof of \autoref{P:n1:t0} only uses \autoref{L:n1form} and \autoref{L:hankel}, which hold for every monic polynomial. So it applies to $\fct'$ and gives $\det G'=\pm\fct'(0)^{(\ell-1)r}(\det G_A)^{\ell}$. Finally, $\fct'(0)=\fct(-\lambda)=\prod_j(-\lambda-\kapp_j)=(-\eta)^{-\ell}\prod_j(t+\eta\kapp_j)$, and $(-\lambda)(-\eta)=t$. Together this gives \autoref{P:n1:ht}.

The unimodularity and non-degeneracy conditions follow now directly.
\end{proof}

\begin{Corollary}\label{L:Aprime}
  Let $\eta=0$ and $t$ arbitrary. Then the CRaFH superalgebra  $\Hhd[1](A)=\HHd[1](A)/(X_1^{\ell})$ of  level $\ell$ is Frobenius with trace $\Tr(X_1^{k}\bb)=\delta_{k,\ell-1}\tr(\bb)$ and Nakayama automorphism
  $\psi(X_1^{k}\bb)=X_1^{k}\psi^{\ell}(\bb)$.
\end{Corollary}

\begin{proof}
  If $\fct=X_1^{\ell}$ and $\eta=0$ then
  $\Tr(X_1^{k}\bb)=\delta_{k,\ell-1}\tr(\bb)$, see \autoref{factor}, which is the asserted trace. The form is  non-degenerate, by \autoref{P:n1}\autoref{P:n1:eta0}. The other statements follow from \autoref{main:bilinearform} noting that here $\ell_0=\ell$.
\end{proof}

\subsection{Non-degeneracy in the  symmetric group case using the RaFH-CRaFH trick}
\label{SS:Trick}
We again consider CRaFH superalgebras $\Hhd$ as in  \autoref{ass}  but
assume additionally $\eta=0$. By \autoref{formform}, the trace map is
then the top degree functional, $\Tr(X^\bd
H_w\bb)=\delta_{w,e}\,\delta_{\bd,\ELL}\,\tr(\bb)$. To be allowed to treat $\Tr$ as
a top functional and to put $\Hhd$ into the framework from
\autoref{extop} we use \autoref{filtannoying}.

We first describe the associated graded algebra of $\Hhd$ with respect
to the filtration from \autoref{filtannoying}. The description uses a
trick, which we call the \emph{RaFH-CRaFH trick}: we use a CRaFH
superalgebra as the Frobenius part of a RaFH superalgebra. This is
useful because the associated graded algebra is the RaFH superalgebra whose Frobenius algebra is the CRaFH superalgebra for $n=1$.

By \autoref{P:n1}\autoref{P:n1:eta0}, the CRaFH superalgebra $\Hhd[1](A)=\HHd[1](A)/(X_1^{\ell})$ with $\eta=0$ is a Frobenius superalgebra over $R$ with trace $\tr_{\Hhd[1](A)}(X_1^k\bb)=\delta_{k,\ell-1}\tr(\bb)$, see \autoref{L:Aprime}. By \autoref{T:FrobTensor} and \autoref{tensorFrob}, $\Hhd[1](A)^{\otimes n}$ is then a Frobenius-tensor superalgebra with trace $\tr_{\Hhd[1](A)}^{\otimes n}$. As $\eta=0$, it is trivially $\eta$-symmetric. Thus,  the RaFH superalgebra $\Hfind[n;\hbar,0,0](\Hhd[1](A)^{\otimes n})$ is defined. To distinguish it from the generators $X_j$ of $\Hhd$, we write $x$ for the generator $X_1$ of $\Hhd[1](A)$, viewed as an element of this Frobenius algebra.

\begin{Proposition}\label{L:gr}Consider the filtered superalgebra $\Hhd=\HHd[n;\hbar, t, 0](\Bo)/(\fct)$ from \autoref{filtannoying}. Then
  \[
   \Phi\colon \Hfind[n;\hbar,0,0](\Hhd[1](A)^{\otimes n})\bijection\gr\Hhd,\qquad
    H_i\mapsto H_i,\quad a^{(j)}\mapsto a^{(j)},\quad x^{(j)}\mapsto x_j:=X_j+F_0
  \]
  for $0\leq i\leq n-1$, $a\in A$,
  defines an isomorphism of $R$-superalgebras. It identifies the trace map of
  \autoref{formlevel1} with the functional induced by $\Tr$ on $\gr\Hhd(\Bo)$ via \autoref{L:antitri}.
\end{Proposition}

\begin{proof}Abbreviate $\Bo'=\Hhd[1](A)^{\otimes n}$. By \autoref{notationclash} it suffices to
check compatibility with the defining relations from \autoref{D:RFHS}. The first one is fine since
$\eta=0$ on both sides, the second and third are clear, and in the fourth we may take
$\bb=a^{(j)}$ or $\bb=x^{(j)}$. Together with the relations of $\Bo'$ this leaves the following.
\begin{enumerate}
\item\label{Tee1} $H_ix_j=x_{s_i(j)}H_i$: for $j\not=i,i+1$ this is $H_iX_j=X_jH_i$, while for
$j=i,i+1$ the relations $H_iX_i=X_{i+1}H_i-t\tau^{(i,i+1)}$ and
$H_iX_{i+1}=X_iH_i+t\tau^{(i+1,i)}$ of \autoref{RaFHrels}, valid as $\gamma_i=t$, have
correction terms in $F_0$, which vanish in $F_1/F_0$.
\item\label{Tee2} $a^{(j)}x^{(j)}=x^{(j)}\psi(a)^{(j)}$, the commutation of different tensor
factors, and $(x^{(j)})^{\ell}=0$: the first two are $\bb X_j=X_j\psi^{(j)}(\bb)$ and
$X_iX_j=X_jX_i$ from \autoref{RaFHrels}, and the last holds since $X_j^{\ell}\in F_{\ell-1}$ by
\autoref{filtannoying}.
\end{enumerate}
Thus $\Phi$ is well-defined. It maps the basis $\set{H_w\bb'}$ of \autoref{P:H0Basis}, with
$\bb'$ running through the basis $\set{x^\bd\bb}$ of $\Bo'$, to the classes of the
$X^{w(\bd)}H_w\bb$, which form a basis of $\gr\Hhd(\Bo)$ by \autoref{filtannoying}. So $\Phi$ is
an isomorphism. Finally, $\Tr$ vanishes on $F_{N-1}$ by \autoref{formform}, so $\bar{\Tr}$ is
determined by its values in the top degree, where it agrees with the trace map of
\autoref{formlevel1}.
\end{proof}

We collect our results  and deduce that in the symmetric group case, that means if $\eta=0$, the  trace map $\Tr$ is always non-degenerate, and compute the Nakayama automorphism.
\begin{Theorem}\label{T:nondeg}
  Consider a CRaFH superalgebra  $\Hhd$ as in \autoref{ass} but assume $\eta=0$. Then the Gram matrix of $\<\ ,\ \>$  is invertible. In particular  $\<\ ,\ \>$ is unimodular and therefore non-degenerate. In particular, $\Hhd(\Bo)$ is a Frobenius
  $R$-superalgebra. Its Nakayama automorphism is $\psi^{\ell_0}$ .
\end{Theorem}

\begin{proof}
  By \autoref{extop} and \autoref{L:antitri}, $\det G=\pm\det\bar G$. By \autoref{L:gr} $\bar G$ this is the Gram matrix of the trace $\Tr$ on $\Hfind[n;\hbar,0,0](\Bo')$ from \autoref{formlevel1}  with $\Bo'=\Hhd[1](A)^{\otimes n}$. The trace of $\Hhd[1](A)$ is unimodular by \autoref{P:n1}\autoref{P:n1:eta0}, hence its $n$-fold tensor power is unimodular as well. Moreover, its Gram matrix is up to signs a Kronecker power. So $\det\bar G\in R^{\times}$ by \autoref{formlevel1}. The Nakayama automorphism is given by \autoref{main:bilinearform}.
\end{proof}

\begin{Remark}\label{R:BKfiltration}
The strategy of this subsection goes back to Brundan and Kleshchev \cite[Appendix~A]{BK:HigherSchurWeyl}. For $A=\kk$ and $\theta=1$, the CRaFH superalgebras with $\hbar=t=1$ and $\eta=0$ are the degenerate cyclotomic Hecke algebras. Then \autoref{filtannoying} and \autoref{L:gr} specialise to their description of the associated graded algebra as $\kk[x_1,\dots,x_n]/(x_1^\ell,\dots,x_n^\ell)\rtimes\kk\Sym_n$. The trace of \autoref{formlevel1} on it is the form of \cite[Lemma~A.1]{BK:HigherSchurWeyl}, and \autoref{T:nondeg} is \cite[Theorem~A.2]{BK:HigherSchurWeyl}, which we obtain in addition with unimodularity. Also the case $x=X_1$, $\eta=0$ in the proof of \autoref{main:bilinearform} follows their argument, using that $X^{\ELL}$ commutes with all $H_i$. For the degenerate case see also  \cite[Theorem~4.6]{MazStr}. The same filtration as ours is used in \cite{Brundan:degenCentre} and \cite[\S4.3]{Savage:AffineWreath} to study centres.
\end{Remark}

\subsection{Non-degeneracy and Gram determinants in case $\eta\not=0$ and $t=0$}
In this subsection we still consider $\Hhd$ as in \autoref{ass} but assume now $t=0$ and $\eta\not=0$ in which case
$\Tr(X^{\bd}H_w\bb)=\delta_{\bd,0}\,\delta_{w,e}\,\tr(\bb)$. Recall that $\eta\not=0$ forces $\psi=\op{id}$ by \autoref{rketasymm} which ensures the relations $\bb X_i=X_i\bb$, $\taui[i,i+1]=\taui[i+1,i]$ and  the symmetry of $\<\ ,\ \>$. We will use this tacitly.

We compute the Gram determinant by an inductive argument in the spirit of \cite{MalleMathas}, where the Ariki--Koike case is treated. As in loc. cit.\ and in the proof of \autoref{main:bilinearform} we consider for $1\leq \tl\leq n$ the subalgebras $\Kr{\tl}$ of $\Hhd$ generated by $\Bo$,
$X_1,\dots,X_{\tl}$ and $H_1,\dots,H_{\tl-1}$ with $\Kr{0}=\Bo$.  By \autoref{C:HdBasis}, $\Kr{\tl}$ is free with basis the $X^{\bd}H_w\bb$ with $d_i=0$ for $i>\tl$ and $w\in\langle s_1,\dots,s_{\tl-1}\rangle$. We will use this to construct an ordered basis, in which the Gram matrix becomes triangular. For this we put a total ordering $\prec$ on $\N^{\ell,n}$ by letting $\bd\prec\ba$ if $d_i<a_i$ for $i$ maximal with $d_i\not=a_i$.

    \begin{Lemma}\label{L:Krdec}
   Let $1\leq \tl\leq n$. Multiplication gives a decomposition of $\Kr{\tl}$
      \begin{equation}\label{Krdec}
    \Kr{\tl}=\bigoplus_{c=0}^{\ell-1}\bigoplus_{j=1}^{\tl}\Kr{\tl-1}X_{\tl}^{c}H_{\tl-1}\cdots H_{j+1}H_j \quad\text{with}\quad \varepsilon_{\tl}\map{\Kr{\tl}}{\Kr{\tl-1}},
    \end{equation}
     into free left $\Kr{\tl-1}$-modules of rank one  with  the $\Kr{\tl-1}$-linear projection map $\varepsilon_{\tl}$ onto the $(c,j)=(0,\tl)$ summand. Moreover, $\Tr(z)=\Tr(\varepsilon_{\tl}(z))$ for all $z\in\Kr{\tl}$.
 \end{Lemma}
 \begin{proof}We deduce the decomposition from \autoref{C:HdBasis} and the relations \autoref{RaFHrels} for $t=0$. Every $w\in\langle s_1,\dots,s_{\tl-1}\rangle$ factors uniquely as $w=w's_{\tl-1}\cdots s_j$ with $w'\in\langle s_1,\dots,s_{\tl-2}\rangle$ and $1\leq j\leq\tl$ and  lengths adding, so $H_w=H_{w'}H_{\tl-1}\cdots H_j$. In particular the  product formula $X^{\bd'}H_{w'}\bb\cdot X_{\tl}^{c}H_{\tl-1}\cdots H_j=\pm X^{\bd'+ce_{\tl}}H_{w}\,\leftindex^{s_j\cdots s_{\tl-1}}\bb$ with $w=w's_{\tl-1}\cdots s_j$  for a basis element $X^{\bd'}H_{w'}\bb$ of $\Kr{\tl-1}$ with a basis element $X_{\tl}^{c}H_{\tl-1}\cdots H_j$ of the proposed free module shows that a basis of that module is mapped to a basis of $\Kr{\tl}$. The trace map $\Tr$ kills the product unless $\bd'=c=0$ and $w=e$. Moreover, $w=e$ forces $w'=e$ and $j=\tl$, by uniqueness of the factorisation. In particular,  $\Tr$  is zero on all summands in \autoref{Krdec} labelled by $(c,j)\not=(0,\tl)$.
\end{proof}

Let $\Kz{\tl}$ be the span of the $H_w\bb$ with $w\in\langle s_1,\dots,s_{\tl-1}\rangle$. The relations \autoref{RaFHrels} not involving the $X_i$ imply that the span of the $H_w\bb$ is multiplicatively closed and so $\Kz{\tl}$ is a subalgebra. \autoref{P:H0Basis} together with \autoref{notationclash} identifies $\Kz{\tl}$ with $\Hfind[\tl]\otimes A^{\otimes(n-\tl)}$. Its trace  $\Tr(H_w\bb)=\delta_{w,e}\tr(\bb)$ is the trace of \autoref{formlevel1}.

\begin{Lemma}\label{L:inverses}
If $\fct(0)\in R^{\times}$ then $X_i$ is invertible in $\Hhd$ for any $1\leq i\leq n$ with   \[
    X_1^{-1}=-\fct(0)^{-1}\bigl(X_1^{\ell-1}+f_{\ell-1}X_1^{\ell-2}+\dots+f_1\bigr)\qquad\text{and}\qquad
    X_{i+1}^{-1}=\hbar^{2}H_i^{-1}X_i^{-1}H_i^{-1}.
  \]
  In particular, $X_i^{\pm1}\in\Kr{i}$  for all $1\leq i\leq n$.
\end{Lemma}
\begin{proof}For $X_1$ we compute $X_1(X_1^{\ell-1}+f_{\ell-1}X_1^{\ell-2}+\dots+f_1)=\fct(X_1)-\fct(0)=-\fct(0)$, and observe  $X_1^{\pm1}\in\Kr{1}$.

The formula for  $X_{i+1}^{-1}$ holds since $\hbar^2X_{i+1}=H_iX_iH_i$ by \autoref{L:HXH} and $H_i$ is invertible with $H_i^{-1}=\hbar^{-2}(H_i-\eta\taui)$, see \autoref{Rk:hinv}\autoref{Hinv1}. Since $H_i$ and $\taui$ lie in $\Kr{i+1}$, so does $H_i^{-1}$ and then $X_{i+1}^{-1}\in\Kr{i+1}$ by induction.
\end{proof}

By \autoref{L:inverses} $X^{\bd}$ is invertible in $\Hhd$ even for $\bd\in\Z^{n}$, and $X^{\ba-\bd}=X^{\ba}(X^{\bd})^{-1}=\prod_iX_i^{a_i-d_i}$ for $\ba,\bd\in\Z^{n}$.

  Let $\ba\in\Z^{n}$ and $u\in\Hhd$ of $\X$-degree zero. We call $\tl\in\set{1,\dots,n}$ \emph{extremal for $X^{\ba}u$} if $a_i=0$ for $i>\tl$ and $u\in\Kz{\tl}$. Then $X^{\ba}u\in\Kr{\tl}$. Note that $n$ is extremal for every $X^{\ba}u$.

\begin{Lemma}\label{MMargument}
  Assume $\fct(0)\in R^{\times}$. Let $\ba\in\Z^n$ and $u\in\Hhd$ of $\X$-degree $0$ such that $\tl$ is extremal for $X^{\mathsf a}u$.  If $1\leq a_{\tl}\leq\ell-1$ then $\Tr(X^{\ba}u)=0$. If $a_{\tl}=0$ then $\Tr(X^{\ba}u)=\Tr(X^{\ba}\varepsilon_{\tl}(u))$ and $\varepsilon_{\tl}(u)\in\Kz{\tl-1}$. Moreover, for arbitrary $\ba,\bd\in\N^{\ell,n}$ with $\bd\prec\ba$ and any $u\in\Hhd$ of $\X$-degree $0$ we have $\Tr(X^{\ba-\bd}u)=0$.
 \end{Lemma}

\begin{proof} Assume first that $1\leq a_{\tl}\leq\ell-1$ and write
$X^{\ba}=X^{\ba'}X_{\tl}^{a_{\tl}}$ with $X^{\ba'}\in\Kr{\tl-1}$
(explicitly, $a_i'=a_i$ for $i\not=\tl$ and $a_{\tl}'=0$) and
  $u=\sum_{j=1}^{\tl}u_jH_{\tl-1}\cdots H_{j}$ with $u_j\in\Kz{\tl-1}$ (take $c=0$ in \autoref{Krdec} as $u\in\Kz{\tl}$). As $X_{\tl}$ commutes with $\Kz{\tl-1}$, we get $X^{\ba}u=\sum_{j=1}^{\tl}(X^{\ba'}u_j)\,X_{\tl}^{a_{\tl}}H_{\tl-1}\cdots H_{j}$, with $X^{\ba'}u_j\in\Kr{\tl-1}$. By  \autoref{Krdec} we have $\Tr\left((X^{\ba'}u_j)\,X_{\tl}^{a_{\tl}}H_{\tl-1}\cdots H_{j}\right)= 0$ for each $j$ because the element lies in the summand for $(a_{\tl},j)\not=(0,\tl)$. It is therefore killed by $\Tr=\Tr\circ\varepsilon_{\tl}$.
   Thus,  $\Tr(X^{\mathsf a}u)=0$ as desired.

Assume next that $a_{\tl}=0$. Then we have $X^{\ba}\in\Kr{\tl-1}$, so
  $\varepsilon_{\tl}(X^{\ba}u)=X^{\ba}\varepsilon_{\tl}(u)$ by $\Kr{\tl-1}$-linearity, and
  $\Tr(X^{\ba}u)=\Tr(\varepsilon_{\tl}(X^{\ba}u))$ and $\varepsilon_{\tl}(u)\in\Kz{\tl-1}$ by \autoref{Krdec}.

  Finally let $\be=\ba-\bd$ and let $s$ be the largest index with $a_s\not=d_s$.
  Then $e_i=0$ for $i>s$ and $1\leq e_s\leq\ell-1$, because $\bd\prec\ba$. Hence we can apply the above second case in order for
   $\tl=n,n-1,\dots,s+1$ to get
  $\Tr(X^{\be}u)=\Tr(X^{\be}u')$ with
  $u'=\varepsilon_{s+1}\cdots\varepsilon_{n}(u)\in\Kz{s}$. Each step is legitimate since after applying $\varepsilon_{\tl+1},\dots,\varepsilon_n$, the index $\tl$ is extremal for $X^{\theta}\varepsilon_{\tl+1}\cdots\varepsilon_n(u)$, and $e_{\tl}=0$ for $\tl>s$. Then we can apply the first case for $\tl=s$, since $s$ is extremal for $X^{\be}u'$ and $1\leq e_s\leq\ell-1$. We obtain
  $\Tr(X^{\be}u')=0$. This finishes the proof.
\end{proof}
\begin{Proposition}\label{P:triangular}
 Let $\fct(0)\in R^{\times}$. For $\ba,\bd\in\N^{\ell,n}$, $x,y\in\Sym$ and
  $\bb,\bc\in\cBo$ we have
  \[
    \<X^{\ba}H_x\bb,\;\bc H_{y^{-1}}X^{-\bd}\>=
    \begin{cases*}
      0 & if $\bd\prec\ba$,\\
      \delta_{x,y}\,\hbar^{2\ell(x)}\,\tr\bigl(\leftindex^{y}(\bb\bc)\bigr) & if $\bd=\ba$.
    \end{cases*}
  \]
\end{Proposition}

\begin{proof}
  By \autoref{main:bilinearform}, applied to $x=X^{-\bd}$, which is a polynomial in the $X_i$ by \autoref{L:inverses}, with $\psi=\op{id}$, and since the $X_i^{\pm1}$ pairwise commute,
  $\Tr(X^{\ba}H_x\bb\bc H_{y^{-1}}X^{-\bd})
  =\Tr(X^{\ba-\bd}\,H_x\bb\bc H_{y^{-1}})$ with
  $H_x\bb\bc H_{y^{-1}}\in\Kz{n}$. If $\bd\prec\ba$ the latter vanishes by
  \autoref{MMargument}. If $\bd=\ba$ then  $X^{\ba-\bd}=1$. Now the
  claim is true by \autoref{pairing}. (This is applicable because $H_x\bb\bc H_{y^{-1}}$ lies in $\Kz{n}$ which we identify compatibly with $\Tr$ with $\Hfind$.)
  \end{proof}

\begin{Remark}\label{P:radical}
  Assume $\fct(0)=0$. Since $t=0$, the algebra $\HHd$ is graded by the $\X$-degree. Its degree zero
  part $\HHd^{(0)}$ is the subalgebra spanned by the $H_w\bb$, which is isomorphic to $\Hfind$ by
  \autoref{notationclash}, and the positive degree part $\HHd^{(\geq1)}$ is a two-sided ideal with
  $\HHd/\HHd^{(\geq1)}\cong\Hfind$. As $\fct(X_1)=X_1\tilde{\mathtt f}(X_1)\in\HHd^{(\geq1)}$, the
  projection $\HHd\to\Hfind$ factors through $\Hhd$ and yields a surjective algebra homomorphism
  \[
    \pi\map{\Hhd}{\Hfind},\qquad \pi(X_i)=0,\quad \pi(H_i)=H_i,\quad \pi(\bb)=\bb .
  \]
  Since $\pi(X^{\bd}H_w\bb)=\delta_{\bd,0}H_w\bb$, we have $\Tr=\Tr_{\Hfind}\circ\pi$ on the basis and
  hence everywhere. Therefore $\ker\pi\subseteq\rad\<\ ,\ \>$, and equality holds because $\pi$ is
  onto and $\Tr_{\Hfind}$ is non-degenerate by \autoref{formlevel1}. Thus
  $\rad\<\ ,\ \>=\ker\pi$, and
  $\<\ ,\ \>$ is degenerate if and only if $\ell\geq2$.
\end{Remark}

\begin{Proposition}\label{T:main}
 Consider $\Hhd$ with $t=0$ and $\eta\not=0$. Then
   the Gram determinant of $\<\ ,\ \>$ with respect to the bases
  of \autoref{basisthmcycl} and \autoref{C:HdBasis}, recalling $r=|\cB|$, is
  \[
   \det G=\pm\mu\hbar^{\binom{n}{2}\op{rk}\Hhd}
           \fct(0)^{n(\ell-1)\op{rk}\Hhd/\ell},
    \quad\text{for some $\mu\in\kk^{\times}$ depending only on $A$, $n$, $\theta$, $(\ell_0,\ell_1)$}.
  \]
  In particular, for $n\geq1$, $\ell\geq2$, the form $\<\ ,\ \>$ is non-degenerate if and only if
  $\fct(0)\not=0$, that is $\ell_0=0$. It is unimodular if and only if
  $\fct(0)\in R^{\times}$, that is all $\kapp_j\in R^{\times}$.
\end{Proposition}

\begin{Remark}\label{rootsinR}
 \autoref{T:main} and its proof are also valid if  $\fct$ is a  $\theta$-cyclotomic polynomial whose roots lie not in $\kk$ but  in
  $R$.  Unimodularity obviously holds for $n=0$. The case  $\ell=1$  is contained in
\autoref{formlevel1}.
\end{Remark}

\begin{proof} We denote by $G_\Bo=(\tr(\bb\bc))_{\bb,\bc}$ the Gram matrix for $\tr_\Bo$ in the basis $\cBo$ of $\Bo$ induced from $\cB$.
Assume first that $\fct(0)\in R^{\times}$. We show that $\<\ ,\ \>$ is unimodular. For this let $v_I=X^{\ba}H_x\bb$ run
through the basis of \autoref{basisthmcycl} and $u_J=\bc H_{y^{-1}}X^{-\bd}$ through the elements
of \autoref{P:triangular}, which exist by \autoref{L:inverses} if $\fct(0)\in R^{\times}$. We order $I=(\ba,x,\bb)$ by $\ba$ and $J=(\bd,y,\bc)$ by
$\bd$ with respect to $\prec$. By \autoref{P:triangular}, the matrix $M=(\Tr(v_Iu_J))_{I,J}$ is block upper triangular with diagonal blocks
$D_{\ba}=\bigl(\delta_{x,y}\hbar^{2\ell(x)}\tr(\bb\bc)\bigr)_{(x,\bb),(y,\bc)}$, where we used $\tr(\leftindex^{y}(\bb\bc))=\tr(\bb\bc)$, as $\tr$ on $\Bo$ is invariant under the superflip maps. Thus,  $D_{\ba}$ is block diagonal in $x$ with blocks $\hbar^{2\ell(x)}G_\Bo$.  Because $G_\Bo$ is up to signs the $n$-fold Kronecker power of $G_A$ and $\det G_A\in\kk^{\times}$, we obtain
\[
\det D_{\ba}=\hbar^{2r^{n}\sum_{x\in\Sym}\ell(x)}(\det G_\Bo)^{n!}\in R^{\times}.
\]
Hence $\det M\in R^{\times}$ and the functionals $\Tr(\,\cdot\,u_J)$ form an $R$-basis of $\op{Hom}_R(\Hhd,R)$. In particular, the map $\Hhd\to\op{Hom}_R(\Hhd,R)$, $y\mapsto\Tr(\,\cdot\,y)$ is a surjection between free $R$-modules of the same finite rank. Hence it is an isomorphism, and $\<\ ,\ \>$ is unimodular.  In particular $\det G\in R^{\times}$. (A posteriori the $u_J$ form an $R$-basis of $\Hhd$.)

 We next prove that $\det G$ has the asserted form by treating the roots of $\fct$ as variables. Note that $\theta=1$, since $\eta\neq0$ forces $\psi=\op{id}$ by \autoref{rketasymm}. Let $\kk'=\kk[Z_1,\dots,Z_\ell]$ with $\mathtt f_0=\prod_{j=1}^{\ell}(X_1-Z_j)$ and let  $R_0=\kk'[\hbar^{\pm1},\eta]$. Consider  $A'=A\otimes_\kk\kk'$. This is a Frobenius superalgebra over $\kk'$ with basis $\cB$. As $\eta\not=0$ in $R_0$ and $A'$ is symmetric, the $R_0$-algebra $\mathcal H_0=\HHd[n;\hbar,0,\eta](A'^{\otimes n})/(\mathtt f_0)$ is defined. Its structure constants and its trace are computed from \autoref{D:Hecke} modulo $\mathtt f_0$ and \autoref{formform}. Hence its Gram matrix $G_0$ with respect to the bases of \autoref{basisthmcycl} and \autoref{C:HdBasis} has entries in $R_0$, and $\det G_0\in R_0$ is a polynomial in the $Z_j$ and $\eta$ with coefficients in $\kk[\hbar^{\pm1}]$.

To identify this polynomial, pass to the domain $R_1=R_0[(Z_1\cdots Z_\ell)^{-1}]\supseteq R_0$. There, $\mathtt f_0(0)=\pm Z_1\cdots Z_\ell$ is a unit, so the first part of the proof, applied to the CRaFH algebra $\mathcal H_0\otimes_{R_0}R_1$, shows that $\det G_0$ is a unit of $R_1=\kk[Z_1^{\pm1},\dots,Z_\ell^{\pm1}][\hbar^{\pm1},\eta]$. The units of this ring are the monomials $\mu\prod_jZ_j^{m_j}\hbar^{m}$ with $\mu\in\kk^{\times}$ and $m_j,m\in\Z$, since for a domain $D$ the units of $D[\eta]$ are $D^{\times}$ and those of $D[x^{\pm1}]$ are $D^{\times}x^{\Z}$. Therefore, $\det G_0=\mu\prod_jZ_j^{m_j}\hbar^{m}$ for some $\mu\in\kk^{\times}$, $m_j,m\in\Z$. Since both sides lie in $R_0$ it is an equality in $R_0$.

We next study the exponents $m_j$. They are nonnegative, because $\det G_0$ lies in $R_0=\kk[Z_1,\dots,Z_\ell][\hbar^{\pm1},\eta]$ and a monomial $\mu\prod_jZ_j^{m_j}\hbar^{m}$ lies in $R_0$ only if all $m_j\geq0$. Let now  $\Sym_\ell$ act on $R_0$ by permuting the variables $Z_j$ while fixing $\hbar,\eta$. Then $\mathtt f_0=\prod_j(X_1-Z_j)$ is $\Sym_\ell$-invariant, since the coefficients of $\mathtt f_0$ are up to signs the elementary symmetric polynomials in the $Z_j$. Let $c^{K}_{IJ}\in R_0$ be the structure constants of $\mathcal H_0$  defined by $v_Iv_J=\sum_Kc^{K}_{IJ}v_K$ for the basis of \autoref{basisthmcycl}. They depend on the $Z_j$, but only through $\mathtt f_0$. Hence a base change through $\sigma\in\Sym_\ell$
 turns $\mathcal H_0$ into the CRaFH algebra over $R_0$ with parameters $\sigma(\hbar)=\hbar$, $\sigma(\eta)=\eta$ and cyclotomic polynomial $\sigma(\mathtt f_0)=\mathtt f_0$ which is  $\mathcal H_0$ itself.  By  \autoref{basisthmcycl} $\sigma\in\Sym_\ell$ sends the structure constants $c^{K}_{IJ}$ to  $\sigma(c^{K}_{IJ})$ and therefore  $\sigma(c^{K}_{IJ})=c^{K}_{IJ}$.  As the values of $\Tr$ on the basis lie in $\kk$, every entry of $G_0$ is $\Sym_\ell$-invariant and therefore also  $\det G_0$. Thus $\det G_0=\mu\prod_jZ_j^{m_j}\hbar^{m}$ must be symmetric in the $Z_j$ which means $m_i=m_j$ for all $i,j$.  Altogether we get
\[
\det G_0=\mu\prod_jZ_j^{m_j}\hbar^{m}=\pm\mu\,\hbar^{m}\,\mathtt f_0(0)^{m_1}\quad\text{in }R_0\text{, with }m_1\geq0,
\]
in which $\eta$ does not occur.

The exponents follow from two gradings of $\mathcal H_0$, both compatible with the relations \autoref{RaFHrels} because $t=0$, and with $\mathtt f_0$, which is homogeneous of degree $\ell$ for the first and of degree $0$ for the second grading. In each, a product of basis elements of degrees $k_1,k_2$ has coefficients of degree $k_1+k_2-k$ at basis elements of degree $k$, and $\Tr$ reads off the coefficients of the $\bb$, which have degree $0$; so $\det G_0$ is homogeneous of degree $\sum_Ik_I+\sum_Jk_J$, where $k_I,k_J$ are the degrees of the row and column basis elements. For $\deg X_i=\deg Z_j=1$ and $\deg H_i=\deg\bb=\deg\hbar=\deg\eta=0$ the basis elements $X^{\be}H_w\bb$ and $\bb H_wX^{\be}$ have degree $|\be|$, so $\det G_0$ has degree $2\,n!\,r^{n}\sum_{\be\in\N^{\ell,n}}|\be|=n!\,n(\ell-1)\ell^{n}r^{n}=\ell m_1$, whence $m_1=n!\,n(\ell-1)\ell^{n-1}r^{n}$. For $\deg H_i=\deg\hbar=\deg\eta=1$ and $\deg X_i=\deg Z_j=\deg\bb=0$ they have degree $\ell(w)$, so $\det G_0$ has degree $2\,\ell^{n}r^{n}\sum_{w\in\Sym}\ell(w)=n!\tbinom n2\ell^{n}r^{n}$ in $(\hbar,\eta)$, which is $m$ as $\eta$ does not occur.

Finally, $Z_j\mapsto\kapp_j$, $\hbar\mapsto\hbar$, $\eta\mapsto\eta$ defines a $\kk$-algebra homomorphism $R_0\to R$, as the $\kapp_j$ lie in $\kk\subseteq R$. It exists for any $\kapp_j$, including $\kapp_j=0$, because it is defined on $R_0$ and not on $R_1$. It maps $\mathtt f_0$ to $\fct$ and, by \autoref{basisthmcycl} and the definition of $\Tr$, the structure constants and trace of $\mathcal H_0$ to those of $\Hhd$, hence $\det G_0$ to $\det G$. This gives the formula for $\det G$. For $\ell\geq2$ we have $m_1>0$, so $\det G$ is nonzero if and only if $\fct(0)\not=0$ and a unit if and only if $\fct(0)\in R^{\times}$, as $R$ is a domain.
\end{proof}

\subsection{The non-degeneracy theorem}
In this subsection we establish non-degeneracy criteria for $\Tr$ with formulas for the determinants of the Gram matrix for $\langle_-,_-\rangle$.

Let $\Hhd$ be as in \autoref{ass} of level
$\ell$ with level vector $(\ell_0,\ell_1)$. By \autoref{basisthmcycl} it is free
over $R$ of rank $n!\ell^{n}r^{n}$, where $r=|\cB|$ and
$N=n(\ell-1)$ is its top $\X$-degree.  The singular locus of the Gram determinant of the form $\langle_-,_-\rangle$ will depend on the roots of $\fct$. To make this precise we introduce the following elements in $R$.
\begin{Definition}
Write $\fct=\sum_{i=0}^{\ell}a_iX_1^{i}$. We define the \emph{$(t,-\eta)$-homogenisation of $\fct$} to be the element
\begin{equation}\label{E:g}
  \gcy=\sum_{i=0}^{\ell}a_it^{i}(-\eta)^{\ell-i}
      =\prod_{\varrho}(t+\eta\varrho)
      =t^{\ell_0}\prod_{i=1}^{\ell_1}
          \bigl(t^{\theta}-(-\eta)^{\theta}\kapp_i^{\theta}\bigr) \in R,
\end{equation}
 where $\varrho$ runs over the
$\ell$ roots of $\fct$ with multiplicity (i.e., $0$ taken $\ell_0$ times
together with the $\omega^{j}\kapp_i$).
\end{Definition}
\begin{Remark}\label{specialg}
By definition, we have  $\gcy=t^{\ell}$ if $\eta=0$,
and $\gcy=(-\eta)^{\ell}\fct(0)$ if $t=0$.
\end{Remark}

For arbitrary $t,\eta\in R$ let $\T\map{\Hhd}R$ be the $R$-linear map given on the basis of \autoref{basisthmcycl} by
\begin{equation}\label{generalT}
  \T(X^{\bd}H_w\bb)=\delta_{w,e}\,(-\eta)^{N-|\bd|}t^{|\bd|}\,\tr(\bb),\quad\text{so that}\quad \T=\nu^{N}\Tr\text{ with } \nu=
  \begin{cases}
  t&\text{if $\eta=0$ or $t\in R^{\times}$},\\
-\eta&\text{if $t=0\neq\eta$.}
  \end{cases}
\end{equation}
  Let $G_{\T}$ and $G=G_{\Tr}$ respectively be
the Gram matrices of $(x,y)\mapsto T(xy)$ and $\<\ ,\ \>$ with respect to the
bases of \autoref{basisthmcycl}. By \autoref{generalT} they satisfy
\begin{equation}\label{GTG}
  \det G_{\T}=\nu^{N\op{rk}\Hhd}\det G,
\end{equation}
whenever $\Tr$ is defined. Our \emph{non-degeneracy theorem} is now the following:
\begin{Theorem}\label{T:summary}
Assume the setup of the beginning of this subsection. Then the following holds
\begin{enumerate}
\item There exists $\mu\in\kk^{\times}$, depending only on $A$, $n$, $\theta$ and
  $(\ell_0,\ell_1)$, such that we have for all $t,\eta\in R$
  \[
    \det G_{\T}=\pm\mu\,\hbar^{\binom{n}{2}\op{rk}\Hhd}\,\gcy^{\,m},
    \qquad\text{where}\qquad m=n(\ell-1)\op{rk}\Hhd/\ell.
  \]
  \item If $\eta=0$, then
  $\det G=\pm\mu\,\hbar^{\binom{n}{2}\op{rk}\Hhd}$, so $\<\ ,\ \>$ is unimodular.
 \item  If $\eta\neq0$ and $\Tr$ is defined, then $\nu$ from \autoref{generalT} is nonzero, $\gcy\nu^{-\ell}\in R$ and the Gram determinant $\det G$ is
 \begin{equation}\label{madeleine}
 \det G=\pm\mu\,\hbar^{\binom{n}{2}\op{rk}\Hhd}(\gcy\nu^{-\ell})^{m}.
 \end{equation}
In particular, the  form $\<\ ,\ \>$ is non-degenerate if and only if
  $m=0$ or $\gcy\neq0$. It is moreover unimodular if and only if $m=0$ or  ($m>0$ and $\gcy\in\nu^{\ell}R^{\times}$).
  \end{enumerate}
\end{Theorem}

\begin{proof}
  Note that $m=n(\ell-1)\,n!\,\ell^{n-1}r^{n}$ with $m=0$ exactly when
  $n=0$ or $\ell=1$ and $\nu^{N\op{rk}\Hhd}=(\nu^{\ell})^{m}$. We prove the formula separately for the cases  $\eta=0$, $t=0\neq\eta$ and $t\eta\neq0$. They are exhaustive, since $R$ is a domain.

 \noindent\emph{Case $\eta=0$.} We argue similarly as in the proof of \autoref{T:main}. Let $R_2=\kk[Z_1,\dots,Z_{\ell_1}][\hbar^{\pm1},t]$ with formal variables $Z_j$ and  $\hbar$ and let $\mathtt f_2=X_1^{\ell_0}\prod_i(X_1^{\theta}-Z_i^{\theta})$. Denote by
 $\mathcal H_2$ the  $R_2$-CRaFH superalgebra corresponding to $\Hhd$ with $\eta=0$, and let  $G_2$ the Gram matrix of the corresponding $R_2$-bilinear form $\<\ ,\ \>$. By  \autoref{T:nondeg},
\begin{equation}\label{G2det}
\det G_2\in R_2^{\times}=\kk^{\times}\hbar^{\Z}, \;\text{even},\;\det G_2\in\alpha\hbar^{\Z}, \text{where $\alpha\in\kk^{\times}$ only depends on  $A$, $n$, $\theta$ and $(\ell_0,\ell_1).$}
\end{equation}
For the latter  note that $\mathcal H_2$,
 its two bases and $\Tr$ are built from
 $A$, $n$, $\theta$, $(\ell_0,\ell_1)$ only. (The charge and $t$ have been
 replaced by the indeterminates $Z_i$ and $t$.) The same holds then for  $G_2$ and $\det G_2$.

To get a formula for  $\det G_2$ we define a grading on  $\mathcal H_2$ by setting
  $\deg H_i=\deg\hbar=\deg t=1$ and $\deg X_i=\deg Z_i=\deg\bb=0$. As
  $\deg(\gamma_i)=\deg t=1$ and $\deg(\mathtt f_2)=0$, the relations
  \autoref{RaFHrels} are indeed homogeneous. Then  $\deg(X^{\bd}H_w\bb)=\ell(w)=\deg(\bb H_wX^{\bd})$ for any $\bd\in\N^{\ell,n}$, $w\in\Sym$
  and $\bb\in\Bo$. Moreover, $\Tr$ is of  degree $0$. Therefore, every entry of $G_2$ is
 homogeneous of degree the sum of the degrees of the labelling basis elements and  $\det G_2$ is homogeneous of degree the sum of the degrees of all basis elements labelling rows and columns, that is
  \[
    2\sum_{\bd\in\N^{\ell,n}}\sum_{\bb\in\cBo}\sum_{w\in\Sym}\ell(w)
      =2\ell^{n}r^{n}\sum_{w\in\Sym}\ell(w)
      =\binom{n}{2}n!\,\ell^{n}r^{n}=\binom{n}{2}\op{rk}\Hhd,
  \]
  where $\cBo$ is the basis of $\Bo$ induced from $\cB$. Hence by \autoref{G2det},
  $\det G_2=\mu\hbar^{\binom{n}{2}\op{rk}\Hhd}$, with $\mu\in\kk^{\times}$
  depending only on $A$, $n$, $\theta$ and $(\ell_0,\ell_1)$. Specialising
  $Z_i\mapsto\kapp_i$, $\hbar\mapsto\hbar$, which is compatible with the bases and $\Tr$ by
  \autoref{basisthmcycl}, gives the desired formula for $\det G$. The one for $\det G_{\T}$ follows then directly with
   \autoref{GTG} and  \autoref{specialg}.

   \noindent\emph{Case $t=0\neq\eta$.} Here $\nu=-\eta$ and $\gcy=\nu^{\ell}\fct(0)$
  by \autoref{specialg}. By \autoref{T:main}, with the exponents rewritten, we have $\det G=\pm\mu'\hbar^{\binom{n}{2}\op{rk}\Hhd}\fct(0)^{m}$,
  with $\mu'\in\kk^{\times}$ depending only on $A$, $n$, $\theta$ and
  $(\ell_0,\ell_1)$. Now \autoref{GTG} gives the formula \autoref{madeleine}, but with $\mu=\mu'$ .

   \noindent\emph{Case $t\eta\neq0$.} Both sides of the formula  \autoref{madeleine}  lie in $R$.
   Since $R'=R[(t\eta)^{-1}]$ is a domain containing $R$ we scalar extend from $R$ to $R'$,  where $t$ and $\eta$ are units and
  $\nu=t$ and $\lambda=t\eta^{-1}$ is defined. Recall the shift isomorphism $\Sh$ with $\fct'$,
  ${\Hhd}'$ as in \autoref{Phidef}, and let  $G'$ be the Gram matrix of
  $\Tr_{{\Hhd}'}$. Then $\Tr=(-\lambda)^{-N}\Tr_{{\Hhd}'}\circ\Sh$  with $\Sh$
  unitriangular on both bases by \autoref{extro}, so $\det G=(-\lambda)^{-N\op{rk}\Hhd}\det G'$.
  As ${\Hhd}'$ has $t=0\neq\eta$, the previous case applies by \autoref{rootsinR} and gives
  $\det G'=\pm\mu'\hbar^{\binom{n}{2}\op{rk}\Hhd}\fct'(0)^{m}$. By
  \autoref{E:g}, $\fct'(0)=\prod_{\varrho}(-\varrho-\lambda)=(-\eta)^{-\ell}\gcy$,
  and $(-\lambda)(-\eta)=\nu$. So rescaling the formula for $\det G'$ using \autoref{GTG} gives \autoref{madeleine}, again with $\mu=\mu'$.

   \noindent\emph{The constant.} We need to show that the cases produce the same constant $\mu=\mu'$. If $\psi\neq\op{id}$, only the case $\eta=0$ occurs, since $\eta\neq0$ forces $\psi=\op{id}$, and there is nothing to compare.
   Over $R_3=\kk[\hbar^{\pm1},t^{\pm1},\eta]$ we
  have $t\eta\neq0$, so the last case applies. Specialising $\eta\mapsto0$ is
  compatible with the bases and $\T$ by \autoref{basisthmcycl} and
  \autoref{generalT}, and sends $\gcy$ to $t^{\ell}$ by
  \autoref{specialg}. It turns that conclusion into the one of the first case,
  so $\mu'=\pm\mu$.

   \noindent\emph{The criteria.} Let $\eta\neq0$ with $\Tr$ defined. Then
  $\nu\neq0$ and $\gcy\nu^{-\ell}\in R$. It is $\fct(0)$ if $t=0$ and has
  the unit $\nu^{\ell}$ as denominator if $t\in R^{\times}$. So
  \autoref{GTG} gives
  $\det G=\pm\mu\hbar^{\binom{n}{2}\op{rk}\Hhd}(\gcy\nu^{-\ell})^{m}$. As $\mu$
  and $\hbar$ are units and $R$ is a domain, the criteria follow for $m>0$. For
  $m=0$ the determinant is a unit.
\end{proof}

\begin{Remark}\label{R:KMsymmetricconj}
  Recall the notation of \autoref{R:KMbasisconj} with. $A$ a symmetric Frobenius algebra and
  $\eta=0$. Then $\theta=1$ divides $\ell_0$, so $\psi^{\ell_0}=\op{id}$ and
  \autoref{main:bilinearform} gives that $\<\ ,\ \>$ is symmetric, while
  \autoref{T:nondeg} gives that it is unimodular. Hence $\Hhd$ is a symmetric
  Frobenius algebra, which proves \cite[Conjecture~3.22]{KleshchevMuth} for
  scalar charge. It is a statement about graded algebras by
  \autoref{Rk:Zgraded}. The
  cellularity statement \cite[Conjecture~3.23]{KleshchevMuth} might be addressed by
  the respectable bases of \autoref{T:Respectable}.
\end{Remark}

\begin{Remark}\label{R:recover}
  This section gives a uniform treatment of trace
  forms on CRaFH algebras that recover many results in the literature.
  Surprisingly, our approach is new, as are the explicit Gram determinants in \autoref{T:summary}.
\begin{enumerate}
\item Taking $A=\kk$ and $\eta=0$, the CRaFH algebra is a degenerate cyclotomic Hecke algebra. \autoref{T:nondeg} shows that these algebras are symmetric algebras with the same trace as \cite[Theorem~A.2]{BK:HigherSchurWeyl}; see \autoref{R:BKfiltration}.
\item Taking $A=\kk$ and $t=0\neq\eta$, and the normalisation of \autoref{ex:aff}, the CRaFH algebra is an  Ariki--Koike algebras, and $\Tr$ corresponds to the trace form of Bremke and Malle \cite{BM:redwds}. \autoref{T:main} recovers the main result of \cite{MalleMathas}, which shows that this form is symmetric when the cyclotomic parameters are invertible, and it shows that the form is degenerate when $\ell\geq2$ and one of the parameters is zero.
\item Taking $A$ arbitrary and $\eta=0$ recovers \cite[Theorem~6.15]{Savage:AffineWreath}; see \autoref{R:SavageRS}. This case includes the cyclotomic wreath Hecke algebras when $A=\kk G$ and $G$ is a cyclic group cyclic $G$, see \cite[Theorem~2.4]{Cui:degYokonuma}, and the cyclotomic Sergeev superalgebras when $A=\mathfrak{Cl}_1$. The latter are Frobenius by \cite[Corollary~15.6.4]{Kleshchevbook}.
\item Taking $A$ symmetric and $\eta\neq0$ proves that the cyclotomic
quotients of the quantum affine wreath algebras  are symmetric from \cite{RossoSavage:QuantumAffineWreath}.
\item  Taking $A=\kk$, $t=1$, $\eta=\xi-1$ and $\hbar^2=\xi$ we obtain the unified cyclotomic Hecke algebras of \cite[Definition~2.1]{HeHuLin:traceForm}, and $\Tr$ is their trace form \cite[Definition~3.5]{HeHuLin:traceForm}. \autoref{T:summary} recovers \cite[Theorem~3.16]{HeHuLin:traceForm}.
\item As noted in \autoref{badh},
\autoref{T:summary} excludes the CRaFH algebras with $\hbar=0$. In
particular, \autoref{T:summary} does not include the statement that the
cyclotomic nilHecke algebras are symmetric~\cite{SVV:center}. An
explicit form in this case is given in
\cite[Definition~4.11, Theorem~5.11]{HuLiang:nilHecke}.
\end{enumerate}
\end{Remark}

\section{Braided diagonalisable RaFH superalgebras}\label{S:Diagonalisable}
This section develops and studies a framework that allows to define
analogues of the (anti)symmetrisers and of (anti)spherical
Hecke modules for braided diagonalisable RaFH superalgebras. From a symmetric
group perspective, this amounts to defining analogues of the trivial and
sign representations for the superalgebras~$\HHd$.

To define spherical and antispherical modules for $\HHd$ we need to diagonalise
the action of the Hecke generators $H_1,\dots,H_{n-1}$. In fact, it is
enough to diagonalise $H_1$. This is not routine because the quadratic relations in
$\HHd$ involve elements of $\Bo$. Many of the usual arguments break down, and we must
develop our treatment from scratch. To this end, we first study
roots of the quadratic relations $H_i^2-\eta\tau^{(i,i+1)}H_i-\hbar^2=0$
in $\HHd$.  In general, these roots are elements of $\Bo$, rather than
being scalars in $\kk$. Such roots do not necessarily exist
and even if they do exist, they are not necessarily unique. Assuming
existence, we pick roots for $H_1$ and use them to inductively
construct roots for all of the $H_i$. Again, some care is needed here
because these roots need to be chosen in a way that is compatible with
the symmetric group action on $\Bo$. A braided diagonalisable RaFH
superalgebra is then a RaFH superalgebra that contains such roots $\xi^+_i$ and
$\xi^-_i$ that satisfy certain $\xi$-braid relations.

  \subsection{Braided diagonalisability}\label{SS:Diagonalisable}
  We start with the main definition.
  \begin{Definition}\label{D:Diagonalisable}
    The RaFH superalgebra $\HHd=\HHd[n;\hbar, t, \eta](\Bo)$ is \emph{diagonalisable} if $\hbar$ is
    invertible and there exist
    $\xi_{1}^\pm\in R[\eta\taui[1,2]]=R[\eta\taui[2,1]]\subset R\otimes_\kk\Bo$, called \emph{roots},  such that
    $(H_1-\xi_{1}^+)(H_1-\xi_{1}^-)=0=(H_1-\xi_{1}^-)(H_1-\xi_{1}^+)$. Moreover,  $\HHd$ is \emph{braided diagonalisable} if it is diagonalisable and
    the following two \emph{$\xi$-braid relations} hold:
    \begin{equation}\label{etabraid}
    \xi_1^\pm\leftindex^{s_2}{{\xi_1^\pm}}\leftindex^{s_1s_2}{{\xi_1^\pm}}
    =\leftindex^{s_1s_2}{{\xi_1^\pm}}\leftindex^{s_2}{{\xi_1^\pm}}\,\xi_1^\pm.
    \end{equation}
    These $\xi$-braid relations make sense if $n\geq3$, and they hold by convention if
  $n\leq2$.
  \end{Definition}

  \begin{Remark}\label{xiinv}  By convention, $\xi^\pm\in R$ if $\eta=0$. If $\eta\not=0$, then $\taui[1,2]=\taui[2,1]$ by $\eta$-symmetricity, \autoref{rketasymm}. When passing to a suitable localised extension of $R$ where $\eta$ is invertible, the roots $\xi^\pm$ become polynomials in $\taui[1,2]=\taui[2,1]$. In any case the roots satisfy
  \begin{equation}\label{xxinvformula}
  \leftindex^{s_1}\xi_1^\pm=\xi_1^\pm,\quad\xi_{1}^+\xi_{1}^-=\xi_{1}^-\xi_{1}^+,
  \end{equation}
 In fact, $(H_1-\xi_{1}^+)(H_1-\xi_{1}^-)=0$ implies
  $0=H_1^2-(\xi^+_1+\xi^-_1)H_1+\xi^+_1\xi^-_1=(H_1-\xi_{1}^-)(H_1-\xi_{1}^+)$.
\end{Remark}
  \begin{Lemma}\label{conddiag}The superalgebra $\HHd$ is diagonalisable
  if and only if $\hbar$ is
    invertible and there exists \[\text{$\quad\xi_{1}^\pm\in R[\eta\tau^{(1,2)}]=R[\eta\tau^{(2,1)}]\quad $ satisfying}
  \quad\xi_{1}^+\xi_{1}^-=-\hbar^2=\xi_{1}^-\xi_{1}^+\quad\text{ and }\quad\xi_{1}^++\xi_{1}^-=\eta\taui[1,2].\]
\end{Lemma}
\begin{proof}
By expanding the condition $(H_1-\xi_{1}^+)(H_1-\xi_{1}^-)=0$ in the first basis from \autoref{C:BasisHd} we get $\xi_{1}^+\xi_{1}^-=-\hbar^2$ and $\eta\taui[1,2]=\xi_{1}^++\leftindex^{s_1}\xi_{1}^-$, and $(H_1-\xi_{1}^-)(H_1-\xi_{1}^+)=0$ gives $\xi_{1}^-\xi_{1}^+=-\hbar^2$ and $\eta\taui[1,2]=\leftindex^{s_1}\xi_{1}^++\xi_{1}^-$.
The statements follow, since the $\xi^\pm$ are $s_1$-invariant by \autoref{xiinv}.
\end{proof}
\subsubsection*{We first deal with the question of \emph{existence of roots}.}
\begin{Lemma}\label{L:roots}
    Consider $\HHd$ with  $\hbar$ invertible and let $x= \left(\eta\taui[1,2]\right)^2+4\hbar^2$.  If $\eta=0$ or if  $2$ is invertible in $\kk$ and $x$ has a square root
  $\alpha\in R[\eta\tau^{(1,2)}]$, then  $\HHd$  is diagonalisable with
  $2\xi^\pm_1=\eta\taui[1,2]\pm\alpha$. If $\eta=0$ then $\xi^\pm=\pm\hbar$ or $\xi^\pm=\mp\hbar$ gives all possible choices for roots.
\end{Lemma}
\begin{proof}
If $\eta=0$, then $\alpha=\pm2\hbar$ are the unique square roots of $x$ in the domain $R$,  and we get $\xi^\pm=\pm\hbar$ or $\xi^\pm=\mp\hbar$. If $\eta\not=0$,  then   $\xi^\pm_1=\tfrac12(\eta\taui[1,2]\pm\alpha)$ satisfy \autoref{conddiag} by construction, since $\alpha$ commutes with $\taui[1,2]=\taui[2,1]$ by definition.
\end{proof}

 \begin{Remark}
   If $\hbar=0$ then, by definition, $\HHd$ is not diagonalisable. Nevertheless, roots exist, for instance $\xi^+=\eta\taui[1,2]$ and $\xi^-=0$, and $\xi^\pm=0$ are roots if and only if $\eta=0$.
    \end{Remark}

\subsubsection*{We next consider the  $\xi$-braid relations in more detail.}
\begin{Remark}\label{obvious}
   If $\eta=0$, then $\xi^\pm_1=\pm\hbar$ or $\xi^\pm_1=\mp\hbar$,  and the $\xi$-braid
    relation is obvious.  \autoref{Ex:Taft}, \autoref{Ex:Clifford}, and \autoref{Alambda} all give rise to braided
  diagonalisable RaFH superalgebras, where $\Bo=A^{\otimes n}$ with $A$ a
  Frobenius superalgebra with a nontrivial Nakayama automorphism. By $\eta$-symmetricity, if $\eta\not=0$ then
    $\tau^{(1,2)}=\tau^{(2,1)}$, so $\xi_1^\pm$ commutes with $\taui[1,2]$
    and with $H_1$ by definition.
\end{Remark}

\begin{Example} The $\xi$-braid relations hold automatically, and $\HHd$ is braided diagonalisable, if $\tau^{(1,2)}=1^{\otimes n}$ and $\xi^\pm_1$ exist.  This
    case includes the classical Hecke algebra where $A=\kk$ with $q\in R$ invertible,
    $\eta=q^{-1}-q$ and $\hbar=1$. Then $\alpha=\pm(q+q^{-1})$ and $\xi^\pm_1=\pm q^{\mp}$ are the roots.
    Similarly for other normalisations, see \autoref{ex:aff}.
\end{Example}

  \begin{Example}
    If $\xi^\pm_1$ exists, then the $\xi$-braid relations hold whenever
    $\Bo$ is supercommutative, so $\HHd$ is braided diagonalisable.
This includes for example the case
    $\Bo=A^{\otimes n}$ where $A$ is the cohomology algebra of a smooth
    projective variety~$X$ viewed as a Frobenius algebra via the Poincar{\'e} pairing.
  \end{Example}

 \begin{Example}
   Let $A=\F_2[x]/(x^2)$ and $\hbar=1$.  Let $\tau=1\otimes x+x\otimes1$. Then $\tau^2=0$ and $(H_i+1)^2=H_i^2+1$. If $\eta\neq0$, no roots exist.  The system $\xi^+\xi^-=1$, $\xi^+ +\xi^-=\tau$ has no solution in $\F_2[\tau]/(\tau^2)$. With $\eta=0$ roots exist and $\HHd$ is
   braided diagonalisable with $\xi_1^\pm=1$.
 \end{Example}
 \begin{Example}
   Let $A=\F_3[x]/(x^3)$.  Let $\tau=1\otimes x^2+x\otimes x+x^2\otimes
   1$. Then $\tau^2=0$ and we may set $\alpha=2\hbar$. Then
   $\xi^+_1=\hbar+\tfrac{\eta}{2}\tau,\ \xi^-_1={-\hbar}+{\tfrac{\eta}{2}}\tau$.
   If we take for instance $\eta=2$ and $\hbar=1$, then $\HHd$ is
   braided diagonalisable with  $\xi^\pm_1=\tau\pm1$.
 \end{Example}

Interesting non-commutative Frobenius algebras arise in
 geometric and algebraic representation theory, often equipped with a
 positive grading.  The algebra $\mathcal{Z}ig_m$ from  \autoref{Ex:ZigZagm}
 is a very special case of such an algebra.
\begin{Lemma}\label{P:diagbraided}
   Let  $A$ be a positively graded symmetric Frobenius algebra
   with a homogeneous trace form of negative degree $-m$, where $m$ is
   the top degree of~$A$. Assume $\hbar$ is invertible, and that $\eta=0$ or that $2\in\kk^\times$.  Then $\HHd$ is diagonalisable.
 \end{Lemma}
 \begin{proof} Note that $(\taui[1,2])^2$ is homogeneous of top degree in $A\otimes A$, hence $(\taui[1,2])^i=0$ for $i\geq 3$. If $2\in\kk^\times$ then
 ${2\hbar}+{\frac{\eta^{2}}{4\hbar}}(\taui[1,2])^2$ is a square root of $x= \left(\eta\taui[1,2]\right)^2+4\hbar^2$. By \autoref{L:roots}  the claim follows.
    \end{proof}

We next provide diagonalisable RaFH-algebra which are in general not braided diagonalisable.
 \begin{Example}\label{tauZigZag}
    Consider the zig-zag algebra $\mathcal{Z}ig_m$, from \autoref{Ex:ZigZagm} as a symmetric $\mathbb{Z}$-graded Frobenius algebra. Then $\tr$ is homogeneous of degree $-2$ and  $\tau=\taui[1,2]$ is of degree $2$. We have
    \begin{equation*}
 \tau=\sum_{i=1}^{m}\left(e_i\otimes l_i+l_i\otimes e_i\right)
    +\sum_{i=1}^{m-1}\left(a_i\otimes\bar a_i+\bar a_i\otimes a_i\right), \quad
   \tau^2=\sum_{x,y\in \cB}xy\otimes x^\vee y^\vee
     =2\sum_{i=1}^{m}l_i\otimes l_i+\sum_{|i-j|=1}l_i\otimes l_j.
     \end{equation*}
     Then indeed $\tau^3=0$ by degree reasons. If $\eta\not=0$ and $2\in\kk^\times$,  then $2\hbar+\frac{\eta^{2}}{4\hbar}\tau^2$ is a square root of $x= \left(\eta\tau\right)^2+4\hbar^2$ and  $\xi^\pm_1=\pm\hbar+\frac{\eta}{2}\taui\pm\frac{\eta^{2}}{8\hbar}(\taui)^{2}$ are roots. Thus, $\HHd(\mathcal{Z}ig_m^{\otimes n})$ is diagonalisable. Let $\xi^\pm(T)=\pm\tfrac{\eta^{2}}{8\hbar}T^{2}+\tfrac{\eta}{2}T\pm\hbar$. The $\xi$-braid relations are equivalent to the vanishing of
\begin{equation}\label{cheese}
\begin{split}
    \xi^\pm(\taui[1,2])\xi^\pm(\taui[1,3])\xi^\pm(\taui[2,3])-\xi^\pm(\taui[2,3])\xi^\pm(\taui[1,3])\xi^\pm(\taui[1,2])\\
     =\pm\hbar (\tfrac{\eta}{2})^2([\taui[1,2],\taui[1,3]]+ [\taui[1,2],\taui[2,3]]+[\taui[1,3],\taui[2,3]]).
     \end{split}
     \end{equation}
They are in fact equivalent, since the terms of degrees $0$,  $2$ and also of degree $6$ cancel.

  If now $m=2$, then $\pm\hbar\frac{\eta^{2}}{4}\not=0$ implies that \autoref{cheese}
 is a nonzero multiple of
  \begin{align*}
   l_{12}\otimes(a_1\otimes\bar a_1-\bar a_1\otimes a_1)+(a_1\otimes\bar a_1-\bar a_1\otimes a_1)\otimes l_{12}+a_1\otimes l_{21}\otimes\bar a_1 +\bar a_1\otimes l_{12}\otimes a_1,
  \end{align*}
 where $l_{ij}=(l_i-l_j)$.  Thus  $\HHd(\mathcal{Z}ig_2^{\otimes n})$ is braided diagonalisable if $\eta=0$, and  $\HHd(\mathcal{Z}ig_2^{\otimes n})$ is  diagonalisable but not braided diagonalisable if $\eta\not=0$, $2\in\kk^\times$ and $n\geq3$. Similarly for any $m\geq 2$.
     \end{Example}

\begin{Example}\label{semisimple}
Assume $\Bo=A^{\otimes n}$, where $A$ is the semisimple Frobenius
superalgebra $A=M(\mathtt{m}|\mathtt{n})$ of matrices from \autoref{exsemisimple}. One
easily checks that $(\taui[1,2])^2=1\otimes1$.
Assume $\eta\not=0$. If $2\in\kk^{\times}$ and $\eta^{2}+4\hbar^{2}$ has a square root $\alpha\in R$ then $\HHd$ is diagonalisable with
$\xi^\pm_1=\tfrac12\left(\eta\taui[1,2]\pm\alpha\right)=\tfrac12\left(\eta s_1\pm\alpha\right)$. The $\xi$-braid relations read then as
\[\tfrac12\left(\eta\taui[1,2]\pm\alpha\right)\tfrac12\left(\eta\taui[1,3]\pm\alpha\right)\tfrac12\left(\eta\taui[2,3]\pm\alpha\right)
-\tfrac12\left(\eta\taui[2,3]\pm\alpha\right)\tfrac12\left(\eta\taui[1,3]\pm\alpha\right)\tfrac12\left(\eta\taui[1,2]\pm\alpha\right)=0.\]
If we expand, only the terms with two $\tau$ factors do not cancel. These vanish if and only if $\eta^2\alpha=0$ or the action of $(12)(13)+(12)(23)+(13)(23)=(132)+(123)+(132)$ agrees with that of   $(23)(13)+(23)(12)+(13)(12)=(123)+(132)+(123)$. Thus for $\mathtt{m}+\mathtt{n}\geq 2$, $\HHd(\Bo)$ is braided diagonalisable if and only if $\eta^2\alpha=0$, that is $\alpha=0$. If for instance $\eta=q^{-1}-q$, then  $\alpha=q+q^{-1}$ and $\eta\alpha=q^{-2}-q^{2}$, so that,  for $\mathtt{m}+\mathtt{n}\geq 2$,  $\HHd(\Bo)$ is only braided diagonalisable if $q^2=-1$.
\end{Example}
\begin{Remark}The $\xi$-braid relations hold if the roots are solutions of an appropriate Yang-Baxter equation. This property is explored, independently of us,  in a nice way in \cite{MinetsYBE}.
\end{Remark}

 Unless otherwise stated, we assume that $\HHd$ is
 braided diagonalisable for the rest of this section.

 \subsection{Compatible choice of roots}
We next define roots for all quadratic relations in a compatible way.

 \begin{Definition}\label{D:xis}
    Suppose that $\HHd$ is braided diagonalisable. For $1\le i\le n-2$, set
    ${\xi_{i+1}^\pm} =\leftindex^{s_is_{i+1}}{{\xi_{i}^\pm}}$.
 \end{Definition}

 \begin{Proposition}\label{P:Hieigenvalues}
   Assume $\HHd$ is braided diagonalisable and let $1\le i<n$. Then $\xi_{i}^\pm\in R[\eta\taui[i,i+1]]$ and
   $(H_i-\xi_{i}^{+})(H_i-\xi_{i}^{-})=0=(H_i-\xi_{i}^{-})(H_i-\xi_{i}^{+})$.  The \emph{general $\xi$-braid relations}  for $1\leq i\le n-2$ hold,
    \begin{equation}\label{genxibraid}
    \xi_i^\pm\leftindex^{s_{i+1}}{{\xi_i^\pm}}\leftindex^{s_is_{i+1}}{{\xi_i^\pm}}
    =\leftindex^{s_is_{i+1}}{{\xi_i^\pm}}\leftindex^{s_{i+1}}{{\xi_i^\pm}}\,\xi_i^\pm.
    \end{equation}
 \end{Proposition}
\begin{Corollary}\label{xxxinvformulas} Formulas as in \autoref{xxinvformula} hold for $\xi_i^\pm$, $1\le i<n$, i.e., when replacing all indices $1$ by $i$.
\end{Corollary}
 Comparing with the quadratic relation in \autoref{D:Hecke} we have:
 \begin{Corollary}\label{C:xii}
    If $\HHd$ is braided diagonalisable then $\xi^+_i\xi^-_i=-\hbar^2$ and
  $\xi^+_i+\xi^-_i=\eta\taui$ hold for $1\le i<n$.
 \end{Corollary}

 \begin{proof}[Proof of \autoref{P:Hieigenvalues}] We prove the three claims for ${\xi_{i}^\pm}$ by induction on~$i$. The first claim is obvious. We prove the other two claims for ${\xi_{i}^\pm}$ by induction on~$i$.
   If $i=1$ there is nothing to do by
   \autoref{D:Diagonalisable}. By induction, ${\xi_{i}^\pm}\in R[\eta\taui[i,i+1]]$, so
   ${\xi_{i+1}^\pm}\in\leftindex^{s_is_{i+1}}R[\eta\taui[i,i+1]]=R[\eta\taui[i+1,i+2]]$.

   For the second claim, since~$\hbar$ is invertible, so is $H_k$ for
   $1\le k<n$ with $H_k^{-1}=\hbar^{-2}(H_k-\eta\taui[k,k+1])$. Thus,
   $H_iH_{i+1}\cdot H_i\cdot(H_iH_{i+1})^{-1}
          = (H_iH_{i+1}H_i)\cdot H_{i+1}^{-1}H_i^{-1}
          = (H_{i+1}H_iH_{i+1})\cdot H_{i+1}^{-1}H_i^{-1}
          = H_{i+1}$,
          since the $H_1,\dots,H_{n-1}$ satisfy the braid relations.
   Therefore,  by induction,
   \begin{align*}
     0 &= H_iH_{i+1}\cdot(H_i-{\xi_{i}^{+}})
        \cdot H_{i+1}^{-1}H_i^{-1} \cdot H_iH_{i+1}
        \cdot(H_i-\xi_{i}^{-})\cdot H_{i+1}^{-1}H_i^{-1}
         = (H_{i+1}-\leftindex^{s_is_{i+1}}{{\xi_{i}^{+}}})
         (H_{i+1}-\leftindex^{s_is_{i+1}}{{\xi_{i}^-}})\\
         &= (H_{i+1}-\xi_{i+1}^+)(H_{i+1}-\xi_{i+1}^-).
   \end{align*}
  By \autoref{xxinvformula}, it remains to prove \autoref{genxibraid}.
   When $i=1$ this holds by \autoref{D:Diagonalisable} and so we  assume that \autoref{genxibraid} holds for $\xi_i^\pm$ and want to show that $\xi_{i+1}^\pm\leftindex^{s_{i+2}}{{\xi_{i+1}^\pm}}\leftindex^{s_{i+1}s_{i+2}}{{\xi_{i+1}^\pm}}
    =\leftindex^{s_{i+1}s_{i+2}}{{\xi_{i+1}^\pm}}\leftindex^{s_{i+2}}{{\xi_{i+1}^\pm}}\xi_{i+1}^\pm$.  We recall that
   $\leftindex^{s_j}{{\xi_i^\pm}}=\xi_i^\pm$ for $j=i$ and $j\geq i+2$ by definition and \autoref{xxinvformula}, and calculate

   \begin{equation*}
\begin{gathered}
     \leftindex^{s_is_{i+1}s_{i+2}}{{\xi_i^\pm}}
        =\leftindex^{s_is_{i+1}}{{\xi_i^\pm}} =\xi_{i+1}^\pm,\qquad
     \leftindex^{s_is_{i+1}s_{i+2}}{{(\leftindex^{s_{i+1}}{{\xi_i^\pm}})}}
        = \leftindex^{s_i(s_{i+2}s_{i+1}s_{i+2})}{{\xi_i^\pm}}
         = \leftindex^{s_{i+2}s_is_{i+1}}{{\xi_i^\pm}}
         = \leftindex^{s_{i+2}}{{\xi_{i+1}^\pm}},\\
     \leftindex^{s_is_{i+1}s_{i+2}}{{(\leftindex^{s_is_{i+1}}{{\xi_i^\pm}})}}
        =\leftindex^{s_is_{i+1}s_{i+2}s_is_{i+1}}{{\xi_i^\pm}}
          =\leftindex^{s_{i+1}s_is_{i+2}s_{i+1}s_{i+2}}{{\xi_i^\pm}}
         =\leftindex^{s_{i+1}s_{i+2}}{{(\leftindex^{s_is_{i+1}}{{\xi_i^\pm}})}}
         =\xi_{i+2}^\pm
         =\leftindex^{s_{i+1}s_{i+2}}{{\xi_{i+1}^\pm}}.
\end{gathered}
\end{equation*}
   Hence, acting on the identity \autoref{genxibraid}  by $s_is_{i+1}s_{i+2}$ gives the required identity for
   $\xi_{i+1}^\pm$.
 \end{proof}

 \begin{Remark}
   By \autoref{P:Hieigenvalues}, $\xi_i^+$ and $\xi_i^-$ are
   roots of the quadratic equation $H_i^2-\eta\taui[i,i+1] H_i-\hbar^2=0$. In
   general, these roots are not unique because $\HHd$ contains zero divisors.
   \autoref{P:Hieigenvalues} should be viewed as fixing a choice of
   eigenvalues for~$H_i$, which are the roots of this quadratic
   equation. For Iwahori--Hecke algebras, these roots depend only on the
   conjugacy classes of the corresponding simple reflections.
   \autoref{P:Hieigenvalues} ensures that this is marginally more
   complicated for RaFH superalgebras, where the roots of conjugate simple
   reflections are conjugate. (Of course, all simple reflections are
   conjugate in type~$A$.)
 \end{Remark}

\subsection{Diagonalisable coefficients}
 This subsection introduces analogues $\xi^\pm_w$ of the elements $\xi^\pm_i$
 for each permutation $w\in W$. In the case of the Iwahori--Hecke
 algebras, when the Frobenius superalgebra is trivial,
 $\xi^\pm_w=(\pm q)^{\mp\ell(w)}$ for any $w\in\Sym$.
 In general, these coefficients are more
 complicated and it is often wise not to unpack their dependence on $\Bo$.  Instead we use skew-derivations of $\Bo$,  in the sense of \autoref{skewDer}, to handle them, see \autoref{D:phi} below.

 There is also a (minor) dependence of the formulas on $\hbar$ which is however only a scaling factor. For simplicity, one may assume  that $\hbar=1$. Under the assumption that $\hbar$ is invertible, this is not a further restriction, since rescaling the generators $H'_i=\hbar^{-1}H_i$ and setting $\eta'=\hbar^{-1}\eta$ and $t'=t$ induces an identification of $\RR$-algebras $\HHd[n;\hbar, t, \eta](\Bo)=\HHd[n;1, t', \eta'](\Bo)$ which is the identity on the $X_i$'s and on $\Bo$. We nevertheless keep $\hbar$ in the formulas below.
  \begin{Definition}\label{D:hw}
   Let $\bi=(i_1,\dots,i_k)$, where $1\le i_j<n$ for $1\le j\le k$.
   Define $\xi^\pm(\bi)\in\Bo$ inductively by
   \[
        \xi^\pm(\bi)=\xi^\pm(i_1,\dots,i_k) = \begin{cases*}
            \xi_{i_1}^\pm\leftindex^{s_{i_1}}{{\xi^\pm(i_2,\dots,i_k)}} & if $k\geq1$,\\
            1& if $k=0$.
        \end{cases*}
   \]
   In particular, if $w=s_{i_1}\dots s_{i_k}$ is a reduced expression
   for $w\in W$, set $\xi^\pm_w=\xi^\pm(i_1,\dots,i_k)\in\Bo$.
 \end{Definition}

An analogue of \autoref{D:hw} in case $\Bo$ is commutative can be found in
 \cite[Section~5.4]{LaiMinets}.

 We need to show that $\xi^\pm_w$ is well-defined, i.e.\ independent of the choice of reduced
 expression for~$w$.

 \begin{Lemma}\label{L:hw}
    Let $w\in W$. The element $\xi^\pm_w$ is independent of the choice of
    reduced expression for $w$.
 \end{Lemma}

 \begin{proof}
   It is enough to show that $\xi^\pm(i,i+1,i)=\xi^\pm(i+1,i,i+1)$ and that
   $\xi^\pm(i,j)=\xi^\pm(j,i)$ when $|i-j|>1$, for all admissible~$i$
   and~$j$. The first condition is equivalent to the identity
   $\xi_i^\pm\leftindex^{s_{i+1}}{{\xi_i^\pm}}\leftindex^{s_is_{i+1}}{{\xi_i^\pm}}
   =\leftindex^{s_is_{i+1}}{{\xi_i^\pm}}\leftindex^{s_{i+1}}{{\xi_i^\pm}}\,\xi_i^\pm$,
   which is the general braid relation from  \autoref{P:Hieigenvalues}.
    The requirement that
   $\xi^\pm(i,j)=\xi^\pm(j,i)$ when $|i-j|>1$ is equivalent to
   $\xipm_i\leftindex^{s_i}{{\xipm_j}}=\xipm_j\leftindex^{s_j}{{\xipm_i}}$. This holds because  $\leftindex^{s_i}{{\xipm_j}}=\xipm_j$ and $\leftindex^{s_j}{{\xipm_i}}=\xipm_i$ by \autoref{xxxinvformulas}, and because $\taui[i,i+1]$ and $\taui[j,j+1]$ commute by \autoref{CommutingTaus}.
 \end{proof}

 \begin{Corollary}\label{C:hdv}
     Suppose that $w=dv\in W$ with $\ell(w)=\ell(d)+\ell(v)$. Then
     $\xi^\pm_w=\xi^\pm_d\leftindex^d{{\xi^\pm_v}}$.
 \end{Corollary}

 \begin{proof}
     Let $d=s_{i_1}\dots s_{i_k}$ and $v=s_{j_1}\dots s_{j_l}$ be
     reduced expressions for $d$ and $v$. Then $w=s_{i_1}\dots
     s_{i_k}s_{j_1}\dots s_{j_l}$ is a reduced expression for~$w$, so
     the result follows in view of \autoref{L:hw}.
 \end{proof}

 \begin{Corollary}\label{C:hwrev}
  Let $w\in\Sym$ with reduced expression $w=s_{i_1}\cdots s_{i_k}$ and let
  $u_m=s_{i_1}\cdots s_{i_m}$ with $u_0=1$. Then
  \[\leftindex^w{{\xipm_{w^{-1}}}}
      =\leftindex^{u_{k-1}}{{\xipm_{i_k}}}\cdots
       \leftindex^{u_{1}}{{\xipm_{i_2}}}\;\xipm_{i_1},\]
 hence  $\leftindex^w{{\xipm_{w^{-1}}}}$ is the twisted product  from  \autoref{D:hw} read backwards. Moreover,
  \[\xi^+_w\,\leftindex^w{{\xi^-_{w^{-1}}}}=(-\hbar^2)^{\ell(w)}
      =\leftindex^w{{\xi^-_{w^{-1}}}}\,\xi^+_w\quad\text{ and }\quad
    \xi^-_w\,\leftindex^w{{\xi^+_{w^{-1}}}}=(-\hbar^2)^{\ell(w)}
      =\leftindex^w{{\xi^+_{w^{-1}}}}\,\xi^-_w.\]
\end{Corollary}

 \begin{proof}
   Since $w^{-1}=s_{i_k}\cdots s_{i_1}$ is a reduced expression,
   $\xipm_{w^{-1}}=\xipm_{i_k}\leftindex^{s_{i_k}}{{\xipm_{i_{k-1}}}}\cdots
   \leftindex^{s_{i_k}\cdots s_{i_2}}{{\xipm_{i_1}}}$ by \autoref{D:hw}. On the other hand,
   $\leftindex^{s_{i_1}\cdots s_{i_k}s_{i_k}\cdots s_{i_{m+1}}}{{\xipm_{i_m}}}
   =\leftindex^{s_{i_1}\cdots s_{i_m}}{{\xipm_{i_m}}}
   =\leftindex^{s_{i_1}\cdots s_{i_{m-1}}}{{\xipm_{i_m}}}$, using
   $\leftindex^{s_{i_m}}{{\xipm_{i_m}}}=\xipm_{i_m}$ by \autoref{xxxinvformulas}. This is the first claim.
   For the second,
   $\xi^+_w\,\leftindex^w{{\xi^-_{w^{-1}}}}
      =\leftindex^{u_0}{{\xi^+_{i_1}}}\cdots\leftindex^{u_{k-2}}{{\xi^+_{i_{k-1}}}}
        \left(\leftindex^{u_{k-1}}{(\xi^+_{i_k}\xi^-_{i_k})}\right)
        \leftindex^{u_{k-2}}{{\xi^-_{i_{k-1}}}}\cdots\leftindex^{u_0}{{\xi^-_{i_1}}}=(-\hbar^2)^{\ell(w)}$ by applying  $\xi^+_i\xi^-_i=-\hbar^2=\xi^-_i\xi^+_i$ from
   \autoref{C:xii} successively. The rest follows similarly.
 \end{proof}

\begin{Remark}\label{R:GH}
  Since $H_i\bb=\leftindex^{s_i}\bb H_i$ in $\HHd$, we have
  $\xipm_wH_w=(\xipm_{i_1}H_{i_1})(\xipm_{i_2}H_{i_2})\cdots(\xipm_{i_k}H_{i_k})$
  for every reduced expression $w=s_{i_1}\cdots s_{i_k}$. Since
  $H_w$ is invertible, \autoref{L:hw} is equivalent to the statement
  that the elements $\hat{H}^\pm_i=\xipm_iH_i\in\Hfind\subseteq\HHd$ satisfy the braid
  relations. This is where the $\xi$-braid relations \autoref{etabraid} come from. Consequently
  $\hat{H}^\pm_w=\xipm_wH_w$ is well-defined and
  $\hat{H}^\pm_{dv}=\hat{H}^\pm_d\hat{H}^\pm_v$ whenever $\ell(dv)=\ell(d)+\ell(v)$.
  \end{Remark}
\subsection{Moving roots past $\Bo$}
 We describe next how elements of $\Bo$ move past the coefficients $\xipm_w$. As for the commutation of $H_i$ with $\Po$ in \autoref{L:crucial}, we have a smash product relation with a correction term controlled by a skew-derivation. The relevant skew-derivations are the following.

\begin{Definition}\label{D:phi}
  For $w\in\Sym$ let $\phi_w=\op{id}-w\colon\Bo\to\Bo$, that is
  $\phi_w(\bb)=\bb-\leftindex^w\bb$ for $\bb\in\Bo$. We abbreviate $\phi_i=\phi_{s_i}$ for $1\le i<n$.
\end{Definition}

\begin{Lemma}\label{L:phi}
  Let $u,v\in\Sym$, $1\le i,j<n$ and $\bb,\bc\in\Bo$.
  \begin{enumerate}
  \item\label{phi:skew} $\phi_u$ is a $u$-skew-derivation of $\Bo$, that is
    $\phi_u(\bb\bc)=\phi_u(\bb)\bc+\leftindex^u\bb\,\phi_u(\bc)$.
  \item\label{phi:equiv} $\phi_{uv}=\phi_u+u\phi_v$ and $v\phi_uv^{-1}=\phi_{vuv^{-1}}$.
    In particular
    \begin{gather*}
      \phi_is_i=s_i\phi_i=-\phi_i,\quad\phi_i^2=2\phi_i,\quad
      \phi_is_j=s_j\phi_i\quad \text{ if } |i-j|>1,\\
      \phi_js_is_j=s_is_j\phi_i,\quad \phi_is_js_i=s_js_i\phi_j\quad \text{ if } j=i+1 .
    \end{gather*}
  \item\label{phi:kernel} $\phi_u(\bb)=0$ if and only if $\leftindex^u\bb=\bb$.
    In particular, $\phi_i$ vanishes on $R[\eta\taui[i,i+1]]$ and hence on $\xipm_i$.
  \item\label{phi:translates} For $1\le i_1,\dots,i_l<n$,
    the element $\phi_{i_1}\circ\dots\circ\phi_{i_l}(\bb)$ is a $\Z$-linear combination of
    the $\leftindex^x\bb$ with $x\in \langle s_{i_1},\dots,s_{i_l}\rangle\subseteq\Sym$.
  \end{enumerate}
\end{Lemma}

\begin{proof}
  For \ref{phi:skew} note that
  $\bb\bc-\leftindex^u(\bb\bc)
   =(\bb-\leftindex^u\bb)\bc+\leftindex^u\bb(\bc-\leftindex^u\bc)$,
  as  $u$ acts on $\Bo$ by an algebra automorphism.  For \ref{phi:equiv},
  $\op{id}-uv=(\op{id}-u)+u(\op{id}-v)$ and $v(\op{id}-u)v^{-1}=\op{id}-vuv^{-1}$ as
  operators on $\Bo$. The displayed identities follow from $s_i^2=1$, from
  $s_is_j=s_js_i$ if $|i-j|>1$, and from the braid relation
  $s_is_js_is_js_i=s_j$ for $j=i+1$ respectively.  For \ref{phi:kernel}  only the second part needs an argument.  If $\eta=0$ then $R[\eta\taui[i,i+1]]=R$, and if
  $\eta\not=0$ then $\leftindex^{s_i}{{\taui[i,i+1]}}=\taui[i+1,i]=\taui[i,i+1]$ by
  $\eta$-symmetricity, see \autoref{rketasymm}. Part \ref{phi:translates} follows by
  expanding the product $(\op{id}-s_{i_1})\circ\dots\circ(\op{id}-s_{i_l})$.
\end{proof}

\begin{Remark}\label{R:phivsdelta}
The $\phi_u$ extend by the same formula to $u$-skew-derivations of $\Po$, which are similar to the Frobenius--Demazure operators from \autoref{D:FrobDem}. If $\psi$ is trivial, then
  $\delta_i=\taui[i,i+1]\partial_i$ on $R[\X]$ by \autoref{R:classical}, whereas
  $\phi_i=(X_i-X_{i+1})\partial_i$ on $R[\X]$ by definition of the divided
  difference operator $\partial_i$. Hence $\taui[i,i+1]\phi_i=(X_i-X_{i+1})\delta_i$
  on $R[\X]$. On $\Bo$, however, $\delta_i$ vanishes while $\phi_i$
  measures the failure of $s_i$-invariance.
\end{Remark}

\begin{Lemma}\label{L:alpha}
  There exist elements $\alpha_1,\dots,\alpha_{n-1}\in\Bo$ such that for all $1\le i<n$
  we have
  \begin{enumerate}
    \item\label{alpha:def} $\alpha_i\in R[(\eta\taui[i,i+1])^2]$ with
    $\xipm_i-\alpha^\pm_i\in \eta\taui[i,i+1]R[(\eta\taui[i,i+1])^2]$, where
    $\alpha^\pm_i:=\pm\alpha_i$.
  \item\label{alpha:compatible} $\alpha_{i+1}=\leftindex^{s_is_{i+1}}\alpha_i$ if $i<n-1$.
  \item\label{alpha:central} $\alpha_i\in Z(\Bo)$, $\leftindex^{s_i}\alpha_i=\alpha_i$, and $\leftindex^{s_j}\alpha_i=\alpha_i$ for $|i-j|>1$.
  \item\label{alpha:twist} $\bb(\xipm_i-\alpha^\pm_i)=(\xipm_i-\alpha^\pm_i)\leftindex^{s_i}\bb$ for all $\bb\in\Bo$.
  \end{enumerate}
\end{Lemma}

\begin{proof}
  By \autoref{D:Diagonalisable}, $\xi^+_1=f(\eta\taui[1,2])$ for some polynomial
  $f\in R[z]$. Write $f(z)=f_0(z^2)+zf_1(z^2)$ with $f_0,f_1\in R[z]$ and set
  $\alpha_1=f_0\left((\eta\taui[1,2])^2\right)$ and then inductively
  $\alpha_{i+1}=\leftindex^{s_is_{i+1}}{{\alpha_i}}$ for $1\le i<n-1$ which gives \ref{alpha:compatible}.

  Since
  $\leftindex^{s_is_{i+1}}{{\taui[i,i+1]}}=\taui[i+1,i+2]$ and
  $\xipm_{i+1}=\leftindex^{s_is_{i+1}}{{\xipm_i}}$ by \autoref{D:xis}, we get
  $\alpha_i=f_0\left((\eta\taui[i,i+1])^2\right)$ and
  $\xi^+_i=f(\eta\taui[i,i+1])$ for all $i$, so
  $\xi^+_i-\alpha_i=\eta\taui[i,i+1]f_1\left((\eta\taui[i,i+1])^2\right)$.
  As $\xi^-_i=\eta\taui[i,i+1]-\xi^+_i$ by \autoref{C:xii}, also
  $\xi^-_i+\alpha_i=\eta\taui[i,i+1]\left(1-f_1\left((\eta\taui[i,i+1])^2\right)\right)$.
  This proves \ref{alpha:def}.  If  $\eta=0$ then $\alpha_i=\xi^+_i=\pm\hbar$ and $\xi^-_i=-\xi^+_i$,
  and therefore $\xi^+_i-\alpha_i=0=\xi^-_i+\alpha_i$, and
  \ref{alpha:central} and \ref{alpha:twist} hold. If $\eta\not=0$
  then $\psi$ is trivial by $\eta$-symmetricity. By \autoref{Nakayama}, see also \autoref{athroughtau},
  $(\taui[i,i+1])^2$, and then also $\alpha_i$, is central and \ref{alpha:central}  holds.
  Moreover, $\bb\eta\taui[i,i+1]g\left((\eta\taui[i,i+1])^2\right)=\eta\taui[i,i+1]g\left((\eta\taui[i,i+1])^2\right)\leftindex^{s_i}\bb$
  for $g\in R[z]$, which is \ref{alpha:twist}.
\end{proof}

We fix elements $\alpha_1,\dots,\alpha_{n-1}$ as in \autoref{L:alpha} for the rest of this section.
We like to get formulas for moving $\xi^\pm_w$ past elements $\bb\in\Bo$. The next result is the analogue of \autoref{L:crucial} for the roots.

\begin{Lemma}\label{L:xicomm}
  Let $1\le i<n$ and $\bb\in\Bo$. Then
  \begin{equation}\label{xicomm}
    \bb\,\xipm_i=\xipm_i\,\leftindex^{s_i}\bb+\alpha^\pm_i\phi_i(\bb),
    \qquad
    \leftindex^{s_i}{\left(\bb\,\xipm_i\right)}=\xipm_i\,\bb-\alpha^\pm_i\phi_i(\bb),
        \qquad
      \xipm_i\,\bb=\leftindex^{s_i}\bb\,\xipm_i+\alpha^\pm_i\phi_i(\bb).
  \end{equation}
\end{Lemma}

\begin{proof}
   The first formula holds, since
$\bb\,\xipm_i=\bb\alpha^\pm_i+\bb(\xipm_i-\alpha^\pm_i)
      =\alpha^\pm_i\bb+(\xipm_i-\alpha^\pm_i)\leftindex^{s_i}\bb
      =\xipm_i\leftindex^{s_i}\bb+\alpha^\pm_i\left(\bb-\leftindex^{s_i}\bb\right),
 $
   by \autoref{L:alpha}\autoref{alpha:central} and \autoref{alpha:twist}. Since $\xipm_i$ and $\alpha_i$ are $s_i$-invariant,  \autoref{xxxinvformulas} and  \autoref{L:alpha}\autoref{alpha:central}, and $s_i\phi_i=-\phi_i$, \autoref{L:phi}\autoref{phi:equiv}, the other formulas follow by applying $s_i$.
\end{proof}
We have more general commutation formulas as follows.
 \begin{Lemma}\label{L:axirec}
   Suppose that $w=s_iv$ with $\ell(w)=\ell(v)+1$, and let $\ba\in\Bo$. Then
   \[
     \ba\,\xipm_w
       =\xipm_i\,\leftindex^{s_i}{\left(\ba\,\xipm_v\right)}
        -\alpha^\pm_i\leftindex^{s_i}{\left(\phi_i(\ba)\,\xipm_v\right)}.
   \]
 \end{Lemma}

 \begin{proof}
   By \autoref{C:hdv}, $\xipm_w=\xipm_i\leftindex^{s_i}{{\xipm_v}}$, and  $\ba\xipm_i\leftindex^{s_i}{{\xipm_v}}
    =\xipm_i\leftindex^{s_i}\ba\,\leftindex^{s_i}{{\xipm_v}}
     +\alpha_i^\pm\phi_i(\ba)\leftindex^{s_i}{{\xipm_v}}$ by the first
   formula in \autoref{xicomm}.
   The claim follows, since $s_i\phi_i=-\phi_i$ by \autoref{L:phi}\autoref{phi:equiv}.
 \end{proof}

\begin{Remark}\label{R:Theta}
  Set $\Theta^\pm_i=\alpha^\pm_i\phi_i$, an $s_i$-skew-derivation of $\Bo$ because
  $\alpha_i$ is central, by the observation preceding \autoref{L:gammacentral}
  applied to $\Bo$. In this notation \autoref{xicomm} turns into
  \[
     \xipm_i\,\bb=\leftindex^{s_i}\bb\,\xipm_i+\Theta^\pm_i(\bb),
  \]
  which is \autoref{crucial} under the dictionary
  $(H_i,\Delta_i,\delta_i,-\gamma_i)\leftrightarrow(\xipm_i,\Theta^\pm_i,\phi_i,\alpha^\pm_i)$.
  Note that $\Theta^\pm_i(\bb)=\bb\,\xipm_i-\xipm_i\leftindex^{s_i}\bb$ does
  not depend on the choice of $\alpha_i$ in \autoref{L:alpha}, and that
  $\Theta^-_i=-\Theta^+_i$. The analogues of the relations \autoref{RaFHDerA} and \autoref{RaFHDerB} hold,   by \autoref{L:phi}\ref{phi:equiv} and \autoref{L:alpha}\ref{alpha:compatible},\ref{alpha:central},  in
  the following simpler form
  \begin{gather*}
    \Theta^\pm_is_j=s_j\Theta^\pm_i,\quad \Theta^\pm_i\Theta^\pm_j=\Theta^\pm_j\Theta^\pm_i\ \ (|i-j|>1),\qquad
    s_is_j\Theta^\pm_i=\Theta^\pm_js_is_j,\quad \Theta^\pm_is_js_i=s_js_i\Theta^\pm_j\ \ (j=i+1),\\
    (\Theta^\pm_i)^2=2\alpha^\pm_i\Theta^\pm_i,\qquad \Theta^\pm_is_i=s_i\Theta^\pm_i=-\Theta^\pm_i.
  \end{gather*}
\end{Remark}

For a given reduced expression $w=s_{i_1}\dots s_{i_k}$  for $w\in\Sym$ and $J=\set{j_1<\dots<j_l}\subseteq\set{1,\dots,k}$ set
 \begin{equation}\label{defy}
   \phi_J=\phi_{i_{j_l}}\circ\dots\circ\phi_{i_{j_1}}
   \quad\text{and}\quad
   y^\pm_{J,w}=y^\pm_{J,1}\leftindex^{s_{i_1}}{{y^\pm_{J,2}}}\dots
   \leftindex^{{s_{i_1}\dots s_{i_{k-1}}}}{{y^\pm_{J,k}}},\quad\text{where}\quad
      y^\pm_{J,m} = \begin{cases*}
          \alpha^\pm_{i_m} & if $m\in J$,\\
          \xi^\pm_{i_m}    & if $m\notin J$.
        \end{cases*}
 \end{equation}
 By convention, $\phi_\emptyset=\op{id}$ and $y^\pm_{\emptyset,w}=1$ for $w=1$. Note that $y^\pm_{\emptyset,w}=\xipm_w$ by \autoref{D:hw}, and $y^\pm_{\set{1,\dots,k},w}$ is central.

 \begin{Proposition}\label{L:axi}
   Suppose $w=s_{i_1}\dots s_{i_k}$ is a reduced expression for
   $w\in\Sym$ and  $\ba\in\Bo$. Then
   \begin{equation}\label{apastxis}
     \ba\xi^\pm_w = \sum_{J\subseteq\set{1,2,\dots,k}}
       (-1)^{|J|}y^\pm_{J,w}\leftindex^w{{\phi_J(\ba)}}.
   \end{equation}
 \end{Proposition}

  \begin{proof}
   We argue by induction on $k=\ell(w)$, the case $k=0$ being trivial. Let $k>0$
   and $v=s_{i_2}\dots s_{i_k}$, so that $w=s_{i_1}v$ with $\ell(w)=\ell(v)+1$. For $K\subseteq\set{2,\dots,k}$ let $y^\pm_{K,v}$ and
   $\phi_K$ be defined with respect to the reduced expression $v=s_{i_2}\dots s_{i_k}$. Then
   we have $y^\pm_{K,w}=\xipm_{i_1}\leftindex^{s_{i_1}}{{y^\pm_{K,v}}}$ and
   $y^\pm_{K\cup\set1,w}=\alpha^\pm_{i_1}\leftindex^{s_{i_1}}{{y^\pm_{K,v}}}$, as well as
   $\phi_{K\cup\set1}=\phi_K\circ\phi_{i_1}$.
   Applying \autoref{L:axirec}, and then the induction hypothesis to the elements $\ba$ and $\phi_{i_1}(\ba)$ of $\Bo$, we obtain
   \begin{align*}
     \ba\xipm_w
       &=\xipm_{i_1}\leftindex^{s_{i_1}}{\left(\sum_{K\subseteq\set{2,\dots,k}}(-1)^{|K|}y^\pm_{K,v}\leftindex^v{{\phi_K(\ba)}}\right)}
         -\alpha^\pm_{i_1}\leftindex^{s_{i_1}}{\left(\sum_{K\subseteq\set{2,\dots,k}}(-1)^{|K|}y^\pm_{K,v}\leftindex^v{{\phi_K\phi_{i_1}(\ba)}}\right)}\\
       &=\sum_{K\subseteq\set{2,\dots,k}}(-1)^{|K|}y^\pm_{K,w}\leftindex^w{{\phi_K(\ba)}}
         +\sum_{K\subseteq\set{2,\dots,k}}(-1)^{|K\cup\set1|}y^\pm_{K\cup\set1,w}\leftindex^w{{\phi_{K\cup\set1}(\ba)}}.
   \end{align*}
   This proves the claim, since every $J\subseteq\set{1,\dots,k}$ is either of the form $J=K$ or $J=K\cup\set1$.
 \end{proof}

 \begin{Remark}\label{R:closed}  Formula \autoref{apastxis} appears more conceptual, but for our purposes less practical, when formulated in terms of the $\Theta_i=\alpha_i\phi_i$ from \autoref{R:Theta}. By \autoref{xicomm}, $\leftindex^{s_i}{(\bb\,\xipm_i)}=(\xipm_i-\Theta_i)(\bb)$ for all $\bb\in\Bo$, hence $\ba\,\xipm_{i_1}\leftindex^{s_{i_1}}{{\xipm_v}}
   =\leftindex^{s_{i_1}}{\left((\xipm_{i_1}-\Theta_{i_1})(\ba)\,\xipm_v\right)}$ for all $\ba\in\Bo$, and induction gives the closed formula
  \begin{equation}\label{apastxisTheta}
    \ba\,\xipm_w=\leftindex^w{\Bigl[(\xipm_{i_k}-\Theta_{i_k})\circ\dots\circ(\xipm_{i_1}-\Theta_{i_1})(\ba)\Bigr]}.
  \end{equation}
  Expanding \autoref{apastxisTheta}, gives $2^k$ terms and each is a composition of left multiplications and skew-derivations. \end{Remark}

Applying the commutation relation $H_i\bb=\leftindex^{s_i}\bb H_i$ from  \autoref{D:Hecke}, we obtain:

 \begin{Corollary}\label{C:axiH}
   Let $w=s_{i_1}\dots s_{i_k}$ be a reduced expression for $w\in\Sym$ and $\ba\in\Bo$. Then
   \[
     \ba\xi^\pm_wH_w = H_w\sum_{J\subseteq\set{1,2,\dots,k}}\left(
     (-1)^{|J|}\,\leftindex^{w^{-1}}{{\!(y^\pm_{J,w})}}\phi_J(\ba) \right)
     = \sum_{J\subseteq\set{1,2,\dots,k}}
    \left( (-1)^{|J|}\,y^\pm_{J,w}\,H_w\,\phi_J(\ba)\right).
   \]
   Moreover, each $\phi_J(\ba)$ is a $\Z$-linear combination of the elements
   $\leftindex^x\ba$ with $x\in\langle s_{i_j}\mid j\in J\rangle$. In particular, $x$ lies in every
   parabolic subgroup of $\Sym$ that contains $w$.
 \end{Corollary}

 \begin{proof}
   The two equalities arise by multiplying \autoref{L:axi} on the right by $H_w$ and then move $H_w$ to the left.
   The last statement is \autoref{L:phi}\ref{phi:translates}.
 \end{proof}

\subsection{Parabolic RaFH modules}

We are now ready to introduce the \emph{parabolic RaFH modules}
$\MJ^\pm_J$. As special cases, they include the \emph{spherical} and the
\emph{antispherical modules} $\MJ^+_S$ and  $\MJ^-_S$, respectively. In
the case of the symmetric groups, these modules correspond to
the induced trivial and sign representations from Young subgroups.

Throughout this subsection we assume that $\HHd$ is braided diagonalisable. Recall from
\autoref{C:xii} that $\xi^+_i\xi^-_i=-\hbar^2$ and
$\xi^+_i+\xi^-_i=\eta\taui[i,i+1]$, so that the eigenvalues of $H_i$ are
$\xi^+_i$ and $\xi^-_i$.

Recall that $S=\set{s_1,\dots,s_{n-1}}$ is the set of Coxeter generators
for $W$. For $J\subseteq S$, let $\Sym_J=\<s\mid s\in J\>$ be the corresponding
\emph{parabolic} (or Young) subgroup of $W$. If $J\subseteq K\subseteq
S$, let
\[
 \Sym_K^J = \set{d\in\Sym_K|\ell(dw)\ge\ell(d)\text{ for all }w\in\Sym_J}.
\]
As is well-known (see, for example,
\cite[Proposition~3.3]{Mathas:ULect}), $\Sym^J_K$ is a complete set of
(minimal length) left coset representatives of $\Sym_J$ in $\Sym_K$.
Similarly, let $\leftindex^J_K W=\set{d^{-1}|d\in\Sym^J_K}$ be the set
of (minimal length) right coset representatives of~$\Sym_J$ in~$\Sym_K$. For
convenience, set $\Sym^J=\Sym^J_S$ and $\leftindex^J W=\leftindex^J_S W$,
which are sets of left and right coset representatives of~$\Sym_J$ in~$W=\Sym_S$,
respectively.

\begin{Definition}\label{D:HJ}
  Given  $K\subseteq S$ we let
  $\HJ_K^\pm = \sum_{w\in\Sym_K}\hbar^{-2\ell(w)}\xi^\pm_w H_w$ and define
  $\MJ^\pm_K = \HHd\HJ_K^\pm$.
\end{Definition}
We call $\MJ^\pm_K$ the \emph{parabolic RaFH module}  associated with $K$.

\begin{Remark}
The $\HJ_K^\pm$ are RaFH analogues of the (anti)symmetrizer in the group algebra of $W$, or of the Kazhdan--Lusztig basis elements in its Hecke algebra, for the longest element in the parabolic subgroup $W_K$.
\end{Remark}

 \begin{Lemma}\label{L:HJK}
   Suppose that $J\subseteq K\subseteq S$. Then
   \[
     \HJpm_K = \sum_{d\in\Sym^J_K}\hbar^{-2\ell(d)}\xi^\pm_d H_d\HJpm_J
      = \sum_{d\in\leftindex^J_K W}\hbar^{-2\ell(d)}\HJpm_J H_d
        \leftindex^{d^{-1}}{{\xi^\pm_d}}.\]
 \end{Lemma}

 \begin{proof}
    Every $w\in\Sym_K$ factorises uniquely as $w=dv$ with $d\in\Sym^J_K$,
    $v\in\Sym_J$ and $\ell(w)=\ell(d)+\ell(v)$. Hence, by \autoref{C:hdv} using the relation $ \leftindex^y\bb H_y=H_y\bb$ for $y\in\Sym$, $\bb\in\Bo$, we get
    \[
      \HJ_K^\pm=\sum_{d\in\Sym^J_K}\sum_{v\in\Sym_J}
        \hbar^{-2(\ell(d)+\ell(v))}\xi^\pm_d\leftindex^d{{\xi^\pm_v}}H_dH_v
      =\sum_{d\in\Sym^J_K}\hbar^{-2\ell(d)}\xi^\pm_dH_d\,\HJpm_J .
    \]
    Factorising instead $w=vd$ with $v\in\Sym_J$, $d\in\leftindex^J_KW$ gives
    $\xi^\pm_{vd}H_{vd}=\xi^\pm_vH_v\,\xi^\pm_dH_d$ and
    $\xi^\pm_dH_d=H_d\;\leftindex^{d^{-1}}{{\xi^\pm_d}}$.
 \end{proof}

 The elements $\HJ^\pm_J$ are eigenvectors for the $H_w$ with $w\in\Sym_J$, with eigenvalues in $\Bo$:
 \begin{Theorem}\label{P:HW}
Let $J\subseteq S$ and  $w\in\Sym_J$. Then
   $H_w\HJ_J^\pm = \leftindex^w{{\xi^\pm_{w^{-1}}}}\HJ_J^\pm$.
 \end{Theorem}

 \begin{proof}
   By \autoref{C:xii} we have $\eta\taui[i,i+1]\xi^\pm_i=(\xi^\pm_i)^2-\hbar^2$,
   and $\leftindex^{s_i}{{\xi^\pm_i}}=\xi^\pm_i$ by \autoref{xxxinvformulas}, so that
   $H_i\xi^\pm_iH_i=\xi^\pm_iH_i^2$. Therefore, by the quadratic relation in
   \autoref{RaFHrels},
   \[
     H_i\HJpm_{\set{s_i}}=H_i+\hbar^{-2}\xi^\pm_iH_i^2
       =\left(1+\hbar^{-2}\eta\taui[i,i+1]\xi^\pm_i\right)H_i+\xi^\pm_i
       =\hbar^{-2}(\xi^\pm_i)^2H_i+\xi^\pm_i=\xi^\pm_i\HJpm_{\set{s_i}} .
   \]
   By the second formula of \autoref{L:HJK}, with $J$ replaced by $\set{s_i}$ and $K$
   by $J$, we have $\HJpm_J=\HJpm_{\set{s_i}}Y$ for some $Y\in\HHd$, so
   $H_i\HJpm_J=\xi^\pm_i\HJpm_J$ whenever $s_i\in J$. Now let
   $w=s_{i_1}\dots s_{i_k}$ be reduced and set $v=s_{i_1}\dots s_{i_{k-1}}$. Then
   $H_w\HJpm_J=H_v\left(\xi^\pm_{i_k}\HJpm_J\right)
     =\leftindex^v{{\xi^\pm_{i_k}}}H_v\HJpm_J$. So by induction,
   $H_w\HJpm_J=\leftindex^{v}{{\xi^\pm_{i_k}}}\cdots
     \leftindex^{s_{i_1}}{{\xi^\pm_{i_2}}}\,\xi^\pm_{i_1}\HJpm_J$, which is the claim
   by \autoref{C:hwrev}.
 \end{proof}

 \begin{Corollary}\label{C:TrivialContents}
 Assume $\HHd$ is braided diagonalisable and let $J\subseteq S$.
    Suppose that $X_1\HJpm_S=\conts\HJpm_S$. Then $X_k\HJpm_S=\conts^\pm_k\HJpm_S$, where
    $\conts^\pm_1=\conts$ and
    \begin{equation}\label{contents}
      \conts^\pm_{k+1}=\hbar^{-2}\left(\xi^\pm_k\leftindex^{s_k}{{\conts^\pm_k}}\xi^\pm_k
        + t\taui[k,k+1]\,\xi^\pm_k\right).
    \end{equation}
 \end{Corollary}

 \begin{proof}
   The case $k=1$ is by assumption. By \autoref{L:HXH} we have
   $\hbar^2X_{k+1}=H_kX_kH_k+tH_k\taui[k+1,k]$. Note that
   $\xi^\pm_k$ commutes with
   $X_k$ by $\eta$-symmetricity, since  $\xi^\pm_k=\pm\hbar\in R$ if $\eta=0$ and by  \autoref{Lfunnyequation} if $\eta\not=0$. Now \autoref{P:HW} gives  $H_kX_kH_k\HJpm_S=H_k\xi^\pm_kc^\pm_k\HJpm_S
       =\leftindex^{s_k}{\left(\xi^\pm_kc^\pm_k\right)}H_k\HJpm_S
       =\xi^\pm_k\leftindex^{s_k}{{c^\pm_k}}\xi^\pm_k\HJpm_S$,
       whereas  $tH_k\taui[k+1,k]\HJpm_S=t\taui[k,k+1]\xi^\pm_k\HJpm_S$. Taking the sum gives the eigenvalues as in \autoref{contents}.
 \end{proof}

 Combining \autoref{C:BasisHd} and \autoref{P:HW} we obtain a basis for
 the parabolic RaFH modules. These bases are a precursor to the
 respectable bases of the cyclotomic RaFH superalgebras constructed in
 \autoref{SS:LevelOneRespectable} and
 \autoref{SS:HigherLevelRespectable} below.

 \begin{Corollary}\label{C:MJBasis}
    Let $J\subseteq S$.  Then $\MJ^\pm_{J}$
      is free as an $R$-module with
     basis \[\set{\bb X^\ba  H_d\HJpm_J|\bb\in\cBo, \ba\in\N^n_{0}, \text{ and }d\in\Sym^J }.\]
 \end{Corollary}

 \begin{Remark}\label{R:LaiMinets}
   Parabolic RaFH modules  for commutative Frobenius algebras were already constructed in \cite{LaiMinets} in the framework of polynomial quantum wreath product algebras.
Motivated by geometric Euler classes, the authors constructed elements
analogous to our $\xi^\pm_w$. The $\xi$-braid relations are automatic in
this case since the Frobenius algebras considered in \cite{LaiMinets} are
commutative. It would be nice to have an interpretation of our elements in
terms of (non-commutative Borel-Moore) homology classes.  In
\cite{LaiMinets}, Schur algebras associated with polynomial quantum
wreath product algebras were introduced. The construction mimics the
classical approach, see for example \cite{Mathas:ULect}, but the Euler classes make the construction surprisingly difficult and technically involved.  See also \cite{Muthandco} for an approach to Schur algebras for wreath product algebras using web categories.
\end{Remark}

\section{Respectable algebras}
  To better understand the representation theory of CRaFH superalgebras, we want to show
  that they  have based structures that share features of
  cellular algebras~\cite{GrahamLehrer:Cellular} and finite based
  stratified algebras~\cite{BrundanStroppel:Semiinfinite}. To do
  this, we introduce \emph{respectable algebras} with their \emph{respectable bases},
  develop some basic theory and finally show that CRaFH algebras are special cases.

   The prototype of a respectable algebra is the matrix algebra from
   \autoref{prototype}. More generally, \autoref{Ex:Cellular} shows that the class of
   respectable algebras includes the \emph{cellular algebras} of Graham
   and Lehrer~\cite[Definition~1.1]{GrahamLehrer:Cellular}, which are
   one of our motivating examples. Respectable algebras generalise the \emph{standardly based} algebras of \cite{DuRui:strat}, which slightly weaken the conditions of cellularity to include algebras that do not have an anti-isomorphism. These are precisely the $\kk$-respectable algebras with $|\hlam|=1$ for all $\lambda\in\Lambda$, see \autoref{Ex:StandardlyBased}. Respectability also makes sense for superalgebras and the construction is similar to the generalised cellular algebras of Li and Shi~\cite{LiShi:QuiverHeckeClifford}.  Respectable algebras  also include, see \autoref{Ex:BasedQH},  finite dimensional  \emph{based quasi-hereditary} algebras of \cite{KlMuth, BrundanStroppel:Semiinfinite}.

  \subsection{The definition of respectable algebras}\label{SS:Respectable}
  Recall that $\kk$ is a commutative integral domain with $1$.

    \begin{Definition}\label{D:Respectable}
    Let $\CA$ be a $\kk$-algebra that is free and of finite rank as a
    $\kk$-module.

    A \emph{respectable datum} for $\CA$ is a tuple
    $(\Lambda,\ge,\plam,\hlam,\olam,\qlam,\Olam,a)$, where
    $(\Lambda,\ge)$ is a poset, and for $\lambda\in\Lambda$ there
    are \emph{indexing sets} $\plam,\hlam,\qlam$ and a $\kk$-algebra
    $\Olam$ with basis $\olam$, and $a$ is an injective map
    \[
      a\colon\coprod_{\lambda\in\Lambda}
      \plam\times\hlam\times\olam\times\qlam\to \CA; (p,h,o,q)\mapsto a_{phoq},
    \]
    such that the three conditions below hold. Here, for each $\lambda\in\Lambda$,  the  $\lambda$-\emph{stratum} $\CA^{(\lambda)}$ is the $\kk$-submodule of $\CA$ spanned by
    $\set{a_{phoq}|p\in\plam,h\in\hlam,o\in\olam, q\in\qlam}$
    and $\CA^{\ge\lambda}=\CA^{(\lambda)}+\CA^{>\lambda}$ with $\CA^{>\lambda}=\sum_{\mu>\lambda}\CA^{(\mu)}$.
  \begin{thmenumerate}[label=\upshape(R\arabic*), ref=\upshape(R\arabic*)]
  \item\label{resp1}  The set
      $\set{a_{phoq}|p\in\plam,h\in\hlam, o\in\olam, q\in\qlam\text{ for
      }\lambda\in\Lambda}$ forms a $\kk$-basis of~$\CA$.
    \item\label{resp2} For $\lambda\in\Lambda$, fix $p\in\plam$, $h\in\hlam$, $o\in\olam$,  $q\in\qlam$. Then for $x\in \CA$ there exist scalars
      $r_{pop'o'}(x)$ and $s_{oqo'q'}(x)$ in $\kk$,
      for $p'\in\plam$, $o'\in\olam$, $q'\in\qlam$, such that, modulo $\CA^{>\lambda}$,
      \[
        xa_{phoq}\equiv\sum_{\substack{p'\in\plam\\o'\in\olam}}r_{pop'o'}(x)a_{p'ho'q}
        \quad\text{and}\quad
        a_{phoq}x\equiv\sum_{\substack{q'\in\qlam\\o'\in\olam}}s_{oqo'q'}a_{pho'q'},
      \]
      where $r_{pop'o'}$ does not depend on $h$ and $q$, and
      $s_{oqo'q'}$ does not depend on $p$ and $h$.
      \item\label{resp3}  Write $oo'=\sum_{o''\in\olam}c_{oo'o''}o''$, for
      $c_{oo'o''}\in\kk$. If $\lambda\in\Lambda$, then
      $\CA^{(\lambda)}$ is an $\Olam$-bimodule via
      \[
          oa_{pho'q} = \sum_{o''\in\olam}c_{oo'o''}a_{pho''q}
             = a_{phoq}o'
      \]
      for all $p\in\plam$, $h\in\hlam$, $o,o'\in\olam$ and $q\in\qlam$. For
      $x\in \CA$ we have, modulo $\CA^{>\lambda}$, the compatibility
      \[
        (xa_{phoq})o'\equiv x(a_{phoq}o')
        \quad\text{and}\quad
        (o'a_{phoq})x\equiv o'(a_{phoq}x).
      \]

    \end{thmenumerate}
    The basis $\set{a_{phoq}}$ is the corresponding \emph{respectable basis} of~$\CA$.
    \end{Definition}
     To improve clarity, we sometimes write
  $a_{phoq}=a_{p,h,o,q}=a^\lambda_{p,h,o,q}$. We sometimes omit the $h$ or
  $o$ and write $a_{poq}$, $a_{phq}$ or $a_{pq}$ if
  $|\hlam|=1$, $|\olam|=1$ or $|\hlam|=|\olam|=1$ respectively.
           \begin{Definition}\label{perfectresp}
    A   \emph{respectable algebra} $\CA$ is an algebra $\CA$ with a respectable datum.
    The algebra $\CA$ is  \emph{perfectly respectable} if it is respectable and there are maps
    \[
      \plam\to \CA; p\mapsto\dot p,\quad
      \hlam\to \CA; h\mapsto\dot h,\quad
      \olam\to \CA; o\mapsto\dot o,\quad
      \qlam\to \CA; q\mapsto\dot q,
      \qquad\text{for }\lambda\in\Lambda,
    \]
    such that $a_{phoq}=\dot p\dot h \dot o \dot q$.
    If $\Olam=\kk$ for all $\lambda\in\Lambda$, we call the datum, and $\mathcal{A}$, (perfectly) \emph{$\kk$-respectable}. We also call an algebra $\kk$-respectable if it admits a $\kk$-respectable datum.
    For $\lambda\in\Lambda$, the tuple
    $(|\plam|,|\hlam|,|\olam|,|\qlam|)$ is the \emph{size} of the
    stratum $\CA^{(\lambda)}$, whose rank equals the product
    of these four numbers. We call $\Olam$ the \emph{fibre algebra} corresponding to $\lambda$, see \autoref{Ex:LerayHirsch}.
\end{Definition}

\begin{Remark}
   The condition of \autoref{resp2} is a two-sided analogue of one of
   the axioms of a cellular basis (see \autoref{Ex:Cellular}), whereas
   \autoref{resp3} is new. As we will see, respectable algebras are
   considerably more general than cellular algebras, even though they
   enjoy many of the same properties. Examples \ref{Ex:Taftresp},
   \ref{R:KoenigXi}, and~\ref{Ex:Kronecker}
   describe respectable algebras that are \textit{not} cellular.
 \end{Remark}

\begin{Remark}\label{R:associativity}
  For $\kk$-respectable algebras, the bases $\olam$ can be omitted and
  the condition~\autoref{resp3} is automatic. Expanding
  $(xa_{phoq})u\equiv x(a_{phoq}u)\pmod{\CA^{>\lambda}}$ and
  $(ua_{phoq})x\equiv u(a_{phoq}x)\pmod{\CA^{>\lambda}}$ shows
  that  the condition is in general equivalent to
  \[
    \sum_{o'}r_{pop'o'}(x)c_{o'uo''}=\sum_{o'}c_{ouo'}r_{po'p'o''}(x)
    \quad\text{and}\quad
    \sum_{o'}c_{uoo'}s_{o'qo''q'}(x)=\sum_{o'}s_{oqo'q'}(x)c_{uo'o''},
  \]
  for $x\in \CA$, $p,p'\in\plam$, $q,q'\in\qlam$ and
  $u,o,o''\in\olam$. In examples,  \autoref{resp3} usually follows from
  associativity of multiplication in~$\CA$ using the embeddings
  $\Olam\to \CA, o\mapsto a_{phoq}$, for chosen $p\in\plam$,
  $h\in\hlam$, $q\in\qlam$.
\end{Remark}

\begin{Definition}
    Instead of working with algebras we could instead work with
    \emph{algebra objects} in the symmetric monoidal category of
    $\Z$-graded or $\Z/2\Z$-graded $\kk$-modules and degree preserving
    maps. \autoref{D:Respectable} now defines \emph{respectable graded
    algebra} and \emph{respectable superalgebra}, respectively. The
    index sets $\plam,\hlam,\olam,\qlam$ then carry degree and parity
    functions, respectively, such that $|a_{phoq}|=|p|+|h|+|o|+|q|$. For
    perfectly respectable datum these functions give the degrees and
    parities of $\dot p,\dot h,\dot o,\dot q$. By abuse of language, we call such algebra objects just \emph{algebras}.
\end{Definition}

Before we start developing the theory of respectable algebras, we give a
series of examples to show their ubiquity.

\subsection{Examples of respectable algebras}
We provide now some explicit examples of respectable algebras, more precisely of algebra objects in vector (super)spaces, depending on the context.
\begin{Example}\label{verystupid}
The algebra $\kk$ is perfectly respectable with all indexing sets being $\set{1}$ and basis  $a_{1,1,1,1}=1$.
\end{Example}
\begin{Example}\label{sumrespectable}
If $\CA$ and $\CA'$ are respectable algebras, then so is $\CA\oplus \CA'$ with respectable datum given by the disjoint union of the respectable data.
\end{Example}
Any finite dimensional algebra is respectable and a fibre algebra in a trivial way:
\begin{Example}\label{stupid} Every finite dimensional superalgebra $B$ is perfectly
respectable with $\Lambda=\{\lambda\}$ by setting  $\Olam=B$ with a fixed basis $\olam$,
$\plam=\hlam=\qlam=\set{1}$ and $a_{1,1,o,1}=1\cdot 1\cdot o\cdot 1=o$.
Useful respectable bases usually have larger indexing sets, ideally with
each fibre algebra, $\Olam$, being a division algebra. As we will
show, if $\Olam$ is a division algebra, then it is the endomorphism ring
of the corresponding simple $B$-module.
\end{Example}

The following in particular applies if $\CO$ is a division algebra.

  \begin{Example}\label{prototype}
    Let $\CO$ be a  $\kk$-algebra with $\kk$-basis
    $\set{o_i|i\in I}$, and let $\CA=\Mat_m(\CO)$ be the algebra of
    $m\times m$-matrices over $\CO$, for $m\ge1$. Let $E_{pq}$ be the
    elementary matrix with $1$ in row~$p$ and column~$q$, and zeros
    elsewhere. It is easy to see that $\CA$ is a perfectly respectable algebra
    with respectable datum
    \begin{equation*}
    \begin{gathered}
    \Lambda=\set{\bullet},\quad \plam[\bullet]=\set{E_{p1}|1\le p\le m},\quad
    \hlam[\bullet]=\set{E_{11}}, \quad \olam[\bullet]=\set{o_i|i\in I},
    \quad\qlam[\bullet]=\set{E_{1q}|1\le q\le m}, \\
     \text{$a_{phoq}=phoq=\dot p\dot h\dot o\dot q,\quad$ with $\quad\dot{p}=p$, $\dot{h}=h$, $\dot{o}=o$, $\dot{q}=q$}.
     \end{gathered}
     \end{equation*}
     In particular, $\Mat_m(\kk)$ is a $\kk$-respectable and $\Mat_m(D)$ is respectable for every division algebra $D$ over~$\kk$, and then $\Delta^\bullet=D^m$ is the unique up to isomorphism simple module, see \autoref{T:Respectable}.

     This respectable basis has a single stratum of size $(m,1,|\olam[]|,m)$.

    Together with \autoref{sumrespectable}, this shows that if $G$ is a
    finite group, then the group algebras $\Q G$, $\R G$ and $\C G$ are
    all respectable. Respectable bases are far from unique. If, for
    example, $M\in\Mat_{m}(\CO)$ is an invertible matrix and $\CO$ is a
    local ring, then $e=ME_{11}M^{-1}$ is a primitive idempotent
    in~$\CA$. It is easy to see that, for any $M$,  there exists
    a respectable basis of $\CA$ with $\hlam[\bullet]=\set{e}$.
  \end{Example}

   \begin{Example}\label{excomm}
   Let $\CA$ be commutative over an algebraically closed field $\kk$. Its semisimplification $\CA/\op{rad}\CA$ is also commutative,  and the
   simple $\CA$-modules are $1$-dimensional. Consider  a composition series
   $\{0\}\subset \CA_1\subset
   \CA_2\subset\ldots\subset \CA_m=\CA$ of~$\CA$, so that $\dim\CA_i/\CA_{i-1}=1$. This is a
   filtration of $\CA$ by two-sided ideals, which implies that $\CA$ is a
   cellular algebra with the identity involution \cite{KX:CellAlg}. Let
   $\set{a_1,a_2,\ldots, a_i}$ be a basis of $\CA_i$ and
   $\Lambda=\{1>2>\cdots>m\}$. Setting $\plam=\qlam=\set{1_{\CA}}$,
   $\olam=\set{1_\kk}$, $\Olam=\kk$, $\hlam=\set{a_\lambda}$, for
   $\lambda\in\Lambda$, defines a $\kk$-respectable basis of~$\CA$, where
   $a_{phoq}=phoq$, because $\CA^{>\lambda}=\CA_{\lambda-1}$ is an ideal of~$\CA$.
 Alternatively we could also take
   $\hlam=\qlam=\set{1}$, $\plam=\set{a_\lambda}$  or
   $\plam=\hlam=\set{1}$, $\qlam=\set{a_\lambda}$, for each
   $\lambda$. This gives a total of $3^m$ respectable bases.
  \end{Example}

  \begin{Example}
    Let $\CA=\Q \mathtt{C}_4$ be the group algebra of the cyclic group  $\mathtt{C}_4=\<\,x\mid x^4=1\,\>$ of order~$4$. Let $\Lambda=\set{1>2>3}$ and set
    $\plam=\qlam=\set{1}$, for $\lambda\in\Lambda$, and
    \[
      \Olam[1]=\Olam[3]=\Q, \quad
      \olam[1]=\olam[3]=\set{1},\quad
      \Olam[2]=\Q(\iota), \quad\text{and}\quad
      \olam[2]=\set{1,\iota},
    \]
    where $\iota=\sqrt{-1}\in\C$. Finally,
    $\hlam[1]=\set{\sum_{k=0}^3x^k}$,
    $\hlam[2]=\set{1}$, and
    $\hlam[3]=\set{\sum_{k=0}^3(-x)^k}$, and
    \[
      a^\lambda_{poq} = \begin{cases*}
          phoq& if $\lambda\in\set{1,3}$,\\
          \tfrac{1-x^2}2 & if $\lambda=2$ and $o=1$,\\
          \tfrac{x-x^3}2 & if $\lambda=2$ and $o=\iota$.
      \end{cases*}
    \]

    Notice that $x^2\cdot\tfrac12(1-x^2)=-\tfrac12(1-x^2)$, so the
    action of~$\iota$ on $\CA^{(2)}$ coincides with multiplication by~$x$,
    showing that \autoref{resp3} is satisfied. This perfectly
    respectable basis of~$\Q \mathtt{C}_4$ realises the Artin--Wedderburn
    decomposition $\Q\mathtt{C}_4\cong\Q\oplus\Mat_1(\Q(\iota))\oplus\Q$ of this
    algebra.
  \end{Example}

\begin{Example}\label{Ex:Qn}
    Consider the (queer) superalgebra $Q(m)$ of $2m\times 2m$-matrices
    $\set[\big]{\begin{psmallmatrix} A&B\\B&A \end{psmallmatrix}|A,B\in\Mat_m(\kk)}$.
    This is a superalgebra with even elements
    $\set[\big]{\begin{psmallmatrix} A&0\\0&A \end{psmallmatrix}}$ and
    odd elements
    $\set[\big]{\begin{psmallmatrix} 0&B\\B&0 \end{psmallmatrix}}$. As $Q(m)\cong\Mat_m(\Cl_1)$, see \autoref{prototype}, the
    superalgebra $Q(m)$ is respectable with
    \[\Lambda=\set{\bullet},\quad \plam[\bullet]=\set{1,\dots,m}=\qlam[\bullet],\quad
    \hlam[\bullet]=\set{1},\quad \olam[\bullet]=\set{1,C},\quad \Olam[\bullet]=\Cl_1,\]
    where $\Cl_1$ is the Clifford superalgebra from \autoref{Ex:Clifford}, and, with $\overline{i} = i+m$,
    \[
         a_{poq} = \begin{cases*}
             E_{pq}+E_{\overline{p}\overline{q}} & if $o=1$,\\
             E_{p\overline{q}}+E_{\overline{p}q} & if $o=C$.
         \end{cases*}
    \]
  In fact, $Q(m)$ is a perfectly respectable superalgebra, with $\dot p=E_{p1}+E_{\overline{p}\overline{1}}$, $\dot h=1$, $\dot 1=E_{11}+E_{\overline{1}\overline{1}}$, $\dot C=E_{1\overline{1}}+E_{\overline{1}1}$ and $\dot q=E_{1q}+E_{\overline{1}\overline{q}}$, so that $a_{phoq}=\dot p\dot h\dot o\dot q$.
\end{Example}

\begin{Example}\label{Ex:CliffordRespectable}
  Let $\kk$ be an algebraically closed field with $2\in\kk^{\times}$. If $n$ is odd, then $\Cl_n\cong\Mat_m(\Cl_1)\cong Q(m)$ as in \autoref{Ex:Qn}, where $m=2^{\frac{n-1}{2}}$ and
  $\Cl_1$ is generated by the odd element $C=\zeta C_1\cdots C_n$ with $C^2=1$
  for a suitable $\zeta\in\kk^\times$. By \autoref{prototype}, $\Cl_n$ is a perfectly
  respectable superalgebra with a single stratum and $\Olam[\bullet]=\Cl_1$.
  This cannot be improved to a $\kk$-respectable datum with homogeneous basis.
  Any such datum has a single stratum, because $\Cl_n$ has no
  proper nonzero homogeneous ideals. The standard module $\Delta^\bullet$ is then a nonzero $\Cl_n$-supermodule, so it has dimension at least $2m$, as the unique simple $\Cl_n$-supermodule has superdimension $m|m$. Hence $|\plam[\bullet]|\ge2m$, similarly $|\qlam[\bullet]|\ge2m$, and so
  $\dim\Cl_n\ge4m^2$, which is a contradiction because $\dim\Cl_n=2m^2$. For $n$ even,
  $\Cl_n$ is the matrix superalgebra $M(2^{\frac n2-1}|2^{\frac
  n2-1})$ and, hence, $\kk$-respectable.
\end{Example}

  \begin{Example}\label{Ex:Cellular}
    Suppose that $\CA$ is a cellular algebra, in the sense of Graham and
    Lehrer~\cite{GrahamLehrer:Cellular}, with weight poset
    $(\Lambda,\ge)$ and cellular basis $\set{c_{st}|s,t\in T(\lambda),
    \lambda\in\Lambda}$. Then $\CA$ is a $\kk$-respectable algebra with respectable datum
    \[(\Lambda,\ge),\quad \plam=\qlam=T(\lambda), \quad \hlam=\set{1_\kk},\quad \olam=\set{1_\kk}, \quad \Olam=\kk, \quad\text{and} \quad a_{st}=c_{st},  \]
   where  $\lambda\in\Lambda$,  $s\in \plam, t\in \qlam$.  A cellular algebra is respectable, but not always perfectly respectable.
   Here we use the convention of \cite{Mathas:ULect}, where the cellular axiom is a congruence modulo the span of the $c_{st}$ with larger index. This is the opposite order to the one used in  \cite{GrahamLehrer:Cellular}.
     \end{Example}

  \begin{Example}
    As a partial converse to \autoref{Ex:Cellular}, let $\CA$ be a
    $\kk$-respectable algebra with weight poset $(\Lambda,\ge)$ with
    $\plam=\qlam$ and $|\hlam|=1$, for $\lambda\in\Lambda$. Since $\CA$ is
    $\kk$-respectable, $|\olam|=1$ and $\Olam=\kk$. Suppose that the
    unique linear map $*\colon \CA\to \CA$ with $a_{pq}^*=a_{qp}$ is
    an anti-isomorphism. Then it is easy to see that $\set{a_{pq}}$ is
    a cellular basis of~$\CA$.
     \end{Example}

\begin{Example}\label{Ex:StandardlyBased}
  A standardly based algebra with defining base $\set{a^\lambda_{ij}}$
  as in \cite[Definition~1.2.1]{DuRui:strat}, is the same as a
  $\kk$-respectable algebra with $\plam=I(\lambda)$, $\hlam=\set{1}$,
  $\qlam=J(\lambda)$, $a_{p1q}=a^\lambda_{pq}$ with \autoref{resp2}
  being \cite[(1.2.2)]{DuRui:strat}. A finite field extension
  $\mathbb{F}\supsetneq\kk$ is respectable by \autoref{stupid}, but not
  standardly based, because its simple module is not absolutely
  irreducible, see \cite[Theorem~2.4.1(c)]{DuRui:strat} and
  \autoref{T:Respectable}\autoref{late2} below.
\end{Example}

\begin{Example}\label{Ex:TL}Many diagram algebras are perfectly respectable. To give just one example, consider the Temperley--Lieb algebra $TL_m(\delta)$ for $\delta\in\kk$, $m\in\N$, see e.g. \cite{TLsurvey}. This algebra has a basis given by non-crossing matchings of $2m$-points, i.e. planar graphs with $2m$ vertices consisting of $m$ marked points on each of two parallel lines. Each point in the graph is the endpoint
of precisely one edge, and the edges can be drawn by non-intersecting arcs
which lie entirely between the lines. Composition of basis vectors is
given by concatenation and replacing each closed component by a factor
$\delta$. Every such diagram can be decomposed uniquely into its top (a
cap diagram), $p$, middle (through strings), $h$ and bottom (the cup
diagram), $q$, as in the following picture:
\tikzset{
  TL/.style    = {scale=0.4, line width=0.6pt, line cap=round,
                  baseline=(current bounding box.center)},
  frame/.style = {draw=black!55, fill=black!3, rounded corners=1pt},
  sep/.style   = {draw=black!45, dashed, line width=0.4pt},
  pt/.style    = {fill, circle, inner sep=0pt, minimum size=1.8pt},
  caparc/.style= {out=-90, in=-90},   %
  cuparc/.style= {out=90,  in=90},    %
  strand/.style= {out=-90, in=90},    %
}
\[
  \begin{tikzpicture}[TL]
    \draw[frame] (0.4,0) rectangle (6.6,1.8);
    \foreach \i in {1,...,6}{\node[pt] at (\i,1.8){}; \node[pt] at (\i,0){};}
    \draw (1,1.8) to[caparc] (2,1.8)   (5,1.8) to[caparc] (6,1.8);
    \draw (3,1.8) to[strand] (5,0)     (4,1.8) to[strand] (6,0);
    \draw (1,0)   to[cuparc] (4,0)     (2,0)   to[cuparc] (3,0);
  \end{tikzpicture}
  \qquad=\qquad
  \begin{tikzpicture}[TL]
    \draw[frame] (0.4,0) rectangle (6.6,5.4);
    \draw[sep] (0.4,3.6) -- (6.6,3.6)   (0.4,1.8) -- (6.6,1.8);
    \foreach \i in {1,...,6}{\node[pt] at (\i,5.4){}; \node[pt] at (\i,0){};}
    \foreach \i in {3,4}{\node[pt] at (\i,3.6){};}
    \foreach \i in {3,4}{\node[pt] at (\i,1.8){};}
    \draw (1,5.4) to[caparc] (2,5.4)   (5,5.4) to[caparc] (6,5.4);  %
    \draw (3,5.4) -- (3,3.6)        (4,5.4) -- (4,3.6);
    \draw (3,3.6) -- (3,1.8)   (4,3.6) -- (4,1.8);  %
    \draw (1,0)   to[cuparc] (4,0)     (2,0)   to[cuparc] (3,0);    %
    \draw (3,1.8) to[strand] (5,0)   (4,1.8) to[strand] (6,0);
    \node[font=\scriptsize, anchor=west] at (6.9,5.0){$p=$ cap diagram};
    \node[font=\scriptsize, anchor=west] at (6.9,2.7){$h=$ through strings };
    \node[font=\scriptsize, anchor=west] at (6.9,0.4){$q=$ cup diagram};
  \end{tikzpicture}
\]
Note that $h$ is the identity diagram $1_\lambda$ on $m-2\lambda$
strands, where $0\le\lambda\le\lfloor m/2\rfloor$ is the number of caps.
To construct a respectable datum take $\Lambda=\set{0<1<\ldots<\lfloor
m/2\rfloor}$, so that diagrams with fewer through strands are in higher
strata. For $\lambda\in\Lambda$ set $\olam=\set{1}$,
$\hlam=\set{1_\lambda}$, and
let $\plam$ and $\qlam$ be the sets of cap diagrams $p$, and  cup
diagrams $q$, respectively obtained as above when $h=1_\lambda$.
Defining $a_{phq}$ to be the concatenated
diagram, this gives a respectable datum. To construct lifts, for each
$\lambda$, fix a cup diagram $q_0\in\qlam$ and let $p_0\in\plam$ be its
reflection. Then $a_{p1_\lambda q_0}a_{p_01_\lambda q}=\delta^\lambda
a_{p1_\lambda q}$.  Hence, if $\delta\in\kk^\times$, then the datum is
perfectly respectable, with $\dot p=a_{p1_\lambda q_0}$, $\dot h=\delta^{-\lambda}1$ and $\dot q=a_{p_01_\lambda q}$.
\end{Example}

\begin{Example}\label{Ex:BasedQH}
  Let $\CA$ be a \emph{finite based quasi-hereditary algebra} in the sense
  of \cite[Definition~5.1]{BrundanStroppel:Semiinfinite}, with its finite set $I$ of
  orthogonal idempotents and poset $\Lambda$ of special idempotents.
    Then $\CA$ is perfectly $\kk$-respectable. Using the
    notation from \cite[\S5]{BrundanStroppel:Semiinfinite}, a  respectable datum is given by
  \[
    (\Lambda,\ge),\quad \plam=Y(\lambda),\quad \hlam=\set{e_\lambda},\quad
    \olam=\set{1_\kk},\quad \Olam=\kk,\quad \qlam=X(\lambda),\quad
    a^{\lambda}_{yx}=y\,e_\lambda\,x=yx.
  \]
   It is perfectly respectable with $\dot p=y$, $\dot h=e_\lambda$ and
   $\dot q=x$, as $y e_\lambda x=yx$. The axioms
   \cite[(QH5),(QH6)]{BrundanStroppel:Semiinfinite} give
   \autoref{resp2}. Based quasi-heredity is strictly stronger than being
   perfectly $\kk$-respectable. For instance, $\kk[x]/(x^2)$ is perfectly $\kk$-respectable, as in \autoref{Ex:gradedlocal}, but not quasi-hereditary, see also \autoref{C:quasihereditary}.
 \end{Example}

The next example shows that there are respectable algebras with
$|\qlam[i]|\ne|\plam[i]|$, for most~$i\in\Lambda$.

   \begin{Example}\label{Ex:Taftresp}
    Fix an integer $\theta\ge2$ and assume $\kk$
    contains a primitive $\theta$th root of unity $\omega\in\kk$ and consider the
    Taft algebra
    $\Taft=\<\,g,y\>/( g^\theta=1, y^\theta=0, yg=\omega gy)$, see  \autoref{Ex:Taft}.  Now
    let $\Lambda=\set{0,1,\dots,2\theta-1}$, with its natural ordering.
    For $i\in\Lambda$ set, $\olam[i]=\set{1_\kk}$ and $\Olam[i]=\kk$.
     \begin{equation}\label{E:Taft}
     \plam[i]=\set{1}, \quad\quad\hlam[i]=\set{h_i},
    \text{ with }h_i=\prod_{j=i+1}^{i+\theta-1}(g-\omega^{-j}),\quad\quad
        \qlam[i] = \begin{cases*}
            \set{1}                    & if $i=0$,\\
            \set{1,\dots,y^{\theta-2}} & if $0<i< \theta$,\\
            \set{y,\dots,y^{\theta-1}} & if $i=\theta$,\\
            \set{y^{\theta-1}}         & if $\theta<i<2\theta$.
        \end{cases*}
     \end{equation}
     Then $\set{(\plam[i],\hlam[i],\olam[i], \qlam[i])|i\in\Lambda}, a$ is a respectable datum giving a perfectly $\kk$-respectable basis of~$\Taft$, which follows
     by observing that $gh_i=\omega^{-i}h_i$ and $yh_i=\omega^{-1}
     h_{i+1}y$, for $0\le i<2\theta$ (indeed, $h_i=h_{i+\theta}$ for
     $0\le i<\theta$). Hence, $\Taft$ is a respectable algebra.  We show
     below that the Taft algebra $\Taft$ is never cellular. In the
     special case when $\theta=2$ and $\op{char}(\kk)\not=2$ this is
     easy to see because the left ideals $\Taft[2](1\pm g)$ are
     indecomposable projective $\Taft[2]$-modules with simple heads and
     non-isomorphic simple socles. However, by \cite[Theorem
     1.1]{KoenigXi:selfinj},  this cannot happen for Frobenius algebras
     that are cellular.
   \end{Example}

The definition of respectable algebras is sufficiently general to
include many examples from the literature, yet strong enough to support
the theory we develop in the next section. In \autoref{D:Standard}, we
introduce standard and costandard modules for respectable algebras, and
develop some basic theory for them that is reminiscent of highest weight
theory. In contrast to the standard and costandard modules of
quasihereditary algebras, our standard and costandard modules do not
necessarily have simple heads or socles, respectively; see
\autoref{T:Respectable}. However, like cellular algebras, the simple
modules of respectable algebras are quotients of the standard modules
indexed by a subset $\Lambda_+\subset\Lambda$. As such, the standard and
costandard modules of respectable algebras behave like the Specht
modules of the symmetric groups. This connection with
the Specht modules is the origin of the term re\emph{spec(h)t}able, and
we hope it earns them some respect.

\subsection{Respectable modules}
This section defines the standard and costandard modules for respectable
algebras, and uses them to realise the irreducible modules of these
algebras as quotients of the standard modules, and submodules of the
costandard modules, when the fibre algebras are division algebras. These
modules generalise the cell modules of a cellular algebra; see
\autoref{Ex:StandardBasedQH}.

  Let $\set{(\plam,\hlam,\olam,\qlam)|\lambda\in\Lambda}$ be the
  respectable data for a $\kk$-superalgebra~$\CA$. Fix $\lambda\in\Lambda$.
Recall from \autoref{D:Respectable}, that $\CA^{>\lambda}$ $\kk$-submodule of $\CA$ spanned by the $a^\mu_{phoq}$ with $\mu>\lambda$.

We start with a useful lemma. To state this, for $\lambda\in\Lambda$ let
$1_\lambda$ be the unit of the fibre $\kk$-algebra $\Olam$. If
$c=\sum_oc_oo\in\Olam$, set $a_{phcq}=\sum_oc_oa_{phoq}=ca_{ph1_\lambda
q}=a_{ph1_\lambda q}c$, where the last two equalities follow by
\autoref{resp3}. By definition, $a_{phcq}, a_{ph1_\lambda q}\in
\CA^{(\lambda)}$ for $c\in\Olam$, $p\in\plam$, $h\in\hlam$, and
$q\in\qlam$.

 \begin{Lemma}\label{L:hlarge}
   Let $\CA$ be a respectable algebra and suppose that
   $\lambda\in\Lambda$. Then $\CA^{\ge\lambda}$ is a two-sided ideal
   of~$\CA$. Moreover, if $|\hlam|>1$, then
   $\CA^{\ge\lambda}\CA^{\ge\lambda}\subseteq\CA^{>\lambda}$.
   Finally, if $|\hlam|=1$ and
   $p,p'\in\plam$, $h\in\hlam$, $o,o'\in\olam$ and~$q,q'\in\qlam$,
   for $\lambda\in\Lambda$, then there exists $c\in\Olam$ such that
   \[
        a_{phoq}a_{p'ho'q'}\equiv a_{phcq'}\pmod{\CA^{>\lambda}}.
   \]
   where $c$ does not depend on $p$, $h$ or $\mathbf{q'}$.
\end{Lemma}

\begin{proof}
 It is immediate from \autoref{resp2} that $\CA^{\ge\lambda}$ is a
 two-sided ideal. To prove the remaining statements suppose that
 $p,p'\in\plam$, $h,h'\in\hlam$, $o,o'\in\olam$ and~$q,q'\in\qlam$.
 Applying \autoref{resp2} twice, there exist scalars
 $r_{pop''o''}(a_{pohq}),s_{o'q'o''q''}(a_{p'o'h'q'})\in\kk$ such that
 \[
       \sum_{\substack{p''\in\plam\\o''\in\olam}}
          r_{pop''o''}(a_{pohq})a_{p''h'o''q'}
        \equiv a_{pohq}a_{p'h'o'q'}
        \equiv \sum_{\substack{q''\in\qlam\\o''\in\olam}}
          s_{o'q'o''q''}(a_{p'o'h'q'})a_{pho''q''}
        \pmod{\CA^{>\lambda}},
 \]
 where $r_{pop''o''}(a_{pohq})$ does not depend on~$h'$ and~$q'$, and
 $s_{o'q'o''q''}(a_{p'h'o'q'})$ does not depend on~$p$ and~$h$. Comparing
 the coefficients on both sides,
 $r_{pop''o''}(a_{pohq})=\delta_{pp''}\delta_{hh'}\delta_{q'q''}
    s_{o'q'o''q''}(a_{p'h'o'q'})$.
 In particular $a_{phoq}a_{p'h'o'q'}\in\CA^{>\lambda}$ if $h\ne h'$.
 Since these coefficients do not depend on~$h$ and~$h'$, if $|\hlam|>1$
 we can compute them by assuming that $h\ne h'$, which forces them to be
 zero by what we have just shown. Hence,
 $\CA^{\ge\lambda}\CA^{\ge\lambda}\subseteq\CA^{>\lambda}$ if
 $|\hlam|>1$. Finally, suppose that $|\hlam|=1$, so that $h=h'$. We have
 shown that a basis element in $\CA^{(\lambda)}$ appears with nonzero
 coefficient in $a_{pohq}a_{p'h'o'q'}$ only if it is of the form
 $a_{pho''q'}$, for some $o''\in\Olam$, in which case it appears with
 coefficient $r_{popo''}(a_{pohq})=s_{o'q'o''q'}(a_{p'ho'q'})\in\kk$,
 which depends only on $p'$, $o'$ ,$o$ and $q$. Hence,
 $a_{phoq}a_{p'ho'q'}\equiv a_{phcq'}\pmod{\CA^{>\lambda}}$, where
 \[
   c = \sum_{o''\in\olam}r_{popo''}(a_{pohq})o''
     = \sum_{o''\in\olam}s_{o'q'o''q'}(a_{p'o'hq'})o''\in\Olam
 \]
 depends only on $p'$, $o'$,$o$ and $q$. This completes the proof.
 \end{proof}

  \begin{Definition}\label{D:Standard}
   Let $\CA$ be a respectable (super)algebra.
   \begin{enumerate}[series=Defstd, label=(\alph*), ref=\autoref{D:Standard}.(\alph*)]
   \item\label{std1} For $\lambda\in\Lambda$, the \emph{standard module}
   $\Delta^\lambda$ is the $(\CA,\Olam)$-bimodule with basis
   $\set{a_{po}|p\in\plam, o\in\olam}$, of parity $|a_{po}|=|p|+|o|$ in the super case, and with $\CA$-action and
   $\Olam$-actions
   \[
        x\cdot a_{po}=\sum_{\substack{p'\in\plam\\o'\in\olam}}r_{pop'o'}(x) a_{p'o'}
        \qquad\text{and}\qquad
        a_{po}\cdot o''=\sum_{o'\in\olam}c_{oo''o'}a_{po'}
   \]
   for $x\in \CA$ and $p\in\plam$ and $o,o''\in\olam$,
   where the scalars $r_{pop'o'}(x),c_{oo''o'}\in\kk$ come from \autoref{D:Respectable}.
   This action of $\CA$ on $\Delta^\lambda$ is well-defined because the
   scalars $r_{pop'o'}(x)$ from \autoref{resp2} depend only on
   $x\in \CA$, $p, p'\in\plam$ and $o,o'\in\olam$, and not on $h\in\hlam$
   or $q\in\qlam$.

   \item \label{std2} Similarly, the \emph{right standard module} $\leftindex{^\lambda}\Delta$ is
   the $(\Olam,\CA)$-bimodule with basis
   $\set{a_{oq}|q\in\qlam, o\in\olam}$, of parity $|a_{oq}|=|o|+|q|$, and actions
   \[
         a_{oq}\cdot x=\sum_{o'\in\olam,\ q'\in\qlam}s_{oqo'q'}(x) a_{o'q'},
         \qquad\text{and}\qquad
         o''\cdot a_{oq}=\sum_{o'\in\olam}c_{o''oo'}a_{o'q}
   \]
   for $x\in \CA$, $q\in\qlam$ and $o\in\olam$.
   (These are actions of $\CA$ by applying
   \autoref{resp2} twice.)

   \item\label{std3} The \emph{costandard module} is the $(\CA,\Olam)$-bimodule
   $\nabla^\lambda=\Hom_\kk(\leftindex{^\lambda}\Delta,\kk)\cong\Hom_{\Olam}(\leftindex^{\lambda}\Delta,\Olam^*)$ with actions
   $(xfo)(y)=f(oyx)$, for $f\in\nabla^\lambda$, $x\in \CA$, $o\in\Olam$ and $y\in\leftindex{^\lambda}\Delta$.
   \item\label{std4} For a left $\Olam$-module $M$, the \emph{standard} and \emph{costandard modules
   with coefficients in $M$} are the left $\CA$-modules
   \[
     \Delta^\lambda(M)=\Delta^\lambda\otimes_{\Olam}M
     \qquad\text{and}\qquad
     \nabla^\lambda(M)=\Hom_{\Olam}(\leftindex{^\lambda}\Delta,M).
   \]
   If $M=L$ is a simple $\Olam$-module, we call them \emph{proper standard} and \emph{proper
   costandard modules}.

   \item Let $\CA$ be an ordinary algebra. Then an $\CA$-module $M$ has a \emph{$\Delta$-filtration} if there exists a
   filtration $M=M_1\supset M_2\supset\dots\supset M_z\supset M_{z+1}=0$, and
   $\lambda_1,\dots,\lambda_z\in\Lambda$, such that
   $M_k/M_{k+1}\cong\Delta^{\lambda_k}$ as $\CA$-modules, for $1\le k\le z$.
   Similarly, $M$ has a $\nabla$-filtration if it has a filtration with
   subquotients $M_k/M_{k+1}\cong\nabla^{\lambda_k}$ as $\CA$-modules,
   for $1\le k\le z$. For superalgebras the subquotients are only required to be isomorphic to
   $\Pi^{\epsilon_k}\Delta^{\lambda_k}$, respectively $\Pi^{\epsilon_k}\nabla^{\lambda_k}$, for some
   $\epsilon_k\in\set{0,1}$. The \emph{multiplicity} of $\Delta^\lambda$
   in a $\Delta$-filtration is the number of
   subquotients isomorphic to $\Delta^\lambda$ or $\Pi\Delta^\lambda$.
   \end{enumerate}
   \end{Definition}

  \begin{Remark}  To see that  \autoref{std1} and
  \autoref{std2} define bimodules apply \autoref{resp2} twice, and then \autoref{resp3}.
  \end{Remark}

  By definition, the standard and costandard modules $\Delta^\lambda$
  and $\nabla^\lambda$ are both $(\CA,\Olam)$-bimodules but, by
  restriction, we usually think of them as left $\CA$-modules.

  \begin{Lemma}\label{L:DeltaFiltration}
    Let $\CA$ be a respectable superalgebra. Then the left regular module $\CA$ has a
    $\Delta$-filtration in which $\Delta^\lambda$ appears with multiplicity $|\hlam||\qlam|$, for
    $\lambda\in\Lambda$.
  \end{Lemma}

  \begin{proof}
    By \autoref{resp2}, $x\CA^{(\mu)}\subseteq\CA^{(\mu)}+\CA^{>\mu}$ for
    $x\in\CA$ and $\mu\in\Lambda$, so $\CA^{\ge\lambda}$ and $\CA^{>\lambda}$
    are left ideals. Fix $\lambda\in\Lambda$, $h\in\hlam$ and $q\in\qlam$. By \autoref{resp1} and
    \autoref{resp2}, the map
    $
      \Delta^\lambda\longrightarrow\CA^{\ge\lambda}/\CA^{>\lambda},
      a_{po}\longmapsto a_{phoq}+\CA^{>\lambda},
    $
    is injective and $\CA$-linear, because the scalars $r_{pop'o'}(x)$ do not depend on $h$
    and $q$. In the super case it is homogeneous of some parity $\epsilon_{hq}=|h|+|q|$, so its image is
    isomorphic to $\Pi^{\epsilon_{hq}}\Delta^\lambda$. By \autoref{resp1}, the quotient
    $\CA^{\ge\lambda}/\CA^{>\lambda}$ is the direct sum of these images, so it is
    isomorphic to $\bigoplus_{h\in\hlam,\ q\in\qlam}\Pi^{\epsilon_{hq}}\Delta^\lambda$.

    Let $\lambda_1,\dots,\lambda_m$ be the elements $\lambda\in\Lambda$ with
    $\CA^{(\lambda)}\ne0$, numbered so that
    $\lambda_i>\lambda_j$ only if $i>j$. Then
    $M_k=\sum_{i\ge k}\CA^{(\lambda_i)}$ contains $\CA^{>\lambda_k}$, so it is a left
    ideal, and $M_k/M_{k+1}\cong\CA^{\ge\lambda_k}/\CA^{>\lambda_k}$ as
    $\CA$-modules. Refining the filtration $\CA=M_1\supset\dots\supset
    M_{m+1}=0$ along the direct sum decompositions gives the claim.
  \end{proof}

  \begin{Example}\label{Ex:StandardBasedQH}
    For cellular algebras, viewed as respectable algebras as in
    \autoref{Ex:Cellular}, the standard modules are the cell modules
    from \cite{GrahamLehrer:Cellular}. Similarly, we recover the
    standard modules of Brundan and the second author from
    \cite{BrundanStroppel:Semiinfinite} in the situation of
    \autoref{Ex:BasedQH}.  The standard module of \autoref{D:Standard}
    has basis indexed by $\plam=Y(\lambda)$, matching the standard basis
    of the module denoted $\Delta(\lambda)=\CA^{\le\lambda}\bar
    e_\lambda$ in \cite[(5.2)]{BrundanStroppel:Semiinfinite}. Taking
    $q=e_\lambda$ and $p'=e_\lambda$ in \autoref{L:Respectable} gives
    $\<a_{e_\lambda},a_{e_\lambda}\>=1$, so $L^\lambda\ne0$
    for every $\lambda\in\Lambda$ and $\Lambda_+=\Lambda$.
    By \autoref{T:Respectable}, below, this implies that $\CA$ is
    quasi-hereditary.
  \end{Example}

\begin{Example}\label{Ex:ProperlyStratified}
  Let $\CA$ be a finite based properly stratified algebra in the sense of
  \cite[Definition~5.21]{BrundanStroppel:Semiinfinite}, with special idempotents $e_\lambda$, sets
  $Y(\lambda)$, $H(\lambda)$, $X(\lambda)$ and local algebras
  $\CA_\lambda=\bar e_\lambda\CA^{\le\lambda}\bar e_\lambda$ with basis $\set{\bar h|h\in H(\lambda)}$.
  Then $\plam=Y(\lambda)$, $\hlam=\set{e_\lambda}$, $\Olam=\CA_\lambda$ with
  $\olam=\set{\bar h|h\in H(\lambda)}$, $\qlam=X(\lambda)$ and $a^\lambda_{ye_\lambda\bar hx}=yhx$ is a
  perfectly respectable datum for $\CA$. Indeed, \autoref{resp1} is
  \cite[(PS4)]{BrundanStroppel:Semiinfinite}, \autoref{resp2} follows from
  \cite[(PS5)]{BrundanStroppel:Semiinfinite} as in \autoref{Ex:BasedQH}, and \autoref{resp3} holds because  the $h\in H(\lambda)$ multiply like the basis $\olam$ of $\CA_\lambda$ modulo $\CA^{>\lambda}$.
  As in \autoref{Ex:StandardBasedQH}, $a_{y\bar h}\mapsto yh\bar e_\lambda$ identifies
  $\Delta^\lambda$ with the standard module $\Delta(\lambda)=\CA^{\le\lambda}\bar e_\lambda$ of
  \cite[(5.6)]{BrundanStroppel:Semiinfinite} as $(\CA,\CA_\lambda)$-bimodules, and
  $\nabla^\lambda\cong(\bar e_\lambda\CA^{\le\lambda})^*$ with the costandard module $\nabla(\lambda)$ of \cite[(5.7)]{BrundanStroppel:Semiinfinite}. As
  $\CA_\lambda$ is local, $\Delta^\lambda(\kk)$ and $\nabla^\lambda(\kk)$ are the
  proper standard and proper costandard modules of \cite[\S5.3]{BrundanStroppel:Semiinfinite} respectively. Note that if
  $\CA_\lambda$ is not self-injective, then
  $\Hom_{\CA_\lambda}(\bar e_\lambda\CA^{\le\lambda},\CA_\lambda)$ differs from $\nabla^\lambda$.
\end{Example}

\begin{Example}\label{Ex:SergeevTwo}
  Let $2\in\kk^\times$ with $\sqrt2\in\kk$ and let $\Hfind[2]=\mathfrak{Cl}_2\rtimes\kk\Sym_2$ be the Sergeev
  superalgebra of \autoref{ex:Sergeev}. Let $e=\tfrac12(1+H_1)$ and $c=\tfrac1{\sqrt2}(C_1+C_2)e$, so that $c^2=e$. Then
  $e\Hfind[2]e=\kk e\oplus\kk c\cong\mathfrak{Cl}_1$, and $\plam[\bullet]=\qlam[\bullet]=\set{1,C_1}$,
  $\hlam[\bullet]=\set e$, $\olam[\bullet]=\set{e,c}$ is a respectable datum with a
  single stratum as in \autoref{Ex:Qn}. The standard module $\Delta^\bullet=\Hfind[2]e=\Hfind[2](1+H_1)$ is the
  simple module of type $\mathtt Q$ labelled by the strict partition $(2)$ in \cite{Sergeev}. It agrees with the parabolic module of  Mori from \cite[\S6]{Mori}.
\end{Example}

  We next establish the existence of a \emph{balanced pairing} between the right
  $\CA$-module $\leftindex{^\lambda}\Delta$ and the left
  $\CA$-module $\Delta^\lambda$. That is, a bilinear map
  $\<\ ,\ \>\colon\leftindex{^\lambda}\Delta\times\Delta^\lambda\to\Olam$ that is
  compatible with the $\Olam$-actions and such that $\<zx,y\>=\<z,xy\>$,
  for all $z\in\leftindex{^\lambda}\Delta$, $x\in \CA$, and
  $y\in\Delta^\lambda$.

  \begin{Proposition}\label{L:Respectable}
     Let $\lambda\in\Lambda$. There is a unique bilinear map
    $\<\ ,\ \>\colon\leftindex{^\lambda}\Delta\times\Delta^\lambda\to\Olam$ such that
     \[
       a_{phoq}a_{p'h'o'q'}\equiv\<a_{oq},a_{p'o'}\>a_{ph1_\lambda q'}
            \pmod{\CA^{>\lambda}},
     \]
    for $p,p'\in\plam$, $h,h'\in\hlam$, $o,o'\in\olam$, and $q,q'\in\qlam$.
    This form is $\Olam$-bilinear and $\CA$-balanced. That is,
     \[
       \<o''z,y\>=o''\<z,y\>,\qquad
       \<z,yo''\>=\<z,y\>o'',\qquad\text{and}\qquad
       \<zx,y\>=\<z,xy\>,
     \]
    for all $z\in\leftindex{^\lambda}\Delta$, $y\in\Delta^\lambda$, $o''\in\Olam$, and $x\in\CA$.
    Further, $\<\ ,\ \>$ is identically zero if $|\hlam|>1$.
  \end{Proposition}

  \begin{proof}
     Let $p,p'\in\plam$, $h,h'\in\hlam$, $q,q'\in\qlam$ and $o,o'\in\olam$.
     By \autoref{L:hlarge},
     $a_{phoq}a_{p'h'o'q'}\equiv ca_{ph1_\lambda q'}
        \pmod{\CA^{>\lambda}}$,
     where $c\in\Olam$ depends only on $o$, $q$,
     $p'$ and~$o'$. Moreover, $c\ne0$ only if $|\hlam|=1$. Hence,
     the form given by setting $\<a_{oq},a_{p'o'}\>=c$ is
     well-defined.

     It remains to show that $\<\ ,\ \>$ is $\Olam$-linear and
     $\CA$-balanced. By construction, the map
     $\<\ ,\ \>$ is defined on a $\kk$-basis of
     $\leftindex{^\lambda}\Delta\times\Delta^\lambda$, it extends uniquely to a
     $\kk$-bilinear form
     $\<\ ,\ \>\colon\leftindex{^\lambda}\Delta\times\Delta^\lambda\to\Olam$.
     Let $x\in \CA$ and set $s_{o''q''}=s_{oqo''q''}(x)$ and
     $r_{p''o''}=r_{p'o'p''o''}(x)$. Using \autoref{D:Respectable}, and
     what we have just proved,
     \begin{align*}
        \<a_{oq}x, a_{p'o'}\>a_{ph1_\lambda q'}
           &=\quad \sum_{\mathclap{o''\in\olam, q''\in\qlam}}
                  s_{o''q''}\<a_{o''q''},a_{p'o'}\>a_{ph1_\lambda q'}
            \equiv \quad\sum_{\mathclap{o''\in\olam, q''\in\qlam}}
                  s_{o''q''}a_{pho''q''}a_{p'ho'q'}\\
          & \equiv(a_{phoq}x)a_{p'ho'q'}=a_{phoq}(xa_{p'ho'q'})
           \equiv\quad\sum_{\mathclap{p''\in\plam, o''\in\olam}}
               r_{p''o''}a_{phoq}a_{p''ho''q'}\\
         &   \equiv\quad\sum_{\mathclap{p''\in\plam, o''\in\olam}}
               r_{p''o''}\<a_{oq},a_{p''o''}\>a_{ph1_\lambda q'}
           =\<a_{oq},xa_{p'o'}\>a_{ph1_\lambda q'},
     \end{align*}
     where all congruences are modulo $\CA^{>\lambda}$. Comparing the
     coefficient of $a_{ph1_\lambda q'}$ on both sides shows that the form is
     balanced. Finally, suppose that under \autoref{resp3}, $o''\cdot a_{oq}$ and $a_{p'o'}\cdot o''$
     correspond to $o''\cdot a_{phoq}$ and $a_{p'ho'q'}\cdot o''$,
     respectively. Then,  the compatibilities
     $(o''a_{phoq})a_{p'ho'q'}\equiv o''(a_{phoq}a_{p'ho'q'})$ and
     $(a_{phoq}a_{p'ho'q'})o''\equiv a_{phoq}(a_{p'ho'q'}o'')$ of \autoref{resp3}
     give $\<o''z,y\>=o''\<z,y\>$ and $\<z,yo''\>=\<z,y\>o''$, because
     $o''\cdot(ca_{ph1_\lambda q'})=(o''c)a_{ph1_\lambda q'}$ and
     $(ca_{ph1_\lambda q'})o''=(co'')a_{ph1_\lambda q'}$ for $c\in\Olam$.
     This completes the proof.
  \end{proof}

  As the form $\<\ ,\ \>$ vanishes when $|\hlam|>1$, we have:

  \begin{Corollary}
    \label{C:hlarge2}
    Assume $\CA$ is a respectable algebra with $|\hlam|\ge2$, for
    $\lambda\in\Lambda$. Then the pairing $\<\ ,\ \>$ on $\leftindex{^\lambda}\Delta\times\Delta^\lambda$ is zero.
    Hence $\Lambda_+\subseteq\set{\lambda\in\Lambda|{|\hlam|}=1}$.
  \end{Corollary}

  \begin{Definition}\label{D:partialx}
    Let $\lambda\in\Lambda$ and set
    $\nabla^\lambda_{\Olam}=\Hom_{\Olam}(\leftindex{^\lambda}\Delta,\Olam)$.
    For $x\in\Delta^\lambda$, let $\partial_x\colon\leftindex{^\lambda}\Delta\to\Olam$ be the map
    \[
        \partial_x(y)=\<y,x\>\qquad\text{for }y\in\leftindex{^\lambda}\Delta.
    \]
  \end{Definition}

  The map $\partial_x$ is left $\Olam$-linear by
  \autoref{L:Respectable}, so it is an element of
  $\nabla^\lambda_{\Olam}$. Observe that $\nabla^\lambda_{\Olam}$ becomes an $(\CA,\Olam)$-bimodule via $(xfo)(y)=f(yx)o$.

  Since the form $\<\ ,\ \>$ is balanced and $\Olam$-bilinear by \autoref{L:Respectable}, we obtain:

  \begin{Corollary}\label{Ahom}
    Let $\lambda\in\Lambda$. Then
    $\partial^\lambda\colon\Delta^\lambda\to\nabla^\lambda_{\Olam}$, $x\mapsto\partial_x$, is an
    $(\CA,\Olam)$-bimodule homomorphism.
  \end{Corollary}

   If $M$ is an $\CA$-module and $L$ is a simple
   $\CA$-module, let $[M:L]$ be the \emph{Jordan-H\"older
   multiplicity} of~$L$ in~$M$.  In the super case, $[M: L]$ counts the
   occurrences of $L$ and its parity shift $\Pi L$ in $M$. Recall from
   \autoref{uEnd} that $\uHom_{\CA}(M,N)$ is the space of inner
   homomorphisms, which is different from the space of all
   $\CA$-module homomorphisms only in the super case.

   \begin{Definition}\label{D:RespectableSimples}
     Suppose that $\kk$ is a field. For
     $\lambda\in\Lambda$, define
     \[
       L^\lambda=\set{\partial_x|x\in\Delta^\lambda}\subseteq\nabla^\lambda_{\Olam}
       \qquad\text{and}\qquad
       \rad\Delta^\lambda=\set{x\in\Delta^\lambda|\partial_x=0}\subseteq\Delta^\lambda.
     \]
     Set $\Lambda_+=\set{\mu\in\Lambda|L^\mu\ne0}$.
     The \emph{standard
     and costandard decomposition matrices} of $\CA$ are the matrices
     \[
         D_\Delta=\bigl([\Delta^\lambda:L^\mu]\bigr)_{\lambda\in\Lambda, \mu\in\Lambda_+}
         \qquad\text{and}\qquad
         D_\nabla=\bigl([\nabla^\lambda:L^\mu]\bigr)_{\lambda\in\Lambda, \mu\in\Lambda_+}.
     \]
     Let $P^\lambda$ be the projective cover of $L^\lambda$ and set
     $C_{\CA}=\Bigl(\dim\uHom_{\CA}(P^\mu,P^\nu)\Bigr)_{\mu,\nu\in\Lambda_+}$,
     the \emph{Cartan matrix} of $\CA$.
      Let
     $H_{\CA}$ be the diagonal matrix  $\op{diag}\left(\tfrac{|\hlam|}{|\olam|}\right)_{\lambda\in\Lambda}$
     and let $\Olam^{\mathrm{op}}$  be the opposite algebra of $\Olam$.
   \end{Definition}

   By \autoref{Ahom}, $L^\lambda$ and $\rad\Delta^\lambda$ are
   $(\CA,\Olam)$-bimodules.

\begin{Example}\label{Ex:KleinFour}
  Let $\op{char}(\kk)=2$ and let $\CA=\kk V_4$ be the group algebra of the Klein four
  group $V_4=\<g,h\mid g^2=h^2=1,\ gh=hg\>$. With $x=g+1$ and $y=h+1$ we have
  $\CA=\kk[x,y]/(x^2,y^2)$. Set $\Lambda=\set{0<1<2}$ with $\hlam[0]=\set{1}$,
  $\hlam[1]=\set{x,y}$, $\hlam[2]=\set{xy}$. Letting all other indexing sets be equal to
  $\set{1}$, $\Olam=\kk$ and $a_{1h11}=h$, we get a perfectly $\kk$-respectable datum for
  $\CA$, with $\dot h=h$ and all other lifts equal to $1$. The only condition to check
  is \autoref{resp2} for $\lambda=1$. It holds, because modulo $\CA^{>1}=\kk xy$ any
  $z\in\CA$ acts on $x$ and on $y$, from both sides, by the same scalar, the coefficient
  of $1$ in $z$. As predicted by \autoref{C:hlarge2}, the form vanishes on the middle stratum,
  and by \autoref{L:DeltaFiltration} the trivial module $\Delta^1$ occurs twice in a
  $\Delta$-filtration of $\CA$. Since also $(xy)^2=0$, we get $\Lambda_+=\set{0}$ with
  $L^0=\kk$.
\end{Example}

The following is our motivation to call the $\Olam$ fibre algebras.
\begin{Example}\label{Ex:LerayHirsch}
  Natural instances of respectable data with $|\hlam|>1$ and non-trivial fibre algebras
  $\Olam$ arise from fibre bundles. Let $Q\subseteq P$ be parabolic subgroups of a connected
  complex reductive group $G$, and $\pi\colon X=G/Q\to Y=G/P$ the fibre bundle with fibre
  $\iota\colon F=\pi^{-1}(eP)\hookrightarrow X$. By the Bruhat decompositions, the cohomology, over $\kk$,  is even. Lift a basis $\mathfrak o$ of $H^*(F)$ to $\tilde o\in H^*(X)$,
  $\iota^*(\tilde o)=o$. Then $\CA=H^*(X)$ is perfectly respectable with
  $\Lambda=\set{0<\dots<\dim_{\mathbb C}Y}$, $\plam[j]=\qlam[j]=\set{1}$, $\hlam[j]$ a basis of
  $H^{2j}(Y)$, $\Olam[j]=H^*(F)$ with basis $\mathfrak o$, and $a_{1ho1}=\pi^*(h)\tilde o$.
  Indeed, \autoref{resp1} is the Leray--Hirsch theorem. For $z\in\CA$ write
  $z\tilde o=\sum c_{h'o'}\pi^*(h')\tilde o'$. If $h\in\hlam[j]$, then $\pi^*(hh')\tilde o'\in
  \CA^{>j}$ for $\deg h'>0$, hence
  $z\pi^*(h)\tilde o\equiv\sum_{o'}c_{1o'}\pi^*(h)\tilde o'$ modulo $\CA^{>j}$, with
  $\iota^*(z)o=\sum_{o'}c_{1o'}o'$ as $\iota^*\pi^*$ kills positive degrees. So $z$ acts on each
  stratum through $\iota^*(z)\in H^*(F)$, independently of $j$ and $h$. This gives \autoref{resp2}
  and, as $\CA$ and $H^*(F)$ are commutative, \autoref{resp3}. In accordance with \autoref{C:hlarge2}, $\Lambda_+=\set{0}$.
\end{Example}

\begin{Example}\label{Ex:gradedlocal}
  Let $\CA=\bigoplus_{d\ge0}\CA_d$ be a finite dimensional $\N_0$-graded algebra with $\CA_0=\kk$, for instance an exterior algebra or the cohomology ring of say a flag variety. Then $\CA$ is perfectly $\kk$-respectable with weight poset $\Lambda=\set{0<1<2<\dots}=\N_0$ by setting  $\plam[d]=\qlam[d]=\set{1}$ and letting $\hlam[d]$ be a basis of $\CA_d$, since $zh\equiv z_0h\equiv hz$ modulo $\CA_{>d}$ for $z\in\CA$ with degree zero part $z_0$ and $h\in\hlam[d]$. Here $|\hlam[d]|=\dim\CA_d$. In accordance with \autoref{C:hlarge2}, only the stratum $d=0$ carries a simple module.
\end{Example}

\begin{Example}\label{R:KoenigXi}
  Let $\kk$ be a field.  A \emph{finite dimensional affine cellular
  algebra} in the sense of K\"onig and Xi
  \cite[Definition~2.1]{KoenigXi:AffineCellular} is respectable. Indeed,
  each layer $J_j/J_{j-1}\cong V_j\otimes B_j\otimes V_j$ of the cell
  chain gives a stratum, ordered as in the cell chain, with
  $\plam[j]=\qlam[j]$ a basis of $V_j$, $\hlam[j]=\set{1}$ and
  $\Olam[j]=B_j$. Here,  \autoref{resp1} holds as
  $\CA=J'_1\oplus\dots\oplus J'_n$ with $J'_j\cong J_j/J_{j-1}$,
  \autoref{resp2} follows from the formula $a(u\otimes b\otimes
  v)=f(a,u)b\otimes v$ with $f(a,u)\in V_j\otimes B_j$ in the proof of
  \cite[Proposition~2.2]{KoenigXi:AffineCellular} and its mirror image,
  and \autoref{resp3} holds for the action of $B_j$ on $V_j\otimes
  B_j\otimes V_j$ by multiplication on the middle factor, which commutes
  with the action of $\CA$ as $\CA$ acts $B_j$-linearly on
  $\Delta_j=V_j\otimes B_j$ and on $\Delta'_j=B_j\otimes V_j$. As
  mentioned in \cite[Section~2.1]{KoenigXi:AffineCellular}, the $\mathbb
  R$-algebra $\mathbb C$ is affine cellular but not cellular. It is
  respectable with $\Olam=\mathbb C$ as in \autoref{stupid}. Note that
  respectable data do not need an involution, and they allow
  non-commutative fibre algebras. Affine cellularity is, however, most useful for infinite dimensional algebras, such as affine Hecke algebras of type $A$, which are excluded by the finite rank assumption in \autoref{D:Respectable}.
\end{Example}

 Recall that a $\kk$-algebra $\CO$ is a \emph{division algebra} if every
 nonzero element of $\CO$ is invertible. If $\CO$ admits a
 $\Z/2\Z$-grading, then $\CO$ is a \emph{division superalgebra} if every
 nonzero homogeneous element of $\CO$ is invertible. A trivially graded
 algebra is a division superalgebra if and only if it is a division
 algebra, so we only speak of division superalgebras.

Over an algebraically closed field $\kk$, the super division algebras
are $\kk$ and the Clifford algebra $\mathfrak{Cl}_1$, see
\cite{Jozefiak}.  So if $\kk$ is algebraically closed and the fibre
algebra $\Olam[\mu]$ is a division superalgebra, then $\Olam[\mu]=\kk$
or $\Olam[\mu]\cong\mathfrak{Cl}_1$. Recall from \autoref{uEnd} the
notation $\uEnd_{\CA}$ for all linear endomorphisms  in the super case.

The next proposition mimics the theory of cellular algebras,
following~\cite{GrahamLehrer:Cellular} and
\cite[Chapter~2]{Mathas:ULect}.  By definition, $\rad\Delta^\lambda$ is
the radical of the bilinear form $\<\ ,\ \>$. This proposition shows
that $\rad\Delta^\lambda$ agrees with the radical of the module
$\Delta^\lambda$ whenever $L^\lambda\ne0$ and $\Olam$ is a division
superalgebra.

Before we state the proposition we record the relevant properties of
division superalgebras.

\begin{Lemma}\label{L:DivisionFrobenius}
  Let $\kk$ be a field and suppose that the fibre algebra $\Olam[\mu]$ is a
  division superalgebra, for $\mu\in\Lambda_+$. Then there exist
  $y\in\leftindex{^\mu}\Delta$ and $x\in\Delta^\mu$ with $\<y,x\>=1_\mu$.
  Moreover, $\Olam[\mu]\cong\Olam[\mu]^*$. Finally, a finite dimensional
  division superalgebra $\CO$ is Frobenius.
\end{Lemma}

\begin{proof}
  The form $\<\ ,\ \>$ is homogeneous, so $\rad\Delta^\mu$ is a proper graded
  submodule of $\Delta^\mu$ since $\mu\in\Lambda_+$.
  Therefore, there exist a homogeneous $x\in \Delta^\mu$
  with $x\notin\rad\Delta^\mu$. By \autoref{L:Respectable},
  $\set{\<y,x\>|y\in\leftindex{^\mu}\Delta}$ is a nonzero graded left ideal of
  $\Olam[\mu]$. It contains $1_\mu$ since $\Olam[\mu]$ is a
  division superalgebra. That is, $\<y,x\>=1_\mu$, for some $y\in\nabla^\mu$,
  proving the first claim.

 To complete the proof we claim that every finite dimensional division
 superalgebra $\CO$ is Frobenius. To see this, let $\tr\colon\CO\to\kk$
 be any linear map such that $\tr(1)\ne0$. If
 $0\ne a=a_{\ov0}+a_{\ov1}\in\CO$, then define $b=a_{\ov0}^{-1}$, if
 $a_{\ov0}\ne0$, and $b=a^{-1}$, if $a=a_{\ov1}$. Then
 $\tr(ab)=\tr(1)\ne0$, so the form $\tr$ is non-degenerate. In
 particular, $\Olam[\mu]$ is a Frobenius $\kk$-superalgebra, and
 $\Olam[\mu]\cong\Olam[\mu]^*$ as $\Olam[\mu]$-modules.
\end{proof}

For $\mu\in\Lambda_+$, let $\Olam[\mu]^{\mathrm{op}}$ be the
\emph{opposite algebra} of $\Olam[\mu]$. If $\kk$ is a field and $M$ is an
$\CA$-supermodule, its \emph{radical} $\rad M$ is the smallest
submodule of~$M$ such that $M/\rad M$ is semisimple.

\begin{Theorem}\label{T:Respectable}Let $\CA$ be a respectable algebra over a field $\kk$. Then the following hold.
\begin{enumerate}
  \item\label{late1}  Let $\mu\in\Lambda$. Then $L^\mu$ is an
         $\CA$-submodule of $\nabla^\mu_{\Olam[\mu]}$ and $\rad\Delta^\mu$ is
         an $\CA$-submodule of~$\Delta^\mu$. Moreover,
         $L^\mu\cong\Delta^\mu/\rad\Delta^\mu$ as
         $\CA$-modules.
  \item\label{late2}  Let $\lambda,\mu\in\Lambda$ and
  suppose that there exist homogeneous $y\in\leftindex{^\mu}\Delta$ and $x\in\Delta^\mu$
  such that $\<y,x\>=1_\mu$. Let $M$ be a submodule of~$\Delta^\lambda$. Then
  $\Delta^\mu=\CA x$ and $\Hom_{\CA}(\Delta^\mu,\Delta^\lambda/M)\ne0$ only if
  $\lambda\ge\mu$. Moreover,
  \[
     \uEnd_{\CA}(L^\mu)\cong\uEnd_{\CA}(\Delta^\mu)\cong\uEnd_{\CA}(\nabla^\mu)\cong\Olam[\mu]^{\mathrm{op}}.
  \]
\end{enumerate}
Suppose now that the fibre algebra $\Olam[\mu]$ is a division
superalgebra for all $\mu\in\Lambda_+$. Then:
\begin{enumerate}[start=3]
  \item\label{late3} Suppose that $\mu\in\Lambda_+$. Then $L^\mu$ is an
  irreducible $\CA$-module. Moreover, $L^\mu$ is isomorphic to
  the socle of $\nabla^\mu$ and to the head of $\Delta^\mu$. In the
  super case, $L^\mu$ is of type $\mathtt{M}$ if $\Olam[\mu]\cong\kk$,
  and of type $\mathtt{Q}$ if~$\Olam[\mu]\cong\mathfrak{Cl}_1$.
  \item\label{late4}  The set $\set{L^\mu|\mu\in\Lambda_+}$ is a
  complete set of pairwise non-isomorphic irreducible
  $\CA$-modules in case of ordinary algebras. In case of
  superalgebras, every irreducible supermodule is isomorphic to $L^\mu$,
  or to $\Pi L^\mu$, for some $\mu\in\Lambda_+$. Moreover,
  $L^\mu\cong\Pi^\epsilon L^\nu$ for some $\epsilon$ only if $\mu=\nu$,
  with $L^\mu\cong\Pi L^\mu$ exactly in type $\mathtt Q$.
  \item\label{late5}  If $\mu\in\Lambda_+$ and $\lambda\in\Lambda$, then
  $[\Delta^\lambda:\Pi^\epsilon L^\mu]\ne0$ only if $\lambda\ge\mu$.
  Moreover, $[\Delta^\mu:L^\mu]=1$.
  \item\label{late6}  Suppose moreover that $\Olam$ is local for every $\lambda\in\Lambda$. Then the Cartan matrix of~$\CA$ factorises as
  \[
     C_{\CA}=\diag(|\olam[\nu]|)_{\nu\in\Lambda_+}
         \cdot D_\Delta^T\cdot  H_{\CA}
         \cdot D_\nabla\cdot \diag(|\olam[\mu]|)_{\mu\in\Lambda_+}.
  \]
\end{enumerate}
\end{Theorem}

\begin{Remark}
If moreover $\Olam=\kk$ and $|\hlam|=1$ for all $\lambda\in\Lambda$, then $\CA$ is standardly based, see \autoref{Ex:StandardlyBased}. Requiring for each $\mu\in\Lambda$ the existence of $x,y$ with $\<y,x\>=1_\mu$  means $\Lambda_+=\Lambda$ and implies that $\CA$ is standardly full-based \cite[Definition~1.3.1]{DuRui:strat} or split quasi-hereditary, \cite[Theorem~4.2.7]{DuRui:strat}.
\end{Remark}

   \begin{proof}[Proof of \autoref{T:Respectable}]
     \autoref{late1}: By \autoref{Ahom}, the linear map
     $\partial^\mu\colon\Delta^\mu\to\nabla^\mu_{\Olam[\mu]}$, given by
     $\partial^\mu(x)=\partial_x$, is an
     $\CA$-module homomorphism. In particular, $L^\mu=\mathop{\rm im}\partial^\mu$ and
     $\rad\Delta^\mu=\ker\partial^\mu$ are $\CA$-modules. Hence,
     $L^\mu\cong\Delta^\mu/\rad\Delta^\mu$ by the first isomorphism
     theorem.

     \autoref{late2}: By assumption, $x,y\in\Delta^\mu$ are homogeneous
     elements with $\<x,y\>=1$. Write
     $y=\sum_{oq\in\olam[\mu]\qlam[\mu]}c_{oq}a_{oq}$, for some
     $c_{oq}\in\kk$. For $p\in\plam[\mu]$ and $h\in\hlam[\mu]$,
     set $y_{ph}=\sum_{oq\in\olam[\mu]\qlam[\mu]}c_{oq}a_{phoq}\in
     \CA^{(\mu)}$. In view of \autoref{L:Respectable},
     \[
         y_{ph}\cdot x
           = \sum_{oq\in\olam[\mu]\qlam[\mu]}c_{oq}a_{p1_\mu}\<a_{oq},x\>
           =a_{p1_\mu}\biggl\<\sum_{oq\in\olam[\mu]\qlam[\mu]}c_{oq}a_{oq},x\biggr\>
           = a_{p1_\mu}\<y,x\> = a_{p1_\mu}.
     \]
     Since $\<o\cdot y,x\>=o\<y,x\>=o$ for $o\in\olam[\mu]$, running the same
     computation with $o\cdot y$ in place of $y$ gives $a_{po}=(o\cdot y)_{ph}\cdot x\in\CA x$, so
     $\CA x=\Delta^\mu$, proving the first claim.

     Now let $\varphi\colon\Delta^\mu\to\Delta^\lambda/M$ be a nonzero
     $\CA$-module homomorphism. By what we have just shown,  $a_{po}=(o\cdot y)_{ph}\cdot x$
     for $p\in\plam[\mu]$, $o\in\olam[\mu]$ and $h\in\hlam[\mu]$. Since
     $\varphi\ne0$, there exist $a_{po}\in\Delta^\mu$ such that
     $\varphi(a_{po})=(o\cdot y)_{ph}\varphi(x)\ne0$, so $\CA^{(\mu)}$ does not annihilate
     $\Delta^\lambda/M$. As $\CA^{\ge\mu}$ is a two-sided ideal, \autoref{resp1} and \autoref{resp2}
     force $\lambda\ge\mu$, proving the second statement in~\autoref{late2}.
     Now suppose that $\mu=\lambda$ and write $\varphi(x)=x'+M$ with $x'\in\Delta^\mu$ and $c=\<y,x'\>$.
     Then $\varphi(a_{po})=(o\cdot y)_{ph}\cdot x'+M=a_{p1_\mu}\<o\cdot y,x'\>+M=a_{po}c+M$.
     So $\varphi$ is right multiplication, $\rho_c\colon m\mapsto mc$, followed by the projection
     to~$\Delta^\mu/M$. Every $\rho_c$ is an endomorphism of $\Delta^\mu$ by \autoref{std1}, and
     $\rho_c\rho_{c'}=\rho_{c'c}$. Taking $M=0$ gives $\End_{\CA}(\Delta^\mu)\cong\Olam[\mu]^{\mathrm{op}}$.
     By symmetry, $\End_{\CA}(\leftindex{^\mu}\Delta)\cong\Olam[\mu]$. Thus, $\End_{\CA}(\nabla^\mu)\cong\Olam[\mu]^{\mathrm{op}}$
     as $\nabla^\mu=\leftindex{^\mu}\Delta^*$ by \autoref{std3}. Taking $M=\rad\Delta^\mu$ shows, by~\autoref{late1}, that every
     endomorphism of $L^\mu$ is induced by some $\rho_c$.  Conversely, $\rho_c$ preserves $\rad\Delta^\mu$,
     as $\<z,mc\>=\<z,m\>c$, and it induces the zero map on $L^\mu$ only if $c=\<y,x\>c=0$. Hence
     $\End_{\CA}(L^\mu)\cong\Olam[\mu]^{\mathrm{op}}$.

     In the super case the same argument applies to homogeneous maps of
     either parity,  giving the statements for the space of inner
     endomorphisms with the usual signs. For instance, $\rho_c$ for
     homogeneous $c$ is  $m\mapsto(-1)^{|c||m|}mc$.

     \autoref{late3}:
     Suppose that $x\in\Delta^\mu$ is homogeneous and $x\notin\rad\Delta^\mu$. Then
     $\<y,x\>\ne0$, for some $y\in\nabla^\mu$.Without loss of
     generality, $\<y,x\>=1$ since $\Olam[\mu]$ is a division superalgebra.
     Therefore, $\Delta^\mu=\CA x$ by \autoref{late2}.
     As this argument applies to any homogeneous element
     $x\in\Delta^\mu$ that is not in $\rad\Delta^\mu$, this shows that
     $\rad\Delta^\mu$ is the unique maximal submodule of
     $\Delta^\mu$. Hence, $L^\mu\cong\Delta^\mu/\rad\Delta^\mu$ is the
     simple head of $\Delta^\mu$ by \autoref{late1}. By the same
     argument, the right $\CA$-module $\leftindex{^\mu}\Delta$ has
     simple head, so its dual $\nabla^\mu$ has simple socle. By
     construction, $\partial^\mu$ is a nonzero homomorphism from
     $\Delta^\mu$ to $\nabla^\mu$ with image $L^\mu$, so $L^\mu$ is the
     socle of $\nabla^\mu$. Further, by \autoref{std3},
     $\nabla^\mu_{\Olam[\mu]}\cong\nabla^\mu$ since
     $\Olam[\mu]\cong\Olam[\mu]^*$ by \autoref{L:DivisionFrobenius}.
     The remaining parts of \autoref{late3} now follow by
     \autoref{late2}.

     \autoref{late4}: By part~\autoref{late2}, the modules $\set{L^\mu|\mu\in\Lambda_+}$ are
     pairwise non-isomorphic. That this gives a complete set of
     irreducible $\CA$-modules follows, in the case of ordinary
     algebras,  easily by induction on the partial
     order on~$\Lambda$. In fact, we can repeat the proof of
     \cite[Theorem~2.16]{Mathas:ULect}, almost word for word, with the only
     change that the scalars there now lie in the division
     superalgebra $\Olam[\mu]$ by Schur's lemma.

     In the super case the same induction produces every irreducible
     supermodule as a quotient of some standard module~$\Delta^\lambda$.
     Hence, every simple module is isomorphic to $L^\mu$ or $\Pi L^\mu$,
     for some $\mu\in\Lambda_+$. These modules are pairwise
     non-isomorphic by~\autoref{late2}, except that $L^\mu\cong\Pi
     L^\mu$ when $\Olam[\mu]$ has an invertible odd element, which is if
     and only if it is of type~$\mathtt{Q}$.

     \autoref{late5}: This follows from parts~\autoref{late2}-\autoref{late3}.

     \autoref{late6}: As in the proof of \autoref{L:DeltaFiltration},
     $\CA$ is filtered by the two-sided ideals $\CA^{\ge\lambda}$. By
     \autoref{resp2} and \autoref{resp3}, the map
     $a_{po}\otimes a_{o'q}\mapsto a_{ph(oo')q}+\CA^{>\lambda}$ gives an isomorphism of
     $\CA$-bimodules
     $\CA^{\ge\lambda}/\CA^{>\lambda}
         \cong(\Delta^\lambda\otimes_{\Olam}\leftindex{^\lambda}\Delta)^{\oplus|\hlam|}$.
     Assume we work with ordinary algebras. Write $P^\mu=\CA e$ for an
     idempotent $e$. Multiplying the filtration with $e$ from the right
     gives a filtration of $P^\mu$ with subquotients
     $(\Delta^\lambda\otimes_{\Olam}\leftindex{^\lambda}\Delta e)^{\oplus|\hlam|}$.
     The left $\Olam$-module $\leftindex{^\lambda}\Delta e$ is a direct summand of
     the free $\Olam$-module $\leftindex{^\lambda}\Delta$, so it is projective and
     hence free of rank $\op{rk}_{\Olam}(\leftindex{^\lambda}\Delta e)/|\olam|$, as
     $\Olam$ is local. This gives a $\Delta$-filtration of $P^\mu$ in
     which $\Delta^\lambda$, for $\lambda\in\Lambda$, appears with
     multiplicity
     \[
       \tfrac{|\hlam|}{|\olam|}\op{rk}_{\Olam} \leftindex{^\lambda}\Delta\otimes_{\CA}P^\mu
        =\tfrac{|\hlam|}{|\olam|}\op{rk}_{\Olam}\Hom_{\CA}(P^\mu,\nabla^\lambda)
        =\tfrac{|\hlam|}{|\olam|}|\olam[\mu]|[\nabla^\lambda:L^\mu],
     \]
     by tensor-hom adjunction as in the proof of
     \cite[Lemma~2.19]{Mathas:ULect}. Therefore,
     \[
       c_{\nu\mu} = \op{rk}\Hom_{\CA}(P^\nu,P^\mu) =|\olam[\nu]|[P^\mu\colon L^\nu]
          =|\olam[\nu]||\olam[\mu]|\sum_{\lambda\in\Lambda} [\Delta^\lambda:L^\nu]
             \tfrac{|\hlam|}{|\olam|}[\nabla^\lambda:L^\mu].
     \]
     Hence,  $C_{\CA}=\diag(|\olam[\nu]|)D_\Delta^T H_{\CA} D_\nabla\diag(|\olam[\mu]|)$, which completes the proof of \autoref{late6} for ordinary algebras.
     In the super case, with $\uHom_{\CA}$ in place of $\Hom_{\CA}$ as in \autoref{D:RespectableSimples}, the argument works verbatim, since every composition factor isomorphic to $L^\mu$ or $\Pi L^\mu$ contributes $\op{rk}\uEnd_{\CA}(L^\mu)=|\olam[\mu]|$ to $\op{rk}\uHom_{\CA}(P^\mu,M)$. This is $1$ in type $\mathtt M$ and $2$ in type $\mathtt Q$.
  \end{proof}

The algebra $\CA$ is \emph{split} if $\uEnd_{\mathcal{A}}(S)\cong\kk$ or
$\uEnd_{\mathcal{A}}(S)\cong\Cl_1$ for every simple supermodule $S$.

\begin{Corollary}\label{C:Split}
  Let $\mathcal{A}$ be a finite dimensional superalgebra over a field
  $\kk$.
  \begin{enumerate}
  \item The algebra $\mathcal{A}$ is split if and only if $\mathcal{A}$ admits a respectable datum with homogeneous basis and $\Olam\in\set{\kk,\Cl_1}$ for all $\lambda\in\Lambda$.
  \item $\mathcal{A}$ is $\kk$-respectable with a homogeneous basis if and only if $\uEnd_{\mathcal{A}}(S)=\kk$ for every simple supermodule $S$.
  \end{enumerate}
  In the non-super case, we assume that $\operatorname{char}\kk\ne2$, and (b) says that $\mathcal{A}$ is $\kk$-respectable if and only if it is split.
\end{Corollary}

\begin{proof}
  By \autoref{T:Respectable}\autoref{late4}, \autoref{late2} and \autoref{L:DivisionFrobenius}, if $\mathcal{A}$ has a respectable datum with homogeneous basis whose strata have algebras $\kk$ or $\Cl_1$, every simple supermodule is isomorphic to some $L^\mu$ or $\Pi L^\mu$, and $\uEnd_{\mathcal{A}}(L^\mu)$ is the opposite of the algebra of the stratum $\mu$. As $\Cl_1^{\mathrm{op}}\cong\Cl_1$ by \autoref{R:Clop}, this gives the "if" parts. Conversely, if $\mathcal{A}$ is split, then $\mathcal{A}/J\cong\prod_i\Mat_{n_i}(D_i)$ with $D_i\in\set{\kk,\Cl_1}$ for $J=\rad\mathcal{A}$, see \cite[Theorem~(2.6) and Remark~(2.7)]{Jozefiak} with \cite{Wall:GradedBrauer}.  Each layer $J^d/J^{d+1}$ is a sum of blocks $E^{(i)}_{p1}(H\otimes D_{ij})E^{(j)}_{1q}$ with $H$ a graded multiplicity space, $D_{ij}=\kk$ if $D_i=D_j=\kk$ and $D_{ij}=\Cl_1$ otherwise. Lifting these bases gives a respectable datum with strata $(d,i,j,\epsilon)$, where $\epsilon$ is the parity of the basis elements of $H$, ordered by $d$, and with algebra $D_{ij}$ for the stratum $(d,i,j,\epsilon)$. In particular, this datum is $\kk$-respectable if all $D_i=\kk$.
\end{proof}

Recall from \autoref{C:hlarge2} that $|\hlam|=1$ if
$\lambda\in\Lambda_+$. Therefore, if $\Lambda_+=\Lambda$, then
$|\hlam|=1$, for all $\lambda\in\Lambda$.

\begin{Corollary}\label{C:quasihereditary}
  Let $\mathcal{A}$ be a respectable (ordinary) algebra over a field~$\kk$ such that $\Olam$ is a division algebra for all $\lambda\in\Lambda$. If $\Lambda=\Lambda_+$, then $\mathcal{A}$ is quasi-hereditary, and for every total order refining $(\Lambda,\ge)$ the ideals $\mathcal{A}^{\ge\lambda}$ form a heredity chain.
\end{Corollary}

\begin{proof}
  Fix $\lambda\in\Lambda$ and put $B=\mathcal{A}/\mathcal{A}^{>\lambda}$ and $J=\mathcal{A}^{\ge\lambda}/\mathcal{A}^{>\lambda}$. As $\lambda\in\Lambda_+$, $|\hlam|=1$ by \autoref{C:hlarge2}. By \autoref{resp2} and \autoref{resp3}, $a_{po}\otimes a_{o'q}\mapsto a_{p,oo',q}+\mathcal{A}^{>\lambda}$ induces an isomorphism of $B$-bimodules $\Delta^\lambda\otimes_{\Olam}\leftindex{^\lambda}\Delta\cong J$, and by \autoref{D:Respectable} the multiplication of $J$ is $(z\otimes w)(z'\otimes w')=z\<w,z'\>\otimes w'$. By \autoref{L:DivisionFrobenius} there are $x\in\Delta^\lambda$ and $y\in\leftindex^{\lambda}\Delta$ with $\<y,x\>=1_\lambda$, so $e=x\otimes y$ is an idempotent with $J=JeJ$, hence $J^2=J$. Moreover $Be=Je=\Delta^\lambda\otimes y\cong\Delta^\lambda$, and $J\cong(\Delta^\lambda)^{\oplus|\qlam|}$ as a left $B$-module, as $\leftindex^{\lambda}\Delta$ is a free left $\Olam$-module with basis $\set{a_{1q}}$. Hence $J$ is projective. Finally, $eBe=eJe=x\Olam\otimes y$ is isomorphic to $\Olam$, a division algebra, so $e\rad(B)e=\rad(eBe)=0$ and $J\rad(B)J=Je\rad(B)eJ=0$. Thus $J$ is a heredity ideal of~$B$.
\end{proof}

\begin{Corollary}\label{C:Semisimple}
    Suppose that $A$ is a respectable algebra over a field~$\kk$. Then
    $A$ is semisimple if and only if $L^\lambda=\Delta^\lambda=\nabla^\lambda$, for all
    $\lambda\in\Lambda$.
\end{Corollary}

  \begin{proof}
  Suppose that $\mathcal{A}$ is semisimple, let $\lambda\in\Lambda$, and put $B=\mathcal{A}/\mathcal{A}^{>\lambda}$ and $I=\mathcal{A}^{\ge\lambda}/\mathcal{A}^{>\lambda}$. By \autoref{resp2}, $\Delta^\lambda$ is isomorphic to a graded left ideal of $B$ contained in $I$, and by \autoref{L:Respectable} the image $L$ of $\rad\Delta^\lambda$ satisfies $IL=0$. As $B$ is semisimple, $L=Be$ for an even idempotent $e$, so $e=e^2\in IL=0$. Hence $\Delta^\lambda\cong L^\lambda$ by \autoref{T:Respectable}\autoref{late1}, and similarly $\nabla^\lambda\cong L^\lambda$.

  Conversely, if all forms are non-degenerate on both sides, then
  $\Lambda=\Lambda_+$ and $|\hlam|=1$ by \autoref{C:hlarge2}, and
  $|\plam|=|\qlam|$. By \autoref{T:Respectable}\autoref{late2}, \autoref{late4} and the (super) Wedderburn
  theorem~\cite[Theorem~12.2.9]{Kleshchevbook},
  \[
    \dim\mathcal{A}/\rad\mathcal{A}=\sum_{\lambda\in\Lambda}\frac{(\dim L^\lambda)^2}{\dim\Olam}=\sum_{\lambda\in\Lambda}|\plam||\hlam||\olam||\qlam|=\dim\mathcal{A},
  \]
  so $\rad\mathcal{A}=0$.
  \end{proof}

   \begin{Example}\label{Ex:TaftII}
     We continue the analysis of the Taft algebra
     $\Taft$ from \autoref{Ex:Taft}. Let $\kk$ be a field and let $\op{char}(\kk)=p\geq0$. Note that the existence of $\omega$ implies  that $p$ is coprime to $\theta$.
     Fix $i\in\Lambda=\set{0,1,\dots,2\theta-1}$. By
     definition, the standard module $\Delta^i$ is $1$-dimensional and
     $y$ and $g$ act on $\Delta^i$ as multiplication by~$0$ and by $\omega^{-i}$
     respectively. Hence, the bilinear form
     $\leftindex^i\Delta\times\Delta^i\to\kk$ is given by
     $\<a_{y^k},a_{1_T}\>=\delta_{k0}{\prod_{j\ne i}(\omega^{-i}-\omega^{-j})}$.
     The assumption on $\op{char}(\kk)$ implies that  the
     $\set{\omega^i|0\le i<\theta}$ are distinct roots of unity, so
     $L^i=\Delta^i\ne0$, for $0\le i<\theta$, and  $L^i=0$, for
     $\theta\le i<2\theta$.
     Hence, $\Lambda_+=\set{0,1,\dots,\theta-1}$. In particular, for $j\in \Lambda_+$, the
     simple $\Taft$-module $L^j$ is the unique $1$-dimensional $\Taft$-module on which $g$ and $y$ act as
     multiplication by~$\omega^{-j}$ and $0$ respectively. It follows that
     $[\Delta^i:L^j]=1$ if $i\equiv j\pmod \theta$, and $[\Delta^i:L^j]=0$
     otherwise. For the costandard modules, $\nabla^0\cong L^0$ and
     $\nabla^i\cong L^{i-\theta+1}$, for $\theta<i<2\theta$, where the upper index of $L$ is read modulo $\theta$, so that $\nabla^{2\theta-1}\cong L^0$. If $1\le i\le
     \theta$, then
     $\op{rk}\nabla^i=\theta-1$ with the bases being given by \autoref{E:Taft}.
     Using  $h_iy^kg = \omega^{k-i}h_iy^k$ we obtain, for $1\le i<\theta$,
     \[
        [\nabla^i:L^j] = \begin{cases*}
           1 & if $j\not\equiv i+1 \pmod \theta$,\\
           0 & otherwise.
        \end{cases*}
        \qquad\text{and}\qquad
        [\nabla^\theta:L^j] = \begin{cases*}
           1 & if $j\ne0$,\\
           0 & otherwise.
        \end{cases*}
     \]
     The calculations show that $e_i=\bigl(\prod_{j\ne i}(\omega^{-i}-\omega^{-j})\bigr)^{-1}h_i$ is
     an idempotent in~$\Taft$, for $0\le i<\theta$.  The $\Taft$-module
      $P^j=\Taft e_j$ has basis $\set{y^ke_j|0\le k<\theta}$. Hence, it has
      simple head $L^j$, simple socle $L^{j-1}$, and $[P^i:L^j]=1$, for $0\le i,j<\theta$.
      In particular,  $P^j$ is indecomposable. The Cartan matrix of
      $\Taft$ is then the $\theta\times \theta$ matrix that has all entries
      equal to~$1$. The factorisation of \autoref{T:Respectable}\autoref{late6} can easily
      be verified. For instance, if $\theta=4$ we get
     \[
     C_{T_4} = \begin{psmallmatrix}
       1 & 1 & 1 & 1\\
       1 & 1 & 1 & 1\\
       1 & 1 & 1 & 1\\
       1 & 1 & 1 & 1\\
     \end{psmallmatrix}
     = \begin{psmallmatrix}
         1 & 0 & 0 & 0 & 1 & 0 & 0 & 0\\
         0 & 1 & 0 & 0 & 0 & 1 & 0 & 0\\
         0 & 0 & 1 & 0 & 0 & 0 & 1 & 0\\
         0 & 0 & 0 & 1 & 0 & 0 & 0 & 1\\
     \end{psmallmatrix}
     \begin{psmallmatrix}
        1 & 0 & 0 & 0 & 0 & 0 & 0 & 0\\
        0 & 1 & 0 & 0 & 0 & 0 & 0 & 0\\
        0 & 0 & 1 & 0 & 0 & 0 & 0 & 0\\
        0 & 0 & 0 & 1 & 0 & 0 & 0 & 0\\
        0 & 0 & 0 & 0 & 1 & 0 & 0 & 0\\
        0 & 0 & 0 & 0 & 0 & 1 & 0 & 0\\
        0 & 0 & 0 & 0 & 0 & 0 & 1 & 0\\
        0 & 0 & 0 & 0 & 0 & 0 & 0 & 1\\
     \end{psmallmatrix}
     \begin{psmallmatrix}
         1 & 0 & 0 & 0\\
         1 & 1 & 0 & 1\\
         1 & 1 & 1 & 0\\
         0 & 1 & 1 & 1\\
         0 & 1 & 1 & 1\\
         0 & 0 & 1 & 0\\
         0 & 0 & 0 & 1\\
         1 & 0 & 0 & 0\\
     \end{psmallmatrix}
     \]
     For any $\theta$, we have $\det C_{T_\theta}=0$, so $T_\theta$ is not a cellular algebra because
     the determinant of the Cartan matrix of a cellular algebra is always a
     positive integer by \cite[Proposition~1.2]{KX:CellularQH}.
   \end{Example}

We consider the smallest example of a zig-zag algebra as in \autoref{Ex:ZigZagm}.

  \begin{Example}\label{Ex:ZigZag}
    Let $\CA=\mathcal{Z}ig_2$, that is the path algebra of the quiver
       \begin{tikzpicture}[>=stealth]
         \node (1) at (0,0) {1};
         \node (2) at (1,0) {2};
         \draw[->, shorten <=1pt, shorten >=1pt]
            (1) to[bend left=20] node[midway, above] {\scriptsize $a$} (2);
         \draw[->, shorten <=1pt, shorten >=1pt]
            (2) to[bend left=20] node[midway, below] {\scriptsize $b$} (1);
       \end{tikzpicture}
    modulo the relations $aba=0=bab$. Then $\CA$ is a $6$-dimensional
    algebra with basis $\cB=\{e_1,e_2,a,b,ab,ba\}$. (We define composition
    so that $ae_1\not=0$). Set  $\Lambda=\set{1>2>3}$. Then $\CA$ has a
    $\kk$-respectable basis $\set{(\plam[i],\hlam[i],\qlam[i])|i\in\Lambda}$
    with respect to $\Lambda$, where
    \begin{align*}
      \plam[1] &= \set{ba},      & \hlam[1] &= \set{e_1}, & \qlam[1] &= \set{e_1}, \\
      \plam[2] &= \set{a, e_1},  & \hlam[2] &= \set{e_1}, & \qlam[2] &= \set{e_1, b}, \\
      \plam[3] &= \set{e_2},     & \hlam[3] &= \set{e_2}, & \qlam[3] &= \set{e_2}.
    \end{align*}
    and $\Olam=\kk$ and $\olam=\set{1_\kk}$, for $\lambda\in\Lambda$.
    Moreover,
    \[
              \Lambda_+ = \set{2,3}\qquad\text{and}\qquad
              D_\Delta =
              \begin{pmatrix}
                  1 & 0 \\
                  1 & 1 \\
                  0 & 1
              \end{pmatrix}=D_\nabla
    \]
  We leave it to the reader to check that  $\CA$ has two
  $1$-dimensional simple modules, and that the decomposition matrices of
  $\CA$ are as given above. The algebra $\CA$ is
  cellular with an anti-isomorphism that fixes~$e_1$ and~$e_2$,
  and swaps $a$ and $b$. The cellularity can also be seen, as in \cite{BrundanStroppel:KhovanovIII},  by realising $\CA$ as
  an idempotent truncation of a quasi-hereditary algebra for a
  block of parabolic category~$\CO$.
\end{Example}

\begin{Example}\label{Ex:Kronecker}
  Let $\mathcal{A}=\kk Q$ be the path algebra of the Kronecker quiver $Q$ with vertices $1,2$
  and arrows $a,b\colon1\rightrightarrows 2$, with basis $e_1,e_2,a,b$ and $a=e_2ae_1$, $b=e_2be_1$. Let
  $\Lambda=\set{1,2,3}$ with $1,2<3$ and $1,2$ incomparable. For $\lambda\in\Lambda$ set $\Olam=\kk$, and
  $\plam[i]=\hlam[i]=\qlam[i]=\set{e_i}$ for $i=1,2$,
  $\plam[3]=\set{e_2}$, $\hlam[3]=\set{a,b}$, and $\qlam[3]=\set{e_1}$,
  with $a_{phq}=phq$. This defines a perfectly $\kk$-respectable datum, since
  $z=\alpha_1e_1+\alpha_2e_2+\beta a+\gamma b$ multiplies $a$ and $b$ by
  $\alpha_2$ from the left and by $\alpha_1$ from the right, and $e_i$ from both sides by $\alpha_i$ modulo $\mathcal{A}^{>i}=\kk a\oplus\kk b$ for $i=1,2$.
  By \autoref{C:hlarge2}, the form vanishes on the arrow stratum, so $\Lambda_+=\set{1,2}$.
  Here $\Delta^3\cong L^2$ and $\nabla^3\cong L^1$, and as $H_{\mathcal{A}}=\op{diag}(1,1,2)$,
  \autoref{T:Respectable}\autoref{late6} shows that $\mathcal{A}$ has the
  non-symmetric Cartan matrix
  $C_{\mathcal{A}}=\begin{psmallmatrix}1&0\\2&1\end{psmallmatrix}$. In
  particular, $\mathcal{A}$ is not cellular by
  \cite[Theorem~3.7]{GrahamLehrer:Cellular}. Via Beilinson's equivalence \cite{Beilinson:CoherentSheaves} describing $\op{Coh}(\mathbb P^1)$,
  $\mathcal{A}\cong\End(\CO\oplus\CO(1))$ is the endomorphism algebra of a tilting bundle $\CO\oplus\CO(1)$ on $\mathbb P^1$. In
  this language, $\hlam[3]$ is a basis of $H^0(\mathbb P^1,\CO(1))=\Hom(\CO,\CO(1))$.
\end{Example}

\subsection{Respectable data compatible with automorphisms}
We want to refine the notion of respectable algebra with a fixed respectable datum in the presence of an automorphism. For this we introduce the notion of
 respectable data compatible with a given automorphism. Our prototype example will later be Frobenius-tensor algebras with the Nakayama automorphisms
 and its  corresponding CRaFH algebras. For any algebra automorphism $\psi$ of order $\theta$ we identify $\mathtt{C}_\theta=\langle \psi^m\mid m\in\Z\rangle$ with generator $\psi$. If $\mathtt{C}_\theta$ acts on a poset $\Lambda$, we denote by $\Omega_\lambda=\mathtt{C}_\theta.\lambda$ the orbit of $\lambda\in\Lambda$.

\begin{Definition}\label{respcompeasy}
 Let $A$ be respectable with
datum $(\Lambda,\ge,\plam,\hlam,\olam,\qlam,\Olam,a)$ and let $\psi$ be an
automorphism of $A$ of order $\theta$. Call the datum \emph{compatible with $\psi$} if
$\mathtt{C}_\theta$ acts on $\Lambda$ by poset automorphisms so that the following holds.
\begin{enumerate}[label=\upshape(cR\arabic*), ref=\upshape(cR\arabic*)]
  \item\label{cR1} For $\lambda\in \Lambda$ we have $\psi(A^{>\lambda})=A^{>\psi.\lambda}$ and
    $\psi(A^{\ge\lambda})=A^{\ge\psi.\lambda}$.
  \item\label{cR2} For $\mathfrak y\in\set{\mathfrak p,\mathfrak h,\mathfrak o,\mathfrak q}$
    there are bijections $\mathfrak y_\lambda\to\mathfrak y_{\psi.\lambda}$,
    $x\mapsto\psi.x$, and units $\beta^{\lambda}_{phoq}\in\kk^{\times}$ with
    \[
      \psi\bigl(a^{\lambda}_{phoq}\bigr)\equiv\beta^{\lambda}_{phoq}\,
        a^{\psi.\lambda}_{\psi.p,\,\psi.h,\,\psi.o,\,\psi.q}
      \pmod{A^{>\psi.\lambda}}.
    \]
  \item\label{cR3} There are units $\beta_o$, $o\in\olam$ such that
    $o\mapsto\beta_o(\psi.o)$ is an algebra isomorphism
    $\Olam\to\Olam[\psi.\lambda]$.
 \item\label{cR4} The equality
    $\beta^{\lambda}_{phoq}/\beta^{\lambda}_{pho'q}=\beta_o/\beta_{o'}$ holds for all $\lambda\in\Lambda$,
    $o,o'\in\olam$.
  \item\label{cR5} $\psi^{|\Omega_\lambda|}.h=h$ for every $h\in\hlam$.
\end{enumerate}
\end{Definition}
\begin{Example}
An analogous definition can be given if  $\psi$ is an \emph{anti}-automorphism of $\CA$, with bijections $\plam\to\qlam[\psi.\lambda]$ and $\qlam\to\plam[\psi.\lambda]$ in \autoref{cR2}.  In particular $\psi=\ast$ might be an anti-involution  of $A$.
Recall that every cellular algebra is respectable by \autoref{Ex:Cellular}. Then every \emph{skew cellular algebra} with anti-involution $\ast$ in  the sense of \cite{MathasRostam} is in the same way respectable, and the respectable datum is compatible with the anti-involution $\ast$ in the sense of \autoref {respcompeasy} with all scalars being  $1$. The  \emph{graded} skew cellular algebras from \cite{MathasRostam}, thus Hecke algebras of type $G(l,p,n)$, arise via  \autoref {respcompeasy}  if we work with algebra objects $A$, equipped with an anti-involution, in the category of graded $\kk$-modules.
\end{Example}

\begin{Remark}
For $\kk$-respectable
data, \autoref{cR3}, \autoref{cR4} are vacuous. For a cellular algebra as in \autoref{Ex:Cellular}, \autoref{respcompeasy} says exactly
that $\psi$ preserves the cell filtration and permutes the cellular basis up to
units, i.e.,\ $\psi(c_{st})\equiv\beta_{st}c_{\psi.s,\psi.t}$ modulo
$A^{>\psi.\lambda}$.
\end{Remark}

\begin{Remark}\label{cardh1}
  Note the quite different meaning of similar notation: for  $x\in\mathfrak y_\lambda$ with
  $\mathfrak y\in\set{\mathfrak p,\mathfrak q,\mathfrak h,\mathfrak o}$ we have
  $\psi.x\in\mathfrak y_{\psi.\lambda}$. On the other hand, $\psi(\dot x)\in A$.  Note that \autoref{pcR3} below is automatic if $|\hlam|=1$ for all $\lambda$.
\end{Remark}
We give a refinement of \autoref{respcompeasy} for perfectly respectable data.

\begin{Definition}\label{respcompatible}
Let $A$ be a perfectly respectable algebra with an automorphism $\psi$ of finite order $\theta$. Fix a respectable datum
$(\Lambda, \ge, \mathfrak p_\lambda,\mathfrak h_\lambda,\mathfrak o_\lambda,\mathfrak q_\lambda, \Olam,a)$ of $A$ with $a^\lambda_{phoq}=\dot p\dot h\dot o\dot q$.   This datum is \emph{compatible with $\psi$} if  there exists an action,  which we denote by
$m\mapsto(\lambda\mapsto\psi^{m}.\lambda)$,  of $\mathtt{C}_\theta$ on $\Lambda$
by poset automorphisms satisfying the following.
  \begin{enumerate}[label=\upshape(pcR\arabic*), ref=\upshape(pcR\arabic*)]
    \item\label{pcR1} For each $\lambda\in\Lambda$ we have
      $\psi(A^{>\lambda})=A^{>\psi.\lambda}$ and
      $\psi(A^{\ge\lambda})=A^{\ge\psi.\lambda}$.
    \item\label{pcR2} For $\mathfrak y\in\set{\mathfrak p,\mathfrak h,\mathfrak o,\mathfrak q}$
      and $\lambda\in\Lambda$ there exist  bijections $\mathfrak y_\lambda\to\mathfrak y_{\psi.\lambda}$,
          $x\mapsto\psi.x$, and units $\beta_x\in\kk^{\times}$, for
          $x\in\plam\cup\hlam\cup\olam\cup\qlam$, with      \[
            \psi(\dot x)\equiv\beta_x\dot{(\psi.x)}\pmod{A^{>\psi.\lambda}},
          \]
      and such that $o\mapsto\beta_o\,(\psi.o)$ is an algebra isomorphism
      $\Olam\to\Olam[\psi.\lambda]$.
    \item \label{pcR3} For $\lambda\in\Lambda$ there exists $\epsilon_\lambda\in\kk^{\times}$
      such that
      \[
        \psi^{|\Omega_\lambda|}(\dot h)\equiv\epsilon_\lambda\dot h\pmod{A^{>\lambda}}
        \qquad\text{for every }h\in\hlam,
      \]
      where $\Omega_\lambda=\mathtt{C}_\theta.\lambda$ is the orbit of $\lambda$.
  \end{enumerate}
  The respectable datum together with the action of $\mathtt{C}_\theta$, the bijections and the scalars $\beta_x$, $\epsilon_\lambda$ form a \emph{respectable datum compatible with the automorphism $\psi$}.
\end{Definition}

\begin{Remark}\label{seeminglyeasier}
For a perfect datum, \autoref{respcompeasy} is implied by
\autoref{respcompatible}. Indeed, multiplying
the four congruences $\psi(\dot x)\equiv\beta_x\dot{(\psi.x)}$ from \autoref{pcR2} gives \autoref{cR2} with
$\beta^{\lambda}_{phoq}=\beta_p\beta_h\beta_o\beta_q$ because $A^{>\psi.\lambda}$ is an ideal. The conditions
\autoref{pcR2}, \autoref{pcR3} imply
\autoref{cR3}, \autoref{cR4}, \autoref{cR5}. Namely,  the second half of
\autoref{pcR2} is \autoref{cR3}, condition \autoref{cR4} follows from the
factorisation just displayed, and for \autoref{cR5} one compares
$\psi^{|\Omega_\lambda|}(\dot h)\equiv\epsilon_\lambda\dot h$ with the
congruence $\psi^{|\Omega_\lambda|}(\dot h)\equiv\beta\,
\dot{(\psi^{|\Omega_\lambda|}.h)}$ obtained by iterating \autoref{pcR2},
multiplies by $\dot p$ on the left and $\dot o\dot q$ on the right, and uses
that the $a^{\lambda}_{ph'oq}$, $h'\in\hlam$, are linearly independent modulo
$A^{>\lambda}$. \autoref{respcompeasy} is in general weaker than
\autoref{respcompatible}, since \autoref{cR2} constrains $\psi$ only on the
strata, while \autoref{pcR2} constrains the lifts.

To see this let $A=\Mat_2(\kk)$ with the perfect datum from
\autoref{prototype}. Thus,
$\Lambda=\set{\bullet}$, $\plam[\bullet]=\set{E_{11},E_{21}}$,
$\hlam[\bullet]=\set{E_{11}}$, $\qlam[\bullet]=\set{E_{11},E_{12}}$,
$\olam[\bullet]=\set{1}$, $\Olam[\bullet]=\kk$, $a_{pq}=E_{pq}$.
Let $\psi$ be conjugation by $\begin{psmallmatrix}0&1\\1&0\end{psmallmatrix}$. Then
$\psi(E_{pq})=E_{\bar p\bar q}$ with $\bar1=2$, $\bar2=1$. Then
\autoref{respcompeasy} holds with all units equal to $1$ and the bijections
$p\mapsto\bar p$, $q\mapsto\bar q$. For this datum
\autoref{pcR2} fails, since it would require
$\psi(\dot h)=E_{22}$ to be a multiple of $\dot h=E_{11}$ modulo $A^{>\bullet}=0$.
\end{Remark}

For $\lambda\in\Lambda$,  the $\mathtt{C}_\theta$-stabiliser of $\lambda$  is the subgroup
  $\langle\psi^{|\Omega_\lambda|}\rangle\subseteq\mathtt{C}_\theta$ of order
  $\theta_{\Omega_\lambda}=\frac{\theta}{|\Omega_\lambda|}$. For an orbit $\Omega$ and a fixed $\lambda\in\Omega$, the map and
  $\sigma_\Omega\colon\Olam\to\Olam$, $o\mapsto\bigl(\prod_{k=0}^{|\Omega_\lambda|-1}\beta_{\psi^{k}.o}\bigr)(\psi^{|\Omega_\lambda|}.o)$, obtained by
  iterating \autoref{pcR2}, is an automorphism of $\Olam$.

Respectable data are usually not compatible with $\psi$. We illustrate this in an example.
\begin{Example}\label{Ex:T2compatible}
  Let $\theta=2$, so $\omega=-1$, and let $2\in\kk^{\times}$.  Then
  $\Taft[2]=\<\,g,y\>/(g^{2}=1, y^{2}=0, yg=\omega gy)$ with $\psi(g)=-g$
  and $\psi(y)=y$.  The respectable datum of \autoref{Ex:Taftresp} has
  $\Lambda=\set{0<1<2<3}$, $h_i=g-\omega^{i+1}$, and basis $a_0=g+1$,
  $a_1=g-1$, $a_2=(g+1)y$, $a_3=(g-1)y$. This respectable datum  is not compatible with $\psi$, because the chain $\Lambda$ has only the trivial poset automorphism,
  so that condition \autoref{pcR1} of \autoref{respcompatible}, that is
  \autoref{cR1}, requires
  $\psi(A^{>i})=A^{>i}$ for all $i$. But $\psi$ permutes the strata (up to
  sign). It swaps for instance the $0$- and the $1$-stratum, but
  $\psi(a_1)=-a_0\notin A^{>0}$.

We obtain however a respectable datum which is compatible with $\psi$ by changing $\Lambda$ to encode the $y$-degree. Namely, setting
\[\Lambda=\set{0<1}, \quad \plam[j]=\set{y^{j}},\quad
  \hlam[j]=\qlam[j]=\set{1},\quad \olam[j]=\set{1,g},\quad \Olam[j]=\kk[g]/(g^{2}-1)\quad\text{with}\quad a_{g^{i}}^{j}=y^{j}g^{i}
  \]
   defines a perfectly respectable
  datum, since $g\cdot y^{j}g^{i}=(-1)^{j}y^{j}g^{i+1}$ and
  $y^{j}g^{i}\cdot g=y^{j}g^{i+1}$, while $y$ maps $A^{(0)}$ into
  $A^{>0}=A^{(1)}$ and kills $A^{(1)}$ from both sides.  The datum is
  compatible with $\psi$ for the trivial action on $\Lambda$. Indeed,
  \autoref{cR1} holds since $\psi(A^{>0})=A^{>0}$.  \autoref{cR2} holds, since $\psi.x=x$ for all indices with $\beta_{y^{j}}=\beta_{1}=1$ and $\beta_{g}=-1$, and \autoref{cR3} is satisfied since $g\mapsto-g$ is an automorphism of $\Olam[j]$.
 Condition \autoref{cR4} holds because
  the datum is perfect, and \autoref{cR5} is empty by
  \autoref{cardh1}.  Here, $\theta_\Omega=2$ and
  $\sigma_\Omega(g)=-g$.  Thus  $\sigma_\Omega$  interchanges the idempotents
  $\tfrac12(1\pm g)$ of $\Olam[j]\cong\kk\times\kk$. The whole of $\psi$ sits
  inside the fibre algebras preventing them from being  $\kk$.
\end{Example}

\section{Respectable bases on $\Hhd[1](A)$}
In this section we construct respectable bases for CRaFH algebras for $n=1$.
\subsection{Ideals in $A$ versus ideals in the RaFH factors $\Hhd[1](A)$ }
Throughout this section we assume the following.
\begin{Setup}\label{setupresp}
 Let  $A$ be a superalgebra of finite rank over $\kk$ with basis $\cB$ and
$\psi\colon A\to A$ an automorphism of finite order $\theta$. Fix a $\theta$-cyclotomic polynomial of degree $\ell=\ell_0+\theta\ell_1$ with level vector $(\ell_0, \ell_1)$:
\begin{equation*}
  f(X)=X^{\ell_0}\prod_{i=\ell_0+1}^{\ell_0+\ell_1}\bigl(X^{\theta}-\kapp_i^{\theta}\bigr)
      =\prod_{k=1}^{\ell'}f_k(X),
  \quad\text{where }\quad\ell'=\ell_0+\ell_1,\quad
  f_k(X)=\begin{cases}X, & 1\le k\le\ell_0,\\ X^{\theta}-\kapp_k^{\theta}, & \ell_0<k\le\ell'.\end{cases}
\end{equation*}
Set $\kapp_i=0$ for $1\le i\le\ell_0$ and $f_0:=1$. We consider  the CRaFH-algebra $\Hhd[1](A)=\kk[X]\Smash A/(f(X))$, where $X$ is even and $aX=X\psi(a)$ for $a\in A$, and construct respectable bases.
\end{Setup}

\begin{Definition}
Corresponding  to the factors $f_k$ for $1\le k\le\ell'$ we have the \emph{RaFH factors}
  \[
    \HH^{\kapp_k}_1(A)=\kk[X]\Smash A/(f_k(X))=\begin{cases}
      A, & \text{if }1\le k\le\ell_0,\\[2pt]
      \kk[X]\Smash A\big/\bigl(X^{\theta}-\kapp_k^{\theta}\bigr), & \text{if }\ell_0<k\le\ell',
    \end{cases}
  \]
  where the smash product is given by $aX=X\psi(a)$ for $a\in A$.
\end{Definition}

Given a 2-sided ideal $J$ in $A$ let $J^{\psi}=\sum_{m}\psi^{m}(J)$ be the smallest
$\psi$-stable ideal containing $J$.

\begin{Lemma}\label{bijideals}
  Let $\kapp_k\in\kk^\times$ and $\HH^{\kapp_k}_1(A)=\kk[X]\Smash A/(X^{\theta}-\kapp_k^{\theta})$. There is an inclusion preserving  bijection
  \begin{align*}
\Phi\colon  \set{ \text{$\psi$-stable two-sided ideal of $A$}}&\cong\set{\text{two-sided ideals of
      $\HH^{\kapp_k}_1(A)$}}\\
      J&\mapsto J\HH^{\kapp_k}_1(A)
      \end{align*}
 For an arbitrary two-sided ideal $J$ of $A$ we have  $\HH^{\kapp_k}_1(A)J\HH^{\kapp_k}_1(A)=J^{\psi}\HH^{\kapp_k}_1(A)$.
\end{Lemma}

\begin{proof}
For any ideal $J$ of $A$ we have $AJ=JA=J$, hence
  $\HH^{\kapp_k}_1(A) J=\sum_jX^{j}AJ=\bigoplus_jX^{j}J$.
  On the other hand,   $J\HH^{\kapp_k}_1(A)=\sum_jJX^{j}A=\bigoplus_jX^{j}\psi^{j}(J)$.
  If $J$ is $\psi$-stable the two agree. Thus, $\Phi$ is well-defined.

The map $\Phi$ obviously preserves inclusions. By \autoref{nequal1}, $\HH^{\kapp_k}_1(A)$ is a free  right
  $A$-module with basis $X^{0},\dots,X^{\theta-1}$.  Since $X^{\theta}=\kapp^{\theta}$ is
  invertible, so is $X$, with $X^{-1}=\kapp^{-\theta}X^{\theta-1}$. Moreover,
   $X^{-j}bX^{j}=\psi^{j}(b)$ for $b\in A$ and
  $j\in\Z$. So conjugation by $X$ is an automorphism of $\HH^{\kapp_k}_1(A)$ restricting
  to $\psi$ on $A$.  Now  any two-sided ideal in $\HH^{\kapp_k}_1(A)$ is of the form
  $\bigoplus_jX^{j}J$ for some ideal $J$ in $A$. It contains $X^{-1}JX=\psi(J)$, and
  comparing the components in $X^{0}A$ gives $\psi(J)\subseteq J$, so
  $J\supseteq\psi(J)\supseteq\dots\supseteq\psi^{\theta}(J)=J$ and $J$ is
  $\psi$-stable. Thus, $\Phi$ is surjective.

  Comparing components once more gives $J\HH^{\kapp_k}_1(A)\cap A=J$, so
  $J\mapsto J\HH^{\kapp_k}_1$ and thus $\Phi$ is injective.

  Let  $I$ the two-sided ideal of $\HH^{\kapp_k}_1(A)$ generated by a 2-sided ideal
  $J$ of $A$. Then $I$ contains $X^{-j}JX^{j}=\psi^{j}(J)$ for every $j$, and therefore  $J^{\psi}\HH^{\kapp_k}_1(A)\subseteq I$. On the other hand,  $J^{\psi}$ is $\psi$-stable and
  contains $J$, so $J^{\psi}\HH^{\kapp_k}_1(A)$ is a two-sided ideal containing $J$. But then $I=J^{\psi}\HH^{\kapp_k}_1(A)$.
  \end{proof}

\begin{Lemma}\label{L:nilp}
  If $J$ is a two-sided ideal of $A$ and $1\le k\le\ell_0$, then
  $X^{k-1}J+X^{k}\Hhd[1](A)$ is a two-sided ideal of $\Hhd[1](A)$.
\end{Lemma}

\begin{proof}
We have  $b X^{k-1}=X^{k-1}\psi^{k-1}(b)$ and
  $X^{k-1}bX=X^{k}\psi(b)$ for $b\in A$ and $k\in\N$, and
  $X^{k}\Hhd[1](A)=u_{k+1}\Hhd[1](A)$ is a two-sided ideal by
  \autoref{L:1}\autoref{H1}. Since $A$ and $X$ generate $\Hhd[1](A)$, the
  claim follows.
\end{proof}

\subsection{The RaFH factors respectable basis} We give a first construction of a respectable datum:
\begin{Theorem}\label{T:1}
Assume \autoref{setupresp}.  Then the  superalgebra $\Hhd[1](A)$ is respectable with respectable datum
 \[(\Lambda,\le)=\set{1<2<\dots<\ell'},\quad\plam=\qlam=\hlam=\set{1},\quad\Olam=\HH^{\kapp_\lambda}_1(A),\]
       \[ \olam=\begin{cases}
          \set{b\mid b\in\cB} & \text{if }1\le\lambda\le\ell_0,\\
          \set{X^{j}b\mid b\in\cB,\ 0\le j\le\theta-1} & \text{if }\ell_0<\lambda\le\ell',
        \end{cases}
\]
\begin{equation}\label{E:a}
 a\colon\bigcup_{\lambda\in\Lambda}\olam\to\Hhd[1](A),\quad
        \olam\ni X^{j}b\mapsto a^{\lambda}_{X^{j}b}=u_\lambda X^{j}b,      \qquad\text{where }u_\lambda=\prod_{i=1}^{\lambda}f_{i-1}.
\end{equation}
It is moreover perfectly respectable with $\dot p=\dot h=\dot q=1$ and $\dot o=a^{\lambda}_{X^{j}b}=u_\lambda X^{j}b$.
\end{Theorem}

\begin{Remark}\label{uextrem}
  Note that by definition $u_1=f_0=1$ and $u_{\ell'+1}=f$.
\end{Remark}
\begin{Remark}
To illustrate the definitions let  $(\ell_0,\ell_1)=(2,3)$.  Then
  \[
    u_1=1,\quad u_2=X,\quad u_3=X^{2},\quad
    u_4=X^{2}(X^{\theta}-\kapp_3^{\theta}),\quad
    u_5=X^{2}(X^{\theta}-\kapp_3^{\theta})(X^{\theta}-\kapp_4^{\theta}).
  \]
  Elements in $a(\olam)$ are all of the form
  $X^{j}b+\sum_{i<j,\ b'\in\cB}\alpha_{ib'}X^{i}b'$ for some $b\in\cB$ with  \emph{leading term}  $\lt(a(\olam))=X^{j}b$.
  \begin{center}
    \begin{tabular}{l|c|c|c|c|c}
      $\lambda$ & $1$ & $2$ & $3$ & $4$ & $5$\\\hline
      $\lt(a(\olam))$ & $1\cdot b$ & $Xb$ & $X^{2}b,\dots,X^{\theta+1}b$
        & $X^{\theta+2}b,\dots,X^{2\theta+1}b$ & $X^{2\theta+2}b,\dots,X^{3\theta+1}b$
    \end{tabular}\;.
  \end{center}
  \end{Remark}
  \begin{Remark}Complications arise for non-trivial $\psi$ since the  $u_\lambda$ from \autoref{E:a} are not central but satisfy
 \begin{equation}\label{E:orbitpoly}
   u_\lambda X=Xu_\lambda, \quad\text{and}\quad b\,u_\lambda=u_\lambda\,\psi^{d}(b)\qquad\text{with }d=\min(\lambda-1,\ell_0)
\end{equation}
for $b\in A$ and $1\le\lambda\le\ell'$, because $u_\lambda=X^{\min(\lambda-1,\ell_0)}p(X^{\theta})$ for some polynomial $p$ in the central element
$X^{\theta}$.
  \end{Remark}

For the proof of \autoref{T:1} we prove two lemmas assuming the  setup from  \autoref{T:1}.

\begin{Lemma}\label{L:1}
  For $k\in\N$ define the right ideals $H^{\ge k}=u_k\Hhd[1](A)$ in $\Hhd[1](A)$.  They
  satisfy:
  \begin{enumerate}
    \item \label{H1} $H^{\ge k}$ is a two-sided ideal of $H$ for all $k$.
    \item \label{H2}  $H=H^{\ge1}\supsetneq H^{\ge2}\supsetneq\dots\supsetneq H^{\ge\ell'+1}=\set{0}$.
    \item  \label{H3} The map $a$ in \autoref{E:a} is injective. The
      $\bigcup_{\lambda\ge k}a(\olam)$, for $\lambda\in\Lambda=\set{1,\dots,\ell'}$,
      form a basis of $H^{\ge k}$.
  \end{enumerate}
\end{Lemma}

\begin{proof}
  \autoref{H1} Note that $Xu_k=u_kX$ and $bu_k=u_k\psi^{i}(b)$ with
  $i=\min(k-1,\ell_0)$.
  Hence the right ideal $u_k\Hhd[1](A)$ is stable under left multiplication by $X$
  and by $b\in A$, and so it is two-sided.

  \autoref{H2} By \autoref{uextrem}, $H^{\ge1}=1\cdot H=H$ and
  $H^{\ge\ell'+1}=u_{\ell'+1}H=fH=0$, and obviously
  $H^{\ge k}\supseteq H^{\ge k+1}$.  Since $\olam\ne\emptyset$ for all
  $\lambda$, the strict inequality follows from   \autoref{H3}.

  \autoref{H3} For $\lambda\in\Lambda$, the leading terms $\lt(a(\olam))$ of the elements
  in $a(\olam)$ are
  \[
    X^{\lambda-1}b,\ b\in\cB\quad\text{if }\lambda\le\ell_0,
    \qquad
    X^{\theta(s-1)+\ell_0+j}b,\ b\in\cB,\ 0\le j<\theta
    \quad\text{if }\ell_0<\lambda=\ell_0+s\le\ell' .
  \]
  Since for $\ell_0+s=\ell'$ we have $ \theta(s-1)+\ell_0+\theta-1=\theta\ell_1-\theta+\ell_0+\theta-1
    =\theta\ell_1+\ell_0-1=\ell-1$,
  these leading terms are precisely the elements $X^{d}b$ with $0\le d<\ell$
  and $b\in\cB$, each occurring exactly once, and these form a basis of $\Hhd[1](A)$ by
  \autoref{nequal1}. Since every element of $a(\olam)$ is its leading term
  plus terms of lower $X$-degree, the family $\bigcup_{\lambda\in\Lambda}a(\olam)$
 is a basis of $\Hhd[1](A)$. In particular, $a$ is injective.

  We claim that $H^{\ge k}$ is spanned by the $a(\olam)$ for $\lambda\ge k$.
  For this it suffices to show that $H^{\ge k}$ is spanned by $u_kX^{j}b$
  with $j\le\ell-k$ if $k\le\ell_0$, and with $j<\theta(\ell'-k+1)$ if
  $k>\ell_0$.  In the first case note that $u_kX^{m}$ can be rewritten, using
  $f$, as a linear combination of $X^{s}$ with $\ell-1\ge s\ge\ell_0\ge k$,
  hence of $u_kX^{j}=X^{j+k-1}$ with $j+k\le\ell$.  In the second case
  consider $u_kX^{m}$ with $m\ge\theta(\ell'-k+1)$, and note that
  $g:=\prod_{i=k}^{\ell'}f_i$ is monic of degree $\theta(\ell'-k+1)$.  By
  definition, $u_kg=f=0$,
  and so $u_kX^{m}$ can be written in terms of lower $X$-degrees.

The listed elements $u_kX^{j}b$
have the same leading terms $X^{d}b$,
  $\deg u_k\le d<\ell$, $b\in\cB$, as the elements of
  $\bigcup_{\lambda\ge k}a(\olam)$. The latter lie in $H^{\ge k}$ because
  $u_\lambda\in u_kH$ for $\lambda\ge k$ and each occurs once in each family. If $u$
  and $v$ are elements of the two families with the same leading term, then
  $u-v\in H^{\ge k}$ has smaller $X$-degree, and expanding $u-v$ in the first
  family, which spans $H^{\ge k}$, shows that it is a combination of elements
  of smaller $X$-degree. Induction on the $X$-degree now shows that both
  families have the same span. Hence $\bigcup_{\lambda\ge k}a(\olam)$ is a
  basis of $H^{\ge k}$.
\end{proof}

\begin{Lemma}\label{L:2}
  Using the setup from \autoref{L:1} let $\lambda\in\Lambda$ and let $\pi\colon H^{\ge\lambda}\to H^{\ge\lambda}/H^{\ge\lambda+1}$
  be the quotient map.  Then
  \[
    \Phi_\lambda\colon\Olam=\HH^{\kapp_\lambda}_1(A)\longrightarrow H^{\ge\lambda}/H^{\ge\lambda+1},
    \qquad X^{j}b\longmapsto\pi\bigl(a^{\lambda}_{X^{j}b}\bigr)=\pi\bigl(u_\lambda X^{j}b\bigr),
  \]
  is an isomorphism of $\kk$-modules. For  $o,o'\in\olam$ with
  $o=X^{j}b$ and $d=\min(\lambda-1,\ell_0)$,
  \begin{equation}\label{E:bimod}
    \Phi_\lambda(oo')=X^{j}\psi^{-d}(b)\cdot\Phi_\lambda(o')=\Phi_\lambda(o)\cdot X^{j}b .
  \end{equation}
  That is, the left and right action of $o$ on the stratum $H^{(\lambda)}$ is the left and  right
  multiplication by $X^{j}\psi^{-d}(b)$ and  $X^{j}b$  respectively. Both actions realise the multiplication of $\Olam$.
\end{Lemma}

\begin{proof}
  By \autoref{L:1}\autoref{H3}, $H^{\ge\lambda}$ has basis
  $\bigcup_{\mu\ge\lambda}a(\olam[\mu])$ and $H^{\ge\lambda+1}$ has basis
  $\bigcup_{\mu\ge\lambda+1}a(\olam[\mu])$, so that $\pi$ maps $a(\olam)$ onto a
  basis of the quotient. Therefore,  $\Phi_\lambda$ is bijective.  For the second
  equality of \autoref{E:bimod} let $o=X^{j}b$ and $o'=X^{j'}b'$.  From \autoref{E:orbitpoly} we get
  $
    a^{\lambda}_{o}\cdot X^{j'}b'=u_\lambda X^{j}b\,X^{j'}b'=u_\lambda X^{j+j'}\psi^{j'}(b)b' .
  $
  In $\Olam$ we have $oo'=X^{j+j'}\psi^{j'}(b)b'$ modulo $X$ if
  $\lambda\le\ell_0$ and modulo $X^{\theta}=\kapp_\lambda^{\theta}$ if
  $\lambda>\ell_0$.  The two agree modulo $H^{\ge\lambda+1}$. For
  $\lambda\le\ell_0$ because $u_\lambda X=u_{\lambda+1}\in H^{\ge\lambda+1}$,
  and for $\lambda>\ell_0$ because
  $u_\lambda(X^{\theta}-\kapp_\lambda^{\theta})=u_\lambda f_\lambda=u_{\lambda+1}\in H^{\ge\lambda+1}$,
  that is $u_\lambda X^{\theta}\equiv\kapp_\lambda^{\theta}u_\lambda$, which is
  exactly the relation of $\Olam$.  The left hand equality is the same
  computation, the twist $\psi^{-d}$ being forced by
  $\psi^{-d}(b)u_\lambda=u_\lambda b$, which holds because
  $u_\lambda$ is $X^{d}$ times a polynomial in $X^{\theta}$ and $X^{\theta}$
  commutes with $A$.
\end{proof}

\begin{proof}[Proof of \autoref{T:1}]
  \autoref{resp1} is
  \autoref{L:1}\autoref{H3} for $k=1$, which also identifies
  $A^{(\lambda)}$ with the span of $a(\olam)$ and $A^{>\lambda}$ with
  $H^{\ge\lambda+1}$.  For \autoref{resp2} note that $\plam$,
  $\hlam$ and $\qlam$ are singletons, so the independence of the scalars
  from $h$ and $q$ is vacuous and the condition says that $A^{(\lambda)}+A^{>\lambda}=H^{\ge\lambda}$
  is stable under left and right multiplication by $H$,
  which is \autoref{L:1}\autoref{H1}.  For \autoref{resp3},
  $1_{\Olam}=X^{0}1\in\olam$ by construction, and by \autoref{L:2} the map
  $\Phi_\lambda$ carries the regular $\Olam$-bimodule onto $A^{(\lambda)}$. Hereby,
  the two actions are multiplication in $H$ by
  $X^{j}\psi^{-d}(b)$ on the left and by $X^{j}b$ on the right.  The required
  congruences $(ba_{o})o'\equiv b(a_{o}o')$ and $(o'a_{o})b\equiv o'(a_{o}b)$
  modulo $A^{>\lambda}$ are therefore instances of the associativity of
  $H$, since left multiplication by one element commutes
  with right multiplication by another.  Finally the basis is perfect with
  $\dot p=\dot h=\dot q=1$ and $\dot o=a^{\lambda}_{o}$.
\end{proof}
\subsection{Orbit space respectable bases}
The datum of \autoref{T:1} is independent of the internal structure of  $A$. All of $A$ sits inside the RaFH factors $\Olam=\HH^{\kapp_\lambda}_1(A)$.
We next construct under the assumption that $\kapp_{\ell_0+1},\dots,\kapp_{\ell'}\in\kk^\times$  a refinement from a given respectable datum for $A$ which is compatible with $\psi$ in the sense of  \autoref{respcompatible}

For the construction we assume \autoref{setupresp} and  fix a perfectly respectable datum
$(\Lambda,\ge,\plam,\hlam,\olam,\qlam,a)$ of $A$ with
$a^\lambda_{phoq}=\dot p\dot h\dot o\dot q$, and we assume that it is
compatible with $\psi$, with action $\lambda\mapsto\psi.\lambda$,
bijections $x\mapsto\psi.x$ and units $\beta_x$, $\epsilon_\lambda$ as
in \autoref{respcompatible}. We write $\dot{(\psi.x)}$ for the lift of
the index $\psi.x$, so that $\psi(\dot x)\equiv\beta_x\dot{(\psi.x)}$
modulo $A^{>\psi.\lambda}$.

We need to understand the orbits of the different $\psi$-actions. The following uses the  $\mathtt{C}_\theta$-action on $\Lambda$:
\begin{Lemma}
The set $\Lambda_\psi$ of orbits for the  $\mathtt{C}_\theta$-action on $\Lambda$ inherits a  poset structure on $\Lambda$ by setting
\[
  \Omega\preceq\Omega'\quad\Longleftrightarrow\quad
  \lambda\le\psi^{m}.\lambda'\ \text{ for some }\lambda\in\Omega,\
  \lambda'\in\Omega'\text{ and some }m\in\Z .
\]
The projection map
$\Lambda\to\Lambda_\psi$ is order preserving.
\end{Lemma}

\begin{proof}
The condition is well-defined since the $\psi^{m}$ are order automorphisms, since
$\lambda,\lambda'$ can be replaced by other representatives of the orbit by changing possibly $m$.
 The relation is obviously reflexive and transitive because
$\lambda\le\psi^{m}.\lambda'$ and $\lambda'\le\psi^{m'}.\lambda''$ give
$\lambda\le\psi^{m}.\lambda'\le\psi^{m+m'}.\lambda''$. For antisymmetry suppose that $\Omega\preceq\Omega'$, $\Omega'\preceq\Omega$. Let $\lambda\in\Omega$ and pick $\mu'\in\Omega'$
with $\lambda\le\mu'$ and $\mu\in\Omega$ with $\mu'\le\mu$. If now $\Omega\ne\Omega'$ then $\mu,\mu'\not=\lambda$, since they are in different orbits and thus $\lambda\lneq\mu'\lneq\mu=\psi^n.\lambda$ for some $n$. Acting by the order automorphism $\psi^{n}$
repeatedly implies $\lambda\lneq\psi^{n}.\lambda\lneq\psi^{2n}.\lambda\lneq\dots\lneq\psi^{\theta n}.\lambda=\lambda$ which is a contradiction.
Thus, $(\Lambda_\psi,\preceq)$ is a poset. The same argument with
$\Omega=\Omega'$ shows that $\lambda<\psi^{n}.\lambda$ for any $n$. That means that
two distinct elements in the same orbit are never comparable in $\Lambda$.
Now $\mu>\lambda$ implies $[\mu]\succ[\lambda]$, and the projection
$\Lambda\to\Lambda_\psi$ is order preserving.
\end{proof}

We denote by  $\Lambda_\psi$ the set of orbits for the  $\mathtt{C}_\theta$-action on $\Lambda$.
For each  orbit $\Omega\in \Lambda_\psi$ we fix a representative $\lambda_\Omega\in\Omega$. Recall that
$r_\Omega=|\Omega|$, $\theta_\Omega=\frac{\theta}{r_\Omega}$ and the
automorphism $\sigma_\Omega=\psi^{r_\Omega}|_{\Olam[\lambda_\Omega]}$ of
$\Olam[\lambda_\Omega]$, which is $o\mapsto\prod_{k=0}^{r_\Omega-1}\beta_{\psi^{k}.o}\,(\psi^{r_\Omega}.o)$
on the basis $\olam[\lambda_\Omega]$ by iterating \autoref{respcompatible}~\autoref{pcR2}.
\begin{Definition}
Assume  \autoref{setupresp}.  For $\ell_0<k\le\ell'$ and an orbit $\Omega\in  \Lambda_\psi$ we define the \emph{orbit superalgebra}
\begin{equation}\label{E:orbitD}
  \mathcal{D}^{k}_\Omega
   =\kk[W]\Smash  \Olam [\lambda_\Omega]/\left(W^{\theta_\Omega}-\kapp_k^{\theta}\right),
   \qquad\text{with smash product relation}\qquad oW=W\sigma_\Omega(o).
\end{equation}
 Since $\kapp_k^{\theta}=(\kapp_k^{r_\Omega})^{\theta_\Omega}$, the superalgebra $\mathcal{D}^{k}_\Omega$ can be viewed as
CRaFH superalgebra $\HH^{(\kapp_k^{r_\Omega})}_1(\Olam[\lambda_\Omega])$ for $n=1$ with Nakayama
automorphism $\sigma_\Omega$ of order $\theta_\Omega$ and $\theta_\Omega$-cyclotomic polynomial $X^{\theta_\Omega}-\kapp_k^\theta$.
\end{Definition}
\begin{Remark}
 If $\psi$ is trivial
then every orbit is a point, $r_\Omega=1$, $\theta_\Omega=1$ and
$\mathcal{D}^{k}_\Omega=\Olam$.
\end{Remark}
We obtain a refinement of \autoref{T:1} using the orbit superalgebras.
\begin{Theorem}\label{T:orbit}Assume \autoref{setupresp}. Assume $\kapp_k$ is invertible for $\ell_0<k\le\ell'$.
  Let $(\Lambda,\ge,\plam,\hlam,\olam,\qlam,\Olam,a)$ be a perfectly respectable datum of $A$ compatible with $\psi$. Define the following datum
    \begin{enumerate}
    \item\label{orbit1} Let  $\Gamma=\bigl(\set{1,\dots,\ell_0}\times\Lambda\bigr)
          \sqcup\bigl(\set{\ell_0+1,\dots,\ell'}\times\Lambda_\psi\bigr)$
       with the lexicographic ordering, that is
      \[(k,\nu)<(k',\nu')\Longleftrightarrow \text{$k<k'$ or ($k=k'$ and $\nu<\nu'$ in $\Lambda$ in case
      $k\le\ell_0$ and $\nu\preceq\nu'$ in case $k>\ell_0$)}.\]
    \item\label{orbit2}  For $(k,\lambda)$ with $k\le\ell_0$ let $(\plam[k,\lambda],\hlam[k,\lambda],\olam[k,\lambda],\qlam[k,\lambda],\Olam[k,\lambda])=(\plam,\hlam,\olam,\qlam,\Olam)$. Then take $\dot h, \dot q, \dot o$ as for $A$ but define $\dot p$ for $p\in\plam[k,\lambda]$ as $u_k\dot p=X^{k-1}\dot p$ with now $p\in\plam$, and set $a^{(k,\lambda)}_{phoq}=X^{k-1}a^{\lambda}_{phoq}$.
    \item\label{orbit3} For $(k,\Omega)$ with $k>\ell_0$ and $\Omega=\Omega_\lambda$ let
      \begin{gather*}
        \plam[(k,\Omega)]=\set{0,\dots,r_\Omega-1}\times\plam,\qquad
        \hlam[(k,\Omega)]=\hlam,\qquad
        \qlam[(k,\Omega)]=\qlam\times\set{0,\dots,r_\Omega-1},\\
        \olam[(k,\Omega)]=\olam\times\set{0,\dots,\theta_\Omega-1},\qquad
        \Olam[(k,\Omega)]=\mathcal{D}^{k}_\Omega
      \end{gather*}
      with lifts $(j,p)\mapsto u_kX^{j}\dot p$, $h\mapsto\dot h$, $(q,j')\mapsto\dot qX^{j'}$ and
      $(o,d)\mapsto\dot oX^{r_\Omega d}$ respectively. Moreover, let
      \begin{equation}\label{E:orbitE}
        a^{(k,\Omega)}_{(j,p),h,(o,d),(q,j')}
          =u_k\,X^{j}\,\dot p\,\dot h\,\dot o\,X^{r_\Omega d}\,\dot q\,X^{j'},
        \qquad 0\le j,j'<r_\Omega,\quad 0\le d<\theta_\Omega .
      \end{equation}
  \end{enumerate}
  Then $(\Gamma,\ge,\plam[k,\lambda],\hlam[k,\lambda],\olam[k,\lambda],\qlam[k,\lambda],\Olam[k,\lambda],a)$
 is a perfectly respectable datum for  $\Hhd[1](A)$.
\end{Theorem}
Before we give the proof we first give a few remarks on the combinatorics of \autoref{T:orbit}.
\begin{Remark} Collapsing $\set{k}\times\Lambda$ and $\set{k}\times\Lambda/\langle\psi\rangle$
each to a point recovers the chain $\set{1<\dots<\ell'}$ from \autoref{T:1}.
\autoref{T:orbit} refines that datum inside each stratum and inside
$\Olam=\HH^{\kapp_k}_1(A)$.

For $k\leq \ell_0$ the only change to the datum of $A$ is that $A$
acts on the left through $\psi^{k-1}$, see \autoref{E:orbitpoly}.
For $k>\ell_0$ an orbit consisting of $r_\Omega$ strata of $A$ fuses into a single
stratum of $\Hhd[1](A)$ whose shape is that of
$\Mat_{r_\Omega}(\mathcal{D}^{k}_\Omega)$ in \autoref{prototype}. The sets
$\mathfrak p$ and $\mathfrak q$ get multiplied by $r_\Omega$, the
multiplicity space $\mathfrak h$ is unchanged and the algebra is the
orbit superalgebra \autoref{E:orbitD} of the stabiliser.
\end{Remark}

\begin{Remark}
As a sanity check for \autoref{T:orbit}  let us calculate the rank of $\Hhd[1](A)$ using the respectable data.
Writing $\mathtt{N}_\mu=|\plam[\mu]||\hlam[\mu]||\olam[\mu]||\qlam[\mu]|$ for the
rank of $A^{(\mu)}$, the bijections of \autoref{respcompatible}\autoref{pcR2} give
$N_\mu=N_\lambda$ for all $\mu\in\Omega$, and since $r_\Omega\theta_\Omega=\theta$
the stratum at $(k,\Omega)$ has rank
$r_\Omega|\plam|\cdot|\hlam|\cdot\theta_\Omega|\olam|\cdot r_\Omega|\qlam|
  =r_\Omega^{2}\theta_\Omega \mathtt{N}_\lambda
  =\theta\,r_\Omega \mathtt{N}_\lambda
  =\theta\sum_{\mu\in\Omega}\mathtt{N}_\mu=\op{rk}\Hhd[1](A).$
\end{Remark}
\begin{proof}[Proof of \autoref{T:orbit}]
Let $H^{\ge k}=u_k\Hhd[1](A)$ be the ideals of \autoref{L:1}, and write
$H^{>\gamma}=\Hhd[1](A)^{>\gamma}$ for $\gamma=(k,\nu)\in\Gamma$. For $k>\ell_0$ write
$\nu=\Omega$, $\lambda=\lambda_\Omega$, $r=r_\Omega$, $\sigma=\sigma_\Omega$,
$\epsilon=\epsilon_\lambda$, $\mathcal D=\mathcal D^{k}_\Omega$, and extend
$a^{\gamma}_{(j,p),h,(o,d),(q,j')}=u_kX^{j}\dot p\dot h\dot oX^{rd}\dot qX^{j'}$ linearly in
the index $(o,d)$ to all $w\in\mathcal D$, using the basis $oW^{d}$ of $\mathcal D$, see \autoref{E:orbitD}. Similarly extend $a^{\lambda}_{phoq}$ to all $w\in\Olam$. We use three facts.
\begin{enumerate}[label=Fact~\arabic*.,leftmargin=*]
\item For $k>\ell_0$, $u_kX=Xu_k$, $bu_k=u_k\psi^{\ell_0}(b)$
  and $u_kX^{\theta}=\kapp_k^{\theta}u_k+u_{k+1}$, by \autoref{E:orbitpoly}.
\item For an index $x$ at $\lambda$ and $m\in\Z$, $\psi^{m}(\dot x)\equiv\beta_x\dot{(\psi^{m}.x)}$
  modulo $A^{>\psi^{m}.\lambda}$ with a unit $\beta_x$, by \autoref{pcR1} and \autoref{pcR2}.
  If $r$ divides $m$, then $\psi^{m}(\dot h)\equiv\epsilon^{m/r}\dot h$ and
  $\psi^{m}(\dot o)\equiv\dot{(\sigma^{m/r}o)}$ modulo $A^{>\lambda}$ by \autoref{pcR3} and the
  definition of $\sigma$.
\item In the first congruence of \autoref{resp2} for $A$ at $\lambda$,
  $c\dot p\dot h\dot o\dot q\equiv\sum_{p',o'}r_{pop'o'}(c)\dot p'\dot h\dot o'\dot q$ modulo
  $A^{>\lambda}$. The factor $\dot q$ may hereby be replaced by any element of $A$ which is congruent
  modulo $A^{>\lambda}$ to a linear combination of the $\dot q'$, $q'\in\qlam$, because the
  scalars are independent of $q$, and $A^{>\lambda}$ is an ideal. Similarly $\dot p\dot h$ may be
  replaced by any element congruent to a linear combination of the
  $\dot p'\dot h'$. By Fact~2, $\psi^{-rd}(\dot q)$ and $\psi^{rd}(\dot p\dot h)$ are such elements.
\end{enumerate}
We now verify the conditions \autoref{resp1}, \autoref{resp2} and \autoref{resp3} using that they hold for $A$.

\noindent\autoref{resp1}: By \autoref{L:1}\autoref{H2} and \autoref{L:2}, $z\mapsto u_kz$
identifies $\HH^{\kapp_k}_1(A)$ with $H^{\ge k}/H^{\ge k+1}$.

\noindent For $k\le\ell_0$ this is $A$,
so the $a^{\gamma}_{phoq}=u_ka^{\nu}_{phoq}$ with $\gamma=(k,\nu)$ form a basis of
$H^{\ge k}/H^{\ge k+1}$.

\noindent Let $k>\ell_0$. Then $\HH^{\kapp_k}_1(A)=\bigoplus_{n<\theta}X^{n}A$ by
\autoref{nequal1}, so $V=H^{\ge k}/H^{\ge k+1}$ has basis the $X^{n}a^{\mu}_{phoq}$ with
$n<\theta$ and $\mu\in\Lambda$. For $\Omega\in\Lambda_\psi$ let
$V^{\succeq\Omega}$ and $V^{\succ\Omega}$ be the spans of those basis elements with
$\Omega_\mu\succeq\Omega$, respectively $\Omega_\mu\succ\Omega$. Let
$A^{\succeq\Omega}=\sum_{\mu\in\Omega}A^{\ge\mu}$ and $A^{\succ\Omega}=\sum_{\mu\in\Omega}A^{>\mu}$.
Since distinct elements of an orbit are incomparable,
\[
  V^{\succeq\Omega}=\bigoplus_{n<\theta}X^{n}A^{\succeq\Omega}\qquad\text{and}\qquad
  V^{\succ\Omega}=\bigoplus_{n<\theta}X^{n}A^{\succ\Omega}.
\]
In $V$, the basis element $a^{\gamma}_{(j,p),h,(o,d),(q,j')}$ of the stratum $\gamma=(k,\Omega)$
becomes $X^{j}\dot p\dot h\dot oX^{rd}\dot qX^{j'}$. Moving all powers of $X$ to the left with
$aX=X\psi(a)$ and applying Fact~2 to each factor gives, for a unit $\beta$ and
$\mu=\psi^{j'}.\lambda$,
\[
  X^{j}\dot p\dot h\dot oX^{rd}\dot qX^{j'}\equiv
  \beta X^{j+rd+j'}a^{\mu}_{\psi^{rd+j'}.p,\ \psi^{rd+j'}.h,\ \psi^{rd+j'}.o,\ \psi^{j'}.q}
  \pmod{V^{\succ\Omega}},
\]
because the error terms of Fact~2 lie in $X^{n}A^{>\mu}\subseteq V^{\succ\Omega}$ for
$\mu\in\Omega$, and $X^{\theta}=\kapp_k^{\theta}$. The indices on the right run exactly once
through the indices of the basis elements $X^{n}a^{\mu}_{phoq}$ with $\mu\in\Omega$. Indeed,
$\mu$ determines $j'$ because $\Omega=\set{\psi^{j'}.\lambda\mid0\le j'<r}$, then $n$
determines $j$ and $d$ because every $n<\theta$ is uniquely $j+rd+j'$ modulo $\theta$ with
$0\le j<r$ and $0\le d<\theta_\Omega$, as $r\theta_\Omega=\theta$, and finally $p$, $h$,
$o$, $q$ are determined by the bijections of Fact~2. Hence the $a^{\gamma}$ with
$\gamma=(k,\Omega)$ form, up to units, a basis of $V^{\succeq\Omega}/V^{\succ\Omega}$. With
respect to a linear extension of $(\Lambda_\psi,\preceq)$, the matrix expressing the
$a^{\gamma}$ with $\gamma\in\set{k}\times\Lambda_\psi$ in the basis of $V$ is therefore block
triangular with invertible diagonal blocks. So these $a^{\gamma}$ form a basis of $V$, and the
$a^{\gamma'}$ with $\gamma'=(k,\Omega')$, $\Omega'\succ\Omega$, span $V^{\succ\Omega}$.

\noindent As $\Gamma$ is ordered lexicographically, all $a^{\gamma}$ form a basis of
$\Hhd[1](A)$, and
\[
  H^{>\gamma}=u_kA^{>\nu}+H^{\ge k+1}\ \text{ if }k\le\ell_0,
  \qquad
  H^{>\gamma}=u_k\sum_{n<\theta}X^{n}A^{\succ\Omega}+H^{\ge k+1}\ \text{ if }k>\ell_0 .
\]
These are two-sided ideals, by \autoref{L:nilp} for $k\le\ell_0$, and for $k>\ell_0$ by
Fact~1, since $A^{\succ\Omega}$ is a $\psi$-stable ideal of $A$ by \autoref{pcR1}. Since
$H^{>\gamma}\supseteq u_kA^{>\lambda}$, congruences modulo $A^{>\lambda}$ may be applied to any
factor of a basis element of the stratum $\gamma$. From now on all congruences are modulo
$H^{>\gamma}$.

\noindent \autoref{resp2}~and~\autoref{resp3}: As $H^{>\gamma}$ is a two-sided ideal, the elements
$x$ satisfying \autoref{resp2}~and~\autoref{resp3} at $\gamma$ form a subalgebra. Thus it suffices to treat $x=X$ and $x=b\in A$.

\noindent For $k\le\ell_0$ we have $u_k=X^{k-1}$, so $Xa^{\gamma}_{phoq}$ and
$a^{\gamma}_{phoq}X$ lie in $H^{\ge k+1}$. Moreover
$ba^{\gamma}_{phoq}=u_k\psi^{k-1}(b)a^{\nu}_{phoq}$ and $a^{\gamma}_{phoq}b=u_ka^{\nu}_{phoq}b$. Hence \autoref{resp2}~and~\autoref{resp3} follow from those for $A$.

\noindent Let $k>\ell_0$. For $c\in A$ let $\rho^{c}_{pp'}(o)=\sum_{o'}r_{pop'o'}(c)o'$ and
$\varsigma^{c}_{qq'}(o)=\sum_{o'}s_{oqo'q'}(c)o'$ in $\Olam$. Then \autoref{resp2} for $A$
reads $ca^{\lambda}_{phoq}\equiv\sum_{p'}a^{\lambda}_{p'h\rho^{c}_{pp'}(o)q}$ and
$a^{\lambda}_{phoq}c\equiv\sum_{q'}a^{\lambda}_{ph\varsigma^{c}_{qq'}(o)q'}$ modulo
$A^{>\lambda}$, whereas \autoref{resp3} for $A$ says that $\rho^{c}_{pp'}$ is right and
$\varsigma^{c}_{qq'}$ is left $\Olam$-linear.

For a linear map $\varphi\colon\Olam\to\Olam$
let $\varphi_L,\varphi_R\colon\mathcal D\to\mathcal D$ be given by $\varphi_L(oW^{d})=\varphi(o)W^{d}$
and $\varphi_R(W^{d}o)=W^{d}\varphi(o)$. Note that
$\varphi_L$ is right $\mathcal D$-linear if $\varphi$ is right $\Olam$-linear, and $\varphi_R$
is left $\mathcal D$-linear if $\varphi$ is left $\Olam$-linear. The $\mathcal D$-bimodule
structure on the stratum is multiplication in the index $w$, so
\autoref{resp3} holds for $x\in\Hhd[1](A)$ as soon as $x$ acts on the index $w$ by a map commuting with the
multiplication by $\mathcal D$ from the other side. We compute the four products.
\begin{enumerate}[label=(\roman*),wide,labelindent=0pt]
\item $Xa^{\gamma}_{(j,p),h,w,(q,j')}=a^{\gamma}_{(j+1,p),h,w,(q,j')}$ if $j<r-1$, and, as
  $W\cdot oW^{d}=\sigma^{-1}(o)W^{d+1}$ in $\mathcal D$, by Facts 1 and 2
  $Xa^{\gamma}_{(r-1,p),h,w,(q,j')}\equiv\beta_p\epsilon^{-1}a^{\gamma}_{(0,\psi^{-r}.p),h,Ww,(q,j')}$,
  with $\beta_p$ from Fact~2 for $m=-r$. Expanding $Ww$ in the basis of $\mathcal D$ gives
  \autoref{resp2} with scalars depending only on $p$ and $w$, and left multiplication by $W$
  gives \autoref{resp3}.

  \item Likewise, $a^{\gamma}_{(j,p),h,w,(q,j')}X=a^{\gamma}_{(j,p),h,w,(q,j'+1)}$ if $j'<r-1$,
  and for $j'=r-1$
  \[
    a^{\gamma}_{(j,p),h,w,(q,r-1)}X\equiv\beta_qa^{\gamma}_{(j,p),h,wW,(\psi^{r}.q,0)},
  \]
  with $\beta_q$ from Fact~2 for $m=r$. This gives \autoref{resp2} and \autoref{resp3}.
\item Let $b\in A$ and $b'=\psi^{\ell_0+j}(b)$. By Facts 1 and 3,
  \begin{align*}
    ba^{\gamma}_{(j,p),h,(o,d),(q,j')}
    &=u_kX^{j}b'\dot p\dot h\dot o\psi^{-rd}(\dot q)X^{rd+j'}
    \equiv\sum_{p',o'}r_{pop'o'}(b')u_kX^{j}\dot p'\dot h\dot o'\psi^{-rd}(\dot q)X^{rd+j'}\\
    &=\sum_{p',o'}r_{pop'o'}(b')a^{\gamma}_{(j,p'),h,(o',d),(q,j')}
  \end{align*}
  with scalars independent of $h$, $q$. Since $b$ acts on the
  index $w$ by $(\rho^{b'}_{pp'})_L$, \autoref{resp2} and \autoref{resp3} hold.
\item Let $b\in A$ and $b'=\psi^{-j'}(b)$. By Facts 2 and 3,
  \[
    a^{\gamma}_{(j,p),h,(o,d),(q,j')}b
    =u_kX^{j+rd}\psi^{rd}(\dot p\dot h)\psi^{rd}(\dot o)\dot qb'X^{j'}
    \equiv\sum_{o',q'}s_{\sigma^{d}(o)qo'q'}(b')u_kX^{j+rd}\psi^{rd}(\dot p\dot h)\dot o'\dot q'X^{j'},
  \]
  and $u_kX^{j+rd}\psi^{rd}(\dot p\dot h)\dot o'\dot q'X^{j'}
  =u_kX^{j}\dot p\dot h\psi^{-rd}(\dot o')X^{rd}\dot q'X^{j'}
  \equiv a^{\gamma}_{(j,p),h,\sigma^{-d}(o')W^{d},(q',j')}$ by Fact~2. The scalars do not
  depend on $p$ and $h$, which is \autoref{resp2}. Since $\sigma^{-d}(o')W^{d}=W^{d}o'$ and
  $oW^{d}=W^{d}\sigma^{d}(o)$, the element $b$ acts on the index $w$ by
  $(\varsigma^{b'}_{qq'})_R$, which gives \autoref{resp3}.
\end{enumerate}
The datum is perfect by construction.
\end{proof}

\subsection{$\kk$-respectability}
\autoref{T:orbit} and \autoref{T:1} are reduction theorems, namely reduction to the RaFH factors $\HH^{\kapp_k}_1(A)$ or the orbit superalgebras
$\mathcal{D}^{k}_\Omega$. They  usually do not provide $\kk$-respectable bases. For $\kk$-respectable $A$, the $\mathcal{D}^{k}_\Omega$ are commutative and split. The following lemma says that a stratum with $\Olam\cong\kk^{m_\lambda}$ can be split into $m_\lambda$ strata each giving rise to  $\Olam=\kk$.

\begin{Theorem}\label{L:splitstratum}
  Let $B$ be a superalgebra. Assume $B$ is perfectly respectable with respectable datum
  $(\Lambda,\ge,\plam,\hlam,\olam,\qlam,\Olam,a)$, and let $\Lambda'\subseteq\Lambda$ be such
  that for every $\lambda\in\Lambda'$ there is an $m_\lambda\ge1$ with
  $\Olam=\bigoplus_{c=1}^{m_\lambda}\kk e^{\lambda}_{c}$ for homogeneous pairwise orthogonal idempotents
  $e^{\lambda}_{c}$ summing to~$1$.  Writing $e^{\lambda}_c=\sum_{o\in\olam}\alpha_{co}\,o$, set
  $\dot e^{\lambda}_c=\sum_{o\in\olam}\alpha_{co}\,\dot o\in B$. Consider the following data.
  \begin{enumerate}
  \item The poset $\Lambda^{\op{split}}=(\Lambda\setminus\Lambda')\sqcup
      \set{(\lambda,c)|\lambda\in\Lambda',\ 1\le c\le m_\lambda}$, with the ordering pulled
      back along the projection $\pi\colon\Lambda^{\op{split}}\to\Lambda$,
      $(\lambda,c)\mapsto\lambda$, that is, $x<y$ if and only if $\pi(x)<\pi(y)$. In
      particular, $(\lambda,c)$ and $(\lambda,c')$ are incomparable for $c\ne c'$.
  \item For $\lambda\in\Lambda\setminus\Lambda'$ let $\plam,\hlam,\olam,\qlam,\Olam,a^\lambda$
      be as in the assumptions.
  \item For $\lambda\in\Lambda'$ and $1\le c\le m_\lambda$ let
  \[
    \plam[(\lambda,c)]=\plam,\quad \hlam[(\lambda,c)]=\hlam,\quad
    \qlam[(\lambda,c)]=\qlam,\quad \olam[(\lambda,c)]=\set{1_\kk},\quad
    \Olam[(\lambda,c)]=\kk,\quad
    a^{(\lambda,c)}_{phq}=\dot p\,\dot h\,\dot e^{\lambda}_c\,\dot q .
  \]
  \end{enumerate}
  This is a perfectly respectable datum, with $\dot p,\dot h,\dot q$ as before and with
  $\dot e^{\lambda}_c$ the lift of the unique element of $\olam[(\lambda,c)]$. In particular,
  if every $\Olam$, $\lambda\in\Lambda$, is isomorphic to $\kk^{m_\lambda}$ for some
  $m_\lambda\ge1$, then $B$ is perfectly $\kk$-respectable.
\end{Theorem}

\begin{proof}
  Since $\sum_{c}B^{(\mu,c)}=B^{(\mu)}$ for $\mu\in\Lambda'$, the ideal $B^{>\nu}$ is the same
  whether computed in $\Lambda$ or in $\Lambda^{\op{split}}$; in particular
  $B^{>(\lambda,c)}=B^{>\lambda}$. Hence at the strata $\nu\in\Lambda\setminus\Lambda'$ the
  conditions \autoref{resp2} and \autoref{resp3} are those of the given datum, while at the new
  strata \autoref{resp3} is automatic by \autoref{R:associativity}, as
  $\Olam[(\lambda,c)]=\kk$.

  Fix $\lambda\in\Lambda'$, write $e_c=e^{\lambda}_c$, and for $y=\sum_oy_o\,o\in\Olam$ set
  $a_{phyq}=\sum_oy_o\,a_{phoq}$, so that $a^{(\lambda,c)}_{phq}=a_{phe_cq}$ because the datum
  is perfect. By \autoref{resp3}, $y\cdot a_{phy'q}=a_{ph(yy')q}=a_{phyq}\cdot y'$. As the
  $e_1,\dots,e_{m_\lambda}$ form a basis of $\Olam$, the $a_{phe_cq}$ form a basis of
  $B^{(\lambda)}$, which gives \autoref{resp1}.

  For \autoref{resp2} let $b\in B$. By \autoref{resp2} at $\lambda$,
  \[
    b\,a_{phe_cq}\equiv\sum_{p'\in\plam,\ o'\in\olam}r_{pe_cp'o'}(b)\,a_{p'ho'q}
    \pmod{B^{>\lambda}} .
  \]
  By \autoref{resp3} the quotient $B^{\ge\lambda}/B^{>\lambda}$ is an $\Olam$-bimodule on which
  the actions of $\Olam$ and of $B$ commute, so we may act with $e_c$ from the right. On the
  left this gives $(b\,a_{phe_cq})e_c\equiv b\,(a_{phe_cq}e_c)=b\,a_{ph(e_ce_c)q}
  =b\,a_{phe_cq}$, and on the right
  $\sum_{p',o'}r_{pe_cp'o'}(b)\,\rho_c(o')\,a_{p'he_cq}$, where $\rho_c\colon\Olam\to\kk$ is
  defined by $o'e_c=\rho_c(o')e_c$; this is well defined since $\Olam$ is commutative with
  basis the orthogonal idempotents $e_c$. Hence
  \[
    b\,a^{(\lambda,c)}_{phq}\equiv\sum_{p'\in\plam}
      \Bigl(\sum_{o'\in\olam}r_{pe_cp'o'}(b)\,\rho_c(o')\Bigr)a^{(\lambda,c)}_{p'hq}
    \pmod{B^{>(\lambda,c)}}
  \]
  with scalars independent of $h$ and $q$. Note that no basis element of a stratum
  $(\lambda,c')$ with $c'\ne c$ occurs, which is what is required, as $(\lambda,c')$ is not
  above $(\lambda,c)$. The right action is analogous, acting with $e_c$ from the left and using
  $e_co'=o'e_c$. Finally, the datum is perfect by construction.
\end{proof}

\subsection{Basic examples}
In the following three examples the action of $\mathtt{C}_\theta$ on $\Lambda$ is
trivial. Hereby $\psi$ can be  trivial or non-trivial.

\begin{Example}\label{Ex:orbittrivial}
  For $A=\kk$ as in \autoref{trivFrobtensor}, the automorphism $\psi$ is the
  identity of order $\theta=1$. Every $f_k$ is linear,
  $\HH^{\kapp_k}_1(\kk)=\mathcal{D}^{k}_\Omega=\kk$, and \autoref{T:1} agrees with
  \autoref{T:orbit}. We have $\Gamma=\set{1<\dots<\ell}$ and all index sets are
  singletons, $\Olam[(k,\bullet)]=\kk$ and $a^{(k,\bullet)}=u_k$. So
  $\Hhd[1](\kk)=\kk[X]/(f)$ is $\kk$-respectable and \autoref{L:splitstratum} has nothing to split. Note that $\Gamma$ records the chosen ordering of the
  entries of the charge.
\end{Example}

\begin{Example}\label{Ex:orbitClifford}
Let $A=\Cl_1$ with $\psi(C)=-C$ as in \autoref{Ex:Clifford} with
$\theta=2$. Assume $2\in\kk^{\times}$. Put $e_\pm=\tfrac12(1\pm C)$, so
that $A=\kk e_+\oplus\kk e_-$ and $\psi(e_\pm)=e_\mp$. Taking $\Lambda=\set{\lambda_+,\lambda_-}$ with no
relations and $\plam[\lambda_\pm]=\set{e_\pm}$,
$\hlam[\lambda_\pm]=\olam[\lambda_\pm]=\qlam[\lambda_\pm]=\set1$,
$\Olam[\lambda_\pm]=\kk$ shows that $A$ is perfectly $\kk$-respectable, because
$A^{>\lambda_\pm}=0$, and $\kk e_\pm$ are ideals.

Let $\mathtt{C}_2$ act on $\Lambda$ by
$\psi.\lambda_\pm=\lambda_\mp$. This is a poset automorphism because
$\Lambda$ is an anti-chain, and $\psi(e_\pm)=e_\mp$ gives
\autoref{respcompatible} with all $\beta_x=1$, while \autoref{pcR3} is empty by
\autoref{cardh1}.

Note that the idempotents $e_\pm$ are not homogeneous. Thus, the respectable datum constructed  is
a datum of algebras and not of superalgebras. There is no homogeneous
substitute. Viewed as superalgebra, $\Cl_1$ has no proper nonzero homogeneous two-sided
ideal and  \autoref{T:orbit} coincides with
\autoref{T:1}.
\end{Example}

\begin{Example}\label{Ex:orbitTaft}
Let $A=\Taft$ as in \autoref{Ex:Taft} with $yg=\omega gy$,
$\psi(g)=\omega^{-1}g$ and $\psi(y)=y$. Assume $\theta\in\kk^{\times}$. Taking
$\Lambda=\set{0<1<\dots<\theta-1}$ to record the $y$-degree, with
$\plam[j]=\set{y^{j}}$, $\hlam[j]=\qlam[j]=\set1$,
$\olam[j]=\set{g^{i}|0\le i<\theta}$, $\Olam[j]=\kk[g]/(g^{\theta}-1)$ and
$a^{j}_{g^{i}}=y^{j}g^{i}$, shows that $A$ is perfectly respectable, because
$A^{>j}=(y^{j+1})$ is an ideal, $g\cdot y^{j}g^{i}=\omega^{j}y^{j}g^{i+1}$ and
$y^{j}g^{i}\cdot g=y^{j}g^{i+1}$ stay in the $j$-th stratum, while $y$ maps it
into $A^{>j}$ from both sides. This extends \autoref{Ex:T2compatible} from
$\theta=2$ to arbitrary $\theta$.

Let $\mathtt{C}_\theta$ act trivially on $\Lambda$. With $\beta_{y^{j}}=1$ and
$\beta_{g^{i}}=\omega^{-i}$ this gives \autoref{respcompatible}, and
$\sigma_\Omega$ is the automorphism $g\mapsto\omega^{-1}g$ of $\Olam[j]$. Every
orbit is a fixed point, so $r_\Omega=1$ and $\theta_\Omega=\theta$. We see that the whole
group is a stabiliser, which is  the other extreme of \autoref{Ex:orbitClifford}. Since $\theta\in\kk^{\times}$ we have $\Olam[j]\cong\kk^{\theta}$ with
idempotents $e_c=\tfrac1\theta\sum_{i}\omega^{-ci}g^{i}$, and
$\sigma_\Omega(e_c)=e_{c+1}$ permutes them in a single cycle. Then
$\mathcal{D}^{k}_\Omega\cong\Mat_\theta(\kk)$, and the whole of $\psi$ sits inside the
fibre algebras and not in the weight poset.
\end{Example}

\section{Respectable CRaFH algebras}
The goal of this section is to show that CRaFH algebras $\Hhd(\Bo)$ with trivial Nakayama automorphism are respectable, provided $A$ is perfectly $\kk$-respectable. To prove this we work in stages, starting with the Frobenius-tensor algebra $\Bo$, then the CRaFH algebras of level $1$, and finally the CRaFH algebras of arbitrary level $\ell\ge1$. Throughout, we assume that $\Hhd(\Bo)$ is braided diagonalisable as in  \autoref{D:Diagonalisable}.

Since the constructions are technically quite demanding, we restrict
ourselves to trivial Nakayama automorphisms. The general case requires
dealing with $\psi$-compatibility and working with orbit space
respectable bases and the necessary tableau combinatorics. This  is
interesting in its own right and will be treated in a separate article.
For simplicity we also work with algebras rather than superalgebras, so
that partitions suffice and strict partitions do not occur. Then, if
$\Bo=A^{\otimes n}$ and the strata of the respectable datum of $A$ are
indexed by a set $I$ with $|I|=N$, the strata of $\Hhd(\Bo)$ are indexed
by $\ell N$-multipartitions of $n$, where $\ell$ is the level of
$\Hhd(\Bo)$. By \autoref{T:Respectable}\autoref{late4}, a subset $\Lambda_+$ of these then label the simple modules.

To show that CRaFH algebras are respectable we need to assume that $A$ is respectable.
This assumption is in reality no restriction by \autoref{stupid}, but we require nice respectable bases to reflect the representation theory. In the examples for $A$ that we most care about nice respectable bases exist. Such examples include  the trivial Frobenius algebra from \autoref{verystupid}, semisimple group algebras of finite
groups, see \autoref{prototype}, truncated polynomial rings, cohomology rings of
sufficiently nice varieties or more generally any finite dimensional
commutative Frobenius algebra via  \autoref{excomm},  zigzag algebras $\mathcal{Z}ig_m$, see \autoref{Ex:ZigZag}, or any idempotent truncation of a based quasi-hereditary algebra via  \autoref{Ex:BasedQH}. The $\kk$-respectability assumption in the main  \autoref{T:HhdRespectable}  allows the weight poset to consist of multipartitions, see \autoref{Rk:Respectable}. In particular,  a classification of simples is given in terms of certain multipartitions.

\subsection{Preliminaries}
We return from the general theory of respectable algebras to our previous setting of RaFH algebras. As in \autoref{D:Hecke} let $\RR$ be a commutative integral domain which is a $\kk$-algebra.  Note that if $A$ is a Frobenius $\kk$-algebra, then $A_R=R\otimes_\kk A$ is again a Frobenius algebra, but now over $R$. Similarly, for (Frobenius-)tensor algebras $\Bo$ and their scalar extension $\Bo_R$.  Moreover,  $\HHd(\Bo)=\HHd(\Bo_R)$. Hence we could have worked with $\kk=R$ from the beginning noting that we stay by \autoref{C:BasisHd} in the category of free $\kk$-modules. From now on we assume $\kk=\RR$.

\begin{Setup} \label{setups9}
We assume that $\RR=\kk$ is a commutative integral domain with $\hbar,t,\eta\in R$. Let $\Bo=A^{\otimes n}$ be a Frobenius-tensor $\kk$-algebra of finite rank with  trivial Nakayama automorphism.
  Suppose that $A$  is perfectly respectable with datum
  $(I, \ge, \mathfrak p_\lambda,\mathfrak h_\lambda,\mathfrak o_\lambda,\mathfrak q_\lambda, a)$  in the sense of \autoref{perfectresp}.

  Without loss of generality we assume that
  $I=\set{1>2>\dots>N}$ is a totally ordered set.  Let
  $\Delta^1,\dots,\Delta^N$ be the standard $A$-modules as in \autoref{D:Standard}.  We consider $I^n = \set{\bi=(i_1,\dots,i_n)|i_k\in I}$ as a poset
  with $\bi\ge\bj$ if $i_k\ge j_k$, for $1\le k\le n$.
  For $a\in A$ and $1\le k\le n$ let  $a^{(k)}=1_A\otimes\cdots\otimes 1_A\otimes  a\otimes 1_A\otimes\cdots\otimes 1_A\in\Bo$,
  where $a$ is in the $k$th tensor factor.
\end{Setup}

\subsection{Respectable Frobenius-tensor algebras}

\begin{Definition}\label{Tensordatum}
 For $\bi\in I^n$ let $\Delta^\bi=\Delta^{i_1}\otimes \dots\otimes \Delta^{i_n}$, and define the following \emph{indexing sets} for $\Bo$
  \begin{equation}
        \plam[\bi]=\plam[i_1]\times\dots\times\plam[i_n],\qquad
        \hlam[\bi]=\hlam[i_1]\times\dots\times\hlam[i_n],\qquad
        \olam[\bi]=\olam[i_1]\times\dots\times\olam[i_n],\qquad
        \qlam[\bi]=\qlam[i_1]\times\dots\times\qlam[i_n].
  \end{equation}
   \end{Definition}

The following is a direct consequence of  the definitions.

\begin{Lemma}\label{P:FrobeniusRespect}
   Assume \autoref{setups9}. Suppose that $A$ is a perfectly $\kk$-respectable algebra and that $n\ge1$.
   \begin{enumerate}
  \item  The
   algebra $\Bo$ is  perfectly respectable with  weight poset $I^n$, index sets
   and  standard modules $\set{\Delta^\bi|\bi\in I^n}$ as in \autoref{Tensordatum}. More precisely,
   $\mathcal{T}_\Ao:=(I^n, \ge, \plam[\bi],\hlam[\bi],\olam[\bi],\qlam[\bi],a)$ with \[\Olam[\bi]=\Olam[{i_1}]\otimes\cdots\otimes\Olam[{i_n}],\quad a^\bi_{phoq}=(\dot p_{1}\dot h_{1}\dot o_{1}\dot q_{1})\otimes\cdots\otimes (\dot p_{n}\dot h_{n}\dot o_{n}\dot q_{n})\in\Bo\] is a perfectly respectable datum. Here
   we use the notation $y=(y_1,\ldots, y_n)\in \frak{y}_\bi$ for $y_j\in\frak{y}_{i_j}$, $\frak{y}\in\set{\plam[],\hlam[],\olam[],\qlam[]}$.
   \item If $\kk$ is a field, then $\set{L^\bi|\bi\in I^n_+}$ is a complete
   set of pairwise non-isomorphic irreducible $\Bo$-modules, where $L^\bi=L^{i_1}\otimes\dots\otimes L^{i_n}$ and $I^n_+=\set{\bi\in I^n|L^{i_k}\ne0\text{ for }1\le k\le n}$.
   \end{enumerate}
\end{Lemma}
\begin{Remark}The same holds in the super setting if we replace $\otimes$ by $\circledast$ and $L^\bi$ by $\Pi^\epsilon L^\bi$, $\epsilon\in\set{0,1}$.
\end{Remark}

\subsection{Tableau combinatorics}\label{SS:Tableaux}
  We recall some standard tableaux  combinatorics that we need to describe our
  respectable bases of the CRaFH algebras.

  A \emph{composition} of $n$ into $\sL$~parts is a sequence
  $\mu=(\mu_1,\dots,\mu_{\sL})\in\N^{\sL}$ of non-negative integers such that
  $|\mu|=\sum_{k=1}^{\sL}\mu_k=n$. The integers $\set{\mu_i|1\le i\le \sL}$ are
  the \emph{parts} of~$\mu$. Let $\Comp[\sL,n]$ be the set of compositions
  of~$n$ into $\sL$~parts. If $\mu=(\mu_1,\dots,\mu_{\sL})$ is a composition,
  let $\Sym_\mu=\Sym_{\mu_1}\times\dots\times\Sym_{\mu_{\sL}}$ be the
  corresponding parabolic subgroup of~$\Sym$.  Explicitly, set
  $\Sym_\mu=\Sym_J$, where
  $J=S\setminus\set{s_{\mu_1},s_{\mu_1+\mu_2},\dots,s_{\mu_1+\dots+\mu_{\sL-1}}}$.

  A \emph{partition} of $n$ is a composition
  $\lambda=(\lambda_1\ge\lambda_2\ge\dots)$ with weakly decreasing parts
  that sum to $|\lambda|=n$. When writing partitions, we often omit zero
  parts and use exponents for the multiplicity of repeated parts. We
  identify a partition with its \emph{diagram} $[\lambda]=\set{(r,c)|1\le
  c\le \lambda_r\text{ for } r\ge1}$, which we think of as a
  left-justified array of boxes in the plane.   A \emph{node} is any element $(r,c)$ of (the diagram of)
  $\lambda$, which is in \emph{row}~$r$ and \emph{column}~$c$.

  Given an integer $\sL\ge1$, an \emph{$\sL$-multipartition} of~$n$ is an
  $\sL$-tuple $\blam=(\lambda^{(1)}|\dots|\lambda^{(\sL)})$ of
  partitions such that $\sum_{k=1}^{\sL}|\lambda^{(k)}|=n$. The
  partition $\lambda^{(k)}$ is the $k$th \emph{component} of $\blam$.
  Similarly, we identify an $\sL$-multipartition $\blam$ with its
  diagram $[\blam]=\set{(m,r,c)|1\le c\le \lambda^{(m)}_r\text{ for }
       r\ge1\text{ and }1\le m\le\sL}$. A triple $(m,r,c)\in\blam$ is a
       \emph{node}.

  \begin{Example}
  The diagrams of
  the partition $(4^2,2,1)$ and the $5$-multipartition $(4,2|2|1|1|1)$ of $11$ are
  \[
        \Diagram{4^2,2,1}\quad \quad\text{and}\quad\quad \Multidiagram{4,2|2|1|1|1}.
  \]
  \end{Example}
We identify partitions and $1$-multipartitions in the obvious way.
\begin{Definition}\label{Defmultipart}
  Let $\Parts[\sL,n]$ be the set of $\sL$-multipartitions of~$n$.
  Then $\Parts[\sL,n]$ is a poset in two ways:
  \begin{enumerate}
 \item The \emph{dominance order} $\gedom$ is defined by setting
   $\blam\gedom\bmu$ for $\blam,\bmu\in\Parts[\sL,n]$ if
  \[
        \sum_{k=1}^{l-1}|\lambda^{(k)}| +\sum_{i=1}^s\lambda^{(l)}_i
           \ge\sum_{k=1}^{l-1}|\mu^{(k)}|+\sum_{i=1}^s\mu^{(l)}_i
           \qquad\text{for $1\le l\le \sL$ and $s\ge1$.}
  \]
  \item  The \emph{lexicographic order} $>_{\mathrm{lex}}$ on $\Parts[\sL,n]$ is the total order where
  $\blam>_{\mathrm{lex}}\bmu$ if there exists $(l,j)$ such that
  $\lambda^{(l)}_j>\mu^{(l)}_j$, and $\lambda^{(k)}_i=\mu^{(k)}_i$
  whenever $k<l$, or $k=l$ and $i<j$.
  \item Define the \emph{lexicographic order} on nodes $(m,r,c)$ and
  $(m',r',c')$ by $(m,r,c)<(m',r',c')$ if $m<m'$, or $m=m'$ and $r<r'$,
  or $m=m'$, $r=r'$ and $c<c'$.
  \end{enumerate}
    Write $\blam\gdom\bmu$ if $\blam\gedom\bmu$ and $\blam\ne\bmu$.  Similarly, write $\blam\ge_{\mathrm{lex}}\bmu$ if $\blam>_{\mathrm{lex}}\bmu$ or $\blam=\bmu$.
  \end{Definition}
  The dominance order is used as a technical tool in the tableau combinatorics. The weight poset of our respectable data is defined in \autoref{D:contorder}. It refines the dominance order and uses the lexicographic order to compare multipartitions with the same component sizes.

  As with partitions, we identify the $\sL$-multipartition $\blam$ with its
 \emph{diagram} which is the set
 \[
    [\blam]=\set{(k,r,c)|1\le k\le \sL\text{ and }1\le c\le \lambda^{(k)}_r\text{ for }r\ge1}.
 \]
 We think of an element as an $\sL$-tuple of diagrams indexed by $1\le k\le \sL$.
 The \emph{node} $(k,r,c)$ of $\blam$
 is in \emph{component}~$k$, in \emph{row}~$r$ and in \emph{column}~$c$.
 \begin{Example}
For the $3$-multipartitions $\blam=(3^2,1\mid4,1^2\mid3)$,  $\bmu=(3^2,1\mid 3^2\mid3)$ the diagrams are
 \[
        \Multidiagram{3^2,1|4,1^2|3}\quad
        \begin{array}[t]{c}
        >_{\mathrm{lex}}\\
        \not\gdom
        \end{array}
        \quad\Multidiagram{3^2,1|3^2|3}.
 \]
 \end{Example}

\begin{Definition}\label{D:Wlam}
 Given $\blam\in\Parts[\sL,n]$, let
 $\Sym_\blam=\Sym_{\lambda^{(1)}}\times\dots\times\Sym_{\lambda^{(\sL)}}$
 be the corresponding parabolic subgroup of~$\Sym$ with associated composition
 $|\blam|=(|\lambda^{(1)}|,\dots,|\lambda^{(\sL)}|)\in \Comp[\sL,n]$. Then $\Sym_\blam$ is a subgroup of
 $\Sym_{|\blam|}$. Let $\Sym^\blam$ and $\Sym^{|\blam|}$ be the sets of
 minimal length coset representatives of $\Sym_\blam$ and
 $\Sym_{|\blam|}$, respectively, in~$\Sym$. Similarly, let
 $\leftindex^\blam\Sym=\set{d^{-1}|d\in\Sym^\blam}$ be the set of
 minimal length right coset representatives of $\Sym_\blam$ in $\Sym$.
\end{Definition}

\begin{Definition}\label{D:tableaux}
 Let $\blam$ be an $\sL$-multipartition and $X$ be a set. A
 \emph{$\blam$-tableau} with values in $X$ is a function
 $\bT\colon\blam\to X$. We think of a $\blam$-tableau as a labelled
 diagram, which allows us to refer to the components, rows and columns
 of a tableau. If~$\bT$ is a $\blam$-tableau, then the \emph{shape}
 of~$\bT$ is $\Shape(\bT)=\blam$.

 A $\blam$-tableau $\bT\colon\blam\to\set{1,2,\dots,n}$ is \emph{row
 standard} if $\bT$ is a bijection and the entries in $\bT$ increase,
 from left to right, along the rows of each component.  Similarly, $\bT$
 is \emph{column standard} if the entries increase from top to bottom
 down the columns of each component. The tableau $\bT$ is
 \emph{standard} if it is both row and column standard.

 Let $\RStd(\blam)$ and $\Std(\blam)$ be the sets of row standard and
 standard $\blam$-tableaux, respectively.  If $\CP$ is a set of
 $\sL$-multipartitions, set $\Std(\CP)=\bigcup_{\blam\in\CP}\Std(\blam)$ and
 let $\Std^2(\CP)=\bigcup_{\blam\in\CP}\Std^2(\blam)$ be the set of pairs
 of standard $\blam$-tableaux of the same shape, for $\blam\in\CP$.
 \end{Definition}

  If $\bT$ is standard and $1\le m\le n$, let $\bT_{\downarrow m}$ be the
 subtableau of $\bT$ that contains $\set{1,2,\dots,m}$. Then
 $\bT_{\downarrow m}$ is a standard $\bmu$-tableau, where
 $\bmu=\Shape(\bT_{\downarrow m})$ is the shape of $\bT_{\downarrow m}$.

 Let $\bS,\bT\in\Std(\Parts)$ be
 standard tableaux, not necessarily of the same shape. Then $\bS$
 \emph{dominates} $\bT$ if
 \[
        \Shape(\bS_{\downarrow m})\gedom \Shape(\bT_{\downarrow m}),
        \qquad\text{ for }1\le m\le n.
 \]
 We abbreviate this as $\bS\gedom\bT$ and write  $\bS\gdom\bT$ if $\bS\gedom\bT$ and $\bS\ne\bT$. The same applies to row standard tableaux, where $\Shape(\bS_{\downarrow m})$ is a multicomposition and $\gedom$ is given by the inequalities of \autoref{Defmultipart}. If
 $\bS\in\RStd(\blam)$ and $\bT=(i,j)\bS\in\RStd(\blam)$,
 for some $1\le i<j\le n$ then
 \begin{equation}\label{E:DominanceStep}
         \bS\gdom\bT
            \quad\Longleftrightarrow\quad \bS^{-1}(i)<\bS^{-1}(j).
 \end{equation}
 Here, $<$ is the lexicographic order on nodes defined
 in \autoref{Defmultipart}(c). That is, $\bS$ dominates $\bT$ if~$i$ is
 in an earlier component of~$\bS$ than~$j$, or $i$ and $j$ are in the
 same component and $i$ is in an earlier row.

 \begin{Example}\label{initial}
 If
 $\blam=(3^2,1\mid3,2\mid2)$, then the following are standard $\blam$-tableaux
 \[
     \Tlam=\Multidiagram[entries=first]{3^2,1|3,2|2} \quad \gedom\quad \bS=\sigma\Tlam=\Multitableau{123,456,{13}|89{10},{11}{12}|7{14}},
 \]
 where $\sigma=(7,13)\in \Sym_{14}$. Moreover, $\Shape(\Tlam_{\downarrow9})=(3^2,1\mid2\mid\emptyset)$, whereas $\Shape(\bS_{\downarrow 9})=(3^2\mid2\mid1)$.
\end{Example}

For each $N$-multipartition $\blam$, the \emph{initial}
 $\blam$-tableau $\Tlam$ is the unique standard $\blam$-tableau such
 that $\Tlam\gedom\bT$, for all $\bT\in\Std(\blam)$. More explicitly,
 $\Tlam$ is the tableau obtained by placing $1,2,\dots,n$  along
 the rows, from left to right, and then top to bottom, starting with the first component
 then the second component etc., see \autoref{initial} for an example.
If $\bS$ is a $\blam$-tableau, let $d_{\bS}\in \Sym$ be the permutation
 of minimal length determined by $\bS= d_{\bS} \Tlam$.

 Recall from above that $\Sym^\blam$ is the set of minimal length coset
 representatives of $\Sym_\blam$ in $\Sym$. It is straightforward to
 check that the map
 \begin{equation}\label{E:RowStd}
    \RStd(\blam)\to\Sym^\blam, \bS\mapsto d_\bS
 \end{equation}
 is a bijection. Compare with \cite[Chapter~3]{Mathas:ULect}.

  In Coxeter group parlance, the dominance order on row
 standard tableaux of the same shape is  the \emph{
 Bruhat order} $\leq$ on~$\Sym$. That is, if $\bS,\bT\in\RStd(\blam)$,
 then $\bS\gedom\bT$ if and only if $d_\bS\le d_\bT$.

\subsection{Parabolic CRaFH modules of level $1$}
 \begin{Definition}\label{leveloneparab}
Let $\blam\in\Parts$. For $J\subseteq S$ with $\Sym_J=\Sym_\blam$, let $\HJ^\pm_\blam=\HJ^\pm_J$, considered as an
element of~$\Hfind$. The corresponding (level one) \emph{parabolic CRaFH
algebra module} is  $\MJ^\pm_\blam=\Hfind\HJ^\pm_\blam$.
\end{Definition}

In view of \autoref{C:HdBasis} and \autoref{C:hwrev},
\autoref{C:MJBasis} implies:
\begin{Proposition}\label{C:StdMJBasis}
  Suppose that $\blam\in\Parts$. The $\Hfind$-module
  $\MJ^\pm_\blam=\Hfind\HJ^\pm_\blam$ is a
  free $\kk$-module with bases
 \begin{equation} \label{naivebases}
      \set{\bb H_{d_\bS}\HJ^\pm_\blam|\bb\in\cBo, \bS\in \RStd(\blam)}
      \qquad\text{and}\qquad
      \set{H_{d_\bS}\bb\HJ^\pm_\blam|\bb\in\cBo, \bS\in \RStd(\blam)}.
  \end{equation}
\end{Proposition}
We will not use these bases in what follows. The reader might however want to keep them in mind as  a rough approximation to the respectable bases, see \autoref{E:HsphqT} below, which we will construct for the
CRaFH algebras. With the goal of adapting the bases \autoref{naivebases} we first define some important elements.

\begin{Definition}
For $\bS,\bT\in \RStd(\blam)$, define
\begin{equation}\label{E:HS}
  H^\pm_\bS = \hbar^{-2\ell(d_\bS)}\xi^\pm_{d_\bS} H_{d_\bS},
  \qquad\text{and}\qquad
  H^{\pm*}_\bT = \hbar^{-2\ell(d_\bT)}H_{d_\bT^{-1}}\xi^{\pm*}_{d_\bT},
\end{equation}
considered as elements of $\Hfind$. Here $\xi^{\pm*}_d=\leftindex^d{{\xi^\pm_{d^{-1}}}}$ denotes the twisted product from
\autoref{D:hw} read in the opposite direction, see \autoref{C:hwrev}.
\end{Definition}
\begin{Remark}
 Since
$\bb H_w=H_w\leftindex^{w^{-1}}{\bb}$ for $\bb\in\Bo$ we have
$H^{\pm*}_\bT=\hbar^{-2\ell(d_\bT)}\xi^\pm_{d_\bT^{-1}}H_{d_\bT^{-1}}$. Note that  $\HJpm_\blam H^{\pm*}_\bT$
is the summand for $d=d_\bT^{-1}$ in the second formula of \autoref{L:HJK}, whereas
$H^\pm_\bS\HJpm_\blam$ is the summand for $d=d_\bS$ in the first formula.
\end{Remark}

\subsection{Respectable bases of CRaFH algebras in level one}
 \label{SS:LevelOneRespectable}
Assuming \autoref{setups9}  we construct in this subsection respectable bases for the CRaFH
 algebra of level~$1$ in terms of a respectable basis of $A$. This basis vastly generalises results of Murphy~\cite{Murphy:basis} for the
 Iwahori--Hecke algebras of the symmetric group.
 Throughout this subsection we assume that $A$ is perfectly $\kk$-respectable, with datum
 $(I,\ge,\plam[i],\hlam[i],\qlam[i],a)$, and we let $N=|I|$. In particular $\Olam[i]=\kk$, and we
 omit the index $o$.

\begin{Remark}
  Every algebra, in particular $\Hfind$, is respectable with the trivial
  datum of \autoref{stupid}, but this datum says nothing about simple
  modules. To obtain a useful respectable datum indexed by
  multipartitions we need to assume that $A$ is perfectly $\kk$-respectable, see
  \autoref{Rk:Respectable}.
\end{Remark}

 Recall from \autoref{D:Diagonalisable} and \autoref{xxinvformula}  that $\HHd$ is braided
 diagonalisable if there is a suitable factorisation
 $(H_1-\xi^+_1)(H_1-\xi^-_1)=0$ of the quadratic relation
 for~$H_1$. We call the CRaFH algebra $\Hhd$
 \emph{braided diagonalisable} if $\HHd$ is braided  diagonalisable.

 As in the introduction of this section, $\Hhd$ is assumed to be braided
 diagonalisable. Throughout this subsection we abuse notation and also
 write $\Tos$ for the image of the structure map $\Tos\to\Bo$ of
 \autoref{D:ST-algebra}, that is, for the subalgebra of $\Bo$ generated
 by the $\taui[i,j]$. This algebra contains the roots $\xipm_i$ by
 \autoref{P:Hieigenvalues} and,  by \autoref{C:xii}, the elements
 $\alpha^\pm_i$ of \autoref{L:alpha} with their inverses
 $(\xipm_i)^{-1}=-\hbar^{-2}\xi^\mp_i$.

\begin{Definition}\label{fancydatum}
Recall the respectable datum and basis for $\Bo$ from \autoref{P:FrobeniusRespect}.
 For $\blam\in\Parts[N,n]$, set
\begin{equation}\label{eqfancydatum}
\begin{gathered}
   \ilam = (i^\blam_1,\dots,i^\blam_n)=(
         \underbrace{1,\dots,1}_{|\lambda^{(1)}|\text{ times}},
         \underbrace{2,\dots,2}_{|\lambda^{(2)}|\text{ times}},
         \dots,
         \underbrace{N,\dots,N}_{|\lambda^{(N)}|\text{ times}}
   )\in I^n,
\end{gathered}
\end{equation}
\[
\plam[\blam]^\pm=\set{ H^\pm_\bS p|\bS\in \Std(\blam), p\in\plam[\ilam]},
\quad\hlam[\blam]^\pm=\set{\HJpm_\blam},
\quad\qlam[\blam]^\pm=\set{hqH^{\pm*}_\bT|\bT\in \Std(\blam), h\in\hlam[\ilam], q\in\qlam[\ilam]}.
\]
\end{Definition}
Our goal, see
\autoref{T:HfindRespectable}, is to show that $\Hfind$ is respectable with respectable basis of the form
\begin{equation}\label{E:HsphqT}
  H^\pm_{\bS,p,h,q,\bT}= H^\pm_\bS p\HJpm_\blam hqH^{\pm*}_\bT,
      \quad\text{ for } \bS,\bT\in\Std(\blam),
      p\in\plam[\ilam], h\in\hlam[\ilam], \text{ and } q\in\qlam[\ilam].
\end{equation}

 \begin{Definition}\label{fancyB}
Using  \autoref{fancydatum} let $\Basis=\bigcup_{\blam\in\Parts}\Basis_\blam\subset \Hfind$, where
$\Basis_\blam=\plam[\blam]^\pm\hlam[\blam]^\pm\qlam[\blam]^\pm$.

Let $\HH^{\pm\ge\blam}_n$ and
 $\HH^{\pm>\blam}_n$ be the $\kk$-spans of $\bigcup_{\bmu\ge\blam}\Basis_\bmu$ and
 $\bigcup_{\bmu>\blam}\Basis_\bmu$, respectively.
 \end{Definition}

Comparing with \autoref{E:HsphqT}, we want to show that $\Basis$ is a respectable basis for a respectable datum given by \autoref{eqfancydatum}.

We start with a combinatorial sanity check.

For $\mu\in\Comp$ we denote by $\binom n\mu=\tfrac{n!}{\mu!}$ the
corresponding multinomial coefficient, where $\mu!=\mu_1!\dots\mu_N!$.

\begin{Lemma}\label{L:rank}
  The algebra $\Hfind$ is free as a $\kk$-module of rank $|\Basis|$.
\end{Lemma}
\begin{proof}
By \autoref{basisthmcycl} we know that $\Hfind$ is free as a $\kk$-module of rank $n!(\op{rk} A)^n$. The claim of the lemma is therefore that  $n!(\op{rk} A)^n=|\Basis|$. To prove this we compute $|\Basis|$.
  For this fix $\mu\in\Comp$. For $\blam\in\Parts$ with $|\blam|=\mu$,  let $\Std_0(\blam)$ be the
  set of standard $\blam$-tableaux such that the integer~$m$ appears in component $k$ if
  and only if $\sum_{j=1}^{k-1}|\lambda^{(j)}|<m\le\sum_{j=1}^k|\lambda^{(j)}|$.
  Then $|\Std(\blam)|=\binom n\mu |\Std_0(\blam)|$. Therefore, by the Robinson--Schensted
  correspondence for $\Sym_\mu$, see e.g.\ \cite[\S4]{Fulton:Tableaux},
  \begin{equation}\label{RobinsonSchensted}
    \sum_{\substack{\blam\in\Parts\\|\blam|=\mu}}|\Std_0(\blam)|^2
      =|\Sym_\mu| = n!\;\binom n\mu^{-1}.
  \end{equation}
  For $1\le k\le N$, set $d_k=|\plam[k]\hlam[k]\qlam[k]|$. Then
  $\op{rk} A=\sum_{k=1}^Nd_k$. Let $d^\mu=\prod_{k=1}^Nd_k^{\mu_k}$.  Then
  \[
    |\Basis| = \sum_{\mu\in\Comp}\sum_{|\blam|=\mu}|\Std(\blam)|^2d^\mu
        = \sum_{\mu\in\Comp}\binom{n}{\mu}^2 d^\mu\sum_{|\blam|=\mu}|\Std_0(\blam)|^2
        = n!\sum_{\mu\in\Comp}\binom n\mu d_1^{\mu_1}\cdots d_N^{\mu_N}
        = n!\Bigl(\sum_{k=1}^Nd_k\Bigr)^n,
  \]
 and thus $|\Basis|=n!(\op{rk} A)^n$ as claimed.  In the calculation we used \autoref{RobinsonSchensted} for  the third equality and the multinomial theorem for the last equality.
\end{proof}
   To construct respectable bases for CRaFH algebras, we first  define a poset of weights.

 \begin{Definition}\label{D:contorder}
  Recall that $I=\set{1>2>\dots>N}$ by \autoref{setups9}. For
  $\bmu,\blam\in\Parts$ write $\bmu>\blam$ if either
  $|\mu^{(1)}|+\dots+|\mu^{(k)}|\ge|\lambda^{(1)}|+\dots+|\lambda^{(k)}|$ for $1\le k\le N$, with at
  least one strict inequality, or $|\mu^{(k)}|=|\lambda^{(k)}|$ for $1\le k\le N$ and
  $\bmu>_{\mathrm{lex}}\blam$. Write $\bmu\ge\blam$ if $\bmu>\blam$ or $\bmu=\blam$.
\end{Definition}

\begin{Remark}\label{R:glweights}
 \autoref{D:contorder} is the usual order on $\mathfrak{gl}_N$-weights, refined lexicographically on multipartitions of the same weight. To see this let  $\operatorname{wt}(\blam)=\sum_{k=1}^N|\lambda^{(k)}|e_k\in\Z^N$ be the \emph{weight} of $\blam\in\Parts$ This can be  viewed as a $\mathfrak{gl}_N$-weight. Then $\bmu>\blam$ in \autoref{D:contorder} if and only if either $\operatorname{wt}(\bmu)-\operatorname{wt}(\blam)$ is a nonzero sum of positive roots $e_k-e_{k'}$, $k<k'$, or $\operatorname{wt}(\bmu)=\operatorname{wt}(\blam)$ and $\bmu>_{\mathrm{lex}}\blam$. \end{Remark}

  We can now state the main result of this section,
 namely a construction of a respectable basis for  CRaFH algebras of level~$1$.
  \begin{Theorem}\label{T:HfindRespectable}
   Let $\Hfind$ be a braided diagonalisable CRaFH algebra of
   level~$\ell=1$. Assume $A$ is perfectly $\kk$-respectable with datum
   $(I,\ge,\plam[i],\hlam[i],\qlam[i],a)$. Then
   \begin{enumerate}
   \item The $\kk$-module
   $\HH^{\pm>\blam}_n$ is a two-sided ideal of $\Hfind$ for every $\blam\in\Parts$.
   \item
     With the ordering from \autoref{D:contorder}, $\Hfind$ is a perfectly $\kk$-respectable algebra with datum
   \[(\Parts,\ge,\plam[\blam]^\pm,
   \hlam[\blam]^\pm,\qlam[\blam]^\pm,a^{\blam,\pm}_{phq}=\dot p\dot h\dot q)\quad\text{for either of the two choices $+$ or $-$.}
  \]
   In particular, $\Basis=\set{\plam[\blam]^\pm
   \hlam[\blam]^\pm\qlam[\blam]^\pm|\blam\in\Parts}$
   is a perfectly respectable basis of $\Hfind$.
   \end{enumerate}
 \end{Theorem}

\begin{Remark}\label{pmchoice}
 \autoref{T:HfindRespectable} gives two respectable bases of
 $\Hfind$, which depend on the choice of root $\xi^\pm_1$ for the
 quadratic relations satisfied by~$H_1$. In the special case of the group algebras of the
 symmetric groups (when $\eta=0$, $\hbar=1$ and $\Bo=\kk$), these bases
 correspond to choosing between the trivial and sign representations.
\end{Remark}

\begin{Remark}\label{R:RightAction} The factorisation of the basis elements according to the index sets is asymmetric with respect to $h\in\hlam[\ilam]$.  This asymmetry seems to be necessary.  To illustrate this consider the symmetric Frobenius algebra $A=\kk[x,y]/(x^2,y^2)$ from \autoref{Ex:KleinFour} with $\cB=\set{1, x, y, xy}$, $\tr(b)=\delta_{b,xy}$ for $b\in\cB$. Then $\cB$ is a respectable basis for the perfectly $\kk$-respectable datum with weight poset  $I=\set{1>2>3}$, $\plam[i]=\qlam[i]=\set{1}$ for all $i$, and
  $\hlam[1]=\set{xy}$, $\hlam[2]=\set{x,y}$, $\hlam[3]=\set{1}$. Let $n=2$, $\eta=0$ and
  $\blam=(\emptyset| 2|\emptyset)$. Then $\ilam=(2,2)$, $\Std(\blam)=\set{\Tlam}$ and
  $H^\pm_{\Tlam}=1=H^{\pm*}_{\Tlam}$, so that $\Basis_\blam$ consists of the four elements
  $\HJpm_\blam h$ with $h\in\set{x\otimes x,x\otimes y,y\otimes x,y\otimes y}$. We claim that the second congruence of \autoref{resp2} fails for the datum which has as the unique non-trivial index set
  $\hlam[\blam]^\pm=\set{\HJpm_\blam h|h\in\hlam[\ilam]}$. Indeed,
  $\HJpm_\blam=1\pm\hbar^{-1}H_1$ and $\HJpm_\blam H_1=\pm\hbar\HJpm_\blam$ by
  \autoref{L:HJ}\autoref{HJ1}, hence $\HJpm_\blam(x\otimes y)H_1=\HJpm_\blam H_1(y\otimes x)=\pm\hbar\,\HJpm_\blam(y\otimes x)$. This element  is not
  congruent modulo $\HH^{\pm>\blam}_n$ to a multiple of the basis element with $h=x\otimes y$.
  This proves the claim.
\end{Remark}

\begin{Remark}\label{Rk:Respectable}
  In   \autoref{T:HfindRespectable} we assume that $A$ is perfectly $\kk$-respectable. This assumption cannot be dropped.
  Take $n=2$, $\eta=0$, $\op{rk} A>1$ and the datum of \autoref{stupid}
  for $A$, so that $I=\set{1}$ and $\Olam[1]=A$. The strata of \autoref{fancydatum} are then
  $\HJpm\Bo$ for $\blam=(2)$ with $\HJpm=1\pm\hbar^{-1}H_1$, and $\Bo$ for $\blam=(1,1)$. By
  \autoref{resp2} one of them is a two-sided ideal. The subalgebra $\Bo$
  is not a two-sided ideal, and neither is
  $\HJpm\Bo$, since for $a\in A$, $a\notin\kk1_{\HH}$, $(a\otimes1)\HJpm=(a\otimes1)\pm\hbar^{-1}H_1(1\otimes a)
    \notin\HJpm\Bo=\set{\bb\pm\hbar^{-1}H_1\bb|\bb\in\Bo}.$
  Moving $\Bo$ past $\HJpm_\blam$ permutes tensor factors, and the terms recombine
  if the coefficients do not depend on the elements to the right of $\HJpm_\blam$.
  Axiom \autoref{resp2} guarantees this for $h$ and $q$, but not for $o$.
\end{Remark}

 Before we launch into the proof of \autoref{T:HfindRespectable},
 we illustrate the basis in an example.
\begin{Example}\label{Ex:HhdBasis}
   Let $A$ be the zigzag algebra $\mathcal{Z}ig_2$ from  \autoref{Ex:ZigZag}, and let
   $\set{(\plam[i],\hlam[i],\qlam[i])|i\in\Lambda}$ be the $\kk$-respectable basis of $A$
   given in this example. Then $N=|I|=3$, where $I=\Lambda=\set{1>2>3}$ as in \autoref{Ex:ZigZag} (called  $\Lambda$ there.) The three strata have $d_1=1$, $d_2=4$ and $d_3=1$ basis elements.

   Let now  $n=2$. We get $\op{rk}\Hfind=2!(\op{rk} A)^2=72$ which should be compared  with the proof of  \autoref{L:rank}.  The respectable basis of
   $\Hfind[2]$ from
   \autoref{T:HfindRespectable} is indexed by  the set $\Parts[3,2]$ of $3$-multipartitions
   of~$2$. The table below lists, for each $\blam\in\Parts[3,2]$, the possible tableaux $\bS$ and $\bT$ in
   column two.  In columns three and five we assume  $c,c'\in\plam[2],$ and $d,d'\in\qlam[2]$, respectively. The weight poset is a chain, as indicated in the table,  except that $(1|0|1)$ is incomparable with both $(0|2|0)$ and $(0|1^2|0)$.

  The bases of $\Hfind$ given by \autoref{T:HfindRespectable} are:
   \[
     \begin{array}{l*4cr}\toprule
       \blam  & \Std(\blam) & \plam[\blam]^\pm & \hlam[\blam]^\pm
              & \qlam[\blam]^\pm &|\Basis_\blam| \\ \midrule
       (2|0|0)   & \multitableau{12||} &H^\pm_\bS(ba\otimes ba)& \HJ^\pm_{(2)}
              & e_{11}H^{\pm*}_\bT&1\\
       (1^2|0|0) & \multitableau{1,2||}&H^\pm_\bS(ba \otimes ba)&1
              & e_{11}H^{\pm*}_\bT&1\\
       (1|1|0)   & \multitableau{1|2|}, \multitableau{2|1|}&H^\pm_\bS(ba \otimes c)
              & 1 & (e_1\otimes d)H^{\pm*}_\bT&16\\
       (1|0|1)   & \multitableau{1||2}, \multitableau{2||1}&H^\pm_\bS(ba \otimes e_2)
              & 1 & e_{12}H^{\pm*}_\bT&4\\
       (0|2|0)   & \multitableau{|12|}&H^\pm_\bS(c\otimes c')& \HJ^\pm_{(2)}
              & (d\otimes d')H^{\pm*}_\bT&16\\
       (0|1^2|0) & \multitableau{|1,2|}&H^\pm_\bS(c\otimes c')& 1
              & (d\otimes d')H^{\pm*}_\bT&16\\
       (0|1|1)   & \multitableau{|1|2}, \multitableau{|2|1}&H^\pm_\bS(c\otimes e_2)
              & 1 & (d\otimes e_2)H^{\pm*}_\bT&16\\
       (0|0|2)   & \multitableau{||12}&H^\pm_\bS e_{22}& \HJ^\pm_{(2)}
              & e_{22}H^{\pm*}_\bT&1\\
       (0|0|1^2) & \multitableau{||1,2}&H^\pm_\bS e_{22}& 1
              & e_{22}H^{\pm*}_\bT&1\\
       \bottomrule
     \end{array}
   \]
   where $e_{ij}=e_i\otimes e_j$ and
   $\HJ^\pm_{(2)}=1+\hbar^{-2}\xi^\pm_1H_1$. The asymmetry from
   \autoref{fancydatum} is not visible, since $h$ involves only $e_1$ and $e_2$. The ranks  in the right-hand column add up to  $72=\op{rk}\Hfind$ as expected.

   Let now $\kk$ be an algebraically closed field. Recall from \autoref{T:Respectable}\autoref{late4} the classification of simple modules. The  top four strata obviously do not contribute to $\Lambda_+$. One can check that the strata indexed by $(0|2|0)$, $(0|0|2)$, $(0|1^2|0)$, $(0|1|1)$, $(0|0|1^2)$ give rise to simple modules if $\op{char}(\kk)\not=2$, whereas this is only true for the last three ones if $\op{char}(\kk)=2$. This should be compared with  \autoref{Ex:ZigZagSimples}.
\end{Example}

The proof of \autoref{T:HfindRespectable} is a downward induction on $\blam$. Its main step shows
that multiplying an element of $\Basis_\blam$ by a generator of $\Hfind$ gives, modulo
$\HH^{\pm>\blam}_n$, a combination of elements of $\Basis_\blam$, see \autoref{L:InductiveStep}.
For the generators in $\Bo$ this follows from \autoref{L:HJ}, and for the $H_m$ from a
straightening procedure based on Garnir tableaux, see \autoref{P:Straighten}. We first collect commutation
rules between $\HJpm_\blam$, the elements $\xipm_vH_v$ with $v\in\Sym_{|\blam|}$, and $\Bo$.

\begin{Lemma}\label{L:HJ}
  Let $\blam\in\Parts$. Let $\Tos_\blam=\kk$ if $\eta=0$, and let
  $\Tos_\blam\subseteq\Tos$ be the subalgebra generated by the $\taui[a,b]$ with
  $(a,b)\in\Sym_{|\blam|}$ if $\eta\ne0$. Then the following hold for $1\le i<n$.
  \begin{enumerate}
  \item\label{HJ1} If $s_i\in\Sym_\blam$, then $\HJpm_\blam H_i=\HJpm_\blam\xipm_i$.
  \item\label{HJ2} If $w\in\Sym_{|\blam|}$, then $\leftindex^w{{\Tos_\blam}}=\Tos_\blam$. Moreover,
    $\xipm_v,(\xipm_v)^{-1}\in\Tos_\blam$ for $v\in\Sym_{|\blam|}$, and $\alpha^\pm_i\in\Tos_\blam$ for
    $s_i\in\Sym_{|\blam|}$.
  \item\label{HJ3} If $\rho\in\Tos_\blam$ and $\bb\in\Bo$, then
    $\rho\bb\in\sum_{x\in\Sym_{|\blam|}}\leftindex^x\bb\Tos_\blam$ and
    $\bb \rho\in\sum_{x\in\Sym_{|\blam|}}\Tos_\blam\leftindex^x\bb$.
  \item\label{HJ4} If $\rho\in\Tos_\blam$ and $w\in\Sym_\blam$, then
    $\HJpm_\blam \rho\xipm_wH_w\in\HJpm_\blam\Tos_\blam$.
  \item\label{slide2} If $v\in\Sym_{|\blam|}$ and $\bb\in\Bo$, then
    \begin{align*}
      \xipm_vH_v\bb&\in\sum_{x\in\Sym_{|\blam|}}\Tos_\blam\leftindex^x\bb\xipm_vH_v&\text{and}&&
       \xipm_vH_v\bb&\in\sum_{x\in\Sym_{|\blam|}}\leftindex^x\bb\xipm_vH_v\Tos_\blam,\\
      \bb\xipm_vH_v&\in\sum_{x\in\Sym_{|\blam|}}\Tos_\blam\xipm_vH_v\leftindex^x\bb&\text{and}&&
       \bb\xipm_vH_v&\in\sum_{x\in\Sym_{|\blam|}}\xipm_vH_v\leftindex^x\bb\Tos_\blam.
    \end{align*}
  \item\label{slide4} If $\bS\in\Std(\blam)$ and $\bb\in\Bo$, then there are $\bb_1,\bb_2\in\Bo$
    with $\bb H^\pm_\bS=H^\pm_\bS\bb_1$ and $H^{\pm*}_\bS\bb=\bb_2H^{\pm*}_\bS$.
  \end{enumerate}
\end{Lemma}

\begin{proof}
  \autoref{HJ1}: Let $J\subseteq S$ with $\Sym_J=\Sym_\blam$, so that $\HJpm_\blam=\HJpm_J$. By
  \autoref{D:HJ}, $\HJpm_{\set{s_i}}=1+\hbar^{-2}\xipm_iH_i$. Since $\xipm_i$ is $s_i$-invariant by
  \autoref{xxxinvformulas}, it commutes with $H_i$, and hence $\HJpm_{\set{s_i}}$ commutes with
  $H_i$ and $\xipm_i$. If $s_i\in J$, then $H_i\HJpm_{\set{s_i}}=\xipm_i\HJpm_{\set{s_i}}$ by the
  proof of \autoref{P:HW}. Altogether,
  \begin{equation}\label{Leop}
    \HJpm_{\set{s_i}}H_i=H_i\HJpm_{\set{s_i}}=\xipm_i\HJpm_{\set{s_i}}=\HJpm_{\set{s_i}}\xipm_i.
  \end{equation}
  On the other hand, the first formula of \autoref{L:HJK}, applied to $\set{s_i}\subseteq J$, gives
  \[
    \HJpm_J=Y\HJpm_{\set{s_i}}\qquad\text{with}\qquad
    Y=\sum_{d\in\Sym^{\set{s_i}}_J}\hbar^{-2\ell(d)}\xipm_dH_d,
  \]
  and multiplying \autoref{Leop} by $Y$ from the left gives $
  \HJpm_JH_i= \HJpm_J\xipm_i$ and, hence, the claim.

  \autoref{HJ2}: If $\eta=0$, then $\xipm_v=(\pm\hbar)^{\ell(v)}$ and $\alpha^\pm_i\in\kk$ by
  \autoref{L:alpha}. Let $\eta\ne0$. Since, by \autoref{D:ST-algebra}, the structure map is
  $\Sym$-equivariant and
  $(w(a),w(b))=w(a,b)w^{-1}\in\Sym_{|\blam|}$ for $w\in\Sym_{|\blam|}$ the first claim follows. All $\xipm_j$ are
  polynomials in $\eta\taui[j,j+1]$, and $(\xipm_j)^{-1}=-\hbar^{-2}\xi^\mp_j$ by \autoref{C:xii}.
  Thus, $\xipm_v$ and
  $(\xipm_v)^{-1}$ lie in $\Tos_\blam$ by \autoref{D:hw} since $v$ has a reduced expression in the simple reflections of $\Sym_{|\blam|}$. Finally,
  $\alpha^\pm_i\in\kk[(\eta\taui[i,i+1])^2]$ by \autoref{L:alpha}.

  \autoref{HJ3}: There is nothing to show if $\eta=0$. Otherwise it suffices to let $\rho$ be a
  monomial in the generators $\taui[a,b]$ of $\Tos_\blam$. Since $\psi$ is trivial,
  \autoref{Nakayama} gives $\taui[a,b]\bb=\leftindex^{(a,b)}{\bb}\taui[a,b]$ for $\bb\in\Bo$, and the
  claim follows by induction on the degree of $\rho$.

  \autoref{HJ4}: As $\xipm_w\in\Tos_\blam$ by \autoref{HJ2}, it suffices to show
  $\HJpm_\blam \rho H_w\in\HJpm_\blam\Tos_\blam$, by induction on $\ell(w)$. the base case $w=e$ is clear.  Let $w=w's_i$ with
  $s_i\in\Sym_\blam$ and $\ell(w)=\ell(w')+1$. By induction, $\HJpm_\blam \rho H_{w'}=\HJpm_\blam \rho'$
  for some $\rho'\in\Tos_\blam$, and by \autoref{HJ1} and \autoref{HJ2},
  $\HJpm_\blam \rho'H_i=\HJpm_\blam H_i\leftindex^{s_i}{{\rho'}}=\HJpm_\blam\xipm_i\leftindex^{s_i}{{\rho'}}
  \in\HJpm_\blam\Tos_\blam$.

  \autoref{slide2}: For $\rho\in\Tos_\blam$ we have $\xipm_vH_v\rho=\xipm_v\leftindex^v{\rho}(\xipm_v)^{-1}\xipm_vH_v$,
  so $\xipm_vH_v\Tos_\blam=\Tos_\blam\xipm_vH_v$ by \autoref{HJ2}. Moreover,
  $\xipm_vH_v\bb=\xipm_v\leftindex^v\bb(\xipm_v)^{-1}\xipm_vH_v$ and
  $\bb\xipm_vH_v=\xipm_vH_v\leftindex^{v^{-1}}{\left((\xipm_v)^{-1}\bb\xipm_v\right)}$. By \autoref{HJ2}, \autoref{HJ3}, the elements $\xipm_v\leftindex^v\bb(\xipm_v)^{-1}$ and
  $\leftindex^{v^{-1}}{\left((\xipm_v)^{-1}\bb\xipm_v\right)}$ lie in
  $\sum_{x\in\Sym_{|\blam|}}\Tos_\blam\leftindex^x\bb$ and in
  $\sum_{x\in\Sym_{|\blam|}}\leftindex^x\bb\Tos_\blam$. Altogether, this gives all four inclusions.

  \autoref{slide4}: This is clear, since $H^\pm_\bS=\hbar^{-2\ell(d_\bS)}\xipm_{d_\bS}H_{d_\bS}$,
  $H^{\pm*}_\bS=\hbar^{-2\ell(d_\bS)}\xipm_{d_\bS^{-1}}H_{d_\bS^{-1}}$ and $\xipm_w$ is invertible in $\Bo$.
\end{proof}

Recall from \autoref{D:tableaux} that a $\blam$-tableau is a function $\bT\colon\blam\to X$, for
some set~$X$, where we identify the multipartition $\blam$ with its diagram, which is a set of nodes. We will use Garnier data.

\begin{Definition}\label{D:GarnirTableau}
  Let $\blam\in\Parts$ and $\alpha=(k,r,c)\in\blam$ with $(k,r+1,c)\in\blam$. The \emph{Garnir
  tableau} $\bG_\alpha=\bG^\blam_\alpha$ is the $\blam$-tableau given by
  \begin{align*}
    \bG_\alpha(l,a,b)&= \begin{cases*}
         \Tlam(k,r,c)+b   & if $l=k$, $a=r$   and $b\ge c$,\\
         \Tlam(k,r,c)+b-1 & if $l=k$, $a=r+1$ and $b\le c$,\\
         \Tlam(l,a,b)     & otherwise.
      \end{cases*}
  \end{align*}
  Its \emph{Garnir belt} is the set of nodes $(k,r,b)$ with $b\ge c$ and $(k,r+1,b)$ with $b\le c$.
  Let $\bnu$ be the $N$-multicomposition obtained from $\blam$ by replacing
  $\lambda^{(k)}_r,\lambda^{(k)}_{r+1}$ by $c-1,\lambda^{(k)}_r+1,\lambda^{(k)}_{r+1}-c$, and let
  $\bmu\in\Parts$ be obtained from $\bnu$ by reordering the parts of $\nu^{(k)}$, and let
  $J\subseteq S$ with $\Sym_J=\Sym_\bnu\cap\Sym_\blam$. We call $(\alpha,\bG_\alpha,\bnu,\bmu,J)$
  a \emph{Garnir datum} for $\blam$.
\end{Definition}

Thus $\bG_\alpha$ coincides with $\Tlam$ outside the Garnir belt, and inside the belt the entries of
$\Tlam$ are entered in increasing order, first along row $r+1$ and then along row $r$.

For
$w\in\Sym$ let $\op{Inv}(w)=\set{(a,b)|1\le a<b\le n,\ w(a)>w(b)}$, so that $\ell(w)=|\op{Inv}(w)|$. Recall that
if $g,d\in\Sym$ with $\op{Inv}(g)\subseteq\op{Inv}(d)$, then $\ell(dg^{-1})=\ell(d)-\ell(g)$.

\begin{Example}\label{Ex:Garnir}
  Let $N=1$, $\blam=(4,3,1)$ and $\alpha=(1,1,2)$. Then $(\alpha, \bG_\alpha, \bnu, \bmu, J)$ is a Garnir datum where  $\bnu=(1,5,1,1)$, $\bmu=(5,1^3)>\blam$,
  $\Sym_J=\<s_2,s_3,s_5\>$. In terms of tableaux,
  \[
    \aTabset{scale=0.8, box font=\small}
    \Tlam=\Tableau{1*2*3*4,*5*67,8},\qquad \bG_\alpha=\Tableau{1*4*5*6,*2*37,8},\qquad
    w\Tlam=\Tableau{1*3*5*6,*2*47,8},\qquad d\Tlam=\Tableau{1467,235,8}.
  \]
   with the Garnir belt shaded and some choice of $w\in\Sym^J_\bnu$. Such an element is determined  by  the set of entries of $w\Tlam$ in row~$2$ of the belt, which can be any $2$-element subset of  $\set{2,\dots,6}$. Thus, $|\Sym^J_\bnu|=10$. Moreover,  $w\Tlam$ is standard unless this set is $\set{2,3}$, that is $w=d_{\bG_\alpha}$.  Here $\ell(d_{\bG_\alpha})=6$, while $\ell(w)\le5$ for $w\ne d_{\bG_\alpha}$. Finally, note that $d\Tlam$ is not
  standard in column~$2$ and we can write $d=ud_{\bG_\alpha}$ for the cycle $u=(5,6,7)$ and have
  $\ell(d)=8=\ell(u)+\ell(d_{\bG_\alpha})$. This should illustrate  \autoref{L:GarnirElement} and its proof.
\end{Example}

The next lemma contains the combinatorics of Garnir tableaux.

\begin{Lemma}\label{L:GarnirElement}
  Let $(\alpha,\bG_\alpha,\bnu,\bmu,J)$ be a Garnir datum for $\blam\in\Parts$.
  \begin{enumerate}
  \item\label{GE1} We have $\bmu>\blam$ and $\bi_\bmu=\ilam$.
  \item\label{GE2} We have $d_{\bG_\alpha}\in\Sym^J_\bnu\subseteq\Sym_{|\blam|}$, and $w\Tlam$ is
    standard with $\ell(w)<\ell(d_{\bG_\alpha})$ for all $w\in\Sym^J_\bnu\setminus\set{d_{\bG_\alpha}}$.
  \item\label{GE3} If $d\in\Sym^\blam$ and $d\Tlam$ is not standard, then $d=ud_{\bG_\beta}$ with
    $\ell(d)=\ell(u)+\ell(d_{\bG_\beta})$ for some Garnir datum $(\beta,\bG_\beta,\dots)$ for $\blam$
    and some $u\in\Sym$.
  \end{enumerate}
\end{Lemma}

\begin{proof}
  \autoref{GE1}: The components of $\bmu$ and $\blam$ have the same sizes, so $\bi_\bmu=\ilam$.
  Let $j\le r$ be minimal with $\lambda^{(k)}_j\le\lambda^{(k)}_r$. Then $\mu^{(k)}_i=\lambda^{(k)}_i$
  for $i<j$ and $\mu^{(k)}_j=\lambda^{(k)}_r+1>\lambda^{(k)}_j$, so $\bmu>_{\mathrm{lex}}\blam$ and
  $\bmu>\blam$ by \autoref{D:contorder}.

  \autoref{GE2}: Let $B$ be the Garnir belt of $\bG_\alpha$. The entries of $\Tlam$ in $B$ are
  consecutive and form the part of size $\lambda^{(k)}_r+1$ of $\bnu$, while every other part of
  $\bnu$ consists of entries of a single row of $\Tlam$. Hence $\Sym_\bnu=\Sym_B\times\Sym'$, where
  $\Sym_B$ permutes the entries in $B$ and $\Sym'\subseteq\Sym_\blam$. As the right descents of an
  element of $\Sym_\bnu$ lie in $\Sym_\bnu$, the set $\Sym^J_\bnu=\Sym_\bnu\cap\Sym^\blam$ consists
  of the $w\in\Sym_B$ for which the rows of $w\Tlam$ increase. In particular
  $\Sym^J_\bnu\subseteq\Sym_{|\blam|}$. For such $w$, the entries of $w\Tlam$ in row $r$ to the
  left of $B$ are smaller, and those in row $r+1$ to the right of $B$ are larger than all entries in
  $B$. So $w\Tlam$ is standard unless every entry of $w\Tlam$ in row $r$ of $B$ is bigger than every entry
  in row $r+1$ of $B$, in which case $w\Tlam=\bG_\alpha$ or equivalently $w=d_{\bG_\alpha}$. Finally,
  $\op{Inv}(d_{\bG_\alpha})$ consists of all pairs $(a,b)$ with $a$ an entry of $\Tlam$ in row $r$ of $B$
  and $b$ one in row $r+1$ of $B$, and $\op{Inv}(w)\subsetneq\op{Inv}(d_{\bG_\alpha})$ for every other
  $w\in\Sym^J_\bnu$, since $w$ preserves the order of the entries in each row of $B$.

  \autoref{GE3}: The rows of $d\Tlam$ increase, but $d\Tlam$ is not standard. Hence there is a
  node $\beta=(k,r,c)\in\blam$ with $(k,r+1,c)\in\blam$ such that the entry of $d\Tlam$ in
  $(k,r,c)$ exceeds the one in $(k,r+1,c)$. As entries of $d\Tlam$ increase along rows, $d(a)>d(b)$ for all
  entries $a$ of $\Tlam$ in row $r$ and $b$ in row $r+1$ of the Garnir belt of $\bG_\beta$. Thus
  $\op{Inv}(d_{\bG_\beta})\subseteq\op{Inv}(d)$, and $u=dd_{\bG_\beta}^{-1}$ satisfies
  $\ell(u)=\ell(d)-\ell(d_{\bG_\beta})$.
\end{proof}

\begin{Remark}\label{R:products}
  For $v,w\in\Sym$, the product $H_vH_w$ is a sum of elements $\rho H_x=H_x\leftindex^{x^{-1}}{\rho}$ with $\rho\in\Tos$ and $\ell(x)\le\ell(v)+\ell(w)$, by \autoref{L:HiHw}.
 Hence $\xipm_vH_v\bb\xipm_wH_w=\xipm_v\leftindex^v(\bb\xipm_w)H_vH_w$ is a sum of elements $\bc H_x=\xipm_xH_x\leftindex^{x^{-1}}{(\xipm_x\bc)}^{-1}$  with $\bc\in\Bo$ and $\ell(x)\le\ell(v)+\ell(w)$.
\end{Remark}

Combining \autoref{L:HJ}, \autoref{L:GarnirElement} and \autoref{R:products}, we can straighten
products with $H_w$ on the left, and obtain a Garnir relation for multiplication on the right.
For $\blam\in\Parts$ we call $y\in\Bo$ \emph{pure} if $y=\pm hq$ with $h\in\hlam[\ilam]$,
$q\in\qlam[\ilam]$. Every $z\in\Sym_{|\blam|}$ fixes $\ilam$, so $\leftindex^zp\in\pm\plam[\ilam]$
for $p\in\plam[\ilam]$, and $\leftindex^zy$ is pure if $y$ is pure.

Every $z\in\Sym_{|\blam|}$ fixes $\ilam$, so $\leftindex^zp\in\pm\plam[\ilam]$
for $p\in\plam[\ilam]$, and $\leftindex^zy$ is pure if so is $y$.
\begin{Proposition}\label{P:Straighten}
  Let $\blam\in\Parts$ and assume that
  \begin{equation}\label{E:IH}
    \HH^{\pm\ge\bmu}_n\text{ is a two-sided ideal of }\Hfind\text{ for all }\bmu>\blam.
  \end{equation}
  Let $\bS\in\Std(\blam)$, $p\in\plam[\ilam]$ and $y\in\Bo$ pure.
  \begin{enumerate}
  \item\label{Straighten1} For $w\in\Sym$ there are $r_{\bU,p'}\in\Tos$, which do not depend on $y$, such that
    \[
      H_wp\HJpm_\blam y\equiv\sum_{\bU\in\Std(\blam),\ p'\in\plam[\ilam]}
        r_{\bU,p'}H^\pm_\bU p'\HJpm_\blam y\pmod{\HH^{\pm>\blam}_n}.
    \]
  \item\label{Straighten2} For $x\in\Sym$ there are $\bc_{\bU,h',q'}\in\Bo$, which do not depend on $\bS$ and $p$, such that
    \[
      H^\pm_\bS p\HJpm_\blam y\xipm_xH_x\equiv
      \sum_{\substack{\bU\in\Std(\blam)\\h'\in\hlam[\ilam],\ q'\in\qlam[\ilam]}}
        H^\pm_{\bS,p,h',q',\bU}\bc_{\bU,h',q'}\pmod{\HH^{\pm>\blam}_n}.
    \]
  \end{enumerate}
\end{Proposition}

\begin{proof}   Suppose that $(\alpha,\bG_\alpha,\bnu,\bmu,J)$ is a Garnir datum for $\blam$, and let $g=d_{\bG_\alpha}$. Let $\sigma\in\Sym_{|\blam|}$ be of minimal length with $\sigma\Sym_\bmu\sigma^{-1}=\Sym_\bnu$. Then $H_\sigma\HJpm_\bmu=\HJpm_\bnu H_\sigma$ and $\bb H_\sigma=H_\sigma\leftindex^{\sigma^{-1}}\bb$ for $\bb\in\Bo$. As $\sigma$ fixes $\ilam=\bi_\bmu$, $H^\pm_{\Tlam[\bmu]}=1=H^{\pm*}_{\Tlam[\bmu]}$ and $\bmu>\blam$ by \autoref{L:GarnirElement}\autoref{GE1}. Then assumption \autoref{E:IH} gives
  \begin{equation}\label{E:Garnirmu}
    p'\HJpm_\bnu y'=H_\sigma\left(\leftindex^{\sigma^{-1}}{{p'}}\HJpm_\bmu\leftindex^{\sigma^{-1}}{{y'}}\right)H_\sigma^{-1}
    \in\HH^{\pm\ge\bmu}_n\subseteq\HH^{\pm>\blam}_n,
  \end{equation}
 for $p'\in\pm\plam[\ilam]$ and pure $y'$. Let
  \[
    Z=\sum_{x\in\leftindex^J_\blam\Sym}\hbar^{-2\ell(x)}\xipm_xH_x\qquad\text{and}\qquad
    Z^*=\sum_{x\in\Sym^J_\blam}\hbar^{-2\ell(x)}\xipm_xH_x.
  \]
  Applying \autoref{L:HJK} to $J\subseteq K$ with $\Sym_K=\Sym_\blam$, respectively $\Sym_K=\Sym_\bnu$, and using $H_x\leftindex^{x^{-1}}{{\xipm_x}}=\xipm_xH_x$ and $\leftindex^J_\bnu\Sym=\set{f^{-1}|f\in\Sym^J_\bnu}$, we obtain
  \[
    \HJpm_\blam=Z^*\HJpm_J=\HJpm_JZ\qquad\text{and}\qquad
    \HJpm_\bnu=\sum_{f\in\Sym^J_\bnu}\hbar^{-2\ell(f)}\xipm_fH_f\HJpm_J
      =\sum_{f\in\Sym^J_\bnu}\hbar^{-2\ell(f)}\HJpm_J\xipm_{f^{-1}}H_{f^{-1}}.
  \]
  Multiplying the two expressions for $\HJpm_\bnu$ by $Z$ on the right, respectively by $Z^*$ on the left, and separating the summand for $g\in\Sym^J_\bnu$, see \autoref{L:GarnirElement}\autoref{GE2}, gives
  \begin{align}
    \xipm_gH_g\HJpm_\blam&=\hbar^{2\ell(g)}\Bigl(\HJpm_\bnu Z
      -\sum_{f\in\Sym^J_\bnu,\ f\ne g}\hbar^{-2\ell(f)}\xipm_fH_f\HJpm_\blam\Bigr),\label{E:Garnirsum}\\
    \HJpm_\blam\xipm_{g^{-1}}H_{g^{-1}}&=\hbar^{2\ell(g)}\Bigl(Z^*\HJpm_\bnu
      -\sum_{f\in\Sym^J_\bnu,\ f\ne g}\hbar^{-2\ell(f)}\HJpm_\blam\xipm_{f^{-1}}H_{f^{-1}}\Bigr).
      \label{E:Garnirsumright}
  \end{align}

  \autoref{Straighten1}: We argue by induction on $\ell(w)$. Write $w=dv$ with $d\in\Sym^\blam$, $v\in\Sym_\blam$ and $\ell(w)=\ell(d)+\ell(v)$. By \autoref{L:HJ}\autoref{HJ2}, \autoref{slide2} and \autoref{P:HW}, $H_vp\HJpm_\blam$ is a $\Tos$-linear combination of elements $p''\HJpm_\blam$ with $p''\in\pm\plam[\ilam]$, and $H_d\rho=\leftindex^d{\rho}H_d$ for $\rho\in\Tos$. So we may assume $w=d$. If $d\Tlam$ is standard, then $H_d=\hbar^{2\ell(d)}(\xipm_d)^{-1}H^\pm_{d\Tlam}$, and we are done. This is always the case if $\blam$ has no Garnir datum, that is, if every component of $\blam$ has at most one row. Otherwise, by \autoref{L:GarnirElement}\autoref{GE3}, $d=ug$ with $g=d_{\bG_\alpha}$ for a Garnir datum and $\ell(d)=\ell(u)+\ell(g)$, so $H_dp\HJpm_\blam y=H_uH_gp\HJpm_\blam y$. By \autoref{L:HJ}\autoref{HJ2} and \autoref{slide2}, $H_gp$ is a sum of elements $\rho p''\xipm_gH_g$ with $\rho\in\Tos$ and $p''\in\pm\plam[\ilam]$, and we insert \autoref{E:Garnirsum}. The term with $\HJpm_\bnu Zy$ lies in $\HH^{\pm>\blam}_n$ by \autoref{L:HJ}\autoref{slide2} and \autoref{E:Garnirmu}. In the other terms, $p''\xipm_fH_f$ is a sum of elements $\rho'\xipm_fH_fp'$ with $\rho'\in\Tos$ and $p'\in\pm\plam[\ilam]$ by \autoref{L:HJ}\autoref{slide2}, and $H_u\rho\rho'\xipm_fH_f$ is a sum of elements $\rho''H_x$ with $\rho''\in\Tos$ and $\ell(x)\le\ell(u)+\ell(f)<\ell(d)$ by \autoref{R:products} and \autoref{L:GarnirElement}\autoref{GE2}. So the claim follows by induction. As $y$ is never changed, the coefficients do not depend on $y$.

  \autoref{Straighten2}: Let $P(y,x)=H^\pm_\bS p\HJpm_\blam y\xipm_xH_x$. We argue by induction on $\ell(x)$. Write $x=vd^{-1}$ with $v\in\Sym_\blam$, $d\in\Sym^\blam$ and $\ell(x)=\ell(v)+\ell(d)$. Then $\xipm_xH_x=\xipm_vH_v\xipm_{d^{-1}}H_{d^{-1}}$ by \autoref{C:hdv}.

  \emph{Case $v\ne1$.} By \autoref{L:HJ}, $\HJpm_\blam y\xipm_vH_v$ is a sum of elements $\HJpm_\blam y'\bc$ with $y'$ pure and $\bc\in\Tos_\blam$. Moving $\bc$ to the right of $\xipm_{d^{-1}}H_{d^{-1}}$ by \autoref{R:products}, we see that $P(y,x)$ is a sum of elements $P(y',d^{-1})\bc'$ with $\bc'\in\Bo$. As $\ell(d)<\ell(x)$, the claim follows by induction.

  \emph{Case $v=1$ and $d\Tlam$ standard.} Then $\xipm_{d^{-1}}H_{d^{-1}}=\hbar^{2\ell(d)}H^{\pm*}_{d\Tlam}$. As $y$ is pure, $y=\pm hq$ with $h\in\hlam[\ilam]$ and $q\in\qlam[\ilam]$, so $P(y,d^{-1})=\pm\hbar^{2\ell(d)}H^\pm_{\bS,p,h,q,d\Tlam}$ is of the required form with $\bU=d\Tlam$.

  \emph{Case $v=1$ and $d\Tlam$ not standard.} By \autoref{L:GarnirElement}\autoref{GE3}, $d=ug$ with $g=d_{\bG_\alpha}$ for a Garnir datum and $\ell(d)=\ell(u)+\ell(g)$, so $P(y,d^{-1})=P(y,g^{-1})\xipm_{u^{-1}}H_{u^{-1}}$. As $g\in\Sym_{|\blam|}$, we can move $y$ to the right of $\xipm_{g^{-1}}H_{g^{-1}}$ by \autoref{L:HJ}\autoref{slide2}, apply \autoref{E:Garnirsumright}, and move the pure elements back to the left of the $\xipm_{f^{-1}}H_{f^{-1}}$. The term with $Z^*\HJpm_\bnu$ lies in $\HH^{\pm>\blam}_n$ by \autoref{E:Garnirmu}. Hence $P(y,g^{-1})$ is congruent to a sum of elements $P(y',f^{-1})\bc$ with $f\in\Sym^J_\bnu\setminus\{g\}$, $y'$ pure and $\bc\in\Tos_\blam$. By \autoref{R:products}, $P(y,d^{-1})$ is then congruent to a sum of elements $P(y',z)\bc'$ with $\ell(z)\le\ell(f)+\ell(u)<\ell(d)$, and the claim follows by induction.

  The elements of $\Bo$ that occur do not depend on $\bS$ and $p$.
\end{proof}

We collect four facts needed for the inductive step. For $\bj\in I^n$ let $\bom_\bj=(1^{m_1}|\dots|1^{m_N})\in\Parts$, where $m_r$ is the number of entries of $\bj$ equal to $r$.
\begin{enumerate}[label=\upshape(F\arabic*), ref=\upshape(F\arabic*)]
\item\label{F1} By \autoref{P:FrobeniusRespect}, the elements $p_\bj h_\bj q_\bj$ with $\bj\in I^n$, $p_\bj\in\plam[\bj]$, $h_\bj\in\hlam[\bj]$ and $q_\bj\in\qlam[\bj]$ form a respectable basis of $\Bo$, where $I^n$ is ordered componentwise. Condition \autoref{resp3} holds here trivially, as $\Olam=\kk$.
\item\label{F2} For $\blam\in\Parts$ let $\Bo^{>\ilam}$ be the span of the basis elements $p_\bj h_\bj q_\bj$ with $\bj>\ilam$. Every $z\in\Sym_\blam$, see \autoref{D:Wlam}, fixes $\ilam$ from \autoref{fancydatum}, hence $\leftindex^z{(\Bo^{>\ilam})}=\Bo^{>\ilam}$.
\item\label{F3} Every basis element $p_\bj h_\bj q_\bj$ of $\Bo$ lies in $\Hfind\Basis_{\bom_\bj}\Hfind$. Indeed, $\Sym_{\bom_\bj}$ is trivial, so $\HJpm_{\bom_\bj}=1=H^\pm_{\Tlam[\bom_\bj]}=H^{\pm*}_{\Tlam[\bom_\bj]}$ by \autoref{D:HJ}, and $H^\pm_{\Tlam[\bom_\bj],p_1,h_1,q_1,\Tlam[\bom_\bj]}=p_1h_1q_1$.
\noindent Choose $x\in\Sym$ such that it  permutes the entries of $\bj$ into the increasing order $\bi_{\bom_\bj}$. As $x$ permutes the tensor factors of $\Bo$, $\leftindex^x{(p_\bj h_\bj q_\bj)}=\pm p_1h_1q_1$ with $p_1\in\plam[\bi_{\bom_\bj}]$, $h_1\in\hlam[\bi_{\bom_\bj}]$, $q_1\in\qlam[\bi_{\bom_\bj}]$. We obtain
\begin{equation*}
  p_\bj h_\bj q_\bj=H_x^{-1}\leftindex^x{(p_\bj h_\bj q_\bj)}H_x
    =\pm H_x^{-1}H^\pm_{\Tlam[\bom_\bj],p_1,h_1,q_1,\Tlam[\bom_\bj]}H_x
    \in\Hfind\Basis_{\bom_\bj}\Hfind.
\end{equation*}
\item\label{F4} $\Hfind$ is generated by $H_1,\dots,H_{n-1}$ and $\Bo$, since $X_1$ is a scalar by the cyclotomic relation of level one, see \autoref{D:CyclotomicHecke}, and $\hbar^2X_{i+1}=H_iX_iH_i+tH_i\taui[i+1,i]$ by \autoref{L:HXH}.
\end{enumerate}

\begin{Lemma}\label{L:InductiveStep}
  Let $\blam\in\Parts$ satisfy \autoref{E:IH}, and let $\bS,\bT\in\Std(\blam)$, $p\in\plam[\ilam]$,
  $h\in\hlam[\ilam]$, $q\in\qlam[\ilam]$, $\bb\in\Bo$ and $1\le m<n$. In the following, congruences
  are modulo $\HH^{\pm>\blam}_n$, and the sums run over $\bU\in\Std(\blam)$, $p'\in\plam[\ilam]$,
  $h'\in\hlam[\ilam]$ and $q'\in\qlam[\ilam]$.
  \begin{enumerate}
  \item\label{E:BoInIdeal} $\Bo^{>\ilam}\subseteq\HH^{\pm>\blam}_n$.
  \item\label{E:sandwich} $\ba\HJpm_\blam\bc\in\HH^{\pm>\blam}_n$ for all $\ba,\bc\in\Bo$ with
    $(\leftindex^z{\ba})\bc\in\Bo^{>\ilam}$ for all $z\in\Sym_\blam$.
  \item\label{E:BoAction}
    $\bb H^\pm_{\bS,p,h,q,\bT}\equiv\sum_{p'}r_{p'}H^\pm_{\bS,p',h,q,\bT}$ and
    $H^\pm_{\bS,p,h,q,\bT}\bb\equiv\sum_{q'}s_{q'}H^\pm_{\bS,p,h,q',\bT}$.
  \item\label{E:HmAction}
    $H_mH^\pm_{\bS,p,h,q,\bT}\equiv\sum_{\bU,p'}r_{\bU,p'}H^\pm_{\bU,p',h,q,\bT}$ and
    $H^\pm_{\bS,p,h,q,\bT}H_m\equiv\sum_{\bU,h',q'}s_{\bU,h',q'}H^\pm_{\bS,p,h',q',\bU}$.
  \end{enumerate}
  Here all coefficients lie in $\kk$, the $r_{p'}$ and $r_{\bU,p'}$ do not depend on $h$, $q$ and
  $\bT$, and the $s_{q'}$ and $s_{\bU,h',q'}$ do not depend on $\bS$ and $p$.
\end{Lemma}

\begin{proof}
  \autoref{E:BoInIdeal}: Let $\bj>\ilam$. Then $j_s\le i^\blam_s$ for all $s$, using the order on $I^n$. Thus,  $m_1+\dots+m_k\ge|\lambda^{(1)}|+\dots+|\lambda^{(k)}|$ for $1\le k\le N$.
  If all these were equalities, then $\bj$ and $\ilam$ would have the same entries, and hence
  $\bj=\ilam$. Thus $\bom_\bj>\blam$ by \autoref{D:contorder}. Then, \autoref{F3} with
  \autoref{E:IH} gives $p_\bj h_\bj q_\bj\in\HH^{\pm>\blam}_n$.

  \autoref{E:sandwich}: By \autoref{D:HJ} and \autoref{C:axiH}, we have
  \[
    \HJpm_\blam=\sum_{w\in\Sym_\blam}\hbar^{-2\ell(w)}\xipm_wH_w\qquad\text{and}\qquad
    \ba\xipm_wH_w=\sum_{K}(-1)^{|K|}y^\pm_{K,w}H_w\phi_{K}(\ba)\quad\text{for }w\in\Sym_\blam.
  \]
  Since each $\phi_{K}(\ba)$ is a $\Z$-linear combination of elements $\leftindex^z{\ba}$ with $z\in\Sym_\blam$,
  $\ba\HJpm_\blam\bc\in\Hfind\Bo^{>\ilam}\subseteq\HH^{\pm>\blam}_n$ by
  \autoref{E:BoInIdeal}.

  \autoref{E:BoAction}: By \autoref{L:HJ}\autoref{slide4}, $\bb H^\pm_\bS=H^\pm_\bS\bb'$ and
  $H^{\pm*}_\bT\bb=\bb''H^{\pm*}_\bT$ for some $\bb',\bb''\in\Bo$. By \autoref{resp2} for $\Bo$, there
  are $r_{p'},s_{q'}\in\kk$, depending only on $\bb',p$ and on $\bb'',h,q$, respectively, such that
  $\ba=\bb'p-\sum_{p'}r_{p'}p'$ and $\bc=hq\bb''-\sum_{q'}s_{q'}hq'$ satisfy $\ba y\in\Bo^{>\ilam}$
  and $p''\bc\in\Bo^{>\ilam}$ for all pure $y$ and all $p''\in\pm\plam[\ilam]$. As
  $\leftindex^z{(hq)}$ is pure, $\leftindex^zp\in\pm\plam[\ilam]$ and
  $\leftindex^z{(\Bo^{>\ilam})}=\Bo^{>\ilam}$ for $z\in\Sym_\blam$, \autoref{E:sandwich} applies to $\ba$ and $hq$, and to $p$ and $\bc$. Since
  $\bb H^\pm_{\bS,p,h,q,\bT}=H^\pm_\bS\bb'p\HJpm_\blam hqH^{\pm*}_\bT$ and
  $H^\pm_{\bS,p,h,q,\bT}\bb=H^\pm_\bS p\HJpm_\blam hq\bb''H^{\pm*}_\bT$, this gives the claim.

  \autoref{E:HmAction}: As $H_mH^\pm_\bS=\hbar^{-2\ell(d_\bS)}\leftindex^{s_m}{{\xipm_{d_\bS}}}H_mH_{d_\bS}$
  is a $\Tos$-linear combination of elements $H_w$ by \autoref{R:products},
  \autoref{P:Straighten}\autoref{Straighten1} shows that $H_mH^\pm_{\bS,p,h,q,\bT}$ is congruent to a
  sum of elements $\rho H^\pm_{\bU,p',h,q,\bT}$ with $\rho\in\Tos$ not depending on $h$, $q$ and $\bT$.
  Similarly, $H^{\pm*}_\bT H_m=\hbar^{-2\ell(d_\bT)}\xipm_{d_\bT^{-1}}H_{d_\bT^{-1}}(\xipm_m)^{-1}\xipm_mH_m$
  is a sum of elements $\xipm_xH_x\bc$ with $\bc\in\Bo$ by \autoref{R:products}, so by
  \autoref{P:Straighten}, $H^\pm_{\bS,p,h,q,\bT}H_m$ is congruent to a sum of elements
  $H^\pm_{\bS,p,h',q',\bU}\bc'$ with $\bc'\in\Bo$ not depending on $\bS$ and $p$. Now
  \autoref{L:InductiveStep}\autoref{E:BoAction} gives the claim. Here $h'\ne h$ occurs in general, see
  \autoref{R:RightAction}, which is why $h$ belongs to the right hand index in \autoref{fancydatum}.
\end{proof}

\begin{proof}[Proof of \autoref{T:HfindRespectable}]
  Using the order of \autoref{D:contorder}, we show by downward induction on $\blam$ that $\HH^{\pm\ge\blam}_n$ is a two-sided ideal and that both congruences of
  \autoref{resp2} hold for $\Basis_\blam$. So let $\blam\in\Parts$ satisfy \autoref{E:IH}, which is
  empty if $\blam$ is maximal. By induction, \autoref{L:InductiveStep}\autoref{E:BoAction} and
  \autoref{E:HmAction} apply and show that  $\HH^{\pm\ge\blam}_n$ is stable under multiplication by the generators
  $H_m$ and $\bb\in\Bo$ of $\Hfind$ from both sides, and both congruences of \autoref{resp2} hold for
  $\Basis_\blam$. Every basis element
  $p_\bj h_\bj q_\bj$ of $\Bo$ lies in $\HH^{\pm\ge\bom_\bj}_n$ by \autoref{F3}, so the
  identity of $\Hfind$ lies in the two-sided ideal $\sum_{\bj\in I^n}\HH^{\pm\ge\bom_\bj}_n$. This ideal
  is therefore $\Hfind$ and $\Basis$ spans $\Hfind$ and  is even a basis by \autoref{L:rank}. This shows
  \autoref{resp1}. By \autoref{R:associativity} and assumption, \autoref{resp3} is automatic.
  \end{proof}

  \subsection{Application: respectable bases built from the Murphy basis of $\kk W$}
To illustrate how to construct new respectable bases out of old ones, we
consider  \autoref{T:HfindRespectable} in the special case when $\Hfind[n;1,0,0](\kk)=\kk W$. Then
$H^+_{\bS,p,h,q,\bT}=H^+_{\bS,1,1,1,\bT}$ from $\eqref{E:HsphqT}$, with
$\bS,\bT\in\Std(\lambda)$ and $\lambda\in\Parts[1,n]$, is the Murphy basis
element $m_{\bS\bT}$ from \cite[Theorem~4.17]{Murphy:basis}.

Let $\kk[\X]^\Sym\subset \kk[\X]$  be the subalgebra of symmetric polynomials and  $\kk[\X]^\Sym_+$  its augmentation ideal.

 Let $\mathcal{C}=\kk[\X]/(\kk[\X]^\Sym_+)$ be the algebra of
 coinvariants. The $\X$-grading on $\kk[\X]$, with $X_i$ in degree~$1$
 for $1\leq i\leq n$, induces a grading
 $\mathcal{C}=\oplus_{d\in\N_0}\mathcal{C}_d$ on $\mathcal{C}$. Let
 $\cB_d$ be a homogeneous $\kk$-basis of $\mathcal{C}_d$.  Consider the
 quotient smash product algebra $\mathcal{C}\Smash\kk\Sym$ of
 $\kk[\X]\Smash\kk\Sym$, with the  smash product relation
 $wX_i=X_{w(i)}w$ for $1\leq i\leq n$ and $w\in\Sym$.

Recall from \autoref{SS:Tableaux} that $\Parts[1,n]$ is the set of
partitions of~$n$. Let $\Lambda=\set{(d,\lambda)|\mathcal{C}_d\neq0,
\lambda\in \Parts[1,n]}$ be the poset with ordering
$(d,\lambda)>(d',\mu)$ if $d>d'$, or $d=d'$ and
$\lambda\vartriangleright\mu$. Following the notation of
\autoref{T:Respectable}, let $\Parts[1,n,+]=(\Parts[1,n])_+$ be the
indexing set for the simple $\kk W$-modules.

\begin{Proposition}\label{P:coinvariantsmash}
  The algebra $\CA=\mathcal{C}\Smash\kk\Sym$ is perfectly
  $\kk$-respectable with weight poset $\Lambda$ and
  \[
    \plam[(d,\lambda)]=\cB_d\times\Std(\lambda),\qquad
    \hlam[(d,\lambda)]=\set{1},\qquad
    \olam[(d,\lambda)]=\set{1},\qquad
    \qlam[(d,\lambda)]=\Std(\lambda),
  \]
  and $a_{(f,\bS),1,1,\bT}=fm_{\bS\bT}$.
  Moreover, $\Lambda_+=\set{(0,\lambda)|\lambda\in\Parts[1,n,+]}$.
  \end{Proposition}

\begin{proof}
   The basis of $\CA$ is perfect by definition, \autoref{resp3} is automatic
   since $|\hlam|=1$ for $\lambda\in\Lambda$, and \autoref{resp1} holds
   because $\set{m_{\bS,\bT}}$ is a perfectly $\kk$-respectable basis of
   $\kk\Sym$.  Hence, it remains to check \autoref{resp2}. Now,
   \[
    \CA^{>(d,\lambda)}=\mathcal{C}_{>d}\otimes\kk\Sym+\mathcal{C}_d\otimes\kk\Sym^{\vartriangleright\lambda},
    \qquad\text{ for }(d,\lambda)\in\Lambda,
   \]
   is a two-sided ideal since the Murphy basis is respectable and
   multiplication with $X_i$ on either side increases the $\X$-degree
   $d$. Now  \autoref{resp2}  follows, since $X_i$ acts by zero modulo
   $\CA^{>(d,\lambda)}$. Since the action of the~$X_i$ is zero, modulo
   higher terms, the bilinear form vanishes on all strata
   $\CA^{(d,\lambda)}$ with $d\ge1$. On the strata $(0,\lambda)$, the
   form on the costandard module coincides with the form of the Murphy
   basis on $\kk \Sym$, so
   $\Lambda_+=\set{(0,\lambda)|\lambda\in\Parts[1,n,+]}$.
\end{proof}

The construction in \autoref{P:coinvariantsmash} generalises to the
restricted rational Cherednik algebras with \emph{coupling parameter}
$c=0$. Let $\kk[\mathbbm{x}]=\kk[x_1,\dots,x_n]$ and
$\kk[\mathbbm{y}]=\kk[y_1,\dots,y_n]$, for indeterminates $x_i$ and
$y_j$. As a vector space, the \emph{rational Cherednik algebra} is
$\mathbf{H}_c=\mathbf{H}_{0,c}
   =\kk[\mathbbm{x}]\otimes\kk\Sym\otimes\kk[\mathbbm{y}]$,
where the algebra structure is given by requiring each tensor factor to
be a subalgebra, and where the defining relations between the tensor
factors are:
\begin{equation}\label{ratChrel}
  wx_kw^{-1}=x_{w(k)},\qquad  wy_kw^{-1}=y_{w(k)},\qquad   [x_k,y_l]=-c(k,l),
  \qquad [x_k,y_k]=c\sum_{j\neq k}(k,j),
\end{equation}
for $w\in\Sym$, $1\leq k$ and $l\leq n$.
The bi-invariants $\mathcal
Z=\kk[\mathbbm{x}]^{\Sym}\otimes\kk[\mathbbm{y}]^{\Sym}$ are central, and the \emph{restricted rational Cherednik algebra} is
$
  \overline{\mathbf H}_c=\mathbf H_{c}/\mathcal Z_+\mathbf H_{c},
$
where $\mathcal Z_+$ is the augmentation ideal of $\mathcal Z$. The
algebra~$\overline{\mathbf{H}}_c$ is free over $\kk$ of rank $(n!)^3$,
with basis $\set{fwg}$, where  $w\in\Sym$ and $f$ and $g$ run through
bases of the coinvariant algebras $\mathcal{C}^{\mathbbm{x}}$ and
$\mathcal{C}^{\mathbbm{y}}$, respectively, which are attached to the
polynomial rings; see \cite[Section~3.1]{BGS:Frobenius}.

For each $d\in\N_0$ fix homogeneous bases $\cB^{\mathbbm{x}}_d$ of
$\mathcal{C}^{\mathbbm{x}}_d$ and $\cB^{\mathbbm{y}}_d$ of
$\mathcal{C}^{\mathbbm{y}}_d$ respectively. The algebra
$\mathcal{C}^{\mathbbm{x}}\otimes\mathcal{C}^{\mathbbm{y}}$ is bigraded
with $(\mathcal{C}^{\mathbbm{x}}\otimes
\mathcal{C}^{\mathbbm{y}})_{d,e}=\mathcal{C}_d^{\mathbbm{x}}\otimes
\mathcal{C}^{\mathbbm{y}}_e$, and conjugation by $\Sym$ preserves the
bidegree. Now let
$\Lambda=\set{(d,e,\lambda)|\mathcal{C}^{\mathbbm{x}}_d\neq0,
\mathcal{C}^{\mathbbm{y}}_e\neq0, \lambda\in \Parts[1,n]}$, ordered by
$(d,e,\lambda)>(d',e',\mu)$ if $d+e>d'+e'$, or $(d,e)=(d',e')$ and
$\lambda\vartriangleright\mu$.

\begin{Proposition}\label{P:restrictedCherednik}
  The algebra $\overline{\mathbf H}_0$ is perfectly $\kk$-respectable with
  weight poset $(\Lambda,\ge)$, index sets
  \[
    \plam[(d,e,\lambda)]=\cB^x_d\times\Std(\lambda),\qquad
    \hlam[(d,e,\lambda)]=\set{1},\qquad
    \olam[(d,e,\lambda)]=\set{1},\qquad
    \qlam[(d,e,\lambda)]=\Std(\lambda)\times\cB^y_e,
  \]
  for $(d,e,\lambda)\in\Lambda$, and $a_{(f,\bS),1,1,(\bT,g)}=fm_{\bS,\bT}g$.
The datum is perfect with $\dot h$ as in the Murphy basis, $\dot p=fd_\bS^{-1}$, $\dot o=1$, $\dot q=d_\bT g$ for $p=(f,\bS)$ and $q=(\bT,g)$.

  Moreover, the form vanishes on all strata $(d,e,\lambda)$ with
  $d+e\ge1$, and on the strata $(0,0,\lambda)$ it is the Murphy form. In
  particular, $\Lambda_+=\set{(0,0,\lambda)|\lambda\in\Parts[1,n,+]}$.
  \end{Proposition}

\begin{proof}
  For axiom \autoref{resp1} and perfectness argue as in \autoref{P:coinvariantsmash}.  Let $(d,e,\lambda)\in\Lambda$. Then
  \begin{equation}\label{MFO2}
    \CA^{>(d,e,\lambda)}=\sum_{d'+e'>d+e}\mathcal{C}^{\mathbbm{x}}_{d'}\otimes \kk\Sym\otimes \mathcal{C}^{\mathbbm{y}}_{e'}+\mathcal{C}^{\mathbbm{x}}_d\kk\Sym^{\vartriangleright\lambda}\mathcal{C}^{\mathbbm{y}}_e.
  \end{equation}
  We first compute the action of the generators $w\in\Sym$ and $x_k,y_k$ for $1\le k\le n$ on the elements $fm_{\bS\bT}g$ with $f\in\cB^\mathbbm{x}_d$ and $g\in\cB^\mathbbm{y}_e$. Obviously, multiplication by $x_k$ from the left or by $y_k$ from the right increases the degree and hence is zero modulo \autoref{MFO2}.
 Using \autoref{ratChrel},  multiplication with $w\in\Sym$ can be moved
 through the $\mathcal{C}^{\mathbbm{x}}$ part while not changing the
 degree, then it acts as desired, since the Murphy basis is respectable.
 Similarly for the action of $w$ from the right.
  Write $m_{\bS\bT}=\sum_{v\in\Sym}c_vv$. We have $y_kfm_{\bS\bT}g=fy_k m_{\bS\bT}g=
    f\sum_{v\in\Sym}c_vv \leftindex^{v^{-1}}y_kg\in
    \mathcal{C}^{\mathbbm{x}}_d \otimes\kk\Sym
    \mathcal{C}^{\mathbbm{y}}_{e+1}$ which is zero modulo
    \autoref{MFO2}.  For degree reasons, the form is identically zero on
    the strata with $d+e\ge1$. For the strata with $d=e=0$, the form is
    just the form of the Murphy basis, so the  last claim follows.
\end{proof}
  By \cite[Section~3.1]{BGS:Frobenius}, $\overline{\mathbf H}_0$ is Frobenius, and thus a perfectly $\kk$-respectable Frobenius algebra.
\subsection{Respectability for higher level CRaFH algebras}
\label{SS:HigherLevelRespectable}
This subsection contains the main application of the paper. Assuming $\Bo=A^{\otimes n}$ is a Frobenius-tensor algebra with $A$ is perfectly $\kk$-respectable, \autoref{T:HhdRespectable} below shows that the CRaFH algebras $\Hhd(\Bo)$  of arbitrary level~$\ell\ge1$ are perfectly $\kk$-respectable with weight poset the set $\lParts$ of $\ell N$-multipartitions of~$n$. Here, $N=|I|$ is the number of strata of the datum of $A$. This result generalises \autoref{T:HfindRespectable} to arbitrary level. The technical core of the proof was already settled there. The key new ingredient is \autoref{D:Xlam}, which mimics \cite[Definition~3.1]{DJM:cyc}.

Recall from \autoref{D:CyclotomicHecke} that $\charge=(\kapp_1,\dots,\kapp_\ell)\in\kk^\ell$ is the charge of $\Hhd$, where $\theta=1$ as $\psi$ is trivial, so that the cyclotomic polynomial is $\prod_{k=1}^\ell(X_1-\kapp_k)$ with the first $\ell_0$ entries of $\charge$ equal to zero.

\subsubsection*{The combinatorics} We first treat the combinatorics behind the main result. For this recall that for cyclotomic Hecke algebras of type~$A$, $\ell$-multipartitions index the cell representations, see \cite{DJM:cyc}. Now \autoref{T:HfindRespectable} and \autoref{T:simples} suggest that the respectable basis of $\Hhd$ should be indexed by $\ell$-tuples of $N$-multipartitions, or alternatively by the set $\lParts$ of $\ell N$-multipartitions of~$n$, which is in bijection with this set. This second realisation allows us to mimic the arguments from the level $1$ case.

\begin{Definition}\label{D:lcontorder}
We consider the set $\lParts$ of $\ell N$-multipartitions of~$n$ as a poset with the partial order of \autoref{D:contorder}, with $N$ replaced by $\ell N$. In particular $\bmu>\blam$ whenever $\bmu$ is obtained from $\blam$ by moving a box to an earlier component.
\end{Definition}

As in level one, to each $\blam\in\lParts$ we associate the Young subgroup $\Sym_\blam$ of \autoref{D:Wlam} and the corresponding parabolic CRaFH module $\MJ^\pm_\blam=\Hhd\HJ^\pm_\blam$ with $\HJ^\pm_\blam$ as in \autoref{leveloneparab}, and the set $\Std(\blam)$ of standard tableaux of \autoref{D:tableaux}.

\begin{Definition}\label{D:chargeparts}
  Suppose that $\blam=(\lambda^{(1)}|\dots|\lambda^{(\ell
  N)})\in\lParts$. Write $\blam=[\blam^{[1]}\|\dots\|\blam^{[\ell]}]$,
  where $\blam^{[l]}=(\lambda^{((l-1)N+1)}|\dots|\lambda^{(lN)})$ is the
  \emph{$\kapp_l$-part} of $\blam$. For $i\in I$, {the components $\lambda^{(lN+i)}$ with $0\le l<\ell$ form} the \emph{$i$-layer} of
  $\blam$. Reading the components in order from left to right, and
  recording {for each box the layer of its component}, define
  \[
    \ilam = (i^\blam_1,\dots,i^\blam_n)=(
    \underbrace{1,\dots,1}_{|\lambda^{(1)}|\text{ times}},\dots,
    \underbrace{N,\dots,N}_{|\lambda^{(N)}|\text{ times}},
    \dots,
    \underbrace{1,\dots,1}_{\mathclap{|\lambda^{(\ell N-N+1)}|\text{ times}}},\dots,
    \underbrace{N,\dots,N}_{|\lambda^{(\ell N)}|\text{ times}}
    )\in I^n.
  \]
  For {$1\le l\le\ell+1$} let $n_l=\sum_{r\le(l-1)N}|\lambda^{(r)}|$ be
  the number of boxes in~$\blam$ before the $\kapp_l$-part. In
  particular, $n_1=0$, and $n_{\ell+1}=n$. {For $1\le l\le\ell$, the} indices $n_l+1,\dots,n_{l+1}$ of $\ilam$ are the \emph{positions of the $\kapp_l$-part}.
\end{Definition}
\begin{Definition}\label{D:Xlam}
  Let $\blam\in\lParts$. Define
  $X_\blam=\displaystyle\prod_{l=2}^{\ell}\prod_{m=1}^{n_l}(X_m-\kapp_l)$.
\end{Definition}
\begin{Remark}
By \autoref{L:XSymmetric}, $X_\blam \HJpm_\blam=\HJpm_\blam X_\blam${, since $\Sym_\blam$ preserves the positions of each $\kapp_l$-part}.
We use this identity below without mention.
\end{Remark}

\begin{Example}\label{Ex:chargeparts}
  Let $n=6$ and $\ell=3$, so that the cyclotomic polynomial is
  $(X_1-\kapp_1)(X_1-\kapp_2)(X_1-\kapp_3)$, and let $N=2$, so that
  $I=\set{1>2}$. We display an $\ell N$-multipartition by separating the
  different $\kapp_l$-parts by~$\|$, and by writing the charge above each
  component. For example, if
  $\blam=(\emptyset|(2,1)|(1)|\emptyset|(1)|(1))=[\emptyset|21\ \|\ 1|\emptyset\ \|\ 1|1]$, then
  \[
      \blam=\LabelledDiagram[
        label={\kapp_1|\kapp_1|\kapp_2|\kapp_2|\kapp_3|\kapp_3},
        tikz after={
          \draw[aTableauMain,thick]([xshift=5mm]A-2-1-2.north)--++(0,-1);
          \draw[aTableauMain,thick]([xshift=-6mm]A-5-1-1.north)--++(0,-1);
        }
      ]{|2,1|1||1|1}.
  \]
  The three boxes of the $\kapp_1$-part of~$\blam$ lie in its component of layer $2$, the $\kapp_2$-part has one box in its component of layer $1$, and the $\kapp_3$-part has one box in each of its components of layer $1$ and $2$. Hence $\ilam=(2,2,2,1,1,2)$ lists the {layers} of the six boxes of~$\blam$. Moreover $n_1=0$, $n_2=3$ and $n_3=4$, so the $\kapp_1$-part has positions $\set{1,2,3}$, the $\kapp_2$-part has position $\set{4}$, and the $\kapp_3$-part has positions $\set{5,6}$. By \autoref{D:Xlam},
  \[
    X_\blam=\prod_{m=1}^{3}(X_m-\kapp_2)\prod_{m=1}^{4}(X_m-\kapp_3).
  \]
\end{Example}

Recall from \autoref{P:FrobeniusRespect} that
$\set{(\plam[\bi],\hlam[\bi],\qlam[\bi])|\bi\in I^n}$ is a
respectable basis of {$\Bo$}. Abusing notation, we use the same
symbols for the basis elements of $\Hhd$ that we used for $\Hfind$.
There should be no ambiguity, however, because the elements of $\Hfind$
are indexed by $\Parts$, whereas the elements of $\Hhd$ are indexed by $\lParts$.
{For $\blam\in\lParts$ let $\ilam$ be as in \autoref{D:chargeparts}.}
Define
$\plam[\blam]^\pm=\set{H^\pm_\bS p|\bS\in\Std(\blam), p\in\plam[\ilam]}$,
$\hlam[\blam]^\pm=\set{\HJ^\pm_\blam X_\blam}$,
$\qlam[\blam]^\pm=\set{hqH^{\pm*}_\bT|\bT\in\Std(\blam), h\in\hlam[\ilam],q\in\qlam[\ilam]}$.
For convenience, set
$H^\pm_{\bS, p,h,q,\bT} = H^\pm_\bS p \HJ^\pm_\blam X_\blam h qH^{\pm*}_\bT$,
for $\bS,\bT\in\Std(\blam)$, $p\in\plam[\ilam]$, $h\in\hlam[\ilam]$,
and $q\in\qlam[\ilam]$. Let
\begin{equation}\label{E:ellBasis}
    \Basis = \bigcup_{\blam\in\lParts}\Basis_\blam,\qquad
    \text{where }\quad
    \Basis_\blam=\set[\big]{H^\pm_{\bS,p,h,q,\bT}|
        \bS,\bT\in\Std(\blam),
        p\in\plam[\ilam],
        h\in\hlam[\ilam],
        q\in\qlam[\ilam]
    }.
\end{equation}
\begin{Remark}
When $\ell=1$, the elements in $\Basis$ coincide with the corresponding
elements from \autoref{E:HsphqT}. In comparison with
\autoref{SS:LevelOneRespectable}, the elements in $\hlam$ contain an
additional factor $X_\blam$. This is consistent with the previous
notation because $X_\blam=1$, that is the empty product,  when $\ell=1$.
\end{Remark}
We start by showing that the set $\Basis$ contains the right number of
elements to be a basis of~$\Hhd$.

\begin{Lemma}\label{L:lrank}
  The algebra $\Hhd$ is free as a $\kk$-module of rank {$|\Basis|$}.
\end{Lemma}

\begin{proof}
  The proof is similar to that of \autoref{L:rank}. For $\mu\in\Comp[N,n]$ let
  \[
    \Parts[\mu,n]=\set[\Big]{\blam\in\lParts|{\mu_i=\sum_{k=0}^{\ell-1}|\lambda^{(i+kN)}|}
         \text{ for }1\le i\le N}.
  \]
  If $\blam\in\Parts[\mu,n]$, then $|\plam[\ilam]\hlam[\ilam]\qlam[\ilam]|=d^\mu$,
where $d^\mu=\prod_{k=1}^N|\plam[k]\hlam[k]\qlam[k]|^{\mu_k}$. If
  $\blam\in\Parts[\mu,n]$, for $k=1,\dots,N$ {view the $k$-layer of $\blam$ as the $\ell$-multipartition
  $\blam_{\langle k\rangle}=(\lambda^{(k)}|\lambda^{(k+N)}|\dots|\lambda^{(k+(\ell-1)N)})\in\Parts[\ell,\mu_k]$}
  and let $\std_\mu(\blam)=\prod_{k=1}^N\Std({\blam_{\langle k\rangle}})$.
  It is well-known and easy to prove that
  $\sum_{\blam\in\Parts[\ell,m]}|\Std^2(\blam)|=\ell^mm!$.
  (For example, take $A=\kk \Z/\ell\Z$ in \autoref{L:rank}.) Therefore,
  \begin{align*}
    |\Basis| & = \sum_{\blam\in\Parts[\ell N,n]}
                    |\Std^2(\blam)|\cdot|\plam[\ilam]\hlam[\ilam]\qlam[\ilam]|
         =\sum_{\mu\in\Comp[N,n]}{\binom n\mu}^2d^\mu
           \sum_{\blam\in\Parts[\mu,n]}|\std_\mu(\blam)|^2
=\sum_{\mu\in\Comp[N,n]}{\binom n\mu}^2d^\mu \ell^{\mu_1+\dots+\mu_N}\mu!\\
        & =\ell^nn!\sum_{\mu\in\Comp[N,n]}{\binom n\mu}d^\mu
         =\ell^nn!\op{rk}{\Bo}.
  \end{align*}
  Now, the result follows by applying \autoref{basisthmcycl}.
\end{proof}

Next we transfer the straightening results of \autoref{SS:LevelOneRespectable} to higher
levels. For $\blam\in\lParts$ let $\HH^{\pm\ge\blam}_n$ and $\HH^{\pm>\blam}_n$ be the $\kk$-spans
of $\bigcup_{\bmu\ge\blam}\Basis_\bmu$ and $\bigcup_{\bmu>\blam}\Basis_\bmu$, respectively.

\begin{Lemma}\label{L:ColumnStandard}
  Let $\blam\in\lParts$ and assume that $\HH^{\pm\ge\bmu}_n$ is a two-sided ideal of $\Hhd$
  for all $\bmu>\blam$. Then
  \begin{thmenumerate}
  \item\label{ellBoInIdeal} $\Bo^{>\ilam}X_\blam\subseteq\HH^{\pm>\blam}_n$.
  \item\label{ellStraighten} \autoref{P:Straighten} holds in
    $\Hhd$ with $y$ replaced by $X_\blam y$, with the same coefficients.
  \item\label{ellInductive} \autoref{L:InductiveStep}\autoref{E:sandwich}-\autoref{E:HmAction} hold in $\Hhd$ for the $H^\pm_{\bS,p,h,q,\bT}$ of
    \autoref{E:ellBasis}, with the same coefficients.
  \end{thmenumerate}
\end{Lemma}

\begin{proof}
  For \autoref{ellBoInIdeal} let $\bj>\ilam$, so that $j_s\le i^\blam_s$ for all $s$ and $\bj\ne\ilam$, as in the proof of \autoref{L:InductiveStep}\autoref{E:BoInIdeal}. Let $w\in\Sym_n$ preserve the positions of each $\kapp_l$-part of $\blam$ and sort the entries of $\bj$ within them, and let $\bj'=w^{-1}\bj$. Let $\bnu\in\lParts$ be the multipartition with a column for each layer $c$ in each $\kapp_l$-part, of length the number of entries $c$ of $\bj$ in the positions of the $\kapp_l$-part, see \autoref{Ex:chargeparts}. Then $\ilam[\bnu]=\bj'$, $\HJpm_\bnu=1$ and $X_\bnu=X_\blam$, as $\bnu$ and $\blam$ have the same number of boxes in each $\kapp_l$-part. Moreover $\bnu>\blam$, because $j_s\le i^\blam_s$ gives $\sum_{r\le(l-1)N+c}|\nu^{(r)}|\ge\sum_{r\le(l-1)N+c}|\lambda^{(r)}|$ for all $l$ and $c$, with a strict inequality for some $l$ and $c$ as $\bj\ne\ilam$. Since $w$ preserves the $\kapp_l$-parts, $X_\blam$ commutes with $H_w$ by \autoref{L:XSymmetric}, so
  \[
    p_\bj h_\bj q_\bj X_\blam=\pm H_wp_{\bj'}h_{\bj'}q_{\bj'}X_\bnu H_w^{-1}
    =\pm H_wH^\pm_{\Tlam[\bnu],p_{\bj'},h_{\bj'},q_{\bj'},\Tlam[\bnu]}H_w^{-1}
    \in\HH^{\pm\ge\bnu}_n\subseteq\HH^{\pm>\blam}_n .
  \]

  For \autoref{ellStraighten} and \autoref{ellInductive}, the proofs in
  \autoref{SS:LevelOneRespectable} only use the relations between $H_1,\dots,H_{n-1}$ and
  $\Bo$, which also hold in $\Hhd$, the assumption on the ideals $\HH^{\pm\ge\bmu}_n$, and
  \autoref{L:InductiveStep}\autoref{E:BoInIdeal}, which is only applied to elements of $\Bo$
  standing next to $\HJpm_\blam X_\blam$, where \autoref{ellBoInIdeal} applies instead. The
  strata $\bmu>\blam$ that occur in these proofs arise from Garnir moves, which change the shape
  of a single component but not the number of boxes in each component, so $X_\bmu=X_\blam$.
  Finally, the only elements that have to be moved past $X_\blam$ are elements of $\Bo$ and the
  $H_w$ with $w$ in the Young subgroup of a single component of $\blam$, such as the Garnir
  elements $\xipm_{g^{-1}}H_{g^{-1}}$. As $\psi$ is trivial and $X_\blam$ is
  symmetric in the $X_i$ with $i$ in one component, $X_\blam$ commutes with these by
  \autoref{L:XSymmetric}. Hence $X_\blam$ is a passenger in all these arguments.
\end{proof}
  \begin{Example}
 We return to \autoref{Ex:chargeparts}, where {$\blam=(\emptyset|(2,1)|(1)|\emptyset|(1)|(1))$}.  As $1>2$ in $I$, we have
 $\bj\ge\ilam$ if and only if $j_s\le i^\blam_s$ as integers for all
 $s$. Let $\bj=(2,1,2,1,1,1)$, so that $\bj>\ilam$. The symmetric group
 $\Sym$ acts on $I^n$ by place permutations, $(w\bj)_s=j_{w^{-1}(s)}$.
 The element $w=s_1$ preserves the positions of each $\kapp_l$-part, and
 $\bj'=w^{-1}\bj=(1,2,2,1,1,1)$ has its entries sorted within the
 positions of each part. Let $\bnu$ be the $\ell N$-multipartition {whose component of
 layer $i$ in the $\kapp_l$-part is $(1^m)$}, where $m$ is the number
 of entries in~$\bj$ that are equal to~$i$ and have positions in $\blam^{[l]}$:
  \[
      \bnu=\LabelledDiagram[
        label={\kapp_1|\kapp_1|\kapp_2|\kapp_2|\kapp_3|\kapp_3},
        tikz after={
          \draw[aTableauMain,thick]([xshift=5mm]A-2-1-1.north)--++(0,-1);
          \draw[aTableauMain,thick]([xshift=-6mm]A-5-1-1.north)--++(0,-1);
        }
      ]{1|1^2|1||1^2|}.
  \]
  Then $\ilam[\bnu]=\bj'$, and $\HJpm_\bnu=1$ as $\Sym_\bnu$ of
  \autoref{D:Wlam} is trivial. {Both $\blam$ and $\bnu$ have} three, one and two boxes in their $\kapp_1$-, $\kapp_2$- and $\kapp_3$-part, so $X_\bnu=X_\blam$. Finally $\bnu>\blam$ in the order of \autoref{D:lcontorder}, since the partial sums of the component sizes are $(1,3,4,4,6,6)$ for $\bnu$ and $(0,3,4,4,5,6)$ for $\blam$.
\end{Example}

We need two technical lemmas. For $0\le m<n$, let $\Hhd[n,m]$ be the
subalgebra of $\Hhd$ generated by $H_1,\dots,H_{m}$ and~{$\Bo$}, so that $\Hhd[n,0]={\Bo}$. {Here and below $\taui[i]=\taui[i,i+1]=\taui[i+1,i]$, as $\psi$ is trivial.}

\begin{Lemma}\label{L:XHij}
  Suppose that $1\le i\le m<n$. Then
  $X_iH_i\dots H_m= (H_i-\eta\taui)\dots(H_m-\eta\taui[m])X_{m+1}+\alpha$,
  for some $\alpha\in\Hhd[n,m]$.
\end{Lemma}

\begin{proof}
   We argue by induction on $m-i$. If $m=i$ then
   $X_iH_i=(H_i-\eta\taui)X_{i+1}-t \taui$ by
   \autoref{D:Hecke}, proving the base case.  Hence, by induction and
   \autoref{D:Hecke}, there exists $\alpha\in\Hhd[n,m]$ such that
   \begin{align*}
      X_iH_i\dots H_mH_{m+1}
        &= \Bigl((H_i-\eta\taui)\dots(H_m-\eta\taui[m])X_{m+1}+\alpha\Bigr)H_{m+1}
         \\
        &=(H_i-\eta\taui)\dots(H_m-\eta\taui[m])
           \Bigl((H_{m+1}-\eta\taui[m+1])X_{m+2}-t\taui[m+1]\Bigr)
            + \alpha H_{m+1}.
   \end{align*}
   Rearranging, this gives the result.
\end{proof}

\begin{Lemma}\label{L:XXlam}
    Suppose that $1\le m\le n$. Then
    $X_{m}\prod_{k=2}^\ell\prod_{i=1}^m(X_i-\kapp_k)
       =\beta\prod_{k=2}^\ell\prod_{i=1}^m(X_i-\kapp_k)$,
for some $\beta\in\Hhd[n,m-1]$.
\end{Lemma}

\begin{proof}
We argue by induction on $m$. If $m=1$ then
$X_1\prod_{k=2}^\ell(X_1-\kapp_k)=\kapp_1\prod_{k=2}^\ell(X_1-\kapp_k)$ by the cyclotomic relation,
which establishes the base case. For the inductive step let
$P_{m+1}=\prod_{k=2}^\ell\prod_{i=1}^{m+1}(X_i-\kapp_k)$, which commutes with $H_m$ by
\autoref{L:XSymmetric} and satisfies $X_mP_{m+1}=\beta P_{m+1}$ with $\beta\in\Hhd[n,m-1]$ by
induction. By \autoref{L:HXH}, with $\taui[m+1,m]=\taui[m]$,
\begin{align*}
   X_{m+1}P_{m+1}
      &= \hbar^{-2}\Bigl(H_mX_mH_m+t H_m\taui[m]\Bigr)P_{m+1}
       = \hbar^{-2}\Bigl(H_mX_mP_{m+1}H_m+t H_m\taui[m]P_{m+1}\Bigr)\\
      &= \hbar^{-2}\Bigl(H_m\beta H_m+t H_m\taui[m]\Bigr)P_{m+1},
\end{align*}
and the last bracket lies in $\Hhd[n,m]$.
\end{proof}

We can now generalise \autoref{T:HfindRespectable} to higher levels.

\begin{Theorem}\label{T:HhdRespectable}
  Suppose that $\Hhd$ is a braided diagonalisable {CRaFH algebra} of
  level $\ell\ge1$, and assume that $A$ is perfectly $\kk$-respectable with datum
  $(I,\ge,\plam[i],\hlam[i],\qlam[i],a)$. Then
  $\Basis=\set{(\plam[\blam]^\pm, \hlam[\blam]^\pm,\qlam[\blam]^\pm)|\blam\in\lParts}$
  is a perfectly $\kk$-respectable basis of $\Hhd$ with weight poset
  $(\lParts,\ge)$.
\end{Theorem}

\begin{proof}
  We argue by downward induction on $\blam$ in the order of \autoref{D:lcontorder}, as in the proof of \autoref{T:HfindRespectable}. So let $\blam\in\lParts$ and assume that $\HH^{\pm\ge\bmu}_n$ is a two-sided ideal of $\Hhd$ for all $\bmu>\blam$. We show that both congruences of \autoref{resp2} hold for $\Basis_\blam$ and all $x\in\Hhd$. In particular, $\HH^{\pm\ge\blam}_n$ is then a two-sided ideal. All congruences are modulo $\HH^{\pm>\blam}_n$. For $x=H_i$ and $x\in\Bo$ the congruences hold by \autoref{L:ColumnStandard}\autoref{ellInductive}, hence they hold for all $x\in\Hhd[n,n-1]$. As $\Hhd$ is generated by $\Hhd[n,n-1]$ and $X_1$ by \autoref{L:HXH}, it remains to consider $x=X_1$.

  Let $\bS,\bT\in\Std(\blam)$, $p\in\plam[\ilam]$, $h\in\hlam[\ilam]$ and $q\in\qlam[\ilam]$. Recall that $X_\blam$ commutes with $\HJpm_\blam$, and that $X_1,\dots,X_n$ commute with $\Bo$ as $\psi$ is trivial. As $H^\pm_{\Tlam}=1=H^{\pm*}_{\Tlam}$, this gives
  \[
    H^\pm_{\bS,p,h,q,\bT}=H^\pm_\bS H^\pm_{\Tlam,p,h,q,\bT}=H^\pm_{\bS,p,h,q,\Tlam}H^{\pm*}_\bT
  \]
  with $H^\pm_{\Tlam,p,h,q,\bT}=X_\blam p\HJpm_\blam hqH^{\pm*}_\bT$ and $H^\pm_{\bS,p,h,q,\Tlam}=H^\pm_\bS p\HJpm_\blam hqX_\blam$.
  Let $1$ lie in component $k$ of $\bS$, and let $m=|\lambda^{(1)}|+\dots+|\lambda^{(k-1)}|$. Then $d_\bS(m+1)=1$, so $d_\bS=du$ with $u=s_1\cdots s_m$, $d(1)=1$ and $\ell(d_\bS)=\ell(d)+m$. As $d(1)=1$, $H_d$ commutes with $X_1$. Hence \autoref{E:HS}, \autoref{C:hdv} and \autoref{L:XHij} give
  \begin{equation}\label{E:X1HS}
    X_1H^\pm_\bS=\hbar^{-2(\ell(d)+m)}\xipm_d\leftindex^d{{\xipm_u}}H_d\Bigl((H_1-\eta\taui[1])\cdots(H_m-\eta\taui[m])X_{m+1}+\alpha\Bigr)
  \end{equation}
  for some $\alpha\in\Hhd[n,m]$, where the product is $1$ and $\alpha=0$ if $m=0$.

  \emph{Claim.} For every nonempty component $k$ of $\blam$, with $m$ as above, there are $\gamma,\gamma'\in\Hhd[n,m]$, depending only on $\blam$ and $k$, such that
  \[
    X_{m+1}X_\blam p\HJpm_\blam hq\equiv\gamma X_\blam p\HJpm_\blam hq
    \qquad\text{and}\qquad
    p\HJpm_\blam hqX_\blam X_{m+1}\equiv p\HJpm_\blam hqX_\blam\gamma'
  \]
  for all $p\in\plam[\ilam]$, $h\in\hlam[\ilam]$ and $q\in\qlam[\ilam]$.

We assume for the moment that this claim holds. By \autoref{E:X1HS} and the claim, $X_1H^\pm_{\bS,p,h,q,\bT}\equiv YH^\pm_{\Tlam,p,h,q,\bT}$ for some $Y\in\Hhd[n,n-1]$ not depending on $h$, $q$ and $\bT$. So the first congruence of \autoref{resp2} for $x=X_1$ follows from the one for $x=Y$. For the second congruence let $1$ lie in component $k$ of $\bT$, and define $m$, $d$ and $u$ for $\bT$ as above. The mirror image $H_m\cdots H_1X_1=X_{m+1}(H_m-\eta\taui[m])\cdots(H_1-\eta\taui[1])+\alpha$ of \autoref{L:XHij}, which has the same proof, gives
  \[
    H^{\pm*}_\bT X_1=\hbar^{-2(\ell(d)+m)}\Bigl(X_{m+1}(H_m-\eta\taui[m])\cdots(H_1-\eta\taui[1])+\alpha\Bigr)H_{d^{-1}}\xi^{\pm*}_{d_\bT}.
  \]
  By the claim, $H^\pm_{\bS,p,h,q,\bT}X_1\equiv H^\pm_{\bS,p,h,q,\Tlam}Y'$ for some $Y'\in\Hhd[n,n-1]$ not depending on $\bS$ and $p$, and we conclude as before.

  It remains to prove the claim. Write $k=k'+(l-1)N$ with $1\le k'\le N$ and $1\le l\le\ell$, so that component $k$ is the component of stratum $k'$ of the $\kapp_l$-part of $\blam$, and let $a=n_l$. Then $a\le m<n_{l+1}$.

  If $l=1$, then $m+1\le n_{l''}$ for $2\le l''\le\ell$. Hence $X_\blam=PR$, where $P=\prod_{l''=2}^\ell\prod_{i=1}^{m+1}(X_i-\kapp_{l''})$ and $R$ is a polynomial in $X_{m+2},\dots,X_n$. By \autoref{L:XXlam} and its mirror image $PX_{m+1}=P\beta'$, which has the same proof, $X_{m+1}X_\blam=\beta X_\blam$ and $X_\blam X_{m+1}=X_\blam\beta'$ with $\beta,\beta'\in\Hhd[n,m]$. So the claim holds with $\gamma=\beta$ and $\gamma'=\beta'$.

  Now let $l\ge2$ and $z_1=s_{a+1}\cdots s_m$. Applying \autoref{L:HXH} $m-a$ times gives $X_{m+1}=\hbar^{-2(m-a)}H_{z_1^{-1}}X_{a+1}H_{z_1}+\alpha'$ with $\alpha'\in\Hhd[n,m]$. As $z_1$ preserves the positions of each $\kapp_{l''}$-part of $\blam$, $X_\blam$ commutes with $H_{z_1}$ and $H_{z_1^{-1}}$ by \autoref{L:XSymmetric}. Let $\bmu\in\lParts$ be the multipartition whose components are columns with $|\mu^{(r)}|=|\lambda^{(r)}|+\delta_{r,k-N}-\delta_{r,k}$. That is, one box of stratum $k'$ moves from the $\kapp_l$-part to the $\kapp_{l-1}$-part. Then $\bmu>\blam$ by \autoref{D:lcontorder}, and $(X_{a+1}-\kapp_l)X_\blam=X_\bmu$ by \autoref{D:Xlam}, since $n_l$ increases by one and the other $n_{l''}$ do not change. Altogether, with $\gamma=\kapp_l\hbar^{-2(m-a)}H_{z_1^{-1}}H_{z_1}+\alpha'\in\Hhd[n,m]$,
  \[
    X_{m+1}X_\blam=\gamma X_\blam+\hbar^{-2(m-a)}H_{z_1^{-1}}X_\bmu H_{z_1}
    \qquad\text{and}\qquad
    X_\blam X_{m+1}=X_\blam\gamma+\hbar^{-2(m-a)}H_{z_1^{-1}}X_\bmu H_{z_1}.
  \]
  Hence the claim holds with $\gamma'=\gamma$ if $X_\bmu H_{z_1}p\HJpm_\blam hq$ and $p\HJpm_\blam hqH_{z_1^{-1}}X_\bmu$ lie in $\HH^{\pm\ge\bmu}_n$, which is a two-sided ideal contained in $\HH^{\pm>\blam}_n$.

  By \autoref{D:HJ}, $\HJpm_\blam=\sum_{v\in\Sym_\blam}\hbar^{-2\ell(v)}\xipm_vH_v$ with $\xipm_v\in\Tos_\blam$ by \autoref{L:HJ}\autoref{HJ2}. Moreover, $\leftindex^xp\in\pm\plam[\ilam]$ and $\leftindex^x{(hq)}$ is pure for $x\in\Sym_{|\blam|}$. So \autoref{L:HJ}\autoref{HJ3} shows that $p\HJpm_\blam hq$ lies in $\sum_\bb\bb\Hhd$ and in $\sum_\bb\Hhd\bb$, where $\bb$ runs over the basis elements $p'h'q'$ of $\Bo$ with $p'\in\plam[\ilam]$, $h'\in\hlam[\ilam]$ and $q'\in\qlam[\ilam]$. As $X_\bmu H_{z_1}\bb=X_\bmu\leftindex^{z_1}\bb H_{z_1}$ and $\bb H_{z_1^{-1}}X_\bmu=H_{z_1^{-1}}X_\bmu\leftindex^{z_1}\bb$, it suffices to show that $X_\bmu\leftindex^{z_1}\bb\in\HH^{\pm\ge\bmu}_n$ for these $\bb$.

  Let $m'=|\lambda^{(1)}|+\dots+|\lambda^{(k-N)}|$ and $z_2=s_{m'+1}\cdots s_a$. Then $z_2z_1$ maps $m+1$ to $m'+1$, and $r$ to $r+1$ for $m'<r\le m$. Hence $z_2z_1\ilam=\ilam[\bmu]$, and $\leftindex^{z_2z_1}\bb=\pm p''h''q''$ with $p''\in\plam[{\ilam[\bmu]}]$, $h''\in\hlam[{\ilam[\bmu]}]$ and $q''\in\qlam[{\ilam[\bmu]}]$. As $z_2$ preserves the positions of each $\kapp_{l''}$-part of $\bmu$, $X_\bmu$ commutes with $H_{z_2}$ by \autoref{L:XSymmetric}. Finally, $\HJpm_\bmu=1$ as $\Sym_\bmu$ is trivial, and $H^\pm_{\Tlam[\bmu]}=1=H^{\pm*}_{\Tlam[\bmu]}$. Therefore
  \[
    X_\bmu\leftindex^{z_1}\bb=H_{z_2}^{-1}X_\bmu\leftindex^{z_2z_1}\bb H_{z_2}
      =\pm H_{z_2}^{-1}H^\pm_{\Tlam[\bmu],p'',h'',q'',\Tlam[\bmu]}H_{z_2}\in\HH^{\pm\ge\bmu}_n.
  \]
  This proves the claim and completes the induction.

  It remains to show that $\Basis$ spans $\Hhd$. For $\bj\in I^n$ let $\bom_\bj\in\lParts$ be the multipartition with $\kapp_\ell$-part $(1^{m_1}|\dots|1^{m_N})$, where $m_c$ is the number of entries of $\bj$ equal to $c$, and with all other components empty. Then $X_{\bom_\bj}=1$, so the argument of \autoref{F3} shows $p_\bj h_\bj q_\bj\in\Hhd\Basis_{\bom_\bj}\Hhd\subseteq\HH^{\pm\ge\bom_\bj}_n$. As in the proof of \autoref{T:HfindRespectable}, the identity of $\Hhd$ therefore lies in the two-sided ideal $\sum_{\bj\in I^n}\HH^{\pm\ge\bom_\bj}_n$, so $\Basis$ spans $\Hhd$. By \autoref{L:lrank}, $\Basis$ is a basis, which is \autoref{resp1}. Condition \autoref{resp3} holds trivially as $\Olam=\kk$, and the datum is perfect as $H^\pm_{\bS,p,h,q,\bT}$ is the product of $H^\pm_\bS p$, $\HJpm_\blam X_\blam$ and $hqH^{\pm*}_\bT$.
\end{proof}

\makeatletter
\newcommand{\where}[1]{\@ifundefined{r@#1}{}{\pageref{#1}}}
\newcommand{\alsowhere}[1]{\@ifundefined{r@#1}{}{, \pageref{#1}}}
\makeatother
\newenvironment{notationtable}[1]
  {\subsection*{#1}\small
   \begin{longtable}{@{}p{0.225\textwidth}@{\hspace{1em}}p{0.635\textwidth}@{\hspace{0.6em}}p{0.085\textwidth}@{}}
   \toprule \emph{Symbol} & \emph{Meaning} & \emph{Page}\\ \midrule \endfirsthead
   \toprule \emph{Symbol} & \emph{Meaning} & \emph{Page}\\ \midrule \endhead
   \bottomrule \endfoot}
  {\end{longtable}}

\section*{Index of notation}\label{S:index}
\addtocontents{toc}{\protect\setcounter{tocdepth}{0}}
\noindent
The symbols of the paper, in the order in which they are introduced.
\begin{notationtable}{\S1.\ Notation and conventions}
$\kk$ & ground ring: a commutative integral domain with $1$ & \where{S:notation}\\
$\iota$ & an element of $\kk$ with $\iota^{2}=-1$, if one exists & \where{S:notation}\\
$Z(A)$ & centre of an algebra $A$ & \where{S:notation}\\
$\kk G$ & group algebra of $G$, with basis $\set{g|g\in G}$ & \where{S:notation}\\
$\Ztwo=\set{\ov0,\ov1}$ & the group of order two carrying all supergradings & \where{S:notation}\\
$\mathtt{C}_m$ & multiplicative cyclic group of order $m$ & \where{S:notation}\\
$n$ & the fixed rank: a non-negative integer & \where{S:notation}\\
$X_1,\dots,X_n$, $\kk[\X]$ & the fixed indeterminates and $\kk[X_1,\dots,X_n]$ & \where{S:notation}\\
$\Sym=\Sym_n$ & the symmetric group on $\set{1,\dots,n}$ (written $W$) & \where{S:notation}\\
$S=\set{s_1,\dots,s_{n-1}}$ & the simple transpositions $s_i=(i,i+1)$ & \where{S:notation}\\
$\ell(w)$ & Coxeter length of $w\in\Sym$ & \where{S:notation}\\
\end{notationtable}

\begin{notationtable}{\S2.\ Frobenius-tensor superalgebras and universal Jucys--Murphy elements}
$\kk G\Smash B$ & smash product of $\kk G$ with a $G$-algebra $B$ & \where{S:Smash}\\
$\nBraid=\nBraid_n$ & the \emph{universal Frobenius monoid} & \where{univtaurel}\\
$\ttau^{(i,j)}$ & generators of $\nBraid$, for $1\le i\ne j\le n$ & \where{univtaurel}\\
$\Tos=\kk\nBraid$ & monoid algebra of $\nBraid$, with $\Sym$ permuting indices & \where{D:universalJM}\\
$\CL_k$ & the $k$th \emph{universal Jucys--Murphy element}; $\CL_1=0$ & \where{D:universalJM}\\
$\Tosym$ & the \emph{symmetric universal Frobenius algebra} & \\
$\op{Br}_n$ & Artin braid group, for comparison with $\nBraid$ & \\
$\op{Vect}_{\Ztwo}$ & category of $\Ztwo$-graded $\kk$-modules & \where{superstuffeasy}\\
$|a|\in\Ztwo$ & the \emph{degree} or \emph{parity} of a homogeneous element & \where{superstuffeasy}\\
$\Pi$ & parity shift, $\Pi M=\kk^{0|1}\otimes M$ & \where{superstuffeasy}\\
$A^{\mathrm{op}}$ & opposite superalgebra, $a\cdot b=(-1)^{|a||b|}ba$ & \where{superstuffeasy}\\
$\Ao=A^{\otimes n}$ & the $n$th tensor power of a superalgebra $A$ & \where{Symaction}\\
$\leftindex^{s_i}\ba$ & the \emph{superflip}: $\Sym$ acts by \emph{superpermutations} & \where{Symaction}\\
$\Bo$ & a \emph{tensor superalgebra}: $A^{\otimes n}$ with a map $\Tos\to\Bo$ & \where{D:ST-algebra}\\
$\tau^{(i,j)}$ & image of $\ttau^{(i,j)}$ in $\Bo$ & \where{D:ST-algebra}\\
$\psi$, $\psi^{(i)}$ & \emph{Nakayama automorphism} and its $i$th tensor leg & \where{D:Nakayama}\\
$(\Bo,\psi)$ & a \emph{Frobenius-tensor superalgebra} & \where{D:Frobtensor}\\
$\eta$-symmetric & $\eta=0$, or $\eta\ne0$ and $\psi=\op{id}$ & \where{D:Frobtensor}\\
$\Bo=\kk^{\otimes n}=\kk$ & the \emph{trivial} Frobenius-tensor superalgebra & \where{trivFrobtensor}\\
$\tr\map A\kk$ & the \emph{trace map}; $\<a,b\>=\tr(ab)$ is non-degenerate & \where{S:Frobeniusalgs}\\
$\psi$ & Nakayama automorphism: $\tr(ab)=(-1)^{|a||b|}\tr(b\psi(a))$ & \where{S:Frobeniusalgs}\\
$\cB$, $r=|\cB|$ & a fixed homogeneous basis of $A$ and its size & \where{S:Frobeniusalgs}\\
$b^\vee$ & the dual basis: $\tr(a^\vee b)=\delta_{ab}$ & \where{S:Frobeniusalgs}\\
$a^{(i)}$ & $a$ embedded in the $i$th tensor leg of $\Ao$ & \where{S:Frobeniusalgs}\\
$\tau^{(i,j)}$ & the \emph{copairing element} $\sum_b b^{(i)}(b^\vee)^{(j)}$ & \where{parts}\\
$\tau^{(i,j)}_{\ov0}$, $\tau^{(i,j)}_{\ov1}$ & its \emph{$0$-part} and \emph{$1$-part}, by the parity of $b$ & \where{parts}\\
$\epsilon_1,\epsilon_2\in\Ztwo$ & parity indices labelling the parts $\tau^{(i,j)}_{\epsilon}$ & \where{Taussuper}\\
$\taui[i]$ & $\tau^{(i,i+1)}=\tau^{(i+1,i)}$, for $A$ symmetric & \\
$\tau$ & $\tau^{(1,2)}\in A\otimes A$, the case $n=2$ & \where{Frobeniusalgex}\\
$\kk G$ & group algebra of a finite group, $\tr(g)=\delta_{g,1_G}$ & \where{Ex:FrobeniusAlgebrasfingroup}\\
$\kk[x]/(x^{m+1})$ & truncated polynomials, $\tr(x^{k})=\delta_{k,m}$; $H^*(\Poo^{m})$ & \where{Ex:FrobeniusAlgebrasHPm}\\
$A=\kk\<y_1,y_2\>/(\dots)$ & the four-dimensional example with $\psi(y_1)=\lambda y_1$ & \where{Alambda}\\
$\Taft$ & the \emph{Taft algebra} on $g,y$ with $yg=\omega gy$ & \where{Ex:Taft}\\
$\omega$ & a fixed primitive $\theta$th root of unity & \where{Ex:Taft}\alsowhere{D:fcyc}\\
$\Cl_1$, $\Cl_n$ & Clifford superalgebra; $\tau^{(i,j)}=1+C_iC_j$ & \where{Ex:Clifford}\\
$M(\tm|\tn)$, $Q_\tn$ & the simple superalgebras of Wedderburn's theorem & \where{exsemisimple}\\
$\mathcal{Z}ig_m$ & the \emph{zig-zag algebra} of type $A_m$ & \where{Ex:ZigZagm}\\
\end{notationtable}

\begin{notationtable}{\S3.\ RaFH superalgebras}
$\RR$ & commutative integral domain over $\kk$ holding $\hbar,t,\eta$ & \\
$\hbar,\;t,\;\eta$ & the three parameters of a RaFH superalgebra & \where{D:Hecke}\\
$B_\RR=\RR\otimes_\kk B$ & scalar extension $\RR\otimes_\kk B$ & \where{defRaFH}\\
$\HHd(\Bo)=\HHd[n;\hbar,t,\eta](\Bo)$ & the \emph{RaFH superalgebra}: Rees affine Frobenius Hecke & \where{D:Hecke}\\
$H_1,\dots,H_{n-1}$ & Hecke generators; $H_i^{2}=\eta\tau^{(i,i+1)}H_i+\hbar^{2}$ & \where{D:Hecke}\\
$\gamma_i$ & $t+\eta X_{i+1}$, central in $\Po$ & \where{D:Hecke}\alsowhere{L:gammacentral}\\
$H_w$ & $H_{i_1}\cdots H_{i_r}$ for a reduced expression of $w$ & \\
$\CL_k$ & Jucys--Murphy elements of $\HHd(\Bo)/(X_1)$ & \\
$\HHd[n;\varphi(\hbar),\varphi(t),\varphi(\eta)](\Bo)_{\RR'}$ & \emph{specialised} RaFH superalgebra along $\varphi\map\RR{\RR'}$ & \where{specialisation}\\
$\Hfind(\Bo)=\Hfind[n;\hbar,\eta](\Bo)$ & the \emph{Rees Frobenius Hecke superalgebra}: no $X_j$ & \where{D:RFHS}\\
$\HSd$ & the finite Hecke algebra $\Hfind(\kk)$ & \where{D:RFHS}\\
$\Hsmash$, $\Hdegen$ & $\Hd[n;1,0,0](\kk)$ and $\Hd[n;1,1,0](\kk)$ & \where{HHdexamples}\\
$\Heaff$ & extended affine (Iwahori--)Hecke algebra & \where{HHdexamples}\\
$\NC$, $\NC_n$ & the Nil Coxeter algebra; $\Hd[n;0,0,0](\kk)$ affinises it & \where{HHdexamples}\\
$H_n(A)$, $\mathcal{Z}ig^{\op{aff}}_n(\Gamma)$ & the Kleshchev--Muth affinisation and affine zigzag algebra & \where{Ex:KMaffinisation}\\
$\Po=\RR[\X]\Smash\Bo$ & $\RR[\X]\Smash\Bo$; written $\Poo$ when unambiguous & \where{DefPo}\\
$\psi_i$, $\psi_\bd$ & $\psi^{(i)}$ and $\psi_1^{d_1}\cdots\psi_n^{d_n}$ & \\
$w$-skew-derivation & $\Delta(xy)=\Delta(x)y+\leftindex^{w}x\,\Delta(y)$ & \where{skewDer}\\
$\delta_i$ & the \emph{Frobenius--Demazure operator} on $\Po$ & \where{D:FrobDem}\\
$\Delta_i$ & the \emph{RaFH derivation} $\Delta_i=-\gamma_i\delta_i$ & \where{D:RaFHDer}\alsowhere{L:crucial}\\
$\partial_i$ & the classical divided difference operator & \where{R:classical}\\
$\pi_i$ & the $K$-theoretic Demazure operator, $\Delta_i=1-\pi_i$ for $t=0$, $\eta=1$& \where{RkDemop}\\
$\mathcal{X}$ & a chosen $\RR$-basis of $\Bo$ & \where{C:BasisHd}\\
$\K$ & algebraic closure of the fraction field of $\RR$ & \\
$\RX^{d,\Sym}$, $\RX^{\Sym}$ & symmetric polynomials in $X_1^{d},\dots,X_n^{d}$; $d=1$ & \\
$Z_{\bar{0}}(C)$ & the \emph{even centre} of a superalgebra $C$ & \\
\end{notationtable}

\begin{notationtable}{\S4.\ Irreducible representations of RaFH superalgebras}
$\uHom_A(M,N)$, $\uEnd_A(M)$ & \emph{inner} morphisms: the super Hom and End spaces & \where{uEnd}\\
type $\mathtt M$, type $\mathtt Q$ & $\uEnd_A(L)\cong\kk$, resp.\ $\cong\Cl_1$ & \where{uEnd}\\
$\leftindex^{\phi}M$ & $M$ with the $A$-action twisted by $\phi$ & \where{uEnd}\\
$L\circledast L'$ & $e(L\otimes L')$ if both are type $\mathtt Q$, else $L\otimes L'$ & \where{types}\\
$\Omega$, $O_i$, $N_\psi$, $L_i$ & simple $A$-modules, their $\psi$-orbits and representatives & \\
$\theta_L$, $\theta_i=\theta_{L_i}$ & orbit size: least $l$ with $\leftindex^{\psi^{l}}L\cong L$ & \\
$\phi_L$ & a chosen isomorphism $\leftindex^{\psi^{\theta_L}}L\to L$ & \where{E:phinorm}\\
$\beta_{L,L'}$, $\beta_L$ & the maps induced by $\tau^{(1,2)}s_1$ & \where{tausi}\alsowhere{L:tauscalar}\\
$\bar\gamma_L$ & the scalar by which $\beta_L$ acts on $L\circledast L$ & \where{L:tauscalar}\\
$\alpha_L,\alpha'_L,\bar\gamma'_L$ & the refinement of $\bar\gamma_L$ in type $\mathtt Q$ & \where{gammaprime}\\
$L^{\mathrm{nil}}$ & $L$ with $X_1$ acting as zero & \where{D:La}\\
$\mathbf L(\zeta)$, $x_\zeta$ & $\bigoplus_{k<\theta_L}\leftindex^{\psi^{k}}L$ and the action of $X_1$ & \where{D:La}\\
$\mathcal E=\mathcal E_1$, $\mathcal E_n$ & $\uEnd_A(\mathbf L(\zeta))^{\mathrm{op}}$ and its $n$th analogue & \where{L:endoalg}\\
$\psi_{\mathcal E}$, $\tau^{(i,j)}_{\mathcal E}$ & the Nakayama automorphism of $\mathcal E_n$ and its copairings & \where{L:endoalg}\alsowhere{psiE}\\
$\rho_i$ & the superflip of the $i$th and $(i+1)$st factors & \\
$\gamma_L$ & the normalisation of $\bar\gamma_L$ in the atomic parameters & \where{gammadef}\\
$\mathcal H_n(L)$ & the \emph{atomic algebra} attached to $L$ & \where{D:block}\\
Hecke / Clifford / wreath type & the three types of atomic algebra & \where{R:blocks}\\
$\Lo$, $\mathcal C(\Lo)$ & a semisimple $\Bo$-module and the subcategory it cuts out & \\
$\bar\Omega$, $\bar O_i$, $\bar N_\psi$, $\bar I$ & simple $A$-modules mod parity and their orbits: the \emph{colours} & \where{D:colours}\\
$\chi(P)$, $|\chi|$ & the \emph{colour} of a simple $\Bo$-module and its multiplicity & \where{D:colours}\\
$\HHd[\mu](\Bo)$, $\Sym^{\mu}$ & the \emph{parabolic subalgebra} attached to a composition $\mu$ & \\
$\mathbf L_\mu$, $\mathcal E_\mu$, $\So_\mu$ & the data and the induced simple module of the classification & \where{T:simples}\\
\end{notationtable}

\begin{notationtable}{\S5.\ CRaFH superalgebras}
$\theta$ & the order of $\psi$ & \where{D:fcyc}\\
$\fcyc(X_1)=\fct(X_1)$ & the \emph{$\theta$-cyclotomic polynomial} & \where{D:fcyc}\\
$\ell_0,\ell_1$; $\kapp_1,\dots,\kapp_{\ell_1}$ & the two exponents and the nonzero cyclotomic parameters & \where{D:fcyc}\\
$\ell=\ell_0+\theta\ell_1$ & the \emph{level} of $\fcyc$ & \where{D:fcyc}\\
$\charge$ & the \emph{charge} $(0,\dots,0,\kapp_1,\dots,\kapp_{\ell_1})$ & \where{D:fcyc}\\
$\ell_\theta=(\ell_0,\ell_1)$ & the \emph{level vector} & \where{D:fcyc}\\
$\Hhd=\Hhd[n;\hbar,t,\eta](\Bo)$ & the \emph{CRaFH superalgebra} $\HHd/\Ic$ of charge $\charge$ & \where{D:CyclotomicHecke}\\
$\Ic$ & the \emph{cyclotomic ideal} $(\fct)$ & \where{D:CyclotomicHecke}\\
$\CA(n)$ & the Hu algebra, and the variant of $\Hfind$ capturing it & \\
$\cB_{\op{cyc}}$ & the basis $\set{X_1^{d}b|d<\ell,\ b\in\cB}$ of $\Hhd[1]$ & \where{nequal1}\\
$\cBo$ & the induced basis $b_1\otimes\dots\otimes b_n$ of $\Bo$ & \\
$X^\bd$ & $X_1^{d_1}\cdots X_n^{d_n}$ for $\bd\in\N_0^{n}$ & \\
$\N^{\ell,n}_0$ & $\set{\bd\in\N^{n}|0\le d_i<\ell}$ & \\
$Q^{m}_c$ & the alternative cyclotomic quotient when $\hbar=0$ & \where{badh2}\\
$\mathtt f_i$, $\mathtt f_Z$ & $\mathtt f_{i+1}=H_i^{-1}\mathtt f_iH_i$, and products $\mathtt f_Z$ & \\
$\Pi^{n,\ell}$, $\Pi^{n,\ell}_+$ & index sets of pairs $(\bd,Z)$ in the basis theorem proof & \where{domweights}\\
$\Poo_i$ & the subalgebra of $\Poo$ on $X_1,\dots,X_i$ & \\
\end{notationtable}

\begin{notationtable}{\S6.\ Trace forms on CRaFH superalgebras}
$\Tr\map{\Hhd}\RR$ & the trace map of a CRaFH superalgebra & \where{defform}\\
$\<\ ,\ \>$ & the induced bilinear form $\<x,y\>=\Tr(xy)$ & \where{defform}\\
$\ELL$ & the top exponent $(\ell-1,\dots,\ell-1)$ & \\
$N=n(\ell-1)$ & the top $\X$-degree of $\Hhd$ & \\
$d$ & $d_1+\dots+d_n$; by convention $0^{0}=1$ & \\
$\Sh$ & the \emph{shifting isomorphism} $X_j\mapsto X_j-t/\eta$ & \where{Phidef}\\
$\psi_{\Hhd}$ & the Nakayama automorphism of $\Hhd$ in level $\ell=1$ & \where{formlevel1}\\
$\X$-filtration & the filtration by the degree in the $X_i$ & \where{Xfiltr}\\
$F_m$, $F'_m$ & the $m$th pieces of the $\X$-filtration of $\Hhd$ & \where{filtannoying}\\
$\mathbf H$, $S$, $S_m$, $\bar S$ & a filtered algebra with basis, and its associated graded & \where{sectop}\\
$T$, $\beta_T$ & a \emph{top functional} and its \emph{top pairing} $T(xy)$ & \where{sectop}\alsowhere{L:antitri}\\
$\alpha_k$, $\alpha_{k,f}$ & the $k$th coefficient map on $\RR[X]$, resp.\ $\RR[X]/(f)$ & \\
$M(f,\phi)$, $M(f,k)$ & the Hankel matrix $(\phi(X^{a+b}))_{0\le a,b<\ell}$ & \\
$G_A$ & the Gram matrix $(\tr(\bb\bc))$ of $A$ & \\
$G$, $G_{\T}$ & the Gram matrices of $\<\ ,\ \>$ and of $\T$ & \\
$\phi_{\Tr}$ & the factor of $\Tr$ on $\RR[X_1]/(\fct)$ & \where{factor}\\
$S$, $S'$ & the ordered bases $\set{X_1^{k}\bb}$ and $\set{\bb X_1^{k}}$ & \where{SS}\\
RaFH-CRaFH trick & feeding $\Hhd[1](A)$ into a RaFH superalgebra as its $A$ & \where{SS:Trick}\\
$\Kr{\tl}$, $\Kz{\tl}$ & subalgebras of $\Hhd$ on $X_1..X_{\tl}$, $H_1..H_{\tl-1}$, $\Bo$ & \\
$\varepsilon_{\tl}$ & the $\Kr{\tl-1}$-linear projection $\Kr{\tl}\to\Kr{\tl-1}$ & \where{L:Krdec}\\
$\prec$ & the total order on $\N^{\ell,n}$ used for triangularity & \\
extremal & $\tl$ with $a_i=0$ for $i>\tl$ and $u\in\Kz{\tl}$ & \\
$\gcy$ & the \emph{$(t,-\eta)$-homogenisation} $\prod_\varrho(t+\eta\varrho)$ & \where{E:g}\\
$\T$, $\nu$, $m$ & the functional $\T=\nu^{N}\Tr$, its scalar, and the exponent & \where{generalT}\alsowhere{T:summary}\\
\end{notationtable}

\begin{notationtable}{\S7.\ Braided diagonalisable RaFH superalgebras}
$\xi^\pm_1$ & \emph{roots} of the quadratic relation for $H_1$ & \where{D:Diagonalisable}\\
$\xi$-braid relations & the braid relations satisfied by the $\xi^\pm_1$ & \where{D:Diagonalisable}\\
$\alpha$ & a square root of $(\eta\taui[1,2])^{2}+4\hbar^{2}$ & \where{L:roots}\\
$\xi^\pm_i$ & $\xi^+_i\xi^-_i=-\hbar^{2}$, $\xi^+_i+\xi^-_i=\eta\taui[i,i+1]$ & \where{D:xis}\alsowhere{C:xii}\\
$\xi^\pm(\bi)$, $\xi^\pm_w$ & twisted products of the $\xi^\pm_i$ along a reduced word & \where{D:hw}\alsowhere{L:hw}\\
$\xi^{\pm*}_d$ & the twisted product read backwards & \where{C:hwrev}\alsowhere{E:HS}\\
$\hat H^\pm_i$, $\hat H^\pm_w$ & $\xi^\pm_iH_i$ and $\xi^\pm_wH_w$: they satisfy the braid relations & \where{R:GH}\\
$\phi_w$, $\phi_i$ & the $w$-skew-derivation $\op{id}-w$ of $\Bo$ & \where{D:phi}\\
$\alpha_i$, $\alpha^\pm_i=\pm\alpha_i$ & central elements of $\Bo$ controlling $\bb\,\xi^\pm_i$ & \where{L:alpha}\alsowhere{L:xicomm}\\
$\Theta^\pm_i=\alpha^\pm_i\phi_i$ & the skew-derivation $\alpha^\pm_i\phi_i$, analogue of $\Delta_i$ & \where{R:Theta}\\
$\phi_J$, $y^\pm_{J,w}$ & the composite $\phi_J$ and the mixed twisted product & \where{defy}\\
$\Sym_J$, $\Sym^J_K$, $\leftindex^J_K\Sym$ & parabolic subgroups of $\Sym$ and coset representatives & \\
$\HJpm_K$ & $\sum_{w\in\Sym_K}\hbar^{-2\ell(w)}\xi^\pm_wH_w$: the (anti)symmetriser & \where{D:HJ}\\
$\MJ^\pm_K$ & the \emph{parabolic RaFH module} $\HHd\HJpm_K$ & \where{D:HJ}\\
$\conts$, $\conts^\pm_k$ & \emph{contents}: the eigenvalues of $X_k$ on $\HJpm_S$ & \where{contents}\\
\end{notationtable}

\begin{notationtable}{\S8.\ Respectable algebras}
$\CA$ & a $\kk$-algebra carrying a respectable datum & \where{D:Respectable}\\
$(\Lambda,\ge)$ & the \emph{weight poset} of a respectable datum & \where{D:Respectable}\\
$\plam,\hlam,\olam,\qlam$ & the four \emph{indexing sets} attached to $\lambda$ & \where{D:Respectable}\\
$\Olam$ & the \emph{fibre algebra} at $\lambda$, with basis $\olam$ & \where{D:Respectable}\alsowhere{perfectresp}\\
$a_{phoq}=a^\lambda_{p,h,o,q}$ & the \emph{respectable basis}; indices dropped when singletons & \where{D:Respectable}\\
$\CA^{(\lambda)}$, $\CA^{\ge\lambda}$, $\CA^{>\lambda}$ & the $\lambda$-\emph{stratum} and the two filtration terms & \where{D:Respectable}\\
\textup{(R1)}, \textup{(R2)}, \textup{(R3)} & the three axioms of a respectable datum & \where{D:Respectable}\\
$r_{pop'o'}(x)$, $s_{oqo'q'}(x)$ & structure constants of the left and right action & \where{D:Respectable}\\
$c_{oo'o''}$ & structure constants of $\Olam$ & \where{D:Respectable}\\
$\dot p,\dot h,\dot o,\dot q$ & the \emph{lifts} of a \emph{perfectly} respectable datum & \where{perfectresp}\\
$\kk$-respectable & $\Olam=\kk$ for every $\lambda\in\Lambda$ & \where{perfectresp}\\
size of a stratum & $(|\plam|,|\hlam|,|\olam|,|\qlam|)$; their product is its rank & \where{perfectresp}\\
$\Mat_m(\CO)$, $E_{pq}$ & the prototype respectable algebra, and matrix units & \where{prototype}\\
$TL_m(\delta)$, $Q(m)$ & Temperley--Lieb and queer superalgebra examples & \where{Ex:TL}\alsowhere{Ex:Qn}\\
$1_\lambda$, $a_{phcq}$ & the unit of $\Olam$, and the basis element it scales & \\
$\Delta^\lambda$ & the \emph{standard module}, an $(\CA,\Olam)$-bimodule & \where{D:Standard}\\
$\nabla_\lambda$ & the \emph{right standard module} & \where{D:Standard}\\
$\nabla^\lambda$ & the \emph{costandard module} $\Hom_\kk(\nabla_\lambda,\kk)$ & \where{D:Standard}\\
$\Delta^\lambda(M)$, $\nabla^\lambda(M)$ & (co)standard modules with coefficients in $M$ & \where{D:Standard}\\
$\Delta$-, $\nabla$-filtration & filtrations with $\Delta$-, resp.\ $\nabla$-subquotients & \where{D:Standard}\\
$\<\ ,\ \>$ & the \emph{balanced pairing} $\nabla_\lambda\times\Delta^\lambda\to\Olam$ & \where{L:Respectable}\\
$\nabla^\lambda_{\Olam}$, $\partial_x$, $\partial^\lambda$ & $\Hom_{\Olam}(\nabla_\lambda,\Olam)$ and the map $x\mapsto\partial_x$ & \where{D:partialx}\alsowhere{Ahom}\\
$L^\lambda$, $\rad\Delta^\lambda$ & the head of $\Delta^\lambda$, and the radical of the pairing & \where{D:RespectableSimples}\\
$\Lambda_+$ & $\set{\mu|L^\mu\ne0}$: it labels the irreducibles & \where{D:RespectableSimples}\alsowhere{T:Respectable}\\
$D_\Delta$, $D_\nabla$ & the standard and costandard decomposition matrices & \where{D:RespectableSimples}\\
$P^\lambda$, $C_{\CA}$, $H_{\CA}$ & projective cover, Cartan matrix, $\diag(|\hlam|/|\olam|)$ & \where{D:RespectableSimples}\\
$[M:L]$ & Jordan--H\"older multiplicity & \\
$\Omega_\lambda$, $r_\Omega=|\Omega_\lambda|$ & the $\mathtt C_\theta$-orbit of $\lambda$, and its size & \where{respcompeasy}\\
$\theta_\Omega=\theta/r_\Omega$ & the order of the stabiliser of $\lambda$ & \\
$\beta^\lambda_{phoq}$, $\beta_x$, $\epsilon_\lambda$ & the units in a datum \emph{compatible with $\psi$} & \where{respcompeasy}\alsowhere{respcompatible}\\
\textup{(cR1)}--\textup{(cR5)}, \textup{(pcR1)}--\textup{(pcR3)} & the compatibility axioms for the two kinds of datum & \where{respcompeasy}\alsowhere{respcompatible}\\
$\sigma_\Omega$ & the induced automorphism of $\Olam[\lambda_\Omega]$ & \\
\end{notationtable}

\begin{notationtable}{\S9.\ Respectable bases on $\Hhd[1](A)$}
$f(X)=\prod_{k=1}^{\ell'}f_k(X)$ & the $\theta$-cyclotomic polynomial, factored as $\prod_kf_k$ & \where{setupresp}\\
$\ell'=\ell_0+\ell_1$ & the number of factors $f_k$ & \where{setupresp}\\
$\HH^{\kapp_k}_1(A)$ & the \emph{RaFH factors} $\kk[X]\Smash A/(f_k)$ & \\
$J^{\psi}=\sum_m\psi^{m}(J)$ & the smallest $\psi$-stable ideal containing $J$ & \\
$u_\lambda=\prod_{i=1}^{\lambda}f_{i-1}$ & the \emph{orbit polynomials}; $u_1=1$, $u_{\ell'+1}=f$ & \where{E:a}\\
$a^\lambda_{X^{j}b}=u_\lambda X^{j}b$ & the respectable basis built from the RaFH factors & \where{T:1}\\
$\lt$ & the \emph{leading term} of a basis element & \\
$H^{\ge k}=u_k\Hhd[1](A)$ & the ideal chain realising the strata & \where{L:1}\\
$\Lambda_\psi$, $\preceq$, $\lambda_\Omega$ & the poset of $\mathtt C_\theta$-orbits, and a representative & \\
$\mathcal D^{k}_\Omega$ & the \emph{orbit superalgebra} attached to $\Omega$ & \where{E:orbitD}\\
$\Gamma$ & the weight poset of the orbit datum & \where{T:orbit}\\
$\mathtt N_\mu$ & the rank of the stratum $A^{(\mu)}$ & \\
$\Lambda^{\op{split}}$, $e^\lambda_i$ & the refined poset splitting a stratum with $\Olam\cong\kk^{m_\lambda}$ & \where{L:splitstratum}\\
\end{notationtable}

\begin{notationtable}{\S10.\ Respectable CRaFH algebras}
$I=\set{1>2>\dots>N}$, $N=|I|$ & the weight poset of the datum of $A$; $I^{n}$ ordered componentwise & \where{setups9}\\
$\Delta^{\bi}$, $\mathcal T_{\Ao}$ & the standard $\Bo$-modules, and the induced datum & \where{Tensordatum}\alsowhere{P:FrobeniusRespect}\\
$\Comp[\sL,n]$ & compositions of $n$ into $\sL$ parts & \where{SS:Tableaux}\\
$\sL$ & the number of components of a multipartition & \where{SS:Tableaux}\\
$[\lambda]$, $(r,c)$ & the diagram of a partition, and its \emph{nodes} & \where{SS:Tableaux}\\
$\blam=(\lambda^{(1)}|\dots|\lambda^{(\sL)})$ & an \emph{$\sL$-multipartition} of $n$ & \where{SS:Tableaux}\\
$\Parts[\sL,n]$ & the set of $\sL$-multipartitions of $n$ & \where{Defmultipart}\\
$\gedom$, $\gdom$ & the \emph{dominance order}, and its strict version & \where{Defmultipart}\\
$>_{\mathrm{lex}}$ & the lexicographic order on multipartitions and nodes & \where{Defmultipart}\\
$\Sym_\blam$, $|\blam|$ & the Young subgroup of $\blam$, and its composition & \where{D:Wlam}\\
$\Sym^\blam$, $\Sym^{|\blam|}$, $\leftindex^\blam\Sym$ & minimal length coset representatives & \where{D:Wlam}\\
$\bT$, $\Shape(\bT)$ & a \emph{$\blam$-tableau}, and its shape & \where{D:tableaux}\\
$\RStd(\blam)$, $\Std(\blam)$ & row standard and standard $\blam$-tableaux & \where{D:tableaux}\\
$\bT_{\downarrow m}$ & the subtableau of $\bT$ on $\set{1,\dots,m}$ & \\
$\Tlam$ & the \emph{initial}, $\gedom$-maximal, $\blam$-tableau & \\
$d_{\bS}$ & the minimal length $d$ with $\bS=d\Tlam$ & \where{E:RowStd}\\
$\HJpm_\blam$, $\MJ^\pm_\blam$ & $\HJpm_J$ for $\Sym_J=\Sym_\blam$, and $\Hfind\HJpm_\blam$ & \where{leveloneparab}\\
$H^\pm_\bS$, $H^{\pm*}_\bT$ & $\hbar^{-2\ell(d_\bS)}\xi^\pm_{d_\bS}H_{d_\bS}$, and its reverse & \where{E:HS}\\
$\ilam$ & the word in $I^{n}$ read off from the components of $\blam$ & \where{fancydatum}\alsowhere{D:chargeparts}\\
$\plam[\blam]^\pm$, $\hlam[\blam]^\pm$, $\qlam[\blam]^\pm$ & the index sets of the datum of $\Hfind$, resp.\ $\Hhd$ & \where{fancydatum}\alsowhere{E:ellBasis}\\
$H^\pm_{\bS,p,h,q,\bT}$ & the respectable basis element $H^\pm_\bS p\HJpm_\blam hqH^{\pm*}_\bT$ & \where{E:HsphqT}\alsowhere{E:ellBasis}\\
$\Basis$, $\Basis_\blam$ & the respectable basis, and its $\blam$-stratum & \where{fancyB}\alsowhere{E:ellBasis}\\
$\HH^{\pm\ge\blam}_n$, $\HH^{\pm>\blam}_n$ & spans of the $\Basis_\bmu$ with $\bmu\ge\blam$, resp.\ $\bmu>\blam$ & \where{fancyB}\\
$\ge$ on $\Parts$ & the weight order, refined by $>_{\mathrm{lex}}$ & \where{D:contorder}\alsowhere{D:lcontorder}\\
$\op{wt}(\blam)$ & $\blam$ read as a $\mathfrak{gl}_N$-weight & \where{R:glweights}\\
$\binom n\mu=n!/\mu!$ & multinomial coefficient & \\
$\Tos_\blam$ & the subalgebra of $\Tos$ attached to $\blam$ & \where{L:HJ}\\
$\bG_\alpha=\bG^\blam_\alpha$ & the \emph{Garnir tableau}, belt and datum at a node $\alpha$ & \where{D:GarnirTableau}\\
$\op{Inv}(w)$ & the inversion set of $w$, of size $\ell(w)$ & \\
pure element & $y=\pm hq$ with $h\in\hlam[\ilam]$, $q\in\qlam[\ilam]$ & \\
$\bom_\bj$ & $(1^{m_1}|\dots|1^{m_N})$: the multiplicities of $\bj$ & \\
$\Bo^{>\ilam}$ & the span of the basis elements of $\Bo$ above $\ilam$ \\
$\mathcal C$, $m_{\bS\bT}$ & the coinvariant algebra, and the Murphy basis of $\kk\Sym$ & \\
$\mathbf H_{0,c}$, $\bar{\mathbf H}_c$, $\mathcal Z$ & rational Cherednik algebra, restricted version, bi-invariants & \\
$\lParts$ & the $\ell N$-multipartitions of $n$: the level $\ell$ weight poset & \where{SS:HigherLevelRespectable}\\
$\blam^{[l]}$ & the \emph{$\kapp_l$-part} of $\blam$ & \where{D:chargeparts}\\
$i$-layer & the components $\lambda^{(lN+i)}$ of $\blam$ & \where{D:chargeparts}\\
$n_l$ & boxes of $\blam$ before the $\kapp_l$-part, and its \emph{positions} & \where{D:chargeparts}\\
$X_\blam$ & $\prod_{l\ge2}\prod_{m\le n_l}(X_m-\kapp_l)$, after Dipper--James--Mathas & \where{D:Xlam}\\
$\Hhd[n,m]$ & the subalgebra of $\Hhd$ on $H_1,\dots,H_m$ and $\Bo$ & \\
\end{notationtable}
\end{document}